\documentclass[a4paper,11pt,twoside,openright]{report}
\usepackage[utf8]{inputenc}
\usepackage[T1]{fontenc}
\usepackage[english]{babel}
\usepackage[hidelinks]{hyperref}
\usepackage{verbatim}
\usepackage{amssymb}
\usepackage{dsfont}
\usepackage{subfigure}
\usepackage{amsmath}
\usepackage{amsthm}
\usepackage{graphicx}
\usepackage{lmodern}
\usepackage{color}
\usepackage{algpseudocode}
\usepackage{enumerate}
\usepackage{cases}
\usepackage{textcomp}
\usepackage{url}
\usepackage{pgfplots}
\usepackage{chngcntr}
\usepackage{remove}
\usepackage{setspace}
\usepackage{mathtools}
\usepackage{epigraph}

\usepackage{etoolbox}

\usepackage[left=1in,top=1in,right=1in,bottom=1in]{geometry}
\usepackage{fancyhdr}
\fancypagestyle{plain}{%
                    \fancyhf{}
                    
                    }
\newcommand{\helv}{%
\fontsize{9}{11}\selectfont}
\renewcommand{\sectionmark}[1]%
{\markboth{\MakeUppercase{\thechapter.\ #1}}{}}
\renewcommand{\subsectionmark}[1]%
{\markright{\MakeUppercase{\thesection.\ #1}}}

\numberwithin{equation}{chapter}
\makeatletter
\@removefromreset{equation}{section}
\makeatother

\theoremstyle{plain}

\newtheorem{lemme}{Lemme}[section]

\newtheorem{thm}[lemme]{Theorem}

\theoremstyle{remark}

\theoremstyle{definition}
\newtheorem{defi}[lemme]{Definition}

\newcommand\numberthis{\addtocounter{equation}{1}\tag{\theequation}}

\newcommand{\phat}{\hat{p}}

\newcommand{\Nb}{\mathbb{N}}
\newcommand{\Cb}{\mathbb{C}}
\newcommand{\Rb}{\mathbb{R}}

\newcommand{\Zb}{\mathbb{Z}}
\newcommand{\Runo}{\mathbb{R}}
\newcommand{\Rdue}{\mathbb{R}^{2}}

\newcommand{\Rd}{\mathbb{R}^{d}}

\newcommand{\complex}{\mathit{i}}

\newcommand{\into}{\int_{\Omega}}

\newcommand{\dx}[1]{\frac{\partial #1}{\partial x}}
\newcommand{\dy}[1]{\frac{\partial #1}{\partial y}}
\newcommand{\du}[1]{\frac{\partial #1}{\partial u}}
\newcommand{\dz}[1]{\frac{\partial #1}{\partial z}}
\newcommand{\dt}[1]{\frac{\partial #1}{\partial t}}
\newcommand{\Dt}[1]{\frac{\text{d} #1}{\text{d} t}}
\newcommand{\dn}[1]{\frac{\partial #1}{\partial \bold{n}}}

\newcommand{\dxx}[1]{\frac{\partial^2 #1}{\partial x^2}}
\newcommand{\dyy}[1]{\frac{\partial^2 #1}{\partial y^2}}

\newcommand{\dtt}[1]{\frac{\partial^2 #1}{\partial t^2}}
\newcommand{\Dtt}[1]{\frac{\text{d}^2 #1}{\text{d} t^2}}

\newcommand{\bigo}[1]{\mathcal{O}\left( #1 \right)}
\newcommand\norm[1]{\left\lVert#1\right\rVert}

\newcommand{\bu}{\bold{u}}
\newcommand{\Rey}{\text{Re}}
\newcommand{\tn}{t_n}
\newcommand{\tnpu}{t_{n+1}}
\newcommand{\bun}{\bu_n}
\newcommand{\bunpu}{\bu_{n+1}}
\newcommand{\pn}{p_n}
\newcommand{\pnpu}{p_{n+1}}
\newcommand{\bus}{\bu^*}

\newcommand{\busnpu}{\bu^*_{n+1}}
\newcommand{\fho}{\hat{f}(\omega)}
\newcommand{\intinf}{\int^\infty_{-\infty}}
\newcommand{\intAd}{\int_{-A/2}^{A/2}}
\newcommand{\ru}{r_1}
\newcommand{\rd}{r_2}
\newcommand{\ophi}{\overline{\phi}}
\newcommand{\ou}{\overline{u}}

\begin{document}

\begin{titlepage}

\begin{center}

\underline{}
\includegraphics[width=0.3\textwidth]{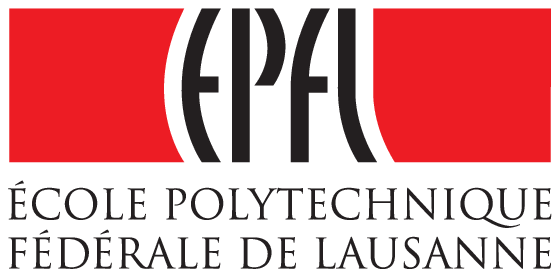}\\[1.4cm]

\textsc{\Large Master Thesis}\\[0.8cm]

\rule{\textwidth}{0.6mm} \\[0.4cm]
{\huge \bfseries Application of Stabilized Explicit\\[0.1cm]Runge-Kutta Methods to the\\[0.1cm]Incompressible Navier-Stokes Equations\\[0.1cm]by means of a Projection Method and a\\[0.2cm]Differential Algebraic Approach}\\[0.4cm]

\rule{\textwidth}{0.6mm} \\[1.3cm]


\begin{minipage}{0.4\textwidth}
\begin{flushleft} \large
\emph{Author:}\\
Giacomo \textsc{Rosilho de Souza}
\end{flushleft}
\end{minipage}
\begin{minipage}{0.4\textwidth}
\begin{flushright} \large
\emph{Supervisor:} \\
Prof. Assyr \textsc{Abdulle}
\end{flushright}
\end{minipage}

\vfill

June 21, 2014

\end{center}

\end{titlepage}

\newpage
\thispagestyle{empty}
\mbox{}
\clearpage



\newpage
\thispagestyle{empty}

\begin{abstract}
In this master thesis we have compared different second order stabilized explicit Runge-Kutta methods when applied to the incompressible Navier-Stokes equations by means of a projection method and a differential algebraic approach. We explored the stability and accuracy properties of the RKC, ROCK2 and PIROCK schemes when coupled with the projection and the differential algebraic approach. PIROCK has shown unexpected instabilities, ROCK2 resulted to be the most efficient and versatile Runge-Kutta method taken into account. The differential algebraic approach sounds computationally costly but it exhibits better accuracy and a larger stability region. These properties make it more efficient than the projection method. The theory presented in the first chapters is supported by numerical experiments.
\end{abstract}
\clearpage

\newpage
\thispagestyle{empty}
\mbox{}
\clearpage

\tableofcontents

\chapter{Introduction}

The aim of this master thesis is to compare the three explicit Runge-Kutta methods RKC, ROCK2 and PIROCK when solving the incompressible Navier-Stokes equations by means of a projection method or a differential algebraic approach. The dimensionless form the Navier-Stokes equations is
\begin{subequations} \label{eq:ns}
\begin{numcases}{}
\dt{\bu} + \left(\bu\cdot \nabla\right)\bu +\nabla p = \nu \nabla^2 \bu & $x\in\Omega$, $t\geq t_0$, \label{eq:ns1} \\
\nabla \cdot \bu = 0 & $x\in\Omega$, $t\geq t_0$, \label{eq:ns2}
\end{numcases}
\text{with Dirichlet boundary condition}
\begin{align} \label{eq:ns3}
\bu(t,x)=\bu_b(t,x) \quad x\in \partial \Omega,\; t\geq t_0
\end{align}
\text{and initial condition}
\begin{align} \label{eq:ns4}
\bu(t_0,x)=\bu_0(x) \quad x\in \Omega.
\end{align}
\end{subequations}
Here $\Omega$ is a bounded domain in $\Rdue$, $t_0\in \Rb$, $\bu=(u,v)$ is the velocity, $p$ is the pressure and $\nu=1/\Rey$ where Re is the Reynolds number. No boundary and initial conditions exists for the pressure. In fact it is given by a hidden constraint and imposing additional conditions would over-determine the system.

The Navier-Stokes equations are the mathematical formulation of many physical phenomena. They may be used to model the weather, ocean currents, water flow in a pipe and air flow around a wing. They help with the design of aircraft and cars, the study of blood flow, the design of power stations, the analysis of pollution, and many other things. Coupled with Maxwell's equations they can be used to model and study magnetohydrodynamics. For this broad range of applications finding effective numerical methods for the approximation of the incompressible Navier-Stokes equations is of significant importance.

The spatial discretization of equation \eqref{eq:ns1} leads to a mildly-stiff system of ordinary differential equations. For this kind of systems it is not worth to use a computationally expensive implicit method. Since the system is not extremely stiff stabilized explicit methods are able to handle the stiffness of the system, in the following we will consider the three methods RKC, ROCK2 and PIROCK. 

The above mentioned Runge-Kutta methods are of order two for ordinary differential equations. Moreover they are adaptive in the number of stages $s$ and have a stability region which increases in size quadratically with $s$. RKC \cite{rkc} and ROCK2 \cite{rock2} are explicit, PIROCK \cite{pirock} is also explicit when applied to the Navier-Stokes equations. The mentioned properties make these methods very attractive for the solution of mildly stiff problems of large dimension, like the one given by equation \eqref{eq:ns1} after spatial discretization.

However one cannot simply integrate equation \eqref{eq:ns1} for two reasons. The first one is that we need the value of the pressure $p$, the second is that the incompressibility constraint \eqref{eq:ns2} has to be taken into account. It turns out that the pressure is just a Lagrange multiplier \cite{chorin}  and its value is so that the velocity is divergence free. Hence the pressure value and the incompressibility constraint are handled together. In practice two different approaches are used to deal with the incompressibility constraint. 

The first one is the use of projection methods. Projection methods have been proposed by Chorin \cite{chorin} and Temam \cite{temam} independently years ago. They use a fractional step approach in which an intermediate velocity (often called virtual velocity) is obtained by solving the momentum equation \eqref{eq:ns1} keeping a constant pressure and disregarding the incompressibility constraint \eqref{eq:ns2}. Then the virtual velocity is projected into the manifold of divergence free fields by solving a Poisson equation for the pressure. Due to they decoupled nature projection methods are much more efficient than fully coupled techniques. The price been
paid, as we will see below, is that it introduces a numerical boundary layer on the velocity field and the velocity–pressure decoupling adversely affects the temporal accuracy of the numerical scheme. After the work of Chorin and Temam many new projection methods appeared \cite{sec_acc,acc_an,classfully,kim-moin,vanKan}, in the following we will concentrate on the method given in \cite{rkcp,gresho}. This is a second order accurate method for the velocity and first order for the pressure. An additional projection for the acceleration gives a second order accurate pressure when needed.

The second approach used to handle the incompressibility constraint is based on the differential algebraic nature of the incompressible Navier-Stokes equations. When discretized in space these equations are differential algebraic of differentiation index 2, the algebraic constraint comes from incompressibility. When Runge-Kutta methods are applied to this kind of equations \cite{ode-ii,acc_an} the algebraic constraint has to be satisfied in each stage and new order conditions arise. In the following we will see that for the particular case of the Navier-Stokes equations the new order conditions affect only the pressure accuracy and moreover they can be circumvented. In fact we will not consider the new conditions but instead other methods that allow to bring the pressure to the same order of accuracy as the velocity \cite{acc_an}.

The goal of this master thesis is to compare RKC, ROCK2 and PIROCK when coupled with the projection method or the differential algebraic approach and see which one of these couplings is the most efficient. RKC and ROCK2 are well suited for diffusion dominated flows because their stability domain is in a neighborhood of the negative real axis. PIROCK instead is more versatile and behaves very well also for advection dominated flows. Unfortunately his partitioned structure has shown an unexpected instability when the projection method has been applied. We did not had time to investigate this issue, hence we focused on the RKC and ROCK2 methods.

For the description of the algorithm we will restrict our attention to homogeneous Dirichlet boundary conditions and assume that there are no external forces. We will make use of two dimensional test cases, which capture much of the computational difficulties for incompressible flow calculations. We will also assume that the mesh spacing is uniform in the $x$ and $y$ directions. These restrictions are not inherent limitations of the method but they have been adopted here for clarity of exposition. 

The present report is organized as follows. In the first section we introduce the stabilized explicit Runge-Kutta methods and the projection method. In the second section we give and introduction to differential algebraic equations, we define the Runge-Kutta methods for semi-explicit differential algebraic equation of index 2 with emphasis on the Navier-Stokes equations and at the end we apply these methods to the stabilized Runge-Kutta methods. The fourth section contains the description of the spatial discretization and the method used for solving the Poisson problem. Finally in section five we show and discuss the results of the numerical experiments.

\chapter{Stabilized explicit Projection method for the Navier-Stokes equations}\label{sec:serkpm}
The first difficulty in solving the incompressible Navier-Stokes equations arises from the coupling of the velocity and the pressure terms. In the 1960s, Chorin \cite{chorin} proposed a projection method in which the approximation of the velocity and the pressure is done through a predictor-corrector procedure legitimated by the Helmholtz-Hodge decomposition. In the first step an intermediate velocity field (virtual velocity) is computed by solving the momentum equation \eqref{eq:ns1} ignoring the incompressibility constraint \eqref{eq:ns2} by taking a constant pressure. In the second step the intermediate velocity is projected into the space of divergence free vector fields and the corrector term is used to obtain the new pressure as well. Due to the decoupled approach this method is much more efficient than fully coupled techniques, this advantage makes the projection method attractive. Many improved methods based on this methodology have been seen in the past years, in this chapter we will present two such methods.

\noindent
Another difficulty arises from the spatial discretization of the incompressible Navier-Stokes equations. In the projection method the computation of the intermediate velocity is done integrating equation \eqref{eq:ns1} discretized in space. This semi-discrete system of equations is mildly stiff thus an integrator with a large stability region has to be used, moreover in order to avoid the solution of non-linear systems in must be explicit. Stabilized explicit Runge-Kutta methods are explicit methods with a large stability region along the negative real axis. Thanks to these schemes one handles the stiffness of the equations without solving non-linear systems.

In this chapter firstly we give a short introduction to Runge-Kutta methods, then we present the three stabilized explicit methods RKC, ROCK2 and PIROCK pointing out the major similarities and differences between them. In the second section we introduce two projection methods. The first one will be used in our numerical experiments in chapter \ref{sec:numexp}, the second one is less general and could not be applied in our case. At the end we point out a minor issue of projection methods in general.

\section{Stabilized Explicit Runge-Kutta methods}\label{sec:serk}
Historically, the goal when constructing a Runge-Kutta formula was to achieve the highest order possible with a given number of stages. In the case where the constructed method is explicit this approach leads to a relatively small stability domain which permits the integration of non stiff systems of equations only. A stiff problem would demand a severe restriction on the time step which prohibits the use of explicit methods. For this reason implicit methods with large stability domains were employed for the solution of stiff equations, requiring the solution of non linear systems at each time step.

Stabilized explicit Runge-Kutta methods are different in the sense that the main goal is to get a stability region which is as large as possible, in a sense that depends on the application. Thanks to this property one can explicitly integrate mildly stiff problems.

In this section we will consider three explicit methods which have a large stability domain along the negative real axis. This kind of schemes are very suitable for mildly stiff problems of large dimension where the eigenvalues are close to the negative real axis, like the ones arising from the spatial discretization of parabolic equations. In the next section we give an introduction to Runge-Kutta methods. In the following we present RKC and ROCK2 which are very similar, then we describe PIROCK, which strongly depends on ROCK2.

\subsection{Introduction to Runge-Kutta methods}
In this section an introduction to Runge-Kutta (RK) methods is given. We will restrict the presentation to the concepts needed to understand the rest of the discussion, for an exhaustive presentation of the subject we refer to \cite{ode-i,ode-ii}.

\subsubsection{Definition of the method}
Let us consider the Cauchy problem
\begin{equation} \label{eq:cauchy}
\begin{cases}
y'(t) = f(t,y(t)) & t>t_0, \\
y(t_0) = y_0,
\end{cases}
\end{equation}
for which to we wish numerically approximate the solution. Here $y:\Runo \rightarrow \Rd$, $f:\Runo \times \Rd \rightarrow \Rd$ is a smooth function and $d\in \Nb$.

Let $\Delta t>0$ be the time step and $t_n=t_0+n\Delta t$ for $n\in\Nb$. Given $y_n$ an approximation of $y(t_n)$ the next approximation $y_{n+1}$ of $y(t_{n+1})$ is given by the following definition.
\begin{defi}
Let $b_i$ and $a_{ij}$ ($i,j=1,\ldots,s$) be real numbers and let $c_i = \sum_{j=1}^s a_{ij}$. An $s$ stage Runge-Kutta method is given by
\begin{align*}
g_i &= y_n + \Delta t \sum_{j=1}^s a_{ij} f(t_n+c_j \Delta t, g_j)  \qquad i=1,\ldots ,s, \\
y_{n+1} &= y_n + \Delta t \sum_{i=1}^s b_i f(t_n + c_i \Delta t,g_i) .
\end{align*}
Setting $k_i:=f(t_n+c_i \Delta t,g_i)$ an equivalent definition is
\begin{align*}
k_i &= f(t_n+c_i \Delta t, y_n + \Delta t  \sum_{j=1}^s a_{ij} k_j)  \qquad i=1,\ldots ,s, \\
y_{n+1} &= y_n + \Delta t \sum_{i=1}^s b_i k_i .
\end{align*}
\end{defi}
For general $a_{i j}$ the computation of $y_{n+1}$ requires the solution of a non-linear system, this should be avoided whenever possible. It can be done setting $a_{ij}=0$ for $j\geq i$, in this case the RK method is called explicit, otherwise it is called implicit.

In the following two sections we will talk about two very important concepts in the framework of the RK methods, stability and accuracy. After we will explain how the time step $\Delta t$ can be chosen dynamically.

\subsubsection{Stability}
Let us consider the special case where $f$ does not depend on $t$ and is linear in $y$, i.e. $f(t,y)=\lambda y$ with $\lambda \in \Cb$ and  $y(t)=y_0 e^{\lambda t}$. Applying the RK method it yields $y_1=R(z)y_0$, where $z=\lambda \Delta t$ and $R(z)$ is given by $a_{ij}$ and $b_i$. $R(z)$ is called the stability function of the RK method, it is rational if the method is implicit and a polynomial if it is explicit. 

Suppose now that $\text{Re}(\lambda)<0$, it follows that $\lim_{t\rightarrow \infty}y(t)=0$ and one would like to preserve this property when $y(t)$ is numerically approximated by the RK method, i.e. $\lim_{n\rightarrow \infty}y_n=0$. Since 
\begin{align}\label{eq:yn}
y_n=R(z)^n y_0
\end{align}
it is necessary that $|R(z)|<1$. If $|R(z)|\leq 1$ then $y_n$ remains bounded and the method is called stable, otherwise $y_n$ starts oscillating and explodes as $n\rightarrow \infty$. The stability domain $\mathcal S$ of a RK method is defined by
\begin{align}
\mathcal{S}=\{ z\in \mathbb C : |R(z)|\leq1 \}.
\end{align}
Let $\lambda \in \mathbb{C}^-$, for some RK methods $\lambda \Delta t\in \mathcal S$ for all $\Delta t>0$ while for others a condition $\Delta t < \delta t(\lambda)$ exists, these methods are called conditionally stable. All the explicit methods are conditionally stable. In the special case $\lambda \in \Rb_-$ the condition becomes $\Delta t<l/|\lambda|$, where $l=-\inf \{ x\in \Rb : x\in \mathcal S \}$.

The concept of stability is extended to non linear functions $f$ in \cite[IV.2]{ode-ii}. A motivation is to preserve the Lyapunov stability of the system, which has significant importance when the computations are done in finite arithmetic. The parameter $\lambda$ is replaced by $\rho(t)$, the spectral radius of $\dy{f}(t,y(t))$.

\subsubsection{Order conditions}
Another important property is the accuracy of the numerical solution, this depends strictly on the order of the method, which is defined by the next definition.
\begin{defi}
A Runge-Kutta method has order $p$ if for all sufficiently regular problems \eqref{eq:cauchy} the local error $\delta y(t)$ satisfies
$$
\delta y(t_{n+1}):=y_{n+1} -y(t_n+\Delta t) = \mathcal{O}(\Delta t^{p+1}) \quad \text{as}\quad \Delta t \rightarrow 0,
$$
where $y(t_n)=y_n$.
\end{defi}
A Runge-Kutta method must satisfy some conditions on $a_{i j},c_i,b_i$ so that a certain order $p$ is achieved. The order conditions for $p\leq 2$ are given in \eqref{eq:ordcond}.
\begin{align}\label{eq:ordcond}
\left.
\begin{aligned}
\sum_{i=1}^s b_i &=1 \quad \text{for } p=1,\\
\sum_{i=1}^s b_i c_i &= \frac 1 2
\end{aligned}
\right\} \quad \text{for } p=2.
\end{align}
For $p=2$ the order conditions arise from the linearisation of $f$, which implies that $R(z)=1+z+\frac{1}{2}z^2 + \bigo{z^3}$ is enough to have second order accuracy.
The number of order conditions grows rapidly, for $p=6$ there is $37$ order conditions and for $p=10$ they are $1205$. For more details on how to find these conditions we refer to \cite[II.2]{ode-i}.

Also the internal stages $g_i$ have an order, it is given by the following definition.
\begin{defi}
We say that the $i$th stage $g_i$ has order $q$ if 
\begin{align}
\sum_{j=1}^s a_{ij}r(c_j)=\int_0^{c_i} r(t) dt \quad \forall r\in \mathbb{P}_{q-1} .
\end{align}
\end{defi}
\noindent
It follows that the $i$th stage has order $q$ if and only if $\sum_{j=1}^s a_{ij}c_j^{k-1} = \frac{1}{k}c_i^k$ for all $k<q$.

\subsubsection{Time step adaptivity with embedded formula}
As we saw in the previous sections the size of $\Delta t$ has an effect on the stability and the accuracy of the solution. The parameter $\rho(t)$ depends on time thus also the stability condition $\rho(t) \Delta t \in \mathcal S$ is time dependent. Likewise since the local error depends on the derivatives of $f$ it depends on $t$, which means that for different $t$s different $\Delta t$s are needed to achieve a required accuracy. These are the reasons for why it is important to dynamically choose $\Delta t$, in the following we will explain how to do that.

Consider two Runge-Kutta methods defined by $a_{ij},b_i$ and $a_{ij},\hat b_i$ of order $p$ and $\hat p$ respectively, where $\hat p<p$. Let $\Delta t_{n+1}=t_{n+1}-t_n$. We integrate from $t_n$ to $t_{n+1}$ computing the stages $k_i$, which are the same for both methods, and we define the error at step $t_{n+1}$ as
\begin{align*}
err_{n+1} := \norm{ \sum_{i=1}^s (b_i-\hat{b}_i) k_i },
\end{align*}
which is the difference between the solutions given by the two RK methods. If this difference increases then probably $\Vert y_{n+1}-y(t_{n+1})\Vert$ is increasing as well. The new time step $\Delta t_{new}$ is given by
\begin{align}\label{eq:dtnew}
\Delta t_{new} = \left(\frac{Tol}{err_{n+1}}\right)^{\frac{1}{\phat+1}} \left(\frac{err_{n}}{err_{n+1}}\right)^{\frac{1}{\phat+1}} \frac{\Delta t_{n+1}}{\Delta t_{n}} \Delta t_{n+1}
\end{align}
and checking $\Delta t_{new}\leq l/|\rho(t_n)|$, where $Tol$ is a user chosen parameter. If $err_{n+1} > Tol$ the $n+1$th step is recomputed defining $\Delta t_{n+1}=\Delta t_{new}$, otherwise $y_{n+2}$ is computed setting $\Delta t_{n+2}=\Delta t_{new}$. In practice there are some safety factors and checks in the implementation. For more details about time step adaptivity we refer to \cite[IV.8]{ode-ii}.

\subsection{The RKC and ROCK2 methods} \label{sec:rkc_rock2}
RKC (Runge-Kutta-Chebyshev)\cite{rkc} and ROCK2 (second order Orthogonal-Runge-Kutta-Che\-byshev)\cite{rock2} are two second order stabilized explicit Runge-Kutta methods intended for the time integration of parabolic partial differential equations. They are presented together in this section since they are very similar and share a lot of properties.

\subsubsection{A brief description}
The optimal stability polynomial of degree $s$ of a second order stabilized explicit method is characterized by
\begin{subequations}
\begin{align}
\overline R_s(z) &= 1 + z + \frac{z^2}{2!} + \sum_{i=3}^s \alpha_{i,s} z^i, \quad \text{with $\alpha_{i,s}\in \Runo$,} \label{eq:opt1}\\
|\overline R_s(z)| &\leq 1 \quad \text{for $z\in[-\overline l_s,0]$ with $\overline l_s$ as large as possible.} \label{eq:opt2}
\end{align}
For every $s$ these polynomials exist and are unique \cite{riha}. An analytic expression exists in terms of elliptic integrals (see \cite{ell_opt}) and the stability region's size is $\overline l_s = 0.821842 s^2$. The practical computation of such polynomials is done numerically. However the realization of these optimal polynomials as Runge-Kutta methods suffer from internal instabilities \cite{dumka} which could be avoided if one uses a recursion formula \cite{dumka,int_stab}. Both RKC and ROCK2 make use of Chebychev polynomials and recursion formulas in order to obtain a stability polynomial which is close the the optimal stability polynomials. Because of the recurrence relations the coefficients $a_{ij},b_i$ of the method are not explicitly given, nonetheless they can be computed recursively.

In order to include a strip of non zero width in the stability region a damping parameter $\eta$ is used. Condition \eqref{eq:opt2} is replaced by 
\begin{align}
|\overline R_s(z)| &\leq \eta \quad \text{for $z\in[-\overline l_s,0]$ with $\overline l_s$ as large as possible,} \label{eq:opt3}
\end{align}
\end{subequations}
where $\eta \in ]0,1[$. The smaller $\eta$ is, the wider the strip will be. On the other hand this will decrease $\overline l_s$ thus a compromise has to be found. Usually $\eta=0.95$.

One of the most remarkable properties of these methods is that they are adaptive in the number of stages $s$ and the stability domain's size $\overline l_s$ increases quadratically with $s$. Such schemes allow unrestricted integration steps as far as stability is concerned by simply taking $s$ large enough. Thus one has to choose the step size $\Delta t$ taking into account only the accuracy requirements. Another good property is that thanks to recursive properties of the stability polynomial they need only a few storage vectors which do not depend on the number of stages $s$.

But RKC and ROCK2 are not the same method, so they have differences. The first one is that the stability polynomial of RKC is available analytically, while it is computed numerically for ROCK2. The second difference is that the stability bound $\overline l_s$ of RKC increases as $0.653 s^2$, while it increases as $0.811 s^2$ for ROCK2. The internal stages of RKC are of order $2$, they are of order $1$ for ROCK2. The last difference is about time step adaptivity. In ROCK2 the local error is estimated using an embedded formula. In RKC an expression of the local error's leading term is analytically computed, this expression is then numerically approximated in order to obtain the estimation of the local error.

We will now define the RKC and ROCK2 methods.

\subsubsection{Definition of the RKC method}\label{sec:rkc}
Here we will define the RKC method following \cite{ode-ii,rkc}. The stability polynomial of RKC is realized by a three term recursion formula defining the internal stages. It uses scaled and shifted Chebychev polynomials. The Chebychev polynomial $T_s(x)$ of degree $s$ is defined by the recursion
\begin{align}\label{eq:cheby}
\begin{aligned}
T_0(x) &= 1, \\
T_1(x) &= x, \\
T_{n+1}(x) &= 2 x T_{n}(x)-T_{n-1}(x),
\end{aligned}
\end{align}
and the stability polynomial of RKC is defined as
\begin{align}\label{eq:pol_rkc}
R_s(z) = a_s+b_s T_s(w_0+w_1 z),
\end{align}
where $w_0=1+\varepsilon/s^2$ and $\varepsilon\ll 1$ is a parameter defining the damping $\eta=a_s+b_s\approx 1-\varepsilon/3$. The other parameters $w_1,a_s,b_s$ are chosen so that the second order conditions
\begin{align}
R_s(0)=1,\quad R_s'(0)=1,\quad R_s''(0)=1
\end{align}
are satisfied, this gives
\begin{align}
w_1=\frac{T'_s(w_0)}{T''_s(w_0)},\quad b_s=\frac{T''_s(w_0)}{\left(T'_s(w_0)\right)^2},\quad a_s=1-b_s T_s(w_0).
\end{align}
Thanks to the recursive properties \eqref{eq:cheby} of Chebychev polynomials in the case of a non linear initial value problem \eqref{eq:cauchy} the scheme defined by 
\begin{align}\label{eq:rkc_rec}
\begin{aligned}
g_0 -y_0 &= 0 \\
g_1 -y_0 &= \kappa_1 \Delta t f(g_0) \\
g_j-y_0 &=\mu_j\left( g_{j-1}-y_0 \right) +\nu_j\left( g_{j-2}-y_0 \right)+\kappa_j\Delta t \left(f(g_{j-1})-a_{j-1} f(g_0)\right)
\end{aligned}
\end{align}
realizes $R_s(z)$ defined in \eqref{eq:pol_rkc} as stability polynomial. In \eqref{eq:rkc_rec} the recursion coefficients are
\begin{align}
b_j&=\frac{T''_j(w_0)}{\left(T'_j(w_0)\right)^2},\quad &\nu_j &=\frac{-b_j}{b_{j-2}},\quad  &\kappa_j&=\frac{2b_j w_1}{b_{j-1}}, \\
a_j &= 1-b_j T_j(w_0), &\mu_j&=\frac{2 b_j w_0}{b_{j-1}}  &
\end{align}
for $j=2,\ldots ,s$ and
\begin{align}
b_0 = b_1 = b_2, \quad \kappa_1=c_1=\frac{c_2}{T'_2(w_0)}.
\end{align}

\begin{figure}[!tbp]
\begin{center}
\subfigure[Stability domain.]{
\includegraphics[trim=0.0cm 0.0cm 0.0cm 0.0cm, clip=true, width=0.45\textwidth]{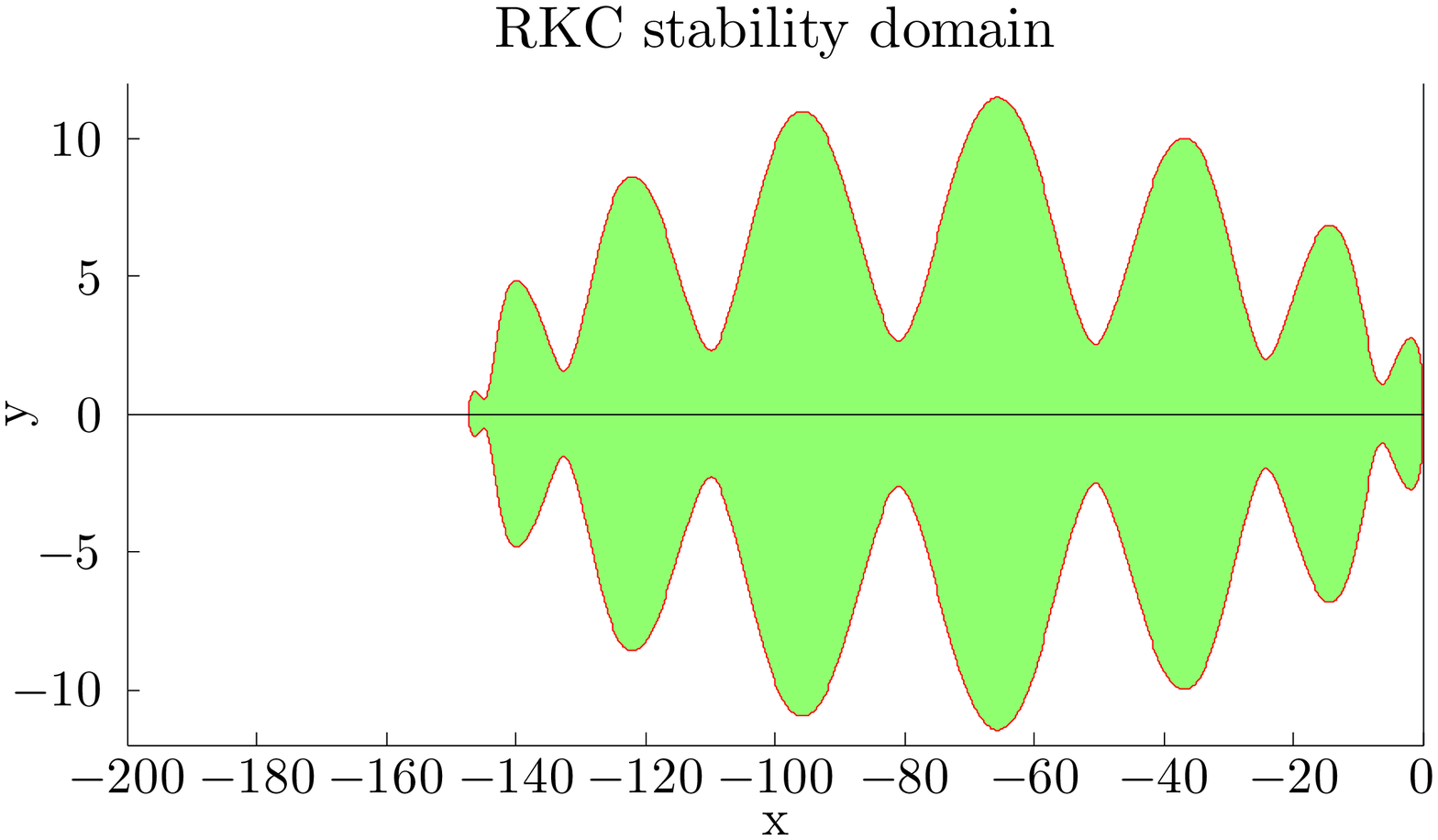} }
\subfigure[Stability polynomial.]{
\includegraphics[trim=0.0cm 0.0cm 0.0cm 0.0cm, clip=true, width=0.45\textwidth]{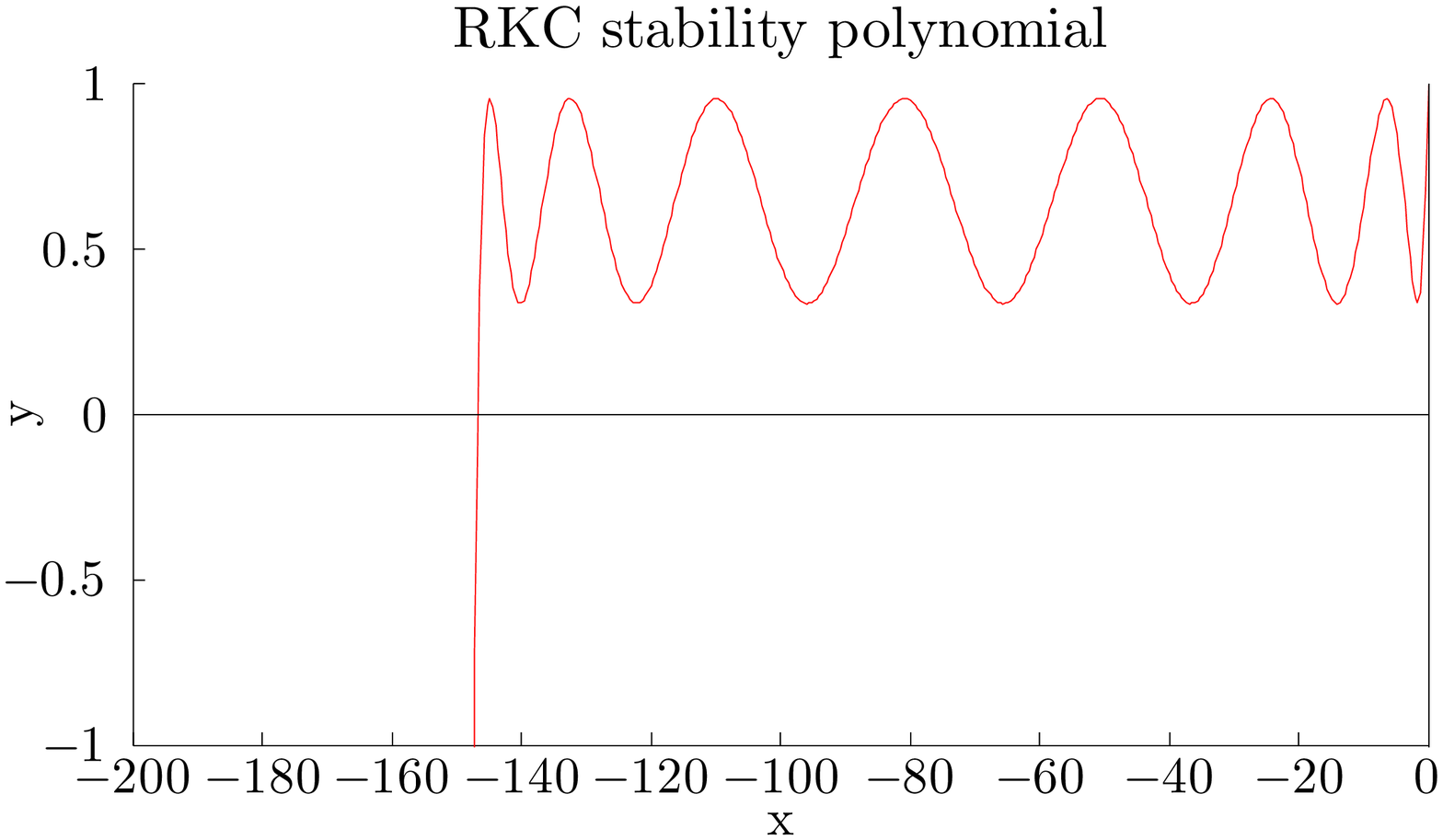} }
\end{center}
\caption{RKC stability domain and polynomial for $s=15$ and damping $\eta\approx 0.95$.}
\label{fig:rkc}
\end{figure}

In Figure \ref{fig:rkc} we show the stability domain and polynomial of RKC with $s=15$ stages and a damping $\eta\approx0.95$.

\subsubsection{Definition of the ROCK2 method} \label{sec:rock2}
The definition of ROCK2 and an exhaustive explanation about its derivation is given in \cite{rock2}, here we will only give its definition. 

The construction of the stability polynomial $R_s(z)$ of ROCK2 strongly relies on the fact that the optimal stability polynomial can be written as $\overline R_s(z)=\overline w(z) \overline P_{s-2}(z)$, where $\overline w(z)$ is a positive second degree polynomial and $\overline P_{s-2}(z)$ is a polynomial of degree $s-2$ with $s-2$ real roots (see \cite{roots}). In ROCK2 $\overline P_{s-2}(z)$ is approximated by a linear combination $P_{s-2}(z)$ of Chebychev orthogonal polynomials. $P_{s-2}(z)$ and the approximation $w(z)$ of $\overline w(z)$ are computed numerically by a fixed point algorithm described in \cite{rock2}. This algorithm computes a set of parameters which recursively define $R_s(z)=w(z) P_{s-2}(z)$, this parameters are stored in table which is used by ROCK2. The recursive scheme is defined by
\begin{align}\label{eq:defrock2}
\begin{aligned}
g_0 &= y_0 \\
g_1 &= y_0 + \Delta t \mu_1 f(g_0) \\
g_j &= \Delta t\mu_j f(g_{j-1})-\nu_jg_{j-1}-\kappa_j g_{j-2} \quad j=2,\ldots ,s-2 \\
g_{s-1} &= g_{s-2}+\Delta t\sigma f(g_{s-2}) \\
g^*_s &= g_{s-1}+\Delta t\sigma f(g_{s-1}) \\
g_s &= g^*_s - \Delta t\sigma \left(1-\frac{\tau}{\sigma^2}\right)\left(f(g_{s-1})-f(g_{s-2})\right)
\end{aligned}
\end{align}
where $\mu_j,\nu_j,\kappa_j,\sigma,\tau$ are computed by the fixed point algorithm and depend on $s$. When applied to the test problem $y'=\lambda y$ we have $g_j=P_{j}(z)y_0$ for $j=0,\ldots ,s-2$ and $g_s=w(z)P_{s-2}(z)y_0$. 

\begin{figure}[!hbtp]
\begin{center}
\subfigure[Stability domain.]{
\includegraphics[trim=0.0cm 0.0cm 0.0cm 0.0cm, clip=true, width=0.45\textwidth]{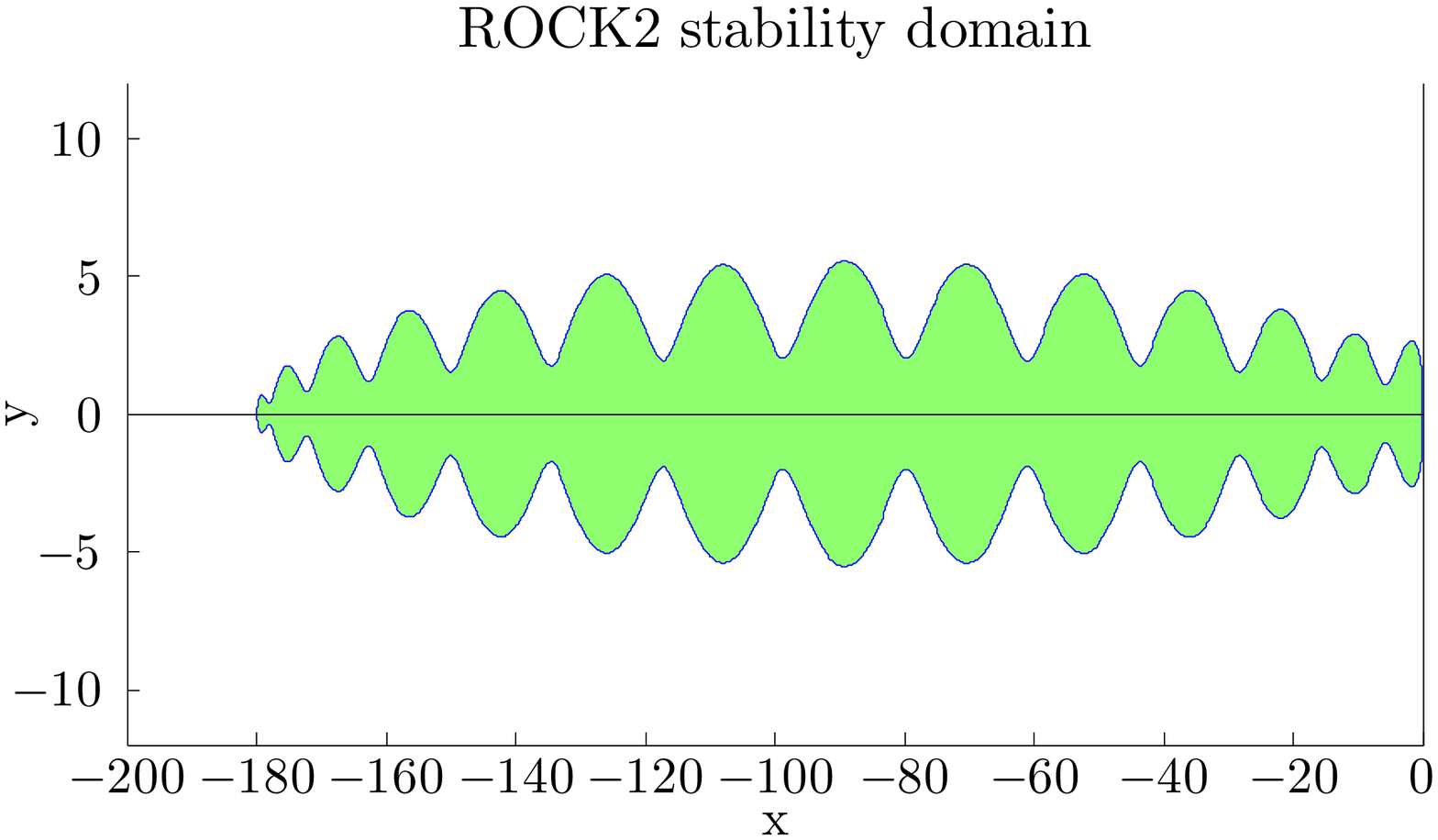} }
\subfigure[Stability polynomial.]{
\includegraphics[trim=0.0cm 0.0cm 0.0cm 0.0cm, clip=true, width=0.45\textwidth]{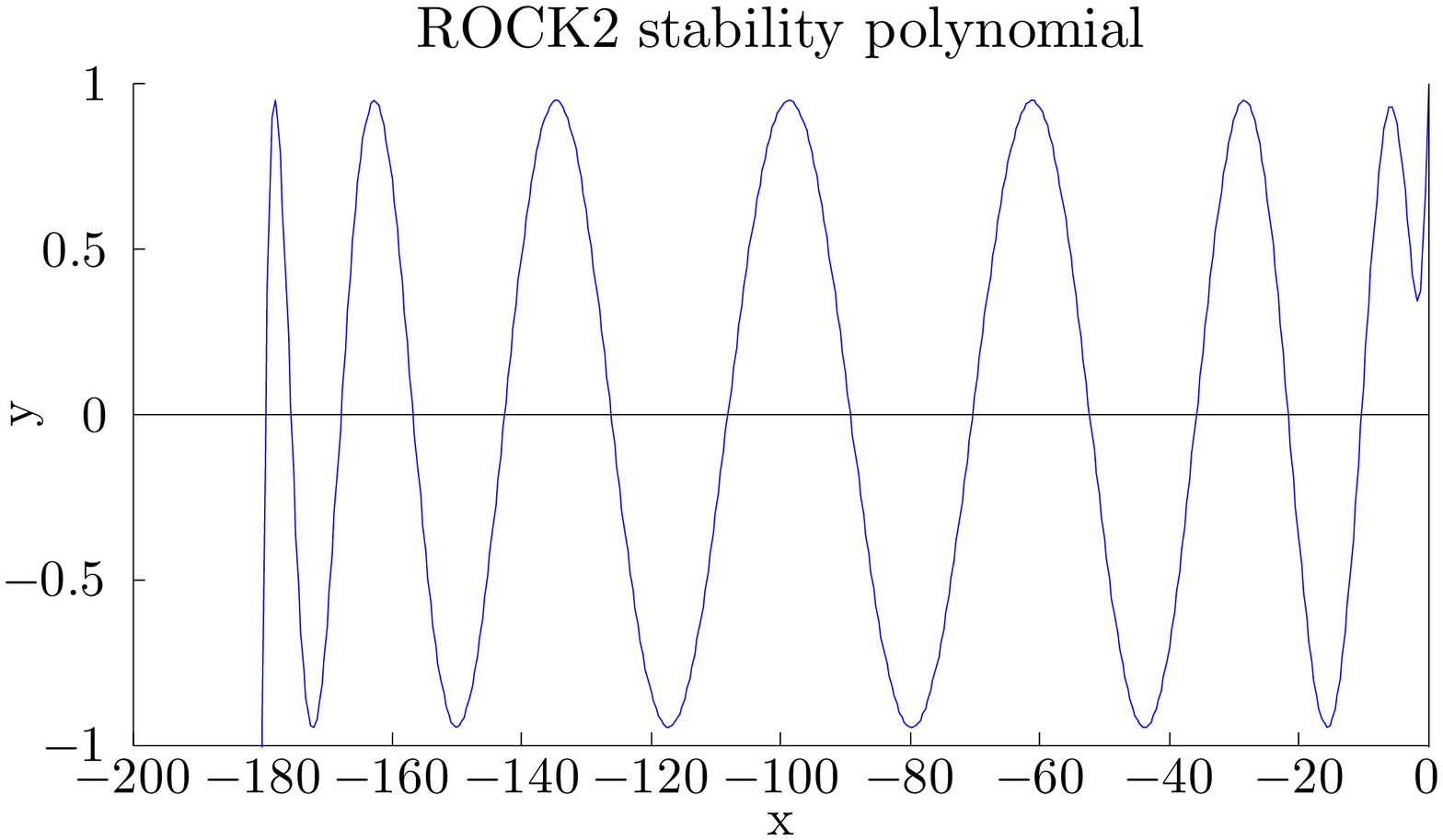} }
\end{center}
\caption{ROCK2 stability domain and polynomial for $s=15$ and damping $\eta=0.95$.}
\label{fig:rock2}
\end{figure}

In Figure \ref{fig:rock2} we show the stability domain and polynomial of ROCK2 with $s=15$ stages and a damping $\eta=0.95$.

\subsubsection{Time step and stages adaptivity}\label{sec:tscrr}
Both RKC and ROCK2 use formula \eqref{eq:dtnew} in order to choose the new time step $\Delta t_{new}$. Once that we have $\Delta t_{new}$ the number of stages $s$ is chosen so that $\Delta t_{new} \leq l_s/\rho(t)$.
The difference between the methods is in how $err_{n+1}$ is computed.

\subsubsection{RKC}
Let us suppose that $y_n=y(t_n)$ and we compute $y_{n+1}$ with the RKC method. The error is defined as
\begin{align}
err_{n+1} :&= \norm{y_{n+1}-y(t_{n+1})} \approx \norm{ \frac{\Delta t^3}{15} y''' (t_n)} + \bigo{\Delta t^4} \label{eq:erkc1}\\
&\approx \frac{1}{15} \left( 12(y_n-y_{n+1})+6\Delta t(f(y_n)+f(y_{n+1})\right) + \bigo{\Delta t^4}. \notag
\end{align} 
The approximation in \eqref{eq:erkc1} is obtained by a Taylor expansion of $y_{n+1}$ and $y(t_{n+1})$ around $t_n$ (both depend on $\Delta t$.).

\subsubsection{ROCK2}
ROCK2 uses an embedded formula in order to estimate $err_{n+1}$. We have that $g^*_s$ and $g_s$ of \eqref{eq:defrock2} are first and second order approximations of $y_{n+1}$ respectively. The error is defined as
\begin{align}
err_{n+1} := \norm{g^*_s-g_s} = \norm{\Delta t\sigma \left(1-\frac{\tau}{\sigma^2}\right)\left(f(g_{s-1})-f(g_{s-2})\right)}.
\end{align}

\subsection{The PIROCK method} \label{sec:pirock}
In this section we will give a short description and define the PIROCK method, for more details we refer to \cite{pirock}.

\subsubsection{A brief description}
PIROCK is a partitioned implicit-explicit orthogonal Runge-Kutta Chebyshev method for the time integration of diffusion-advection-reaction problems of the kind
\begin{equation} \label{eq:ode}
y' = F(y) = F_D(y) + F_A(y) + F_R(y), \quad y(0) = y_0 ,
\end{equation}
where $F_D, F_A, F_R$ are the space discretization of the diffusion, advection and reaction operators respectively. For each operator a well suited Runge-Kutta method is chosen and the coupling is done so that order conditions for additive Runge-Kutta methods are achieved. A symmetric diffusion operator has eigenvalues which lies in an interval of size $[-\bigo{\Delta x^{-2}}, 0]$, where $\Delta x$ is the grid size. For this mildly-stiff problem the ROCK2 \cite{rock2} explicit stabilized second order Chebyshev method is used. The advection operator has imaginary eigenvalues lying in $[-\complex\bigo{\Delta x^{-1}}, \complex \bigo{\Delta x^{-1}}]$ thus an explicit third order Runge-Kutta method is used. The reaction term can be very stiff thus an implicit second order integrator is used. This requires the solution of non-linear systems but since the reaction is local these systems are small and can be solved efficiently by LU decomposition and quasi-Newton methods. The case of a non-symmetric diffusion operator is handled by decomposing $F_D$ in its symmetric and asymmetric parts $F_{D_s}$ and $F_{D_a}$ respectively. Then $F_D$ is replaced by $F_{D_s}$ and $F_{D_a}$ is absorbed in $F_A$ with some minor modifications in the method.

The resulting method is second order accurate and is adaptive in time and in the number of stages $s$ (see section \ref{sec:rock2}). The number of evaluations of $F_A$ and $F_R$ is independent of the number of stages $s$ needed to handle the stiffness of $F_D$. Moreover PIROCK is stable even for non symmetric diffusion operators with eigenvalues located in a sector $S_{\pi/4}=\{ -\rho e^{\complex \tau} : \rho\geq~0, -\pi/4\leq \tau \leq \pi/4\}$. This property also makes PIROCK well suited for advection dominated flows.

\subsubsection{Definition of the PIROCK method}
PIROCK is defined by the following algorithm for $s\geq 3$:
\begin{align*}
&\text{Diffusion integration - ROCK2 method} \\
&K_1 = y_0 + \alpha \mu_1 \Delta t F_D(y_0) \displaybreak[0]\\
&K_j = \alpha \mu_j \Delta t F_D(K_{j-1}) - \nu_j K_{j-1} - \kappa_j K_{j-2} \quad j=2,\ldots ,s-2+l \;\; \text{($l=1$ or $2$)} \displaybreak[0]\\
&K^*_{s-1} = K_{s-2} + \sigma_\alpha \Delta t F_D(K_{s-2}) \displaybreak[0]\\
&K^*_{s} = K^*_{s-1} + \sigma_\alpha \Delta t F_D(K^*_{s-1})\displaybreak[0]\\
&\text{Advection-Reaction integration and coupling} \\
&K = K_{s-2+l} \displaybreak[0]\\
&K_{s+1} = K + \gamma \Delta t F_R(K_{s+1}) \displaybreak[0]\\
&K_{s+2} = K + \beta \Delta t F_D(K_{s+1}) + \Delta t F_A(K_{s+1}) + (1-2\gamma)\Delta t F_R(K_{s+1}) + \gamma \Delta t F_R(K_{s+2}) \displaybreak[0]\\
&K_{s+3} = K +(1-2\gamma)\Delta t F_A(K_{s+1})+ (1-\gamma)\Delta t F_R(K_{s+1}) \displaybreak[0]\\
&K_{s+4} = K + \frac{1}{3}\Delta t F_A(K_{s+1}) \displaybreak[0]\\
&K_{s+5} = K + \frac{2}{3}\beta \Delta t F_D(K_{s+1}) + \frac{2}{3} \Delta t J_R^{-1}F_A(K_{s+4}) + \left(\frac{2}{3}-\gamma\right) \Delta t F_R(K_{s+1}) \\ &\qquad \quad \;\, {} +\frac{2}{3}\gamma \Delta t F_R(K_{s+2}) \displaybreak[0]\\
&\text{Computation of $y_1$}\\
&y_1 = K_s^* -\sigma_\alpha\left(1-\frac{\tau_\alpha}{\sigma_\alpha^2}\right)\Delta t\left(F_D(K_{s-1}^*)-F_D(K_{s-2})\right)+\frac{1}{4}\Delta t F_A(K_{s+1}) \\
&\qquad {}+\frac{3}{4}\Delta t F_A(K_{s+5})+\frac{1}{2}\Delta t F_R(K_{s+1}) + \frac{1}{2} \Delta t F_R(K_{s+2}) \\
&\qquad {}+ \frac{J_R^{-l}}{2-4\gamma}\Delta t\left(F_D(K_{s+3})-F_D(K_{s+1})\right) \numberthis \label{eq:y1}
\end{align*}
where $\gamma=1-\sqrt 2 /2$, $\beta = 1-2\alpha P'_{s-2+l}(0)$, $J_R=I-\gamma \Delta t \dy{F_R}(K_s)$ and $\alpha=1$ for $l=2$ or $\alpha=1/(2P'_{s-1}(0))$ for $l=1$. For $l=2$ we obtain the original ROCK2 method. For $l=1$ we have a larger damping thus the stability domain will be shorter but with a wider extension on the imaginary axis.

\subsubsection{Time step and stages adaptivity}
Three embedded methods $y_{e,D}$, $y_{e,A}$, $y_{e,R}$ are used for time step adaptivity. We obtain three error estimators 
\begin{align*}
err_D = y_1-y_{e,D}, \quad err_A = y_1-y_{e,A},\quad err_R = y_1-y_{e,R}.
\end{align*}
The error is estimated as
\begin{align*}
err_{n+1} = \max\left(||err_{D}||,||err_A||^{2/3},||err_R||\right).
\end{align*}
The factor $2/3$ is due to the fact that $err_A=\bigo{\Delta t^3}$ while $err_D$ and $err_R$ are $\bigo{\Delta t^2}$. The new time step $\Delta t_{new}$ is chosen with formula \eqref{eq:dtnew} and the new number of stages $s$ and $l=1,2$ are chosen in order to keep the method stable.

\section{Projection method based on stabilized Runge-Kutta methods} \label{sec:proj}
As we already said projection methods use a fractional-step technique. Firstly they compute an intermediate virtual velocity and then it is projected into the manifold of divergence-free fields by solving a Poisson equation, the pressure is updated during the projection step as well. The velocity-pressure decoupling of the method makes it cheap but affects the temporal accuracy of the numerical scheme. It has been observed numerically and analytically that second order accuracy for the velocity is readily obtained, while second order accuracy for the pressure is more involved, see \cite{sec_acc,rkcp} and chapter \ref{sec:erkinse}. Moreover the virtual velocity has no physical meaning and artificial boundary conditions should be used. If the same boundary conditions for the virtual and the physical velocity are used the tangential component will loose accuracy near the boundary. The methods presented in sections \ref{sec:method1} and \ref{sec:method2} do not fix this issue thus the solution is affected by a numerical pollution within a boundary layer. The size of the boundary layer has been studied in \cite{gresho} for a particular projection method and results to be $\bigo{\sqrt{\Delta t /\Rey}}$. In section \ref{sec:method3} we propose a method that will be used in section \ref{sec:numexp} in order to motivate the employment of different boundary conditions for the virtual velocity. 

In this chapter we will present two projection methods for the solution of the incompressible Navier-Stokes equations \eqref{eq:ns}. The first one is very versatile and can be applied to any second or higher order time integrator. The second one has a limitation which rely on the advection term, in our case we could use it only for the Stokes equation. At the end we will present a fix which will be used in section \ref{sec:numexp} in order to point out a minor issue very common in projection methods \cite{gresho}. We stress on the fact that this fix is not applicable to real life problems.

\subsection{Method 1}  \label{sec:method1}
In this section we will present the projection method (PM1) which will be used in our numerical experiments of section \ref{sec:numexp}.

This numerical scheme is locally second order accurate for the velocity and first order accurate for the pressure. When necessary an additional projection allows to obtain for the pressure the same order of accuracy as for the velocity. The advantage of PM1 is that the projection step is independent from the integration step thus any second order time integrator can be used. We could apply this method to the three second order schemes RKC, ROCK2 and PIROCK in exactly the same way. 

In the following we will define the method, then we do a local error analysis and give a remark about time step adaptivity. Finally we define a variant which is also used in the numerical experiments. In the definition of PM1 and in the local error analysis we follow \cite{rkcp} but PM1 is also present in \cite{gresho}. Numerical experiments are shown in section \ref{sec:numexp}. In section \ref{sec:numexp} we refer to this method as PM1.

\subsubsection{Definition of the method}\label{sec:defpm1}
Our presentation differs from the one given in \cite{rkcp} in that here equations \eqref{eq:ns} are continuous in space.
The following numerical scheme returns a solution $\bunpu$ and $\pnpu$ at time $\tnpu$ given $\bun$ and $\pn$ at $t_n$.

\begin{enumerate}
\item Compute $\busnpu$ a second order approximation of $\bus(\tnpu)$ the exact solution of
\begin{numcases}{}
\dt{\bus} + \left( \bus \cdot \nabla \right) \bus + \nabla \pn = \nu \nabla^2 \bus & $\tn\leq t\leq\tnpu$, \label{eq:pc}\\
\bus(\tn) = \bun , \notag \\
\left. \bus \right|_\Gamma = \bu_b .\notag
\end{numcases}

\item Solve 
\begin{numcases}{}
\nabla^2 \phi_1 = \nabla \cdot \busnpu & in $\Omega$,\\
\dn{\phi_1}=0 & on $\partial \Omega$,
\end{numcases}
for $\phi_1$.
\item Update $\bunpu=\busnpu-\nabla\phi_1$ and $\pnpu=\pn+\frac{2}{\Delta t}\phi_1$.
\item If at time $t_m$ a second order approximation for the pressure $p_m$ is needed an additional projection for the acceleration is performed.

Set $\bold{F}_m = -\left(\bu_m\cdot\nabla\right)\bu_m - \nabla p_m + \nu \nabla^2\bu_m$, solve
\begin{subequations}  \label{eq:m1}
\begin{numcases}{}
\nabla^2 \phi_2 = \nabla \cdot \bold{F}_m & in $\Omega$,\\
\dn{\phi_2}=0 & on $\partial \Omega$,
\end{numcases}
\end{subequations}
for $\phi_2$ and update $p_m = p_m + \phi_2$.
\end{enumerate}
The motivation of the first projection comes from the constraint $\nabla \cdot \bu=0$ while the second projection arises from the hidden constraint $\nabla \cdot \dt{\bu}=0$.

With this approach the normal component of the velocity at the boundary satisfies the boundary conditions since
\begin{align}
\bold{n}\cdot \bunpu = \bold{n}\cdot \busnpu - \dn{\phi_1} = \bold{n}\cdot \busnpu = \bold{n}\cdot \bu_b(\tnpu) .
\end{align}
But $\frac{\partial \phi_1}{\partial \tau}\neq 0$ thus the tangential component of $\bunpu$ will be affected by a numerical pollution, indeed
\begin{align}
\bold{\tau}\cdot \bunpu = \bold{\tau}\cdot \busnpu - \frac{\partial \phi_1}{\partial \tau} = \bold{\tau}\cdot \bu_b(\tnpu) - \frac{\partial \phi_1}{\partial \tau} \neq \bold{\tau}\cdot \bu_b(\tnpu) .
\end{align}

The theoretical background of the projection \( \bunpu=\busnpu-\nabla\phi_1\) is the decomposition theorem of Helmholtz-Hodge \cite[1.16]{arfken-weber}. The theorem states that a vector field \(\bu\) on a simply connected domain can be uniquely decomposed into a divergence free part \(\bu_{div}\) and an irrotational part \(\bu_{irrot}\). Thus 
\begin{align}\label{eq:hhdec}
\bu = \bu_{div}+\bu_{irrot}=\bu_{div}+\nabla \phi
\end{align}
for some potential function \(\phi\). Taking the divergence of equation \eqref{eq:hhdec} yelds \(\nabla\cdot\bu=\nabla^
2\phi\). If the vector field \(\bu\) is known the potential \(\phi\) is found by solving the Poisson equation and the divergence free part of \(\bu\) can be extracted with the relation \(\bu_{div}=\bu-\nabla\phi\).

\subsubsection{Local error analysis}
In this section we will give the local error analysis of the above projection method as it is done in \cite{rkcp}. We fix $x\in \Omega$ and suppose that the exact solution is given at $\tn$. In the following the explicit $x$-dependence will be left out of notation. 

Suppose $\bun=\bu(\tn)$ and $p_n=~p(\tn)$, doing a formal Taylor expansion of equations \eqref{eq:ns1} and \eqref{eq:pc} around $t_n$ we get
\begin{align}\label{eq:lea1}
\bus(\tnpu)-\bu(\tnpu)=\frac{\Delta t^2}{2}\nabla \dt{p}(\tn) + \bigo{\Delta t^3}.
\end{align}
Since in \eqref{eq:pc} a second order accurate method is used we have 
\begin{align}\label{eq:lea2}
\busnpu-\bus\left( \tnpu\right)=\bigo{\Delta t^3}.
\end{align}
Summing equations \eqref{eq:lea1} and \eqref{eq:lea2} we get
\begin{align} \label{eq:t1}
\busnpu = \bu(\tnpu) +  \frac{\Delta t^2}{2}\nabla \dt{p}(\tn) + \bigo{\Delta t^3}.
\end{align}
From equation \eqref{eq:t1} one sees that setting 
\begin{align} \label{eq:t2}
\bunpu = \busnpu-\frac{\Delta t^2}{2}\nabla \dt{p}(\tn)
\end{align}
we obtain $\bunpu-\bu(\tnpu)=\bigo{\Delta t^3}$, i.e. a second order accurate method for the velocity. In order to compute the correction term we take the divergence of equation \eqref{eq:t1} which leads to
\begin{align}
\nabla \cdot \busnpu = \frac{\Delta t^2}{2}\nabla^2 \dt{p}(\tn)+ \bigo{\Delta t^3},
\end{align}
thus solving the equation
\begin{numcases}{}
\nabla^2 \phi_1 = \nabla \cdot \busnpu & in $\Omega$,\\
\dn{\phi_1}=0 & on $\partial \Omega$,
\end{numcases}
gives $\phi_1=\frac{\Delta t^2}{2} \dt{p}(\tn)+ \bigo{\Delta t^3}$ up to a constant and setting $\bunpu=\busnpu-~\nabla \phi_1$ is equivalent to equation \eqref{eq:t2}. Doing a Taylor expansion of the pressure $p(t)$ around $\tn$ we get
\begin{align}
\pnpu = \pn + \Delta t \dt{p}(\tn) + \bigo{\Delta t^2}
\end{align}
thus setting $\pnpu=\pn+\frac{2}{\Delta t} \phi_1$ leads to a first order approximation of the pressure.

A second projection on the acceleration gives a second order accurate pressure. Supposing that the approximation of the velocity is $p$th order accurate, $\bu_m-\bu(t_m)=\bigo{\Delta t^{p+1}}$, we can write
\begin{align}
\bold{F}_m &= -\left(\bu_m\cdot\nabla\right)\bu_m+\nu\nabla^2\bu_m-\nabla p_m \notag \\
&= -\left( \bu(t_m)\cdot\nabla\right)\bu(t_m) +\nu\nabla^2\bu(t_m)-\nabla p(t_m) +\nabla p(t_m) -\nabla p_m +\bigo{\Delta t^{p+1}} \notag \\ 
&= \dt{\bu}(t_n)+\nabla p(t_m) -\nabla p_m +\bigo{\Delta t^{p+1}} \label{eq:t3}
\end{align}
and taking the divergence of \eqref{eq:t3} gives
\begin{align}
\nabla \cdot \bold{F}_m= \nabla^2\left( p(t_m)-p_m\right) + \bigo{\Delta t^{p+1}}
\end{align}
which implies that $\phi_2$ of equation \eqref{eq:m1} satisfies $\phi_2 = p(t_m)-p_m+\bigo{\Delta t^{p+1}}$ up to a constant. It follows that setting $p_m = p_m + \phi_2$ brings the pressure to the same order of accuracy of the velocity. This is done only for output purposes since it is not necessary for maintaining second order velocity. In order to motivate this last assertion let us neglect the advection term, since $\busnpu$ has second order accuracy we can write
\begin{align}
\busnpu = \bun + \Delta t \nu\nabla^2 \bun - \Delta t \nabla p_n + \frac{1}{2}\Delta t^2\dtt{\bus}(t_n) +\bigo{\Delta t^3} \notag
\end{align}
and formally
\begin{align}
\bunpu &= \left( I-\nabla (\nabla^2)^{-1}\nabla \cdot\right)\busnpu \\
&= \left( I-\nabla (\nabla^2)^{-1}\nabla \cdot\right)\left(\bun + \Delta t \nu\nabla^2 \bun + \frac{1}{2}\Delta t^2\dtt{\bus}(t_n)\right) +\bigo{\Delta t^3}.
\end{align}
Thus the approximation of the pressure $p_n$ affects only the second and higher order terms of $\bunpu$.

Because of equation \eqref{eq:t2} this projection method cannot have convergence order higher than $2$, even if for computing $\busnpu$ one uses a third or higher order numerical integrator.

\subsubsection{Time step adaptivity}\label{sec:projtsc}
The estimation of the error $err_{n+1}$ is done before the projection, thus using $\bus_{n+1}$ and not $\bu_{n+1}$. This means that the real local error $\norm{\bu_{n+1}-\bu(t_{n+1})}$ is overestimated because we are taking into account also the non divergence free part of $err_{n+1}$.

\subsubsection{Variant} \label{sec:pm1var}
As it is done in \cite{rkcp} we have implemented a variant of PM1 in RKC and ROCK2. It consists in projecting each stage $g_i$ during the computation of $\busnpu$ in \eqref{eq:pc}. In section \ref{sec:numexp} we will refer to this method as PM1V. We will define this variant for the RKC method because the recursive formulas are simpler. In order to apply the variant to ROCK2 one just replaces the formulas. 

As we said each stage $g_i$ is projected, we call $g_i^*$ the stage before being projected, $\phi_i$ the solution of
\begin{numcases}{}
\nabla^2 \phi_i = \nabla \cdot g^*_i  & in $\Omega$,\notag \\
\dn{\phi_i}=0  & on $\partial \Omega$, \notag
\end{numcases}
and $g_i=g_i^*-\nabla \phi_i$ the projected stage. Let
\begin{align*}
f(\bu) = -(\bu\cdot \nabla)\bu-\nabla p_n +\nu\nabla^2\bu,
\end{align*}
steps 1, 2, 3 in \ref{sec:defpm1} are replaced by the following algorithm:
\begin{subequations}\label{eq:pm1v_algo}
\begin{flalign}
\text{1. } g_0 &=\bun ,\\
\text{2. } g_1^* &= g_0+\kappa_1\Delta t f(g_0) ,& \\
\text{3. } g_1 &= g_1^*-\nabla \phi_1 ,&
\end{flalign}
for $j=2,\ldots ,s$ do
\begin{flalign}
\text{4. } g^*_j &= g_0+\mu_j\left(g_{j-1}-g_0\right)+\nu_j\left(g_{j-2}-g_0\right)+\kappa_j \Delta t \left(f(g_{j-1})-a_{j-1}f(g_0)\right), \\
\text{5. } g_j &= g_j^*-\nabla \phi_j ,&
\end{flalign}
and finally set
\begin{flalign}
\text{6. } \bunpu &=g_{s} \text{ and } \pnpu=\pn+\frac{2}{\Delta t}\phi_{s}. &
\end{flalign}
\end{subequations}

As you will see in chapter \ref{sec:erkinse} this variant is closer to a Runge-Kutta method for differential algebraic equations than to a projection method. Nonetheless it exists a subtle difference which will be explained in section \ref{sec:rec_dae}.

\subsection{Method 2} \label{sec:method2}
The following method (PM2) has been presented in \cite{classfully} showing good convergence results up to the boundary. In \cite{classfully} the advection term is discretized explicitly at the beginning of the time step, therefore it can be put into a forcing term. It follows that the analysis can be done on the Stokes equations
\begin{subequations} \label{eq:s}
\begin{numcases}{}
\dt{\bu} +\nabla p = \nu \nabla^2 \bu
& in $\Omega$,\label{eq:s1} \\
\nabla \cdot \bu = 0 & in $\Omega$, \label{eq:s2}
\end{numcases}
\end{subequations}
without loss of generality.
However in the case of RKC and ROCK2 the diffusion and advection are advanced together in time, hence the latter cannot be embedded into the forcing term. And in PIROCK neither since we advance the solution in time using the diffusion and then we use these results in order to advance the advection. Thereby the following scheme when applied to RKC, ROCK2 or PIROCK holds only for the Stokes equations. In our numerical experiments the results of this method are very similar to the ones of PM1. In section \ref{sec:numexp} we will present only the results obtained with PM1.

\subsubsection{Definition of the method}
The method is very similar to the previous one, the main difference is the consistent pressure update. Consider the equation 
\begin{align}\label{eq:sm1}
\dt{\bus} +\nabla q = \nu \nabla^2 \bus,
\end{align}
where $q$ is an approximation of the pressure $p$ which as to be defined. As in section \ref{sec:method1} equation \eqref{eq:sm1} is advanced in time from $t_n$ to $t_{n+1}$ and a virtual velocity $\busnpu$ is computed. Afterward the divergence free velocity $\bunpu$ is given as in section \ref{sec:method1} by the Poisson equation for $\phi$. Plugging the velocity $\bunpu = \busnpu-\nabla \phi$ into equation \eqref{eq:sm1} gives
\begin{align} \label{eq:sm2}
\dt{\bunpu} +\nabla \left( q + \dt{\phi}-\nu\nabla^2 \phi \right) = \nu \nabla^2 \bunpu
\end{align}
and comparing equation \eqref{eq:s1} with equation \eqref{eq:sm2} it is clear that $\pnpu$ has to be computed by
\begin{align}
\pnpu=q + \dt{\phi}-\nu\nabla^2 \phi .
\end{align}
In \cite{classfully} it has been shown that setting 
\begin{align}
q(t_n+\epsilon)=p_n+\epsilon \dt{p}(t_n)
\end{align}
second order accuracy is achieved for velocity and pressure. Observe that the boundary conditions are not satisfied in the tangential component of the velocity.

\subsection{Errors at the boundary}\label{sec:method3}
As we already said in section \ref{sec:method1} and \ref{sec:method2} the projection methods PM1 and PM2 give a solution which do not satisfy the boundary conditions exactly. We will see in section \ref{sec:forced} that this issue do not affects the convergence order of the methods but on the other hand we will also see that the errors are coming mostly from the boundary. Inspired by the numerical experiments done with PM1 we propose the following fix (PM3) which needs the derivatives of the exact solution at the boundary. Since it requires an additional information about the exact solution it can be used only to prove, numerically, that a better choice for the virtual velocity boundary conditions could decrease the error.

In sections \ref{sec:method1} and \ref{sec:method2} we saw that the tangential component of the velocity do not satisfies the boundary conditions since we set $\bunpu=\busnpu-~\nabla \phi$ where $\frac{\partial \phi}{\partial \bold{\tau}}\neq 0$ and $\left. \busnpu \right|_\Gamma=\left.\bunpu^b\right|_\Gamma$. The fix is based on the fact that if $\dn{\phi}=0$ it follows that
\begin{align}
\dn{}\frac{\partial \phi}{\partial \bold{\tau}}=\frac{\partial}{\partial \mathbf{\tau}}\dn{\phi}=0,
\end{align}
which implies 
\begin{align}\label{eq:pm3pb}
0 &=\dn{}\left( \bold{\tau}\cdot\nabla \phi\right) = \dn{}\left( \bold{\tau} \cdot\left( \busnpu-\bunpu \right) \right).
\end{align}
Equation \eqref{eq:pm3pb} suggests the following boundary condition for the virtual velocity
\begin{align} \label{eq:pm3b}
\dn{}\left( \bold{\tau}\cdot \busnpu \right)=\dn{}\left( \bold{\tau}\cdot\bunpu^b  \right) \quad \mbox{on } \partial \Omega.
\end{align}

One could approximate $\dn{}\left( \bold{\tau}\cdot\bunpu^b \right)$ by a finite difference formula using the inner points of the mesh but a problem arises when it is embedded in the spatial discretization of the diffusion. At this point the truncation errors become non negligible since they are divided by $\Delta x^2$, where $\Delta x$ is the grid size. We tried to make an approximation with a third order finite difference scheme but the boundary conditions are no more satisfied because of truncation errors. Moreover the errors propagate inside the domain as the solution advance in time.

Nevertheless we can use this fix to show that the boundary conditions for $\bus$ impacts the solution accuracy. Numerical experiments with PM1 and PM3 will be discussed in section \ref{sec:numexp}.

\chapter{Stabilized explicit Runge-Kutta methods for differential algebraic equations with application to the Navier-Stokes equations}\label{sec:erkinse}
Due to the incompressibility constraint \eqref{eq:ns2} the Navier-Stokes equations \eqref{eq:ns} are not ordinary differential equations (ODE) and one cannot simply apply a Runge-Kutta method for ODEs. Instead they are part of a class of equations called differential algebraic equations (DAE) of index 2. In this chapter we will consider Runge-Kutta methods for DAEs of index 2 with focus on the Navier-Stokes equations \cite{acc_an}. We will see that for general DAEs of index 2 new order conditions arise for the differential and algebraic variables (velocity and pressure) but for the specific case of Navier-Stokes most of them vanish. Moreover one can use workarounds in order to obtain high accuracy without satisfying the new order conditions.

In the first section of this chapter we follow \cite[VII.6]{ode-ii} giving a short introduction to differential algebraic equations and we show that the semi-discrete form of the Navier-Stokes equations has differentiation index 2. In the second section we define Runge-Kutta methods for DAEs of index 2 and apply them to the Navier-Stokes equations as it is done in \cite{acc_an}. At the end we will apply the methods defined in the second section to the stabilized Runge-Kutta methods RKC and ROCK2.

\section[Differential algebraic equations and Navier-Stokes differentiation index]{Differential algebraic equations and the Navier-Stokes differentiation index}
In this section we will first give an introduction to differential algebraic equations and then prove that the semi-discrete form of the incompressible Navier-Stokes equations has differentiation index 2.
 
\subsection{Introduction to differential algebraic equations}\label{sec:intro_dae}
A differential algebraic equation is a system of differential equations 
\begin{align} \label{eq:dae}
F(u',u,t)=0
\end{align}
where $F:\Rd\times\Rd\times\Rb\rightarrow \Rd$ is supposed to be sufficiently differentiable and $u:[a,b]\rightarrow \Rd$. The difference from an ordinary differential equation is that usually the system is not solvable for all the components of $u'$.

If the system can be rewritten in the form
\begin{subequations} \label{eq:sedae}
\begin{align} 
y'&=f(y,z,t) \label{eq:f}\\
0 &= g(y,z,t). \label{eq:g}
\end{align}
\end{subequations}
it is called a semi-explicit DAE. The $y$ variable is called \textit{differential}, $z$ is called \textit{algebraic}. In the following we will restrict ourselves to semi-explicit DAEs.

Assuming that $\dz g$ is invertible in a neighborhood of the solution we can transform the above DAE in an ODE by differentiating equation \eqref{eq:g} with respect to $t$. Indeed
\begin{align*}
0 &= \frac{d g}{dt}(y,z,t) = \dy{g}(y,z,t)f(y,z,t)+\dz{g}(y,z,t)z'+\dt{g}(y,z,t),
\end{align*}
which gives
\begin{align}\label{eq:daeind1}
y' &= f(y,z,t) \\
z' &= - \dz{g}(y,z,t)^{-1}\left(\dy{g}(y,z,t)f(y,z,t)+\dt{g}(y,z,t)\right) .
\end{align}
The following definition (taken from \cite[VII.1]{ode-ii}) generalizes this idea.
\begin{defi} \label{def:index}
Equation \eqref{eq:dae} has differentiation index $m$ if $m$ is the minimal number of analytical differentiations
\begin{equation} \label{eq:index}
F(u',u,t)=0,\quad \frac{d F(u',u,t)}{dt}=0,\quad \ldots ,\quad \frac{d^m F(u',u,t)}{dt^m}=0
\end{equation}
such that equations \eqref{eq:index} allow us to extract by algebraic manipulations an explicit ODE $u'=~\varphi(u)$ (which is called the \textit{underlying ODE}).
\end{defi}

Let us take a look at systems of index 1 and 2.

\noindent
\emph{Systems of Index 1.} From \eqref{eq:daeind1} it follows that in the case where $\dz{g}$ is invertible in the neighborhood of the solution the semi-explicit DAE \eqref{eq:sedae} is of index 1.

\noindent
\emph{Systems of Index 2.} Let us consider the system  
\begin{subequations} \label{eq:sedae2}
\begin{align} 
y'&=f(y,z,t) \label{eq:f2}\\
0 &= g(y,t). \label{eq:g2}
\end{align}
\end{subequations}
Differentiating one time the constraint \eqref{eq:g2} gives
\begin{align}\label{eq:hc}
0 &= \Dt{g}(y,t) = \dy{g}(y,t)f(y,z,t)+\dt{g}(y,t),
\end{align}
equation \eqref{eq:hc} is a hidden constraint of the system.
Differentiating \eqref{eq:g2} twice gives
\begin{align*}
0 = \frac{\mbox{d}^2 g}{\mbox{d}t^2}(y,t) &= \dyy{g}(y,t)\left( f(y,z,t),f(y,z,t)\right) +\frac{\partial^2 g}{\partial t\partial y}(y,t)f(y,z,t)\\
&\phantom{=}\phantom{.}+ \dy{g}(y,t)\left(\dy{f}(y,z,t)f(y,z,t) +\dz{f}(y,z,t)z'+\dt{f}(y,z,t) \right)+\dtt{g}(y,t),
\end{align*}
an expression for $z'$ is found under the assumption that $\dy{g}\dz{f}$ is invertible in a neighborhood of the solution. In this case the DAE system \eqref{eq:sedae2} has index 2.

\subsection[Semi-discrete Navier-Stokes equations differentiation index]{Semi-discrete Incompressible Navier-Stokes equations differentiation index}\label{sec:indinse}
The purpose of this section is to show that the spatially discretized incompressible Navier-Stokes equations have differentiation index 2. In the DAE context the velocity $\bu$ is the differential variable $y$ and the pressure $p$ is the algebraic variable $z$.

Let us discretize equations \eqref{eq:ns} in space, for this we follow \cite{acc_an}. Most of the spatial discretization techniques will give a semi-discrete problem which can be written
\begin{subequations}\label{eq:ddae}
\begin{align}
\Omega u'(t) & = -C(u(t)) + \nu D u(t) - G p(t) + \rd(u(t),t), \label{eq:ddae1}\\
M u(t) &= \ru(t), \label{eq:ddae2}
\end{align}
\end{subequations}
where $u(t)\in\Rb^{N_u}$, $p(t) \in \Rb^{N_p}$ are the unknowns vectors for the velocities and pressure respectively, in the remainder their $t$-dependence will be left out of notation. $M$, $C$, $D$ and $G$ represent the discrete divergence, advection, diffusion and gradient operators respectively, $\Omega$ is a diagonal invertible matrix. $\ru$ is a vector with boundary conditions for the continuity equation, $\rd$ contains the boundary conditions and forcing terms of the momentum equation. 

Equation \eqref{eq:ddae} is a DAE where the differential variable $y$ is $u$ and the algebraic variable $z$ is $p$. Identifying \eqref{eq:sedae2} with \eqref{eq:ddae} gives 
\begin{subequations}\label{eq:id}
\begin{align}
f(u,p,t) &= F(u,t) - G p , \label{eq:id1}\\
g(u,t) &= M u - \ru (t). \label{eq:id2} 
\end{align}
\end{subequations}
where $F(u,t) = -C(u)+\nu D u + \rd(t)$ and $\Omega^{-1}$ has been absorbed in the definition of $C$, $D$, $G$ and $\rd$. In section \ref{sec:intro_dae} we saw that if
\begin{align}
\du{g}(u,t)\frac{\partial f}{\partial p}(u,p,t)  = -M G
\end{align}
is invertible in a neighborhood of the solution then the DAE system has index 2. The matrix $L=M G$ is the discretized Laplacian operator. Let us derive the expression for $p'$. Deriving equation \eqref{eq:id2} with respect to time we get
\begin{align}\label{eq:dg}
0=\Dt{g}(u,t)=M u' - \ru'(t)
\end{align}
and inserting $u'=f(u,p,t)$ into \eqref{eq:dg} we obtain the hidden constraint for the pressure
\begin{align}\label{eq:lp}
L p = M F(u,t)-\ru'(t).
\end{align}
Differentiating \eqref{eq:lp} gives
\begin{align}\label{eq:p1}
L p' = M\Dt{F}(u,p,t)-\ru''(t).
\end{align}
Because of the Dirichlet boundary conditions used in the matrices $M$ and $G$ any solution $q\in \Rb^{N_p}$ of $L q = f$ will satisfy Neumann boundary conditions, where $f\in\Rb^{N_p}$. It follows that $L$ is singular since any solution $q$ plus a constant stills a solution. If one changes the last row of $L$ with ones the matrix becomes regular since the mean value of $q$ is imposed. Using this regularized matrix one can solve $\eqref{eq:p1}$ and find an ODE for the algebraic variable $p$. Observe that in the case of a time varying mesh in equation \eqref{eq:lp} appears an additional term depending on $M'$, nonetheless this does not change the following result.

With two differentiations with respect to time we found an expression for $p'$, hence the differentiation index of the incompressible Navier-Stokes equations is 2.

\section{Runge-Kutta methods for differential algebraic equations of index 2} \label{sec:rkdae2}
We saw that given a semi-explicit differential algebraic equation of index 2 one can compute the underlying ODE $u'=\varphi(u)$ (see definition \ref{def:index}). Hence we could solve the system \eqref{eq:sedae2} using a common ODE solver applied to $u'=\varphi(u)$. This is the bad approach. The process used to obtain the underlying ODE is called index reduction by differentiation. Differentiating the constraint \eqref{eq:g2} effectively lowers the index of the system but upon discretization lower index systems satisfy only the derived constraints as $\Dt{g}(y,t)=0$ and $\Dtt{g}(y,t)=0$, not the original constraint $g(y,t)=0$. This phenomenon produced by the local error is called \textit{drift-off}, see \cite[VII.2]{ode-ii} for more details. In the following we will consider Runge-Kutta methods that satisfy the constraint $g(y,t)=0$, the solutions will satisfy all the derived lower index systems as well.

We showed in section \ref{sec:indinse} that the incompressible Navier-Stokes equations are DAE of index 2. In this section we will firstly define the Runge-Kutta methods for DAEs of index 2 and apply them to the Navier-Stokes equations. Then we will talk about the new order conditions arising from the algebraic constraint and show how they can be circumvented. At the end a short remark about time step adaptivity is given.
\subsection{Definition of the method}\label{sec:defrkdae}
An explicit Runge-Kutta method for a semi-explicit index 2 system \eqref{eq:sedae2} is defined by
\begin{subequations}\label{eq:rksedae}
\begin{align}
Y_i &= y_n + \Delta t \sum_{j=1}^{i-1} a_{ij}f(Y_j,Z_j,t_j) &\quad i=1,\ldots s, \\
0 &= g(Y_i,t_i),  \qquad &i=1,\ldots s, \\
y_{n+1}&=y_n + \Delta t \sum_{i=1}^s b_i f(Y_i,Z_i,t_i), &\\
0 &= g(y_{n+1},t_{n+1}). &
\end{align}
\end{subequations}
where $t_i=t_n+\Delta t c_i$. In order to apply method \eqref{eq:rksedae} to the Navier-Stokes equations \eqref{eq:ns} we consider their space discretization, hence the functions $f$ and $g$  of \eqref{eq:rksedae} are given by \eqref{eq:id}. Applying method \eqref{eq:rksedae} to equations \eqref{eq:id} gives
\begin{subequations}\label{eq:rkinse}
\begin{align}
U_i &= u_n + \Delta t \sum_{j=1}^{i-1} a_{ij}(F_j-G P_j)   &\quad i=1,\ldots s, \label{eq:rkinse1}\\
M U_i &= \ru(t_i),  \qquad &i=1,\ldots  s, \label{eq:rkinse2}\\
u_{n+1}&=u_n + \Delta t \sum_{i=1}^s b_i (F_i-G P_i), &\\
M u_{n+1} &= \ru(t_{n+1}). &
\end{align}
\end{subequations}
where $P_i$, $U_i$ and $u_n$ are approximations to $p(t_i)$, $u(t_i)$ and $u(t_n)$ respectively and $F_j=~F(U_j,t_j)$. 
In order to simplify the notation we set $a_{s+1,j}=b_j$ and $c_{s+1}=1$, then method \eqref{eq:rkinse} can be rewritten as
\begin{subequations}\label{eq:rkinse3}
\begin{align}
U_i &= u_n + \Delta t \sum_{j=1}^{i-1} a_{ij}(F_j-G P_j)   &\quad i=1,\ldots s+1, \label{eq:rkinse4}\\
M U_i &= \ru(t_i),  \qquad &i=1,\ldots  s+1, \\
u_{n+1}&=U_{s+1}.
\end{align}
\end{subequations}

For solving this system we eliminate the pressure from equation \eqref{eq:rkinse4}, this is done using the hidden constraint \eqref{eq:lp}. Replacing the pressure $p$ from \eqref{eq:lp} into \eqref{eq:rkinse4} gives
\begin{align}\label{eq:rkddae3}
U_i = u_n + \Delta t \sum_{j=1}^{i-1} a_{ij} F_j - \Delta t G L^{-1} \sum_{j=1}^{i-1}a_{ij} \left( M F_j -\ru'(t_j)\right).
\end{align}
If the internal stage order of the Runge-Kutta method is $q$ we have that
\begin{align}\label{eq:eint}
\Delta t \sum_{j=1}^{i-1} a_{ij}\ru'(t_j) = \int_{t_n}^{t_i}\ru'(t) dt + \bigo{\Delta t^{q+1}} = \ru(t_i)-\ru(t_n) + \bigo{\Delta t^{q+1}}
\end{align}
which suggests to insert the exact integral \eqref{eq:eint} into \eqref{eq:rkddae3}. Finally \eqref{eq:rkddae3} becomes
\begin{align}\label{eq:minse1}
U_i &= u_n + \Delta t \sum_{j=1}^{i-1} a_{ij} (I-G L^{-1}M)F_j + G L^{-1}(\ru(t_i)-\ru(t_n)) \qquad i=1,\ldots ,s+1.
\end{align}
We have $M(I-G L^{-1} M)=0$ (even if a row of $L$ has been modified, see \cite{rkcp}), thus supposing that $u_n$ is consistent, i.e. $M u_n = \ru(t_n)$, it follows $M U_i = \ru(t_i)$ for all the intermediate stages including $u_{n+1}$. In general this is not true for \eqref{eq:rkddae3}.

Setting a pressure-like variable $\phi_i$ for $i=2,\ldots ,s+1$ as the solution of
\begin{align}\label{eq:defphii}
L \phi_i = \frac{1}{c_i}\sum_{j=1}^{i-1} a_{ij} M F_j - \frac{\ru(t_i)-\ru(t_n)}{c_i \Delta t}
\end{align}
and $\phi_1=0$ we can write
\begin{align} \label{eq:phi1}
U_i &= u_n + \Delta t \sum_{j=1}^{i-1} a_{ij} F_j - c_i \Delta t G\phi_i  \qquad i=1,\ldots ,s+1,
\end{align}
in such a way $\phi_{s+1}$ is a first order approximation of the pressure. Indeed, let $\overline U_i$ be the stages of \eqref{eq:rkddae3} but with $u_n$ replaced by $u(t_n)$, since $\norm{u_n-u(t_n)}=\bigo{\Delta t^p}$ with $p\geq 1$ we have $\overline U_i -U_i=\bigo{\Delta t^p}$ and $\overline U_i-u(t_i)=\bigo{\Delta t^{q+1}}$ where $q$ is the stages order. It follows
\begin{align}
\begin{aligned}
L \phi_{s+1} &= \sum_{i=1}^{s} b_i M F_i - \frac{\ru(t_{n+1})-\ru(t_n)}{\Delta t} \\
&= \sum_{i=1}^{s} b_i M \overline F_i - \frac{\ru(t_{n+1})-\ru(t_n)}{\Delta t} + \bigo{\Delta t^p}\\
&= \sum_{i=1}^{s} b_i M F(u(t_i),t_i) - \frac{\ru(t_{n+1})-\ru(t_n)}{\Delta t} + \bigo{\Delta t^{\min(q+1,p)}}\\
&= M F(u(t_{n+1})) - \ru'(t_{n+1}) + \bigo{\Delta t} \\
&= L p(t_{n+1})+ \bigo{\Delta t}.
\end{aligned}
\end{align}
Let us set the vectors $\phi=(\phi_2,\ldots ,\phi_{s+1})$ and $\psi=(\psi_2,\ldots ,\psi_{s+1})$ defined as
\begin{align}\label{eq:lincomb}
\psi = \tilde A^{-1} \text{diag}(c_2,\ldots ,c_{s+1}) \phi,
\end{align}
where $\tilde A_{i j}=a_{i+1,j}$. Comparing equations \eqref{eq:rkinse4} and \eqref{eq:phi1} we see that $\psi_i=P_i$ for $i=2,\ldots ,s$. Setting $p_{n+1}=\psi_{s+1}$ gives a second or 
higher order approximation of the pressure if $\tilde A$ satisfies the order conditions explained in section \ref{sec:daeoc}. To conclude we write down the algorithm used in practice: for $i=1,\ldots ,s+1$ do
\begin{subequations}\label{eq:rkdaealgo}
\begin{flalign}
\mbox{1. } &U^*_i = u_n + \Delta t \sum_{j=1}^{i-1}a_{ij} F_j,&\\
\mbox{2. } &L \phi_i = \frac{1}{c_i \Delta t}\left(M U_i^*-\ru(t_i)\right),\label{eq:rkdaealgo2}\\
\mbox{3. } &U_i = U^*_i-c_i\Delta t G \phi_i,
\end{flalign}
\end{subequations}
and set $u_{n+1}=U_{s+1}$. Observe that at each stage the solution of a Poisson problem has to be computed in \eqref{eq:rkdaealgo2}.

\subsection{Order conditions}\label{sec:daeoc}
The classical order conditions for ODEs are not enough to guarantee the expected order of accuracy of both the velocity and the pressure. This is because here $f$ depends also on $p$ and there is the algebraic constraint, clearly these two facts were not taken into account when developing order conditions for ODEs. The following theorem taken from \cite{bra-hair} gives the order of convergence of the global error in function of the local error, allowing us to concentrate on local error only.
\begin{thm}
Supposing that the initial values $u(t_0)$, $p(t_0)$ are consistent, $a_{i+1,i}\neq 0$, $b_s\neq 0$ and that the local error $\delta u(t)$ satisfies
\begin{align}
\delta u(t) = \bigo{\Delta t^{p+1}}
\end{align}
the method is convergent of order $p$, i.e. 
\begin{align}
u_n-u(t_n) = \bigo{\Delta t^p} \text{ for } t_n-t_0 = n \Delta t \leq T
\end{align}
with $T$ finite.
\end{thm}
The vector $\phi$ used in equation \eqref{eq:lincomb} depends only on $u_n$ and not on $p_n$, thus $p_{n+1}=\psi_{s+1}$ depends only on $u_n$ as well. Consequently the global error is $p_n-p(t_n)=\bigo{\Delta t^q}$ provided that the local error is $\delta p(t)=\bigo{\Delta t^q}$ and $u_n-u(t_n)=\bigo{\Delta t^q}$. This is because the pressure has an instantaneous character, its value is such that the velocity is divergence free and it is independent of the pressure at previous time. In the following we will concentrate on the order conditions for the local error. 

As for the ODEs the order conditions are found by Taylor expansion of the exact and the numerical solution and comparing the coefficients of the differentials. For index 2 DAEs one additional condition for the velocity and two conditions for the pressure appear for order two. For order three there is four new conditions for the velocity and four for the pressure. As you can see the number of new conditions grows rapidly. Hopefully in the case of semi-explicit index 2 DAE some differentials vanish making some conditions trivially satisfied. For the even more specific case of the incompressible Navier-Stokes equations there is no additional order conditions for the velocity at least up to order five (see \cite{acc_an}). For the pressure all the conditions remain, for order two they are
\begin{align}
\sum_{ijk=1}^s b_i \omega_{ij}\omega_{jk}c_{k+1}^2 = 2 \quad \text{ and }\quad  c_s=1,
\end{align}
where $\omega_{ij}$ are the coefficients of $\tilde A^{-1}$. 

\subsection{Circumventing order conditions} \label{sec:circordcond}
The method presented in \ref{sec:defrkdae} needs to satisfy a considerable number of order conditions for the pressure so that the desired accuracy is achieved. If a Runge-Kutta method is built from scratch especially for the Navier-Stokes equations these conditions are taken into account and the method has the desired order. But if one wants to apply an existing Runge-Kutta method to the Navier-Stokes equations probably it has to modify the coefficients $a_{ij}$ so that the conditions for the pressure are satisfied. Moreover it is very likely that at least one stage must be added, therefore one more Poisson problem has to be solved. For these reasons in the following two sections we will show two methods taken from \cite{acc_an} which give higher order pressures without satisfying the additional order conditions. In the second method the stage order of the methods limits the accuracy of the pressure, in the last section we show how this can be avoided. These methods differ from the one in \ref{sec:defrkdae} only in the way the pressure is computed, the velocity is exactly the same and is given by algorithm \eqref{eq:rkdaealgo}.

\subsubsection{Approach 1}
This method (AP1) integrates the velocity as described in \eqref{eq:rkdaealgo} without taking care of order conditions. If the Runge-Kutta method has order $p$ then the velocity will converge with order $p$ since no additional order conditions appear. Given the $p$th order accurate velocity $u_{n+1}$ a $p$th order accurate pressure is computed solving 
\begin{align} \label{eq:ap1p}
L p_{n+1} = M F_{n+1} - \ru'(t_{n+1}).
\end{align}
Indeed, using \eqref{eq:lp}
\begin{align*}
L p_{n+1} &=  M F_{n+1} - \ru'(t_{n+1}) =  M F(u(t_{n+1}),t_{n+1}) - \ru'(t_{n+1}) + \bigo{\Delta t^p} \\&= L p(t_{n+1})+\bigo{\Delta t^p}. 
\end{align*}
This approach is very simple but on the other hand one more Poisson problem must be solved. Moreover it requires that $\ru$ can be differentiated, something that is not required in the computation of $\bu$ and $\phi$. In some practical computations, for example involving a prescribed turbulent inflow, $\ru'$ might not be available. The approach presented in the next section does not need the derivative of $\ru$.

In section \ref{sec:numexp} we will refer to this method with AP1.

\subsubsection{Approach 2}
In section \ref{sec:defrkdae} we saw that, under the assumptions that order conditions are satisfied, a higher order pressure is given by a linear combination of the pressure-like variables $\phi_i$ (see \eqref{eq:lincomb}). Here we will give a different linear combination of the pressure-like variables which gives a $r$th order pressure assuming that the Runge-Kutta method has at least $r$ stages and the order of these is at least $r$. This approach (AP2) is more involved than Approach 1 but it does not require any additional Poisson problem nor the availability of $\ru'$.

Let us denote the average of $p(t)$ in $[t_n,t]$ by $\ophi(t)$, we have
\begin{align}
\ophi(t) = \frac{1}{t-t_n}\int_{t_n}^{t} p(s) ds.
\end{align}
For notation purposes we set $\ophi_1=0$ and $\ophi_i:=\ophi(t_i)$ for $i=2,\ldots ,s+1$, yielding 
\begin{align}\label{eq:ophi}
\ophi_i = \frac{1}{c_i \Delta t}\int_{t_n}^{t_i} p(t) dt .
\end{align}
The goal is to find an accurate point value $p_{n+1}$ from the average values $\ophi_i$, this process is called \textit{reconstruction}.

Let $P(t)$ be the primitive function of $p(t)$. We construct $H(t)$ the polynomial interpolating $P(t)$ at the points $t_{k_1},\ldots ,t_{k_m}$, where the $t_{k_j}$ are distinct, $k_1=1$ and $k_m=s+1$. Then the derivative $h(t)$ of $H(t)$ has the same integral of $p(t)$, indeed
\begin{align}
c_{k_j} \Delta t \ophi_{k_j} = \int_{t_n}^{t_{k_j}} p(t) dt = P(t_{k_j})-P(t_n) = H(t_{k_j})-H(t_n) = \int_{t_n}^{t_{k_j}} h(t) dt.
\end{align}
Using the Lagrange basis 
\begin{align}
l_j(t)=\prod_{i=1,i\neq j}^m \frac{t-t_{k_i}}{t_{k_j}-t_{k_i}}
\end{align}
we have $H(t)=\sum_{j=1}^m P(t_{k_j})l_j(t)$ and $\sum_{j=1}^m l_j(t)=1$ for all $t \in \Rb$. It follows 
\begin{align}
H(t)-P(t_n) = \sum_{j=1}^m (P(t_{k_j})-P(t_n)) l_j(t)=\sum_{j=2}^m c_{k_j} \Delta t \ophi_{k_j} l_j(t).
\end{align}
and
\begin{align}
h(t) = \sum_{j=2}^m c_{k_j} \Delta t \ophi_{k_j} l_j'(t).
\end{align}
Since $H(t)$ interpolates $P(t)$ in $m$ points it follows from Rolle's theorem that $h(t)$ interpolates $p(t)$ in $m-1$ points. Thus $h(t_{n+1})$ is a $m-1$ order approximation of $p(t_{n+1})$ (see \cite[Ch.1]{int-pic}). We set
\begin{align}
p_{n+1}:=h(t_{n+1}) = p(t_{n+1}) + \bigo{\Delta t^{m-1}}.
\end{align}

However the average values $\ophi_i$ are not known in practice and must be approximated. For a $p$th order method with $q$th order internal stages we have, using \eqref{eq:lp} and $\overline U_i$ as in section \ref{sec:defrkdae},
\begin{align} \label{eq:apphi}
\begin{aligned}
L\ophi_i &= \frac{1}{c_i \Delta t}\int_{t_n}^{t_i}L p(t) dt \\
&= \frac{1}{c_i \Delta t}M\int_{t_n}^{t_i} F(u(t),t) dt -\frac{\ru(t_i)-\ru(t_n)}{c_i \Delta t} \\
&= \frac{1}{c_i} \sum_{j=1}^{i-1} a_{ij} M F(u(t_j),t_j)  -\frac{\ru(t_i)-\ru(t_n)}{c_i \Delta t} + \bigo{\Delta t^q} \\
&= \frac{1}{c_i} \sum_{j=1}^{i-1} a_{ij} M F(\overline U_j,t_j)  -\frac{\ru(t_i)-\ru(t_n)}{c_i \Delta t} + \bigo{\Delta t^q}\\
&= \frac{1}{c_i} \sum_{j=1}^{i-1} a_{ij} M F(U_j,t_j)  -\frac{\ru(t_i)-\ru(t_n)}{c_i \Delta t} + \bigo{\Delta t^q}+\bigo{\Delta t^p} \\
&= L\phi_i + \bigo{\Delta t^{\min(p,q)}} =  L\phi_i + \bigo{\Delta t^q}
\end{aligned}
\end{align}
Observe that for $i=s+1$ it holds $p=q$ and if $i\leq s$ then $q\leq p$. Let define $\tilde h(t)$ the approximate counterpart of $h(t)$ as $\tilde h(t) = \sum_{j=2}^m c_{k_j} \Delta t \phi_{k_j} l_j'(t)$. Using \eqref{eq:apphi} we have
\begin{align*}
\tilde h(t) &= \sum_{j=2}^m c_{k_j} \Delta t \phi_{k_j} l_j'(t) \\
&= \sum_{j=2}^m c_{k_j} \Delta t \ophi_{k_j} l_j'(t) + \sum_{j=2}^m \bigo{\Delta t^{q+1}} l_j'(t) \\
&= h(t) + \bigo{\Delta t^q} 
\end{align*}
since $l_j'(t)=\bigo{\Delta t^{-1}}$. Setting $p_{n+1}=\tilde h(t_{n+1})$ we obtain
\begin{align*}
p_{n+1}:=\tilde h(t_{n+1}) = p(t_{n+1}) + \bigo{\Delta t^{\min(m-1,q)}}.
\end{align*}
Finally, to obtain a $r$th order accurate pressure the Runge-Kutta method must have $r$ stages and they must have $r$th order accuracy at least.

In section \ref{sec:numexp} we will refer to this method with AP2.

\subsubsection{Adapting Approach 2 to order one internal stages}\label{sec:ap2ord1}
In the previous section we saw that one can obtain a $r$th order approximation of the pressure by using $r$th order approximations of the $\ophi_i$s. It follows that if only first order approximations of the $\ophi_i$s are available one cannot obtain a second order accurate pressure. In this section we will show how to obtain a second order approximation to $\ophi_j$ using three first order approximations of $\ophi_i$, $\ophi_j$, $\ophi_k$ with $i$, $j$, $k$ distinct. Consequently a second order pressure can be obtained even if only first order approximations to the $\ophi_i$s are available.

Even if in the following we restrict ourselves to $r=2$ with some effort this procedure can be generalized to any $r$.

Using \eqref{eq:lp} one can rewrite $f(u,p,t)$ in \eqref{eq:id1} as $f(u,t)=F(u,t)-G L^{-1} (M F(u,t)-\ru'(t))$. Let us consider equation \eqref{eq:minse1} and write explicitly his dependence on $\Delta t$, we have
\begin{align}
\begin{aligned}
U_i(\Delta t) &= u_n + \Delta t \sum_{j=1}^{i-1} a_{ij} (I-G L^{-1}M)F(U_j(\Delta t),t_n+c_j\Delta t) \\
&\phantom{=} {}+ G L^{-1}(\ru(t_n+c_i\Delta t)-\ru(t_n)),
\end{aligned}
\end{align}
from which it follows $U_i(0)=u_n$ and  
\begin{align}\label{eq:uiun}
\begin{aligned}
U_i(\Delta t)-u_n &= \Delta t U_i'(0) + \bigo{\Delta t^2} \\
&= \Delta t \sum_{j=1}^{i-1} a_{ij} (I-G L^{-1}M)F(u_n,t_n) + c_i\Delta t G L^{-1}\ru'(t_n) + \bigo{\Delta t^2}\\
&= c_i \Delta t \left( F(u_n,t_n) - G L^{-1}\left(M F(u_n,t_n) - \ru'(t_n)\right) \right)+ \bigo{\Delta t^2}\\
&= c_i \Delta t f(u_n,t_n) + \bigo{\Delta t^2}
\end{aligned}
\end{align}
Using \eqref{eq:uiun} we want to compute the error's leading term of $\ophi_i-\phi_i$. Remember that we have (see \eqref{eq:apphi} and \eqref{eq:defphii})
\begin{align}
L\ophi_i &= \frac{1}{c_i \Delta t}M\int_{t_n}^{t_i} F(u(t),t) dt -\frac{\ru(t_i)-\ru(t_n)}{c_i \Delta t}, \label{eq:rap1}\\
L \phi_i &= \frac{1}{c_i}\sum_{j=1}^{i-1} a_{ij} M F_j - \frac{\ru(t_i)-\ru(t_n)}{c_i \Delta t}.\label{eq:rap2}
\end{align}
Let us estimate $L\ophi_i-L\phi_i$. We call $\ou(t)$ the solution of the Navier-Stokes equations but with initial condition $\ou(t_n)=u_n$, it gives
\begin{align}
\int_{t_n}^{t_i} F(u(t),t) dt &= \int_{t_n}^{t_i} F(\ou(t),t) + \bigo{\Delta t^2} dt \notag \\
&= \int_{t_n}^{t_i} F(\ou(t_n),t_n) + \Dt{F}(\ou(t_n),t_n)(t-t_n) +\bigo{\Delta t^2} dt \notag \\
&= c_i \Delta t F(u_n,t_n) + \frac{1}{2}c_i^2\Delta t^2\Dt{F}(u_n,t_n) +\bigo{\Delta t^3}. \label{eq:es1}
\end{align}
And using \eqref{eq:uiun}
\begin{align}
\sum_{j=1}^{i-1} a_{ij} F_j &= \sum_{j=1}^{i-1} a_{ij}\left( F(u_n,t_n)+\du{F}(u_n,t_n)\left(U_j-u_n\right)+\dt{F}(u_n,t_n)(t_j-t_n)+\bigo{\Delta t^2} \right) \notag\\
&= \sum_{j=1}^{i-1} a_{ij}\left( F(u_n,t_n)+ c_j \Delta t\left(\du{F}(u_n,t_n) f(u_n,t_n)+\dt{F}(u_n,t_n)\right)+\bigo{\Delta t^2} \right) \notag\\
&= c_i F(u_n,t_n) + \Delta t\left(\sum_{j=1}^{i-1}a_{ij}c_j\right)\Dt{F}(u_n,t_n) +\bigo{\Delta t^2}.\label{eq:es2}
\end{align}
Using \eqref{eq:rap1},\eqref{eq:rap2},\eqref{eq:es1} and \eqref{eq:es2} we can estimate $\ophi_i-\phi_i$ as 
\begin{align}\label{eq:phiap}
\begin{aligned}
\ophi_i-\phi_i &= \Delta t\left(\frac{1}{2}c_i-\frac{1}{c_i}\sum_{j=1}^{i-1}a_{ij}c_j\right) L^{-1} M \Dt{F}(u_n,t_n) +\bigo{\Delta t^2} \\
&= \Delta t e_i E_n + \bigo{\Delta t^2},
\end{aligned}
\end{align}
for $i=2,\ldots ,s+1$, where 
\begin{align}
e_i &= \frac{1}{2}c_i-\frac{1}{c_i}\sum_{j=1}^{i-1}a_{ij}c_j,\\
E_n &=  L^{-1} M \Dt{F}(u_n,t_n).
\end{align}
Observe that if the internal stages are of order two, i.e. $\sum_{j=1}^{i-1}a_{ij}c_j=\frac{1}{2}c_i^2$, it follows that $\ophi_i-\phi_i=\bigo{\Delta t^2}$. 

Remember that if we have two second order approximations of $\ophi_i$ we can compute a second order approximation of the pressure. For $i=s+1$ we have $\ophi_{s+1}-\phi_{s+1}=\bigo{\Delta t^2}$ because of the second order conditions of the method. It rests to find one more second order approximation to one of the $\ophi_i$s, for this we will use \eqref{eq:phiap}.

Let $i,j,k\in \{2,\ldots ,s+1\}$ be distinct with $j\neq s+1 $, we have
\begin{align} \label{eq:tophi}
\begin{aligned}
\ophi_i &=\ophi(t_i) = \ophi(t_j) + (t_i-t_j)\ophi'(t_j) + \bigo{\Delta t^2},\\
\ophi_k &=\ophi(t_k) = \ophi(t_j) + (t_k-t_j)\ophi'(t_j) + \bigo{\Delta t^2}.
\end{aligned}
\end{align}
Let $\alpha,\beta,\gamma \in \Rb$, using \eqref{eq:tophi}
\begin{align*}
\alpha \ophi_i + \beta \ophi_j + \gamma \ophi_k &= (\alpha+\beta+\gamma)\ophi_j + (\alpha(c_i-c_j)+\gamma(c_k-c_j))\Delta t \ophi_j'(t_j) + \bigo{\Delta t^2}
\end{align*}
and using \eqref{eq:phiap}
\begin{align*}
\alpha \ophi_i + \beta \ophi_j + \gamma \ophi_k = \alpha \phi_i + \beta \phi_j + \gamma \phi_k + (\alpha e_i+\beta e_j +\gamma e_k)\Delta t E_n + \bigo{\Delta t^2}.
\end{align*}
If $\alpha,\beta,\gamma$ are chosen so that 
\begin{align*}
\left\{
\begin{aligned}
\alpha+\beta+\gamma\neq 0 ,\\
\alpha(c_i-c_j)+\gamma(c_k-c_j) = 0 ,\\
\alpha e_i+\beta e_j +\gamma e_k = 0 ,
\end{aligned}
\right.
\end{align*}
then
\begin{align}
\ophi_j = \frac{\alpha \phi_i + \beta \phi_j + \gamma \phi_k}{\alpha+\beta+\gamma}+\bigo{\Delta t^2}.
\end{align}
This is achieved setting
\begin{align*}
\alpha &= \frac{e_j}{c_j-c_i},\\
\beta &= \frac{e_i}{c_i-c_j}-\frac{e_k}{c_k-c_j},\\
\gamma &= \frac{e_j}{c_k-c_j}.
\end{align*}

In chapter \ref{sec:numexp} we will refer to this method with AP2W (to be read as Approach 2 workaround).

\subsection{Time step adaptivity} \label{sec:dae_tsa}
In section \ref{sec:projtsc} we saw that for the projection method the local error is overestimated because also its non divergence-free component is taken into account. Oppositely if a Runge-Kutta method for DAEs is used each stage is projected into the manifold of divergence-free fields. Then the local error computed with an embedded formula uses only divergence-free stages and a better estimation is given.

\section{Stabilized Runge-Kutta methods for the Navier-Stokes equations: differential algebraic approach}
Applying the Runge-Kutta method for differential algebraic equations (see section \ref{sec:defrkdae}) to RKC or ROCK2 requiring that the order conditions are satisfied is not straightforward. Firstly because the coefficients $a_{ij}$ are not explicitly given and, more important, they depend on the number of stages $s$. One can recursively compute them but the new order conditions depends on $w_{ij}$, which is $s$ dependent as well. Second, all the $a_{ij}$ coefficients for $i>j$ are non zero thus one stage must be added in order to gain degrees of freedom and satisfy the additional conditions. Finally, even if the conditions are satisfied one Poisson problem has been added (because of the new stage) and at this point Approach 1 is better since one does not have to modify the method.

Approach 1 can be applied successfully to RKC and ROCK2 without any difficulty when $\ru'$ is available. Approach 2 can be applied to RKC obtaining a second order pressure. This is not true for ROCK2 since it has internal stages of order one preventing higher order pressures, thus the workaround for order one internal stages explained in section \ref{sec:ap2ord1} has to be used.

In the following we will explain how a Runge-Kutta method for differential algebraic equations is applied to a method which is realized by means of recursive formulas like RKC and ROCK2. Then we show how Approach 1 and 2 are applied to RKC and ROCK2. Finally we observe that in RKC time step adaptivity cannot be enabled when the solution is projected at each stage as in algorithm \eqref{eq:rkdaealgo}.

\subsection{Recursive formulas and Runge-Kutta methods for differential algebraic equations of index 2} \label{sec:rec_dae}
In this section we will see how a Runge-Kutta method for differential algebraic equations given in \eqref{eq:rkdaealgo} can be applied to a stabilized Runge-Kutta method which uses recursive formulas instead of the coefficients $a_{ij}$, $b_i$.
Again we will explain the method for the RKC scheme for simplicity, for ROCK2 one just replaces the recursive formulas.

We write down the recursive formulas \eqref{eq:rkc_rec} of RKC when applied to Navier-Stokes without taking into account the pressure (as in \eqref{eq:rkdaealgo}):
\begin{align}\label{eq:rkc_rec_ns}
\begin{aligned}
U_1 &= u_n,\\
U_2 &= u_n+\kappa_1 \Delta t F(u_n), \\
U_i &= u_n+\mu_{i-1}\left( U_{i-1}-u_n \right) +\nu_{i-1}\left( U_{i-2}-u_n \right)+\kappa_{i-1}\Delta t \left(F(U_{i-1})-a_{i-2} F(u_n)\right).
\end{aligned}
\end{align}
for $i=3,\ldots ,s+1$. Let $a_{ij}$ be the coefficients of RKC, if the stages $U_i$ are not projected formulas \eqref{eq:rkc_rec_ns} are equivalent to
\begin{align}\label{eq:rkcaij}
U_i = u_n + \Delta t \sum_{j=1}^{i-1} a_{ij} F(U_j),
\end{align}
for $i=3$ for example we have
\begin{align*}
U_3 &=u_n + \mu_2\kappa_1\Delta t F(u_n)+\kappa_2\Delta t (F(U_2)-a_1 F(u_n)) \\ 
&= u_n + \Delta t \left( (\mu_2\kappa_1-\kappa_2 a_1)F(U_1) + \kappa_2 F(U_2)\right),
\end{align*}
i.e. $a_{3,1}=\mu_2\kappa_1-\kappa_2 a_1$ and $a_{3,2}=\kappa_2$. Consider now the method given in \eqref{eq:rkdaealgo} where the following predictor-corrector procedure is used:
\begin{subequations}\label{eq:predcorr}
\begin{align}
U^*_i &= u_n + \Delta t \sum_{j=1}^{i-1}a_{ij} F_j \label{eq:predcorr1} \\
U_i &= U^*_i-c_i\Delta t G \phi_i,
\end{align} 
\end{subequations}
where $\phi_i$ is given in \eqref{eq:rkdaealgo2}. A naïve application to RKC would be
\begin{subequations}\label{eq:rkc_rec_proj}
\begin{flalign}
\text{1. } U_1^{\phantom{*}} &= u_n,\\
\text{2. } U_2^* &= u_n+\kappa_1 \Delta t F(u_n),& \\
\text{3. } U_2^{\phantom{*}} &= U_2^*-c_2 \Delta t G \phi_2 &
\end{flalign}
and for $i=3,\ldots ,s+1$ do
\begin{flalign}
\text{4. } U_i^* &= u_n+\mu_{i-1}\left( U_{i-1}-u_n \right) +\nu_{i-1}\left( U_{i-2}-u_n \right)+\kappa_{i-1}\Delta t \left(F(U_{i-1})-a_{i-2} F(u_n)\right) &\label{eq:rkc_rec_proj4}\\
\text{5. } U_i^{\phantom{*}} &= U_i^* -c_i \Delta t G \phi_i, &
\end{flalign}
\end{subequations}
but \eqref{eq:rkc_rec_proj4} is not equivalent to \eqref{eq:predcorr1}, in fact for $i=3$ we obtain
\begin{align*}
U_3^* &=u_n + \mu_2(U_2^*-c_2\Delta t G \phi_2-u_n) +\kappa_2\Delta t (F(U_2)-a_1 F(u_n)) \\ 
&=u_n + \mu_2(\kappa_1\Delta t F(u_n)-c_2\Delta t G \phi_2) +\kappa_2\Delta t (F(U_2)-a_1 F(U_1)) \\ 
&= u_n + \Delta t \left( (\mu_2\kappa_1-\kappa_2 a_1)F(U_1) + \kappa_2 F(U_2)\right)-\mu_2 c_2 \Delta t \phi_2 \\
&= u_n + \Delta t \sum_{j=1}^{2}a_{3,j} F_j -\mu_2 c_2 \Delta t \phi_2 ,
\end{align*}
which is wrong since the last term $\mu_2 c_2 \Delta t \phi_2$ does not appear in \eqref{eq:predcorr1}. This additional term comes out because in \eqref{eq:rkc_rec_proj4} we used $U_i$ instead of $U_i^*$. The right realization of of RKC as a Runge-Kutta method for DAEs is given by the following algorithm:
\begin{subequations}\label{eq:rkc_rec_projf}
\begin{flalign}
\text{1. } U_1^{\phantom{*}} &= u_n,&\\
\text{2. } U_2^* &= u_n+\kappa_1 \Delta t F(u_n),& \\
\text{3. } U_2^{\phantom{*}} &= U_2^*-c_2 \Delta t G \phi_2 &
\end{flalign}
and for $i=3,\ldots ,s+1$ do
\begin{flalign}
\text{4. } U_i^* &= u_n+\mu_{i-1}\left( U_{i-1}^*-u_n \right) +\nu_{i-1}\left( U_{i-2}^*-u_n \right)+\kappa_{i-1}\Delta t \left(F(U_{i-1})-a_{i-2} F(u_n)\right) ,&\\
\text{5. }  U_i^{\phantom{*}} &= U_i^* -c_i \Delta t G \phi_i .&
\end{flalign}
\end{subequations}
Observe that the evaluation of $F$ is done on the projected stage $U_i$, otherwise the non projected stage $U_i^*$ is used. Formulation \eqref{eq:rkc_rec_projf} is equivalent to \eqref{eq:rkdaealgo}. 

Algorithm \eqref{eq:pm1v_algo} given in section \ref{sec:pm1var} is equivalent to algorithm \eqref{eq:rkc_rec_proj} but written with a different notation. The only difference is the different pressure update but the results do not change. We confirm this assertion by the numerical experiments of section \ref{sec:errb_ff}.

In section \ref{sec:numexp} we will see that algorithm \eqref{eq:pm1v_algo} (i.e. \eqref{eq:rkc_rec_proj}) behave very similarly to \eqref{eq:rkc_rec_projf}. But when we use Approach 2 of section \ref{sec:circordcond} (also its variant for order one internal stages) it is very important that the realization of the Runge-Kutta method with recursive formulas is consistent with its realization using the $a_{ij}$ coefficients. 


\subsection{Application to RKC and ROCK2}
In order to advance the solution in time one does not have to compute the pressure at each time step so approaches 1 and 2 are used only when one computes the pressure.

It is straightforward to apply Approach 1 to RKC and ROCK2 by simply using algorithm \eqref{eq:rkc_rec_projf} and solving the last Poisson problem \eqref{eq:ap1p}.

In order to apply Approach 2 to RKC we need at least two stages of order two, in our implementation we use $U_s$ and $U_{s+1}$. Since $U_1$ and $U_2$ have stage order less than two, when the second order pressure is computed the minimal number of stages used in \eqref{eq:rkc_rec_projf} is $s=3$, so that $U_1$, $U_2$, $U_s$ and $U_{s+1}$ are distinct. When the pressure is not computed $s=2$ is allowed.

ROCK2 has order one internal stages thus to apply Approach 2 we must use the method explained at the end of section \ref{sec:ap2ord1}. In order to do that we need to compute the coefficients $a_{ij},c_i,b_i$ of ROCK2. Since the method is explicit it follows that $a_{ij}=0$ if $i\leq j$. From \eqref{eq:defrock2} it is not difficult to find the other coefficients recursively, we get
\begin{align*}
a_{ij}&=
\begin{cases}
0 & \mbox{if } i\leq j ,\\
\mu_i & \mbox{if } i=j+1 \mbox{ and } i\leq s-1, \\
-\nu_{i-1}a_{i-1,j} -\kappa_{i-1}a_{i-2,j} & \mbox{if } j+2\leq i \leq s-1,\\
a_{s-1,j} & \mbox{if } i=s \mbox{ and } j\leq s-2, \\
\sigma & \mbox{if } i=s \mbox{ and } j=s-1,
\end{cases} \\
b_i &= 
\begin{cases}
a_{s,i} & \mbox{if } i\leq s-2, \\
2\sigma-\tau/\sigma & \mbox{if } i=s-1, \\
\tau/\sigma & \mbox{if } i=s,
\end{cases}
\end{align*}
and $c_i = \sum_{j=1}^s a_{ij}$. The workaround given in \ref{sec:ap2ord1} for order one internal stages needs three stages $U_i$, $U_j$, $U_k$ of order one in order to compute a second order approximation of $\ophi_j$ with $j\neq s+1$. In our implementation we use $U_2$, $U_3$, $U_4$ and we get a second order approximation of $\ophi_3$. Since ROCK2 is a second order method the last stage $U_{s+1}$ has order two, so we can use $\phi_{s+1}$ as second order approximation of $\ophi_{s+1}$. Hence we need the stages $U_2$, $U_3$, $U_4$ for $\ophi_3$ and $U_{s+1}$ for $\ophi_{s+1}$. The stages $U_4$ and $U_{s+1}$ can coincide meaning that the minimal number of stages is $s=3$, as in the original ROCK2 method.

When using Approach 2 in RKC the computation of the pressure can affect also the velocity since the minimal number of stages changes, this happens only when very small $\Delta t$s are used and it is not true for ROCK2.

\subsection{Time step adaptivity issue in RKC}\label{sec:tsairkc}
We saw in section \ref{sec:tscrr} that in RKC the estimation of the local error $err_{n+1}$ is done by an approximation of $y'''(t_n)$. This approximation is computed under the assumption that the equation being solved is an ODE, so there is no dependence on $p(t)$. In the projection method explained in section \ref{sec:method1} the pressure is kept constant during the integration, in such a way $y'''(t_n)$ does not depends on $p(t)$ and the approximation holds. On the other hand when the differential algebraic approach of section \ref{sec:defrkdae} is used the velocity $\bu$ is projected at each stage. This is equivalent to advance the pressure $p$ in time since its value is so that the velocity is divergence free. So, the approximation of $y'''(t_n)$ given in \eqref{eq:erkc1} is not valid anymore. One can give another approximation taking into account the dependence on $p(t)$ but this is not done here. This discussion is valid also for the variant of the projection method given in \ref{sec:pm1var}.

Observe that an embedded formula for time step adaptivity works also when the velocity is projected at each stage. This is because the theory of embedded formulas uses only the order conditions of the velocity and these are the same for ODEs and the Navier-Stokes equations, hence projecting the stages or not does not affects the correctness of the time step adaptivity procedure.

\chapter{Space discretization and Poisson solver}
In this chapter we will describe the method used for the spatial discretization and explain how the Poisson problem is solved in our implementation. In chapter \ref{sec:erkinse} we made use of matrices for the discretized differential operators in order to describe the theory, however these matrices are not explicitly built in the code. Instead the operators are implemented as functions which given the operand return the matrix multiplication. The boundary conditions and forcing terms are embedded in these functions as well. Consequently the Poisson problem is not solved by means of matrices neither.

For the spatial discretization we use a grid where the unknowns are staggered and in order to solve the Poisson problem we use a fast cosine transform well suited for this particular grid.

\section{The MAC method} \label{sec:mac}
The Marker and Cell (MAC) method has been proposed in 1965 in the Los Alamos laboratory \cite{mac0}. Recently, it has been shown in \cite{mac2} that the MAC method is very competitive and particularly well suited for high Reynolds numbers and free surface problems. 

In this chapter we will motivate the usage of the MAC method and apply it to the incompressible Navier-Stokes equations.

\subsection{Motivating the use of the MAC method's staggered grid}
The MAC method uses a staggered grid, in this arrangement the scalar variables (pressure, density,...) are stored in the cell centers and the momentum variables (velocity) are located at the cell faces.     This is different from a collocated grid, where all the variables are stored at the same position. A staggered grid avoids the odd-even decoupling of pressure and velocity (see equation \eqref{eq:oddeven}). Odd-even decoupling introduces a discretization error which leads to checkerboard patterns in the solution.

\begin{figure}[!htbp]
\begin{center}
\includegraphics[trim=0.0cm 0.0cm 0.1cm 3.7cm, clip=true, width=0.45\textwidth]{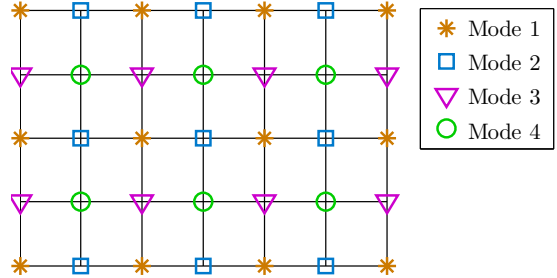}
\end{center}
\caption{A non staggered grid with four independent pressure modes.}
\label{fig:non_staggered}
\end{figure}

Let us motivate the previous considerations with an example. Suppose that $\Omega=[0,1]\times[0,1]$ is divided into a number of uniform cells as in Figure \ref{fig:non_staggered} and the equation to be solved is
\begin{align} \label{eq:ex}
\dt{\bu} + \nabla p = \bold{f} .
\end{align}
In a collocation grid the discretization of equation \eqref{eq:ex} with $\bu=~(u,v)^\top$ and $\bold f = (f^1,f^2)^\top$ leads to
\begin{align}\label{eq:oddeven}
\begin{aligned}
\dt{u_{i,j}} + \frac{p_{i+1,j}-p_{i-1,j}}{2\Delta x} = f^1_{i,j} \\
\dt{v_{i,j}} + \frac{p_{i,j+1}-p_{i,j-1}}{2\Delta x} = f^2_{i,j}
\end{aligned}
\end{align}
where $\Delta x$ is the grid size. With this choice of positioning four independent pressure modes arises, as shown in Figure \ref{fig:non_staggered}. Hence it is easy to build a function $\tilde{p}$ non constant but for which the discretized gradient vanish, making $\tilde{p}$ invisible to the momentum equation \eqref{eq:ex}. In the staggered grid the pressure unknowns are all coupled together. 

\subsection{Application of the MAC method to the Navier-Stokes equations}
\begin{figure}[!hbtp]
\begin{center}
\includegraphics[trim=0.0cm 0.3cm 0.0cm 3.0cm, clip=true, width=0.5\textwidth]{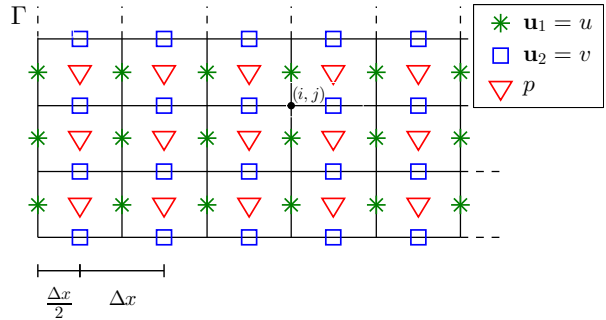}
\end{center}
\caption{Example of a staggered grid and variables positioning.}
\label{fig:staggered}
\end{figure}

In the presentation of the MAC method we will focus on the incompressible Navier-Stokes equations with Dirichlet boundary conditions. As in a collocation grid, in a staggered grid finite differences are employed but in this approach the velocity components $u$, $v$ and the pressure $p$ are not discretized at the same points of the domain. Figure \ref{fig:staggered} displays an example of a staggered grid and the position of the variables.

Let $\Omega=[0,1]\times[0,1]$, $N$ be a positive integer, $\Delta x = 1/N$ be the grid size and $x_{i j}=(i \Delta x, j \Delta x)$ for $i,j=0,1,\ldots ,N$. The component $u$ of the velocity is discretized at the points $x_{i, j-1/2}$, we call $u_{i,j-1/2}$ the numerical approximation of $u(x_{i, j-1/2})$:
\begin{align*}
u_{i,j-1/2} \approx u(x_{i, j-1/2}) \quad \text{for } i=1,\ldots ,N-1 \text{ and } j=1,\ldots ,N.
\end{align*}
Similarly for the component $v$ of the velocity and the pressure $p$ we have
\begin{align*}
v_{i-1/2,j} \approx v(x_{i-1/2,j}) \quad \text{for } i=1,\ldots ,N \text{ and } j=1,\ldots ,N-1
\end{align*}
and
\begin{align*}
p_{i-1/2,j-1/2} \approx p(x_{i-1/2, j-1/2}) \quad \text{for } i=1,\ldots ,N \text{ and } j=1,\ldots ,N.
\end{align*}
The first component of equation \eqref{eq:ns1} expressed in Cartesian coordinates is
\begin{align} \label{eq:nsc}
\dt{u} = - u\dx{u} - v\dy{u} - \dx{p} + \nu\left(\dxx{u}+\dyy{u} \right)
\end{align}
and after discretization in space on a staggered grid by centered finite differences we obtain
\begin{align}
\begin{aligned}
\dt{u_{i,j-1/2}} &= - u_{i,j-1/2} \frac{u_{i+1,j-1/2}-u_{i-1,j-1/2}}{2\Delta x} - v_{i,j-1/2} \frac{u_{i,j+1/2}-u_{i,j-3/2}}{2\Delta x} \\
&\phantom{=}\;\, -\frac{p_{i+1/2,j-1/2}-p_{i-1/2,j-1/2}}{\Delta x} \\
&\phantom{=}\;\, + \nu\frac{u_{i+1,j-1/2}+u_{i-1,j-1/2}-2 u_{i,j-1/2}}{\Delta x^2} + \nu\frac{u_{i,j+1/2}+u_{i,j-3/2}-2 u_{i,j-1/2}}{\Delta x^2}
\end{aligned}
\end{align}
where $v_{i,j-1/2}$ is computed by a simple mean without affecting the second order convergence of the overall scheme. Looking at Figure \ref{fig:staggered} we see that all the pressure unknowns are coupled together and only one pressure mode is allowed by this positioning of the unknowns. 

Near a boundary centered finite differences cannot be used since the grid spacing is not uniform, see Figure \ref{fig:staggered}. A non centered finite difference scheme can be easily found using Taylor series, for $x\in \Rb$ and a function $f\in \mathcal{C}^4\left([x-\Delta x/2, x+2 \Delta x]\right)$ we found the schemes given be the relations
\begin{subequations}
\begin{align} \label{eq:fd1}
\left| f'(x)-\frac{f(x+\Delta x) + 3 f(x)-4f(x-\Delta x/2)}{3 \Delta x} \right| &\leq \frac{1}{4}\Delta x^2 \|f''\|_{\infty}
\end{align}
and
\begin{align} \label{eq:fd2}
\left| f''(x)-\frac{16 f(x-\Delta x/2)- 25 f(x) + 10 f(x+\Delta x)-f(x+2\Delta x)}{5 \Delta x^2} \right| \leq \frac{9}{40}\Delta x^2 \|f^{iv}\|_{\infty},
\end{align}
\end{subequations}
where $\|\cdot\|_{\infty}$ is the maximum norm on the interval $[x-\Delta x/2, x+2 \Delta x]$. When we apply \eqref{eq:fd1} and \eqref{eq:fd2} to \eqref{eq:nsc} at a point $x_{i,1/2}$ we need the value $u_{i,0}$ which is given by the Dirichlet boundary conditions.

For the second component of equation \eqref{eq:ns1} we obtain a similar expression. We have thus obtained a system of $2 N (N-1)$ ODEs. 

The equation \eqref{eq:ns2} is discretized at the pressure points $x_{i-1/2,j-1/2}$ obtaining
\begin{align} \label{eq:pp}
\left( \nabla \cdot \bu \right)_{i-1/2,j-1/2} = \frac{u_{i,j-1/2}-u_{i-1,j-1/2}}{\Delta x}+\frac{v_{i-1/2,j}-v_{i-1/2,j-1}}{\Delta x}.
\end{align}
For both the projection method in section \ref{sec:proj} and the differential algebraic approach in section \ref{sec:rkdae2} we need to solve a Poisson problem with Neumann boundary conditions (see for example step 2 of section \ref{sec:method1} and equation \eqref{eq:defphii}) discretized at the pressure points $x_{i-1/2,j-1/2}$. We decided to solve the Poisson problem with the discrete cosine transform. This is a direct method which automatically imposes Neumann boundary conditions.

\section{Solving the Poisson problem by means of the discrete cosine transform}
In this chapter we will present the method used in order to solve the Poisson equation on $\Omega = [0,1]\times[0,1]$ with Neumann boundary conditions.

We decided to solve the Poisson equation using a discrete cosine transform (DCT). The DCT is a special case of the discrete Fourier transform (DFT) which automatically imposes Neumann boundary conditions. Moreover, like the DFT, the differential equation can be easily solved by algebraic manipulations. Many variants of the DCT exists in the literature (see \cite{dct_dst} for a survey) and the choice depends on the type of discretization of the domain. In the case of staggered grids the most appropriate DCT is called DCT-II, it will be presented later on.

In the following we will derive the DFT and its inverse from the Fourier transforms. Then the DCT-II and its inverse (called DCT-III) will be derived from the DFT. After we define the two dimensional DCT (2d-DCT) and explain how the Poisson problem can be solved with the 2d-DCT. Finally two fast algorithms for the computation of DCT-II and DCT-III will be derived and compared by numerical experiments.

For simplicity in the following we will refer to DCT-II with DCT and to DCT-III with IDCT.

\subsection{The discrete Fourier transform}
There is more than one way to obtain the discrete Fourier transform. In this section we will derive it from the Fourier transforms as in \cite{dft}. But it can be obtained also from Fourier series, trigonometric polynomials and using the delta function (see \cite{dft}).

\subsubsection{Deriving the DFT from the Fourier transform} \label{sec:derdft}
Let $f\in L^1(\Runo)$, the Fourier transform of $f$ is defined by
\begin{align}\label{eq:four}
\fho =  \intinf f(x) e^{-2\pi\complex \omega x} dx
\end{align}
and if $|\hat f| \in L^1(\Runo)$ the inverse is defined by
\begin{align}\label{eq:ifour}
f(x) =  \intinf \fho e^{2\pi\complex \omega x} d\omega .
\end{align}
In the following we will suppose that $f$ and $\hat f$ are sufficiently regular. For more details about the Fourier transform we refer to \cite{fourier-analysis}. $\omega$ is called the frequency variable while $x$ is called the spatial variable. Equation \eqref{eq:ifour} can be seen as an assembling of $f$ as a combination of modes with frequencies $\omega$ in $\Runo$, where each frequency has a weight $\fho$.

Assuming that $f(x)=0$ for $|x|>A/2>0$ equation \eqref{eq:four} reads
\begin{align}
\fho =  \intAd f(x) e^{-2\pi\complex \omega x} dx .
\end{align}
We wish approximate this integral numerically. Let $N\in \Nb$ be even, we define the $N+1$ grid points $x_n=n\Delta x$, where $\Delta x=A/N$ and $n=-N/2,\ldots ,N/2$. Using the composite trapezoidal rule we get
\begin{align}
\fho = \intAd g(x) dx \approx \frac{\Delta x}{2}\left( g\left(-\frac{A}{2}\right) + 2\sum_{n=-N/2+1}^{N/2-1}g\left(x_n\right) + g\left(\frac{A}{2}\right) \right),
\end{align}
where $g(x)=f(x) e^{-2\pi\complex \omega x}$. Under the assumption that $g\left(-\frac{A}{2}\right)=g\left(\frac{A}{2}\right)$ it yields
\begin{align}\label{eq:df1}
\fho = \frac{A}{N}\sum_{n=-N/2+1}^{N/2}f(x_n) e^{-2\pi\complex\omega x_n}.
\end{align}
The approximation \eqref{eq:df1} can be evaluated at any frequency $\omega$ but we would like to approximate also the integral in \eqref{eq:ifour} so must we choose a discretization for the frequency domain. Imagine a function $f$ which modes have an integer number of periods in the domain $[-A/2,A/2]$, so that it fits perfectly in the interval. Consider the wave with one period, it has wave-length $A$ and frequency $1/A$ units. The other waves will have frequency $k/A$ for $k\in \Nb$. It follows that the frequency unit has to be chosen as $\Delta \omega = 1/A$ and the grid points $\omega_k=k\Delta \omega$. In order to uniquely determine $\hat f (\omega_k)$ from $f(x_n)$ and vice-versa we choose also $N$ points $\omega_k$ in the frequency domain. The length of the latter is $\Omega = N\Delta\omega=N/A$ and $\omega_k$ is defined for $k=-N/2+1,\ldots .,N/2$. The following defines the DCT and its inverse IDFT.

\begin{defi}\label{def:dft}
Given the set of $N$ sampled values $f_n:=f(x_n)$ the DFT consists in the $N$ coefficients
\begin{align} \label{eq:dft}
D(f_n)_k=F_k = \frac{1}{N}\sum_{n=-N/2+1}^{N/2}f_n e^{-\frac{2\pi\complex}{N} k n}.
\end{align}
Approximations to $\hat f(\omega_k)$ are given by $A F_k$. The IDFT is another sequence of $N$ complex numbers given by 
\begin{align}\label{eq:idft}
D^{-1}(F_k)_n=f_n = \sum_{k=-N/2+1}^{N/2}F_k e^{\frac{2\pi\complex}{N} k n}.
\end{align}
\end{defi}
It remains to show that the DFT and the IDFT are really inverse one of the other, i.e. $D^{-1}(D(f_n)_k)_n=f_n$ and $D(D^{-1}(F_k)_n)_k=F_k$. This is easy to verify using the identity \eqref{eq:orthdft}.
\begin{align}\label{eq:orthdft}
\frac{1}{N}\sum_{n=-N/2+1}^{N/2}e^{-\frac{2\pi\complex}{N} k n}e^{\frac{2\pi\complex}{N} j n} = \left\{
\begin{array}{ll}
1 & \text{if } j-k\mod N \equiv 0 , \\
0 & \text{otherwise.}
\end{array}
\right.
\end{align}

\subsection{The discrete cosine transform}
The DCT is derived from the DFT using a sequence $f_n$ of sampled points which is even and real. Different sequences give rise of different DCTs, here the sequence $f_n$ is sampled on a staggered grid so that the arising DCT is appropriate for this kind of spatial discretization. Firstly we derive the one dimensional DCT, then we define the two dimensional DCT.

\subsubsection{Deriving the DCT from the DFT} \label{sec:derdct}
Let $N \in \Nb$, $\Delta x=1/N$ and $x_{n+1/2}=(n+1/2)\Delta x$ be a sequence of points on a one dimensional staggered grid. Let $f\in C^0([0,1])$ and $f_{n+1/2}=f(x_{n+1/2})$ its discretization for $n=0,\ldots,N-1$. We define $\left(g_n\right)_{n=-2N+1}^{2N} \in \Rb$ as the sequence of $4 N$ real numbers given by
\begin{align*}
g_n=
\begin{cases}
0 & \text{if $n$ is even,} \\
f_{n/2} & \text{if $n\geq 0$ and $n$ is odd,} \\
f_{-n/2} & \text{if $n\leq 0$ and $n$ is odd.}
\end{cases}
\end{align*}
Applying the DFT given by \eqref{eq:dft} to the above sequence we obtain
\begin{align*}
F_k &= \frac{1}{4 N}\sum_{n=-2N+1}^{2N}g_n e^{-\frac{2\pi\complex}{4N} k n}= \frac{1}{4N} \sum_{n=-2N+1}^{2N}g_n \left( \cos\left( \frac{2\pi}{4N}kn\right)-\complex \sin\left( \frac{2\pi}{4N}kn\right) \right).
\end{align*}
Using $g_n=g_{-n}$ and $g_n=0$ for $n$ even it follows
\begin{align}
F_k &= \frac{1}{2N} \sum_{n=1}^{2N-1}g_n \cos\left( \frac{2\pi}{4N}kn\right) = \frac{1}{2N} \sum_{n=0}^{N-1}g_{2n+1} \cos\left( \frac{2\pi}{4N}k(2n+1)\right) \label{eq:pdct}
\end{align}
Inserting $f_{n+1/2}$ in equation \eqref{eq:pdct} the coefficients $F_k$ of the DCT are given by
\begin{align}\label{eq:dcts}
F_k=\frac{1}{2N} \sum_{n=0}^{N-1}f_{n+1/2} \cos\left( \frac{2n+1}{2N}k\pi\right).
\end{align}
It is easy to see that $F_k=F_{-k}$ and $F_k=-F_{2N-k}$. Thanks to these symmetries the inverse DCT is given by
\begin{align}
g_n &= \sum_{k=-2N+1}^{2N}F_k e^{\frac{2\pi\complex}{4N} k n} = \sum_{k=-2N+1}^{2N}F_k \left( \cos\left( \frac{2\pi}{4N}kn\right)+\complex \sin\left( \frac{2\pi}{4N}kn\right) \right) \nonumber \\
&= F_0 + 2\sum_{k=1}^{2N-1} F_k  \cos\left( \frac{2\pi}{4N}kn\right)+F_{2N}\cos(n\pi) \nonumber \\
&= F_0 \left(1-(-1)^n\right) + 2\sum_{k=1}^{N-1} F_k \left(1-(-1)^n\right) \cos\left( \frac{2\pi}{4N}kn\right) \nonumber 
\end{align}
As expected $g_n=0$ if $n$ is even. For odd indexed $g_n$ we have the coefficients $f_{n+1/2}$ of the IDCT:
\begin{align}\label{eq:idcts}
f_{n+1/2}=g_{2n+1} = 4\left(  \frac{1}{2}F_0 + \sum_{k=1}^{N-1}F_k \cos\left( \frac{2n+1}{2N}k\pi \right) \right).
\end{align}
One sees that a sequence $f_{n+1/2}$ given by \eqref{eq:idcts} is symmetric in $f_0$ and $f_N$, i.e. $f_{-1/2}=f_{1/2}$ and $f_{N-1/2}=f_{N+1/2}$.

It is very common to merge together the scaling factors, the final form of the DCT is given by the following definition.
\begin{defi}
Given a sequence $f_{n+1/2}$ for $n=0,\ldots ,N-1$ the DCT coefficients are given by 
\begin{align}\label{eq:dct}
F_k=\sum_{n=0}^{N-1}f_{n+1/2} \cos\left( \frac{2n+1}{2N}k\pi\right)
\end{align}
for $k=0,\ldots ,N-1$. The IDCT is given by
\begin{align}
f_{n+1/2}&= \frac{2}{N}\left(  \frac{1}{2}F_0 + \sum_{k=1}^{N-1}F_k \cos\left( \frac{2n+1}{2N}k\pi \right) \right) =\frac{2}{N} \sideset{}{'}\sum_{k=0}^{N-1}F_k \cos\left( \frac{2n+1}{2N}k\pi \right) \label{eq:idct}
\end{align}
for $n=0,\ldots ,N-1$. Here $\sum'$ is a sum in which the first term is weighted by one half.
\end{defi}

\subsubsection{The two dimensional DCT}
With the same techniques used in \ref{sec:derdft} and \ref{sec:derdct} one can derive the two dimensional DFT (2d-DFT) from the two dimensional Fourier transform and the two dimensional DCT (2d-DCT) from the 2d-DFT. Here we will skip the derivation of the 2d-DCT giving only his definition.

\begin{defi}
Given a two dimensional sequence $f_{m+1/2,n+1/2}$ for $m,n=0,\ldots ,N-1$ the 2d-DCT coefficients are given by 
\begin{align}\label{eq:2ddct}
F_{j,k}=\sum_{m=0}^{N-1}\sum_{n=0}^{N-1}f_{m+1/2,n+1/2} \cos\left( \frac{2m+1}{2N}j\pi\right)\cos\left( \frac{2n+1}{2N}k\pi\right)
\end{align}
for $j,k=0,\ldots ,N-1$. The 2d-IDCT is given by
\begin{align}\label{eq:2didct}
f_{m+1/2,n+1/2}=\frac{4}{N^2} \sideset{}{'}\sum_{j=0}^{N-1}\sideset{}{'}\sum_{k=0}^{N-1}F_{j,k} \cos\left( \frac{2m+1}{2N}j\pi \right)  \cos\left( \frac{2n+1}{2N}k\pi \right)
\end{align}
for $m,n=0,\ldots ,N-1$.
\end{defi}

\subsection{Solving the Poisson problem}
In this section we show how the Poisson problem with Neumann boundary conditions can be solved using the 2d-DCT.

Consider the Poisson problem
\begin{numcases}{}
\nabla^2 u = f  & in $\Omega$,\nonumber\\
\dn{u}=0 & on $\partial\Omega$.\nonumber
\end{numcases}
Let $u_{m+1/2,n+1/2}$ and $f_{m+1/2,n+1/2}$ be the $u$ and $f$ functions sampled at the staggered points $x_{m+1/2,n+1/2}$ for $m,n=0,\ldots ,N-1$. We note the 2d-DCT coefficients of $u$ and $f$ by $U_{j,k}$ and $F_{j,k}$ respectively. The five points Laplacian in $x_{m+1/2,n+1/2}$ gives
\begin{align*}
f_{m+1/2,n+1/2} &= \frac{u_{m+1/2,n+3/2}+u_{m+1/2,n-1/2}-2 u_{m+1/2,n+1/2}}{\Delta x^2}  \\
&\phantom{=}\phantom{a}+\frac{u_{m+3/2,n+1/2}+u_{m-1/2,n+1/2}-2 u_{m+1/2,n+1/2}}{\Delta x^2} 
\end{align*}
and applying the 2d-DCT we get
\begin{align*}
&\sideset{}{'}\sum_{j=0}^{N-1}\sideset{}{'}\sum_{k=0}^{N-1}\Delta x^2 F_{j,k} \cos\left( \frac{2m+1}{2N}j\pi \right)  \cos\left( \frac{2n+1}{2N}k\pi \right) \\
&= \sideset{}{'}\sum_{j=0}^{N-1}\sideset{}{'}\sum_{k=0}^{N-1} U_{j,k}\left( 2\cos\left(\frac{k\pi}{N}\right)-2 \right) \cos\left( \frac{2m+1}{2N}j\pi \right)  \cos\left( \frac{2n+1}{2N}k\pi \right) \\
&+\sideset{}{'}\sum_{j=0}^{N-1}\sideset{}{'}\sum_{k=0}^{N-1} U_{j,k}\left( 2\cos\left(\frac{j\pi}{N}\right)-2 \right) \cos\left( \frac{2m+1}{2N}j\pi \right)  \cos\left( \frac{2n+1}{2N}k\pi \right).
\end{align*}
Making a comparison term by term gives the equation for $U_{j,k}$:
\begin{align}\label{eq:ujk}
U_{j,k}\left( 2\cos\left(\frac{j\pi}{N}\right)+2\cos\left(\frac{k\pi}{N}\right)-4 \right) = \Delta x^2 F_{j,k}.
\end{align}
The procedure for solving the Poisson equation with the 2d-DCT is:
\begin{enumerate}
\item compute $F_{j,k}$ the 2d-DCT coefficients of $f_{m+1/2,n+1/2}$,
\item compute $U_{j,k}$ with equation \eqref{eq:ujk},
\item compute $u_{m+1/2,n+1/2}$ with the 2d-IDCT.
\end{enumerate}
Notice that a solution of equation \eqref{eq:ujk} for $j=k=0$ fails to exists unless $\Delta x^2 F_{0,0}=0$. Hopefully it is the case, we have 
\begin{align*}
0&=\int_{\partial \Omega}\dn{u}=\into \nabla^2 u= \into f = \sum_{m=0}^{N-1}\sum_{n=0}^{N-1}\left(\Delta x^2 f_{m+1/2,n+1/2} + \bigo{\Delta x^3}\right) \\
&=\Delta x^2 F_{0,0} + \bigo{\Delta x}
\end{align*}
thus $\Delta x^2 F_{0,0}=\bigo{\Delta x}$ as $\Delta x \rightarrow 0$. Consequently it is valid to suppose $\Delta x^2 F_{0,0}=0$. In this way any value of $U_{0,0}$ will satisfy equation \eqref{eq:ujk}, we set $U_{0,0}=0$ which means
\begin{align*}
\into u = \bigo{\Delta x},
\end{align*}
i.e. setting $U_{0,0}=0$ imposes the mean value of $u$.

Observe that the coefficients $u_{m+1/2,n+1/2}$ given by the 2d-IDCT defined in \eqref{eq:2didct} will satisfy the Neumann boundary conditions. For instance we have $u_{1/2,n+1/2}=u_{-1/2,n+1/2}$ for $n=0,\ldots ,N-1$.

\subsection{Fast algorithms for the discrete cosine transform}
In this section we will explain how the 2d-DCT is computed in our implementation. The 2d-DCT in \eqref{eq:2ddct} can be written 
\begin{align*}
F_{j,k}&=\sum_{m=0}^{N-1}\left( \sum_{n=0}^{N-1}f_{m+1/2,n+1/2} \cos\left( \frac{2n+1}{2N}j\pi\right)\right)\cos\left( \frac{2m+1}{2N}k\pi\right) \\
&=\sum_{m=0}^{N-1}\tilde{F}_{j,m}\cos\left( \frac{2m+1}{2N}k\pi\right)
\end{align*}
where $\tilde{F}_{j,m}=\sum_{n=0}^{N-1}f_{m+1/2,n+1/2} \cos\left( \frac{2n+1}{2N}j\pi\right)$. Hence the 2d-DCT reduces to $2N$ DCTs. In the following we will present 2 algorithms for the computation of the DCT.

\subsubsection{Iterative algorithm}
This algorithm has been taken from \cite{rec-dct}. It is an iterative algorithm for the computations of the DCT and the IDCT. Its computational complexity is $\bigo{N^2}$ but the operations involved are mostly add and multiply. So the algorithm is much faster than a naïve computation of \eqref{eq:dct} and \eqref{eq:idct}, which needs a lot of cosine computations.

\textbf{DCT}
Let $f_n \in \Rb$ for $n=0,\ldots ,N-1$ be a real sequence with $N$ even. We set $w_n^k=f_n+(-1)^k f_{N-1-n}$, equation \eqref{eq:dct} can be written
\begin{align*}
F_k &= \sum_{n=0}^{N/2-1} f_n \cos\left(\frac{2n+1}{2N}k\pi\right)+f_{N-1-n} \cos\left(\frac{2(N-1-n)+1}{2N}k\pi\right) \\
&= \sum_{n=0}^{N/2-1} w_n^k \cos\left(\frac{2n+1}{2N}k\pi\right) \\
&= \sum_{n=0}^{N/2-1} w_{N/2-1-n}^k \left( \cos\left(\frac{k\pi}{2}\right)\cos\left(\frac{2n+1}{2N}k\pi\right) +\sin\left(\frac{k\pi}{2}\right)\sin\left(\frac{2n+1}{2N}k\pi\right) \right)
\end{align*}
which yields
\begin{align*}
F_k=
\begin{cases}
G_{N/2-1}(k) (-1)^{k/2} & \text{if $k$ is even,} \\
H_{N/2-1}(k) (-1)^{(k-1)/2} & \text{if $k$ is odd,}
\end{cases}
\end{align*}
where 
\begin{align}\label{eq:ghj}
G_j(k) = \sum_{n=0}^j w_{j-n}^k \cos\left(\frac{2n+1}{2N}k\pi\right), \qquad H_j(k) = \sum_{n=0}^j w_{j-n}^k \sin\left(\frac{2n+1}{2N}k\pi\right).
\end{align}
Using recursive properties of sine and cosine we can compute $G_j(k)$ and $H_j(k)$ in a recursive way as well. Setting $\theta_k=k\pi/N$ we have
\begin{align*}
G_j(k) = \cos\left(\frac{\theta_k}{2}\right)(w_j^k-w_{j-1}^k)+2\cos(\theta_k)G_{j-1}(k)-G_{j-2}(k),\\
H_j(k) = \sin\left(\frac{\theta_k}{2}\right)(w_j^k+w_{j-1}^k)+2\cos(\theta_k)H_{j-1}(k)-H_{j-2}(k).
\end{align*}

\textbf{IDCT}
When computing the IDCT we pre multiply $F_0$ by one half, so that we can write
\begin{align*}
f_n &= \frac{2}{N} \sum_{k=0}^{N-1}F_k \cos\left( \frac{2n+1}{2N}k\pi \right)
\end{align*}
and
\begin{subequations}
\begin{align}
f_n +f_{N-1-n} &= \frac{2}{N} \sum_{k=0}^{N-1}F_k \left(1+(-1)^k \right)\cos\left( \frac{2n+1}{2N}k\pi \right) \notag\\
&= \frac{4}{N} \sum_{k=0}^{N/2-1}F_{2k}\cos\left( \frac{2n+1}{N}k\pi \right) = (-1)^n \frac{4}{N} P_{N/2-1}(n) \label{eq:idcti1}
\end{align}
where $P_j(n)=\sum_{k=0}^j F_{2(j-k)} \sin\left(\frac{(2n+1)}{N}(k+1)\pi\right)$. Similarly
\begin{align}\label{eq:idcti2}
f_n - f_{N-1-n} &= (-1)^n \frac{4}{N} Q_{N/2-1}(n)
\end{align}
\end{subequations}
where $Q_j(n)=\sum_{k=0}^j F_{2(j-k)+1} \sin\left(\frac{(2n+1)}{2N}(2k+1)\pi\right)$. Using \eqref{eq:idcti1} and \eqref{eq:idcti2} we have
\begin{align*}
f_n &= (-1)^n \frac{2}{N}\left( P_{N/2-1}(n)+Q_{N/2-1}(n) \right), \\
f_{N-1-n} &= (-1)^n \frac{2}{N}\left( P_{N/2-1}(n)-Q_{N/2-1}(n) \right) .
\end{align*}
As for $G_j$, $H_j$ also for $P_j$, $Q_j$ it exists a recursive definition. Setting $\theta_n=(2n+1)\pi/N$ we can show that
\begin{align*}
P_j(n) &= \sin(\theta_n)F_{2j}+2\cos(\theta_n)P_{j-1}(n)-P_{j-2}(n), \\
Q_j(n) &= \sin\left(\frac{\theta_n}{2}\right)(F_{2j+1}+F_{2j-1})+2\cos(\theta_n)Q_{j-1}(n)-Q_{j-2}(n).
\end{align*}

\subsubsection{Recursive algorithm}
The following algorithm has been presented in \cite{cordic-dct}. It is recursive and thanks to this property its computational complexity is $\bigo{N\log_2(N)}$. 

\textbf{DCT}
Let $f_n \in \Rb$ for $n=0,\ldots ,N-1$ be a real sequence with $N=2^m$, $m\in\Nb$. We set
\begin{align*}
f^L_n&=f_{2n}+f_{2n+1}, \\
f^H_n&=f_{2n}-f_{2n+1}
\end{align*}
for $n=0,\ldots ,N/2-1$. The DCT is given by
\begin{align*}
F_k &= \sum_{n=0}^{N/2-1}f_{2n}\cos\left(\frac{4n+1}{2N}k\pi\right)+f_{2n+1}\cos\left(\frac{4n+3}{2N}k\pi\right) \\
&=\frac{1}{2}\sum_{n=0}^{N/2-1}(f^L_n+f^H_n)\cos\left(\frac{4n+1}{2N}k\pi\right)+(f^L_n-f^H_n)\cos\left(\frac{4n+3}{2N}k\pi\right)\\
&=\cos\left(\frac{k\pi}{2N}\right)\sum_{n=0}^{N/2-1}f^L_n\cos\left(\frac{2n+1}{N}k\pi\right)+\sin\left(\frac{k\pi}{2N}\right)\sum_{n=0}^{N/2-1}f^H_n\sin\left(\frac{2n+1}{N}k\pi\right).
\end{align*}
Using $\cos\left(\frac{2n+1}{N}\left(\frac{N}{2}-k\right)\pi\right)=(-1)^n\sin\left(\frac{2n+1}{N}k\pi\right)$ it follows
\begin{align*}
F_k &=\cos\left(\frac{k\pi}{2N}\right)\sum_{n=0}^{N/2-1}f^L_n\cos\left(\frac{2n+1}{N}k\pi\right)\\
&\phantom{=}\phantom{1}+\sin\left(\frac{k\pi}{2N}\right)\sum_{n=0}^{N/2-1}f^H_n(-1)^n\cos\left(\frac{2n+1}{N}\left(\frac{N}{2}-k\right)\pi\right)
\end{align*}
and it holds also
\begin{align*}
F_{N-k} &=-\sin\left(\frac{k\pi}{2N}\right)\sum_{n=0}^{N/2-1}f^L_n\cos\left(\frac{2n+1}{N}k\pi\right)\\
&\phantom{=}\phantom{1}+\cos\left(\frac{k\pi}{2N}\right)\sum_{n=0}^{N/2-1}f^H_n(-1)^n\cos\left(\frac{2n+1}{N}\left(\frac{N}{2}-k\right)\pi\right).
\end{align*}
Setting
\begin{align*}
A_k &= \sum_{n=0}^{N/2-1}f^L_n\cos\left(\frac{2n+1}{N}k\pi\right) \\
B_k &= \sum_{n=0}^{N/2-1}f^H_n(-1)^n\cos\left(\frac{2n+1}{N}k\pi\right)
\end{align*}
we get
\begin{subequations} \label{eq:fk}
\begin{align}
&F_0 = A_0, \\
&F_{N/2} = \frac{1}{\sqrt{2}} B_0, \\
&F_k = \cos\left(\frac{k\pi}{2N}\right)A_k+\sin\left(\frac{k\pi}{2N}\right)B_{N/2-k}, \label{eq:fk1}\\
&F_{N-k} = -\sin\left(\frac{k\pi}{2N}\right)A_k+\cos\left(\frac{k\pi}{2N}\right)B_{N/2-k}\label{eq:fk2}
\end{align}
\end{subequations}
for $k=1,\ldots ,N/2-1$. The computation of a $N$ point DCT has been reduced to two $N/2$ point DCTs and $N/2$ rotations. This recursion gives a computational complexity of $\bigo{N\log_2(N)}$. 

\textbf{IDCT}
In \cite{cordic-dct} the recursive algorithm for IDCT if found using the signal flow of the DCT. Here we will derive the IDCT algorithm using the identity
\begin{align*}
\sum_{k=0}^{N/2-1}\cos\left(\frac{2n+1}{N}k\pi\right)\cos\left(\frac{2m+1}{N}k\pi\right) = \frac{N}{4}\delta_{m,n} + \frac{1}{2}.
\end{align*}
From \eqref{eq:fk}(c,d) we have
\begin{align}\label{eq:ak}
A_k = \cos\left(\frac{k\pi}{2N}\right)F_k -\sin\left(\frac{k\pi}{2N}\right) F_{N-k}.
\end{align}
Multiplying the right hand side of \eqref{eq:ak} by $\cos\left(\frac{2m+1}{N}k\pi\right)$ and summing it over $k$ gives
\begin{align*}
\sum_{k=0}^{N/2-1}A_k &\cos\left(\frac{2m+1}{N}k\pi\right) = \sum_{n=0}^{N/2-1}f^L_n\sum_{k=0}^{N/2-1}\cos\left(\frac{2n+1}{N}k\pi\right) \cos\left(\frac{2m+1}{N}k\pi\right) \\
&= \frac{N}{4}f^L_m+\frac{1}{2}A_0 = \frac{N}{4}f^L_m+\frac{1}{2}F_0,
\end{align*}
doing the same on the left hand side we obtain
\begin{align*}
&\sum_{k=0}^{N/2-1}\left(\cos\left(\frac{k\pi}{2N}\right)F_k -\sin\left(\frac{k\pi}{2N}\right) F_{N-k}\right) \cos\left(\frac{2m+1}{N}k\pi\right)  \\
&= F_0 + \sum_{k=1}^{N/2-1}\left(\cos\left(\frac{k\pi}{2N}\right)F_k -\sin\left(\frac{k\pi}{2N}\right) F_{N-k}\right) \cos\left(\frac{2m+1}{N}k\pi\right).
\end{align*}
Setting $w_k=\cos\left(\frac{k\pi}{2N}\right)F_k -\sin\left(\frac{k\pi}{2N}\right) F_{N-k}$ for $k=0,\ldots ,N/2-1$ we have
\begin{align} \label{eq:xl}
f^L_m &=\frac{4}{N}\sideset{}{'}\sum_{k=0}^{N/2-1} w_k \cos\left(\frac{2m+1}{N}k\pi\right)
\end{align}
and in a very similar fashion we obtain
\begin{align}\label{eq:xh}
(-1)^m f^H_m &=\frac{4}{N}\sideset{}{'}\sum_{k=0}^{N/2-1} v_k \cos\left(\frac{2m+1}{N}k\pi\right)
\end{align}
where $v_k=\frac{\sqrt{2}}{2}\left( \cos\left(\frac{k\pi}{2N}\right)\left(F_{N/2+k}+F_{N/2-k}\right)+\sin\left(\frac{k\pi}{2N}\right)\left(F_{N/2+k}-F_{N/2-k}\right) \right)$.
Using \eqref{eq:xl} and \eqref{eq:xh} we have
\begin{align*}
f_{2m} &= \frac{1}{2}\left(C_m +(-1)^m D_m\right) \\
f_{2m+1} &= \frac{1}{2}\left(C_m -(-1)^m D_m\right) 
\end{align*}
for $n=0,\ldots ,N/2-1$, where
\begin{align}
C_m = \frac{2}{N/2}\sideset{}{'}\sum_{k=0}^{N/2-1} w_k \cos\left(\frac{2m+1}{N}k\pi\right) \\
D_m = \frac{2}{N/2}\sideset{}{'}\sum_{k=0}^{N/2-1} v_k \cos\left(\frac{2m+1}{N}k\pi\right).
\end{align}
In this way the $N$ point IDCT has been decomposed in two $N/2$ IDCTs. The DCT and IDCT recursion stops when $N=2$, at this point the transform is computed using \eqref{eq:dct} or \eqref{eq:idct}.

\subsubsection{Comparison of the two algorithms}
\begin{figure}[!hbtp]
\begin{center}
\subfigure[CPU time of $N$ points DCT and IDCT.]{
\includegraphics[trim=0.0cm 0.0cm 0.0cm 0.0cm, clip=true, height=0.18\textheight]{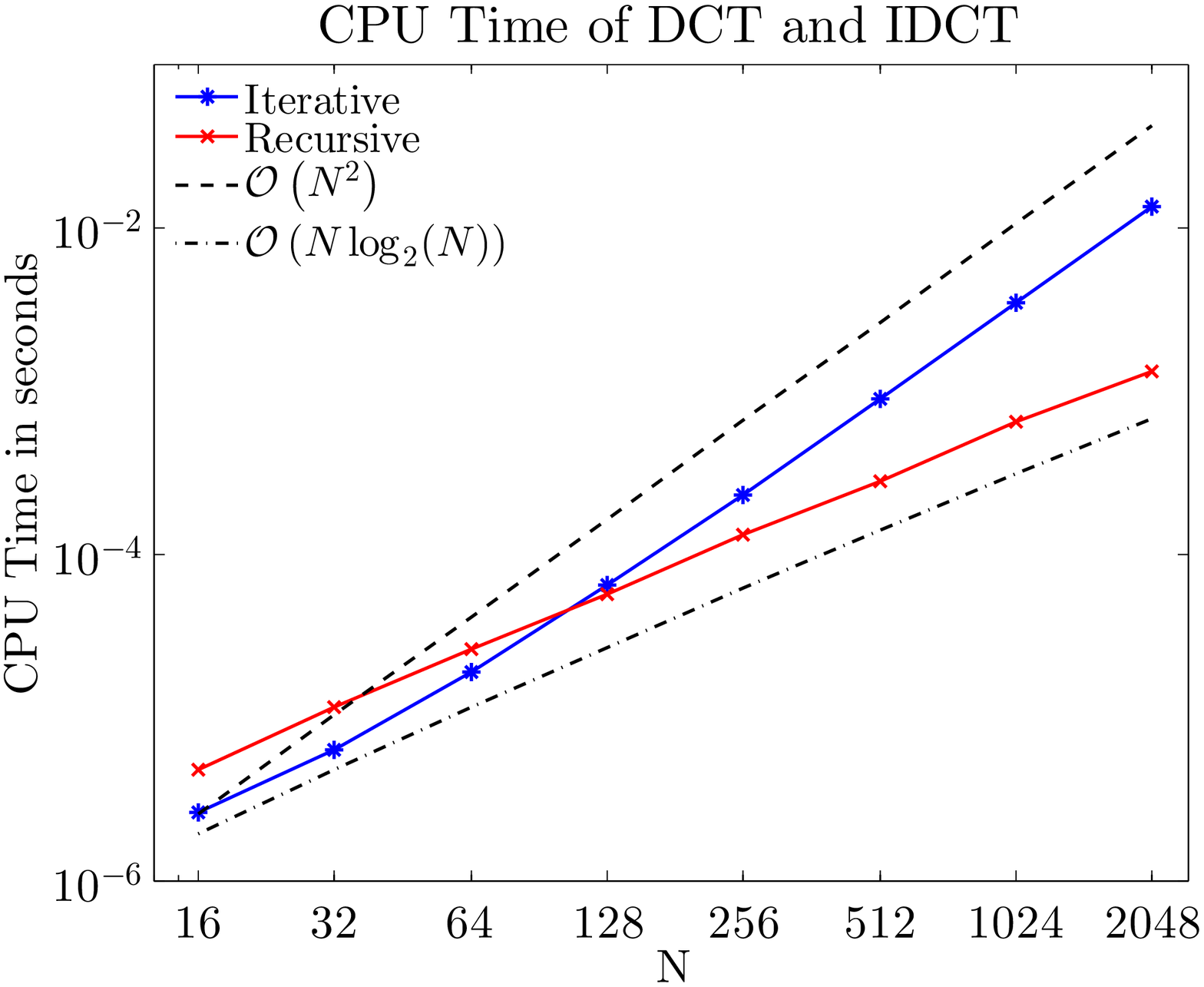} \label{fig:dct}}
\subfigure[CPU time of Poisson solver on a $N\times~N$ grid.]{
\includegraphics[trim=0.0cm 0.0cm 0.0cm 0.0cm, clip=true, height =0.18\textheight]{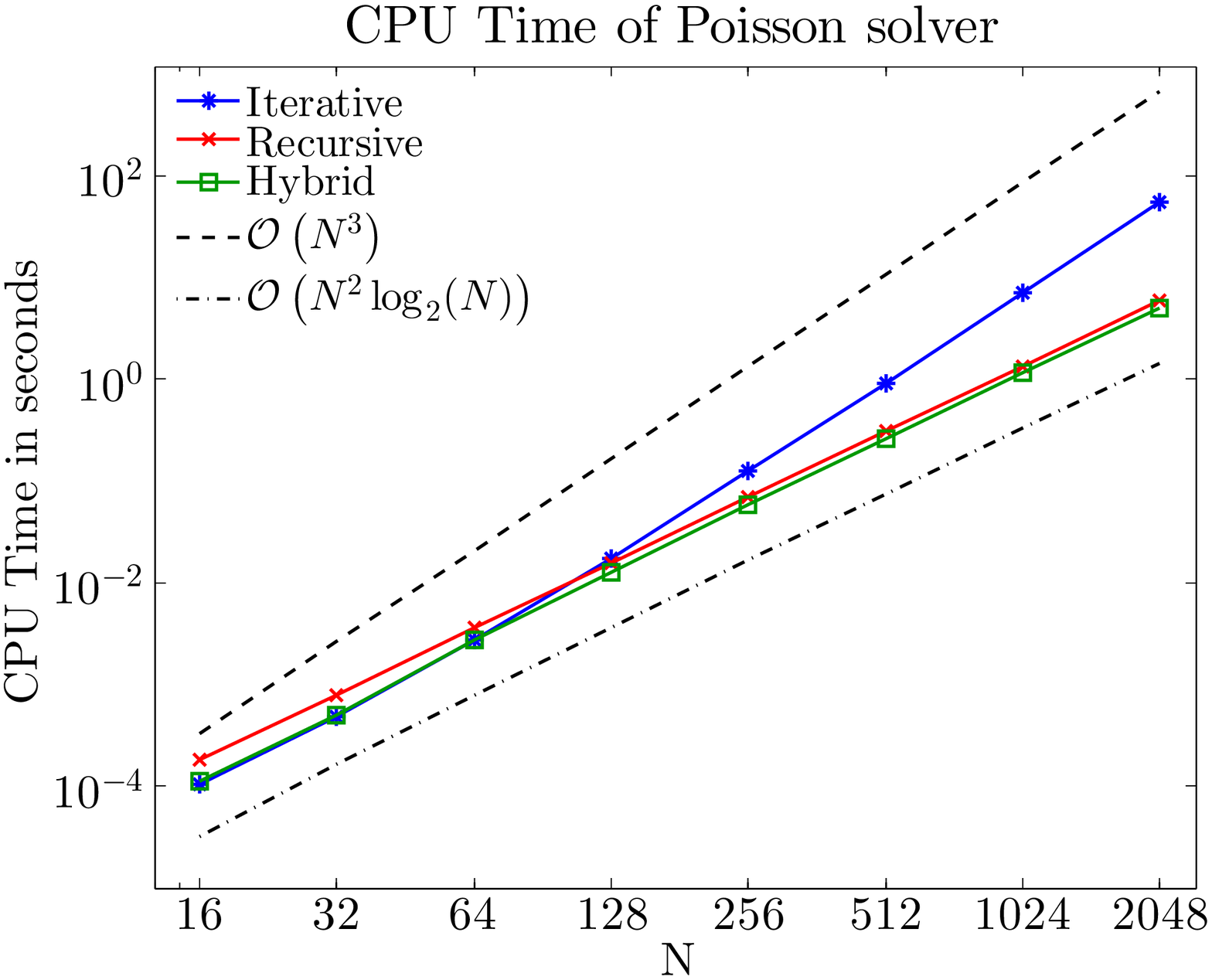} \label{fig:pois}}
\end{center}
\caption{Comparison of CPU time of different DCT algorithms.}
\label{fig:dctcomp}
\end{figure}

In Figure \ref{fig:dct} we show the CPU time of the recursive and the iterative algorithms. The plot suggests to use the iterative method if $N\leq 64$ and the recursive one otherwise. We have implemented an hybrid method which is like the recursive one but the recursion stops when $N=64$, at this point the iterative method is called. In Figure \ref{fig:pois} we compare the three approaches when solving the Poisson problem.

\newpage
\thispagestyle{empty}
\mbox{}
\clearpage

\chapter{Numerical experiments} \label{sec:numexp}
In this chapter we will compare the three schemes RKC, ROCK2 and PIROCK when coupled with the different approaches for solving the incompressible Navier-Stokes equations explained in chapters \ref{sec:serkpm} and \ref{sec:erkinse}. Each numerical experiment is executed several times but using different methods or different test problems, for this reason we will describe each experiment once for all and tell which methods are used before showing the results. After we will consider three test problems, for each one of these we show and comment the results of the previously described experiments.

\section{The codes used and description of the numerical experiments}
In this section we will firstly explain how we name the different couplings between the Runge-Kutta methods and the different approaches for solving the Navier-Stokes equations, we also give a few remarks about these couplings. Then we will describe the numerical experiments and tell which methods are used to perform the experiment.

\subsection{Names for the different methods and remarks}\label{sec:desc_name}
For each numerical method we give a name composed by three parts. 

The first part is the name of the numerical integrator, it can be RKC, ROCK2 or PIROCK. The second component tells which approach has been used to solve the incompressible Navier-Stokes equations. It is PM1 if the projection method of section \ref{sec:method1} is used, PM1V if we use the variant of PM1, i.e. the method explained at the end of section \ref{sec:pm1var}, and PM3 if the method of section \ref{sec:method3} is used. If we use approach 1 or 2 of section \ref{sec:circordcond} the name will be AP1 or AP2 respectively, if the workaround of Approach 2 for order one stages is used it will be AP2W. The third component of the name tells if the second order pressure is computed at each time step or only in the last one, it will be CP1 or CP0 respectively (to be read Compute Pressure true or false). 

For example the name 'ROCK2, AP1, CP0' means that we use ROCK2 coupled with AP1 and the pressure is computed only in the last time step.

Now we give a few remarks about the possible couplings. After the first experiments PIROCK, PM1, CP0 has shown stability issues which we did not had time to investigate, for this reason we did not implemented PM1V, AP1 and AP2 in PIROCK. If RKC is used time step adaptivity can be enabled only when coupled with PM1, otherwise the local error estimation gives wrong results, the reason is given in \ref{sec:tsairkc}. If PM1 or PM1V are used a first order pressure is computed after each time step, disregarding about CP0 or CP1. If AP1, AP2 or AP2W are used with CP0 no pressure is computed until the last time step. Observe that AP1, AP2 and AP2W in general give the same results for the velocity. The only exception is when in RKC a small time step is used and the pressure is computed, in this case the minimal allowed number of stages differs between AP1 and AP2, giving different velocities. Remember that with ROCK2 the approach AP2 cannot be used and we use AP2W, with RKC we use AP2. CP0 or CP1 do not change the way the pressure is computed but just the frequency, the method for the pressure computation is given by PM1, PM1V, AP1, AP2 or AP2W.

\subsection{Convergence order}\label{sec:desc_conv}
The convergence order test is done to check that the effective convergence order matches the theoretical one. In numerical experiments the time and space errors are mixed, in order to split these errors the following method is used.

Following \cite{classfully} the numerical solution of a temporally $p$th order and spatially $r$th order accurate scheme can be expressed as
\begin{align} \label{eq:err}
\eta_{i j}^n= \eta(x_i,y_j,t_n) + \beta_{i j}^n \left( \Delta t \right)^p + \gamma_{i j}^n\left(\Delta x\right)^r + \epsilon
\end{align}
where $\eta_{i j}^n$ and $\eta(x_i,y_j,t_n)$ are the numerical and analytical solutions respectively. $\beta_{i j}^n \left( \Delta t \right)^p$ and $\gamma_{i j}^n\left(\Delta x\right)^r$ are the errors corresponding to the time and spatial discretization and $\epsilon$ is the round-off error. Thanks to expression \eqref{eq:err} we can carry out time and space convergence tests separately. For time convergence we compute a numerical reference solution $\tilde{\eta}_{i j}^n$ using a very small time step so that we can write
\begin{align}
\tilde{\eta}_{i j}^n \approx \eta(x_i,y_j,t_n) + \gamma_{i j}^n\left(\Delta x\right)^r + \epsilon,\notag
\end{align}
then various numerical solutions $\eta_{i j}^n$ are computed with different time steps on the same grid giving the relation
\begin{align}
\|\eta^n-\tilde{\eta}^n\| = \bigo{\Delta t^p},\notag
\end{align}
which gives the temporal convergence order. Similarly for the spatial convergence we compute a reference solution $\hat{\eta}_{i j}^n$ with a very fine grid, yielding
\begin{align}
\hat{\eta}_{i j}^n \approx \eta(x_i,y_j,t_n) + \beta_{i j}^n \left( \Delta t \right)^p + \epsilon. \notag
\end{align}
Computing various numerical solutions $\eta_{i j}^n$ with a fixed time step and different grid size we obtain
\begin{align}
\|\eta^n-\hat{\eta}^n\| = \bigo{\Delta x^r}, \notag
\end{align}
which gives the spatial convergence order.

In our experiments the reference solution $\tilde{\eta}$ for time convergence is noted $\bu_{ref}$ and it is computed with $\Delta x = 1/128$ and a small time step $\Delta t = 2^{-16}$. $\bu_{ref}$ is compared to other solutions noted $\bu_{\Delta t}$ computed with the same $\Delta x$ and $\Delta t = 2^{-m}$, where $m=0,1,\ldots,15$. The reference solution $\hat\eta$ for space convergence is also noted $\bu_{ref}$ and it is computed with $\Delta x=1/512$ and $\Delta t = 2^{-14}$. It is compared to various solutions $\bu_{\Delta x}$ computed with the same $\Delta t$ and different $\Delta x$s. Similarly for the pressure we note as $p_{ref}$ the reference solutions and $p_{\Delta t}$, $p_{\Delta x}$ the other solutions. When we use PM1 or PM1V first and second order pressures are available, in the figures we note them as $p^1$ and $p^2$ respectively. All the tests are done with $\Rey=100$ and the errors are measured with the infinity norm at time $t=1$.

These convergence tests are not done for CP1 but only for CP0. Also, the results of AP1, AP2 and AP2W are put in the same figure since the velocities are the same, only the pressure differs. We will now list the methods used in the convergence tests:
\begin{itemize}
\item RKC with PM1, PM1V, AP1, AP2 and CP0,
\item ROCK2 with PM1, PM1V, AP1, AP2W and CP0,
\item PIROCK with PM1 and CP0.
\end{itemize}
 
\subsection{Stability tests}\label{sec:desc_stab}
The stability tests are done to investigate the size of the stability domain of the different methods. Two different tests are done, in order to confirm the generality of the results. As we already said PIROCK has stability issues which did not permit to perform this test, hence we did it for RKC and ROCK2. 

The first test is done by taking a set of stage numbers and searching for the maximal stable time step when solving the Navier-Stokes equations with $\Delta x = 1/128$ and $\Rey=5$. The goal is to check if the stability domain size's is the same as the one for ODEs. We know that in the ODE context the stability domain grows as $0.653 s^2$ and $0.811 s^2$ along the negative real axis for RKC and ROCK2 respectively, where $s$ it the number of stages. So the theoretical maximal stable time step is given by $\Delta t = 0.653 s^2/\rho$ and  $\Delta t = 0.811 s^2/\rho$ for RKC and ROCK2 respectively, where $\rho$ is the spectral radius of the diffusion operator. In the code we estimate $\rho$ using Gershgorin's theorem.

The first test is done for $\Rey=5$, in order to show that the results do not depend on $\Rey$ we did a second test. For this test we neglected the advection term since we are interested on the Reynolds number effect. Fixing the time step size at $\Delta t =10^{-2}$ we take a set of Reynolds numbers and for each one of these we search for the minimal number of stages $s$ which gives a stable solution. The theoretical minimal $s$ is given by $s=\sqrt{\Delta t \rho /0.653}$ and $s=\sqrt{\Delta t \rho /0.811}$ for RKC and ROCK2 respectively, where $\rho$ depends on $\Rey$.

In both the tests we integrate from $t=0$ to $t=1$. CP0 and CP1 do not affect the result thus these tests are done for CP0. Also, AP1, AP2 and AP2W give the same results since the velocity is the same for both. We will again list the methods used in these test:
\begin{itemize}
\item RKC with PM1, PM1V, AP1, AP2 and CP0,
\item ROCK2 with PM1, PM1V, AP1, AP2W and CP0.
\end{itemize}

\subsection{Errors at the boundaries, codes profiling, accuracy improvement by projections}\label{sec:desc_errb}
This experiment has two goals. The first one is to motivate the employment of different boundary conditions for the virtual velocity $\bus$ of PM1. The second goal is to show that projecting the velocity after each stage gives a much higher accuracy being just a bit slower. 

Using the ROCK2 scheme we will compare the virtual and physical velocities $u^*$ and $u$ given by PM1, PM1V and PM3 (method of section \ref{sec:method3}) in the first two time steps when $\Delta t=10^{-1}$. The tests are done without the advection term thus the PIROCK scheme gives the same results of ROCK2 and the ones given by RKC are very similar. After we will integrate until the end of the interval using PM1, PM1V, AP1, AP2W showing the code profile and the velocity and pressure errors. In particular we will see if the time spent in projecting each stage (PM1V, AP1, AP2W) is worthwhile compared to the increased accuracy with respect to PM1. Finally the codes used for this test are
\begin{itemize}
\item ROCK2 with PM1, PM1V, PM3, AP1, AP2W and CP0.
\end{itemize}

\subsection{Numerical efficiency}\label{sec:desc_eff}
These tests compare the different methods in order to find out which one gives the most accurate solution in less time. Following \cite[II.10]{ode-ii} we have compared the different methods with time step adaptivity enabled and different relative and absolute tolerances 
\begin{align*}
rtol = atol = 10^{-m},
\end{align*}
where $m=2,3,\ldots,12$. The tests are done on an Intel\textregistered Core\textsuperscript{\texttrademark} i7-2600S processor. The Fortran code has been compiled with the GNU compiler and the \verb=-Ofast= option. We compare the solutions given by the different tolerances against a reference solution. We used a $128\times 128$ grid, $\Rey=100$, the starting time is $t=0$ and the end time is $t=1$, the starting time step for ROCK2 and PIROCK is $\Delta t=10^{-3}$, for RKC it is automatically chosen. The reference solution has been computed with the fourth order RK4 method (see \cite[II.1]{ode-i}) using $\Delta t=10^{-7}$ and compensated summation. In the figures the symbols represent the different tolerances.

For RKC we did the test only for PM1 since for the others methods we cannot enable time step adaptivity. For PIROCK also we did the test only for PM1 since the others are not implemented. For ROCK2 all the possible couplings are tested. Moreover this test is done for CP0 and CP1. To summarize the methods used are:
\begin{itemize}
\item RKC with PM1 and CP0, CP1,
\item ROCK2 with PM1, PM1V, AP1, AP2W and CP0, CP1,
\item PIROCK with PM1 and CP0, CP1.
\end{itemize}

\subsection{Different Reynolds numbers behavior}\label{sec:desc_rey}
In this test we are interested in how the different methods behave with different Reynolds numbers. Using a fixed $\Delta x=1/128$ and time step adaptivity we consider different Reynolds numbers. For each Reynolds number we integrate from $t=0$ to $t=1$ and look at the following quantities: the velocity errors against the exact solution, the computational time, the average number of stages per time step, the total number of stages used, the number of time steps and the number of rejected time steps.

PIROCK has stability issues thus this test is done for RKC and ROCK2 only. Again for RKC we use only PM1 because of time step adaptivity. For large Reynolds numbers RKC and ROCK2 showed instabilities arising from the advection term, for this reason it has been neglected. The test is done under the same environment of the efficiency tests (same processor, compiler and options, see \ref{sec:desc_eff}). The methods used are:
\begin{itemize}
\item RKC with PM1 and CP0,
\item ROCK2 with PM1, PM1V, AP1, AP2W and CP0.
\end{itemize}

\subsection{Comparing our results with an established reference}\label{sec:desc_comp}
As a last test we compare our results with a well established reference in order to check the correctness of the implementation. We have compared the solution of the methods that allow time step adaptivity (except CP1) with tabular results given in \cite{ghia} at the stationary point. The methods used are:
\begin{itemize}
\item RKC with PM1 and CP0,
\item ROCK2 with PM1, PM1V, AP1, AP2W and CP0,
\item PIROCK with PM1 and CP0.
\end{itemize}

In the following sections we consider three test problems and show the results of the numerical experiments above described.

\section{The forced flow} \label{sec:forced}
The forced flow problem has been taken from \cite{rkcp}. It has the following exact solution for the incompressible Navier-Stokes equations:
\begin{align} \label{eq:ff}
\begin{aligned}
u(t,x,y) &= -\cos(t)\sin(\pi x)^2\sin(2\pi y), \\
v(t,x,y) &= \cos(t)\sin(2\pi x)\sin(\pi y) ,\\
p(t,x,y) &= -\frac{\sin(t)}{4}\left(2+\cos(\pi x)\right) \left( 2+\cos(\pi y)\right)  \\
&\phantom{=} \phantom{1}+\frac{\pi^2}{2} \cos(t)\left( \cos(\pi x) + \cos(\pi y)+\cos(\pi x)\cos(\pi y)\right) ,
\end{aligned}
\end{align}
with the appropriated forcing terms added to \eqref{eq:ns} to ensure that \eqref{eq:ff} is the exact solution. We have solved the problem in the domain $\Omega=[0,1]\times[0,1]$ with homogeneous Dirichlet boundary conditions for $u$ and $v$.

\subsection{Convergence order}\label{sec:conv_ff}
\begin{figure}[!hbtp]
\begin{center}
\subfigure[All the PM1 and PM1V methods.]{
\includegraphics[trim=0.0cm 0.0cm 0.0cm 0.0cm, clip=true, height=0.2\textheight]{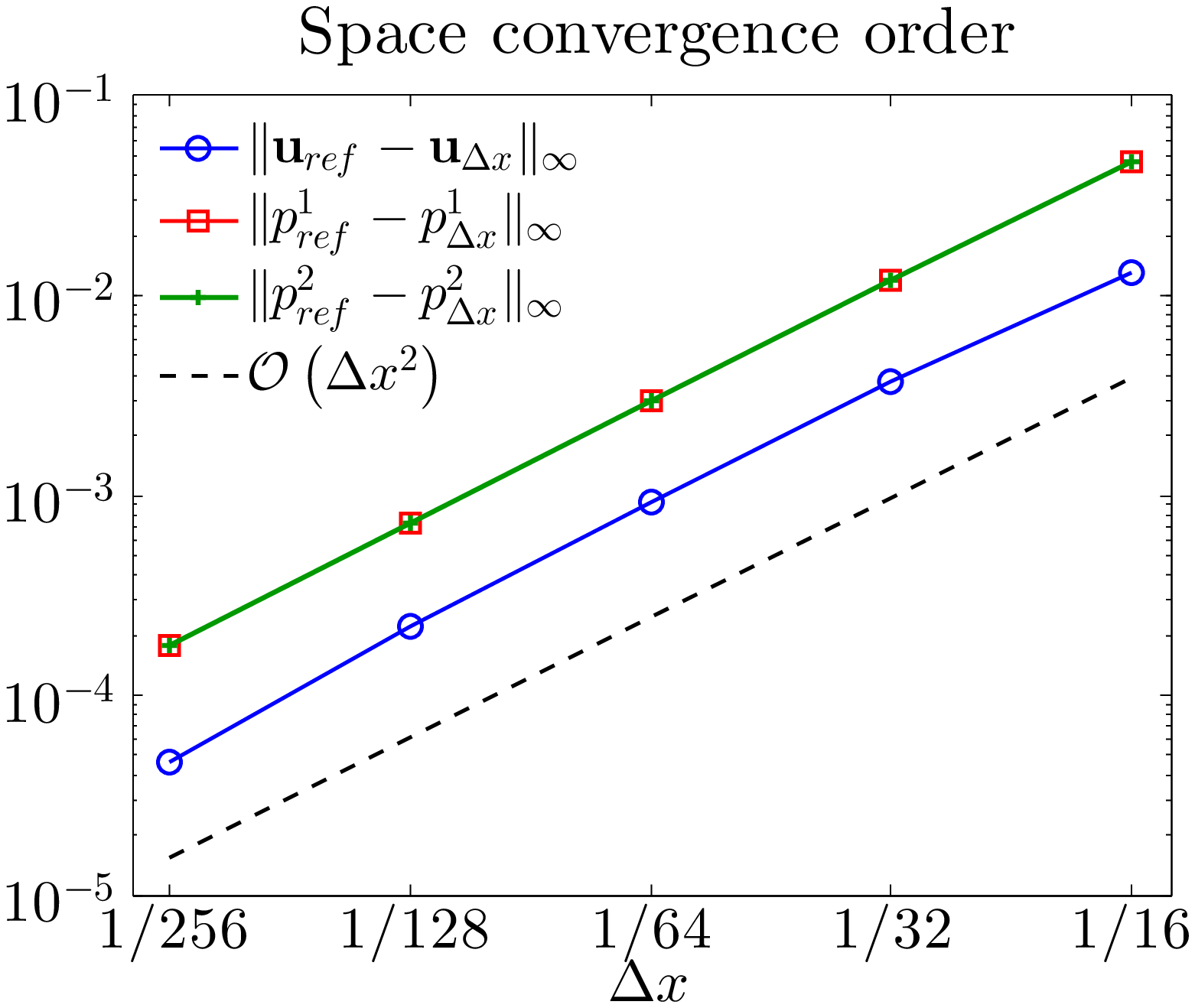} \label{fig:space_Proj}}
\subfigure[All the AP1, AP2 and AP2W methods.]{
\includegraphics[trim=0.0cm 0.0cm 0.0cm 0.0cm, clip=true, height=0.2\textheight]{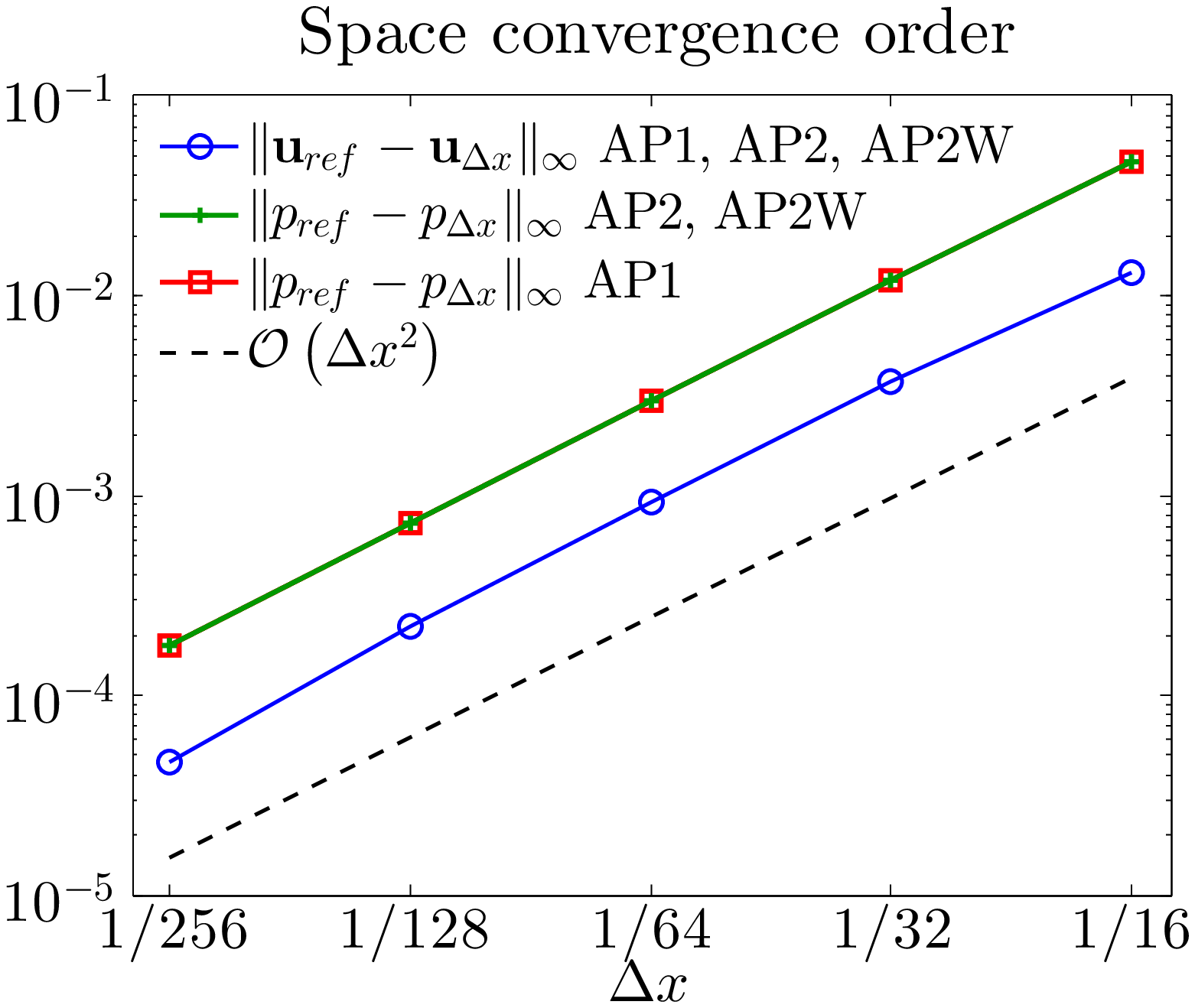} \label{fig:space_ff_DAE}}
\end{center}
\caption{Space convergence results of the forced flow problem.}
\label{fig:space_ff}
\end{figure}

\begin{figure}[!hbtp]
\begin{center}
\subfigure[PIROCK, PM1, CP0.]{
\includegraphics[trim=0.0cm 0.0cm 0.0cm 0.0cm, clip=true, height=0.2\textheight]{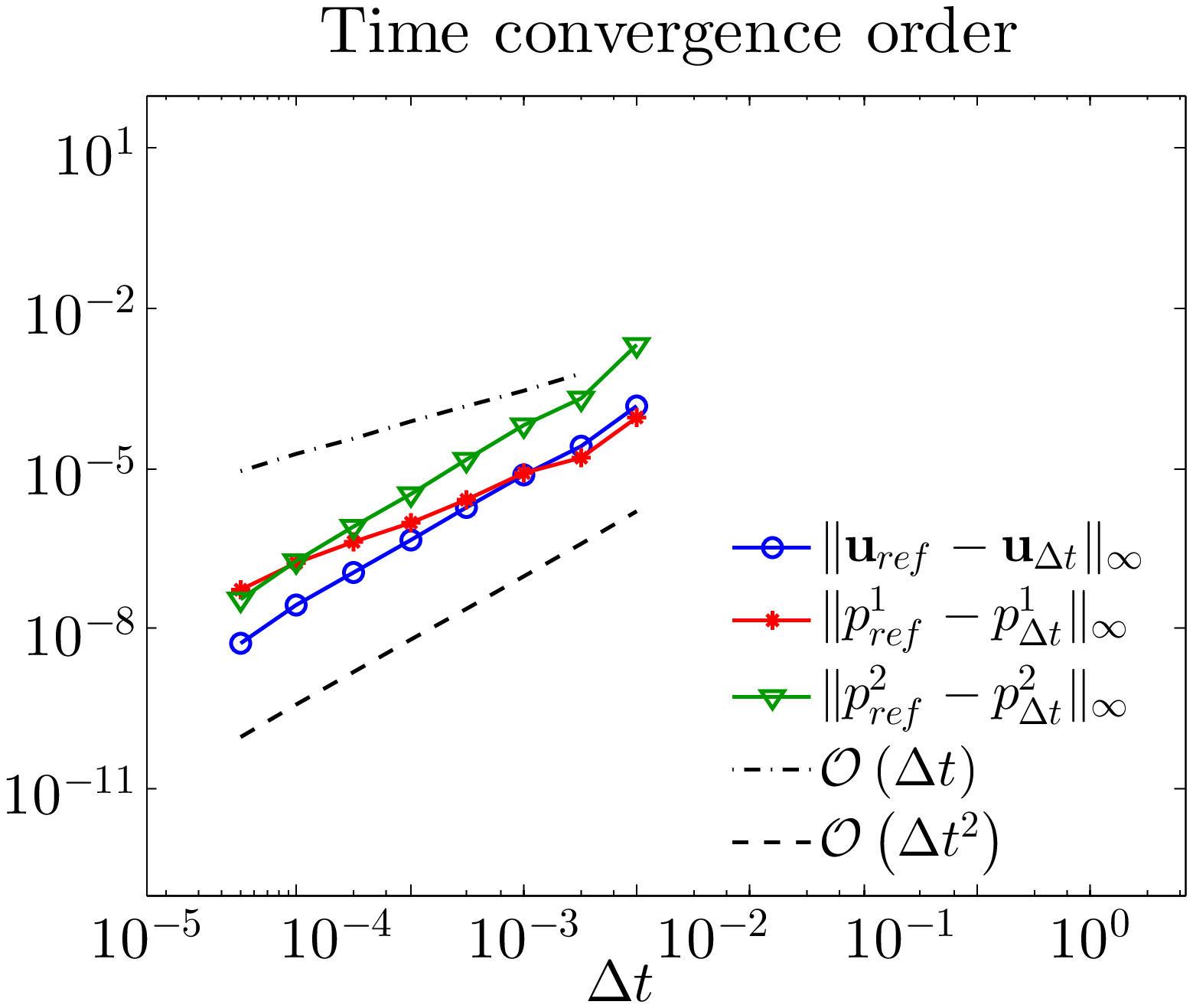} \label{fig:time_ff_pirock_PM1}} \\
\subfigure[RKC, PM1, CP0.]{
\includegraphics[trim=0.0cm 0.0cm 0.0cm 0.0cm, clip=true, height=0.2\textheight]{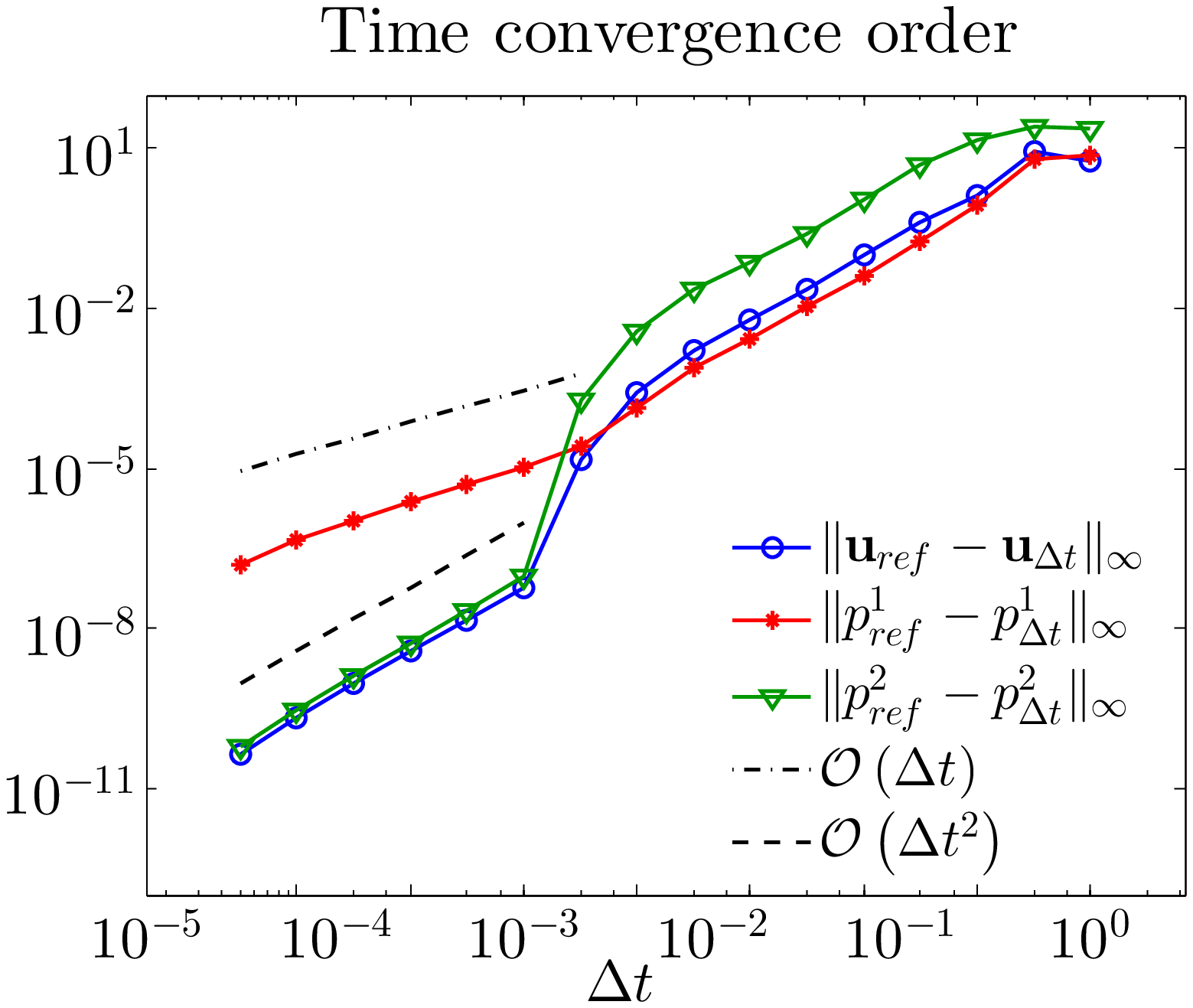} \label{fig:time_ff_rkc_PM1}}
\subfigure[ROCK2, PM1, CP0.]{
\includegraphics[trim=0.0cm 0.0cm 0.0cm 0.0cm, clip=true, height=0.2\textheight]{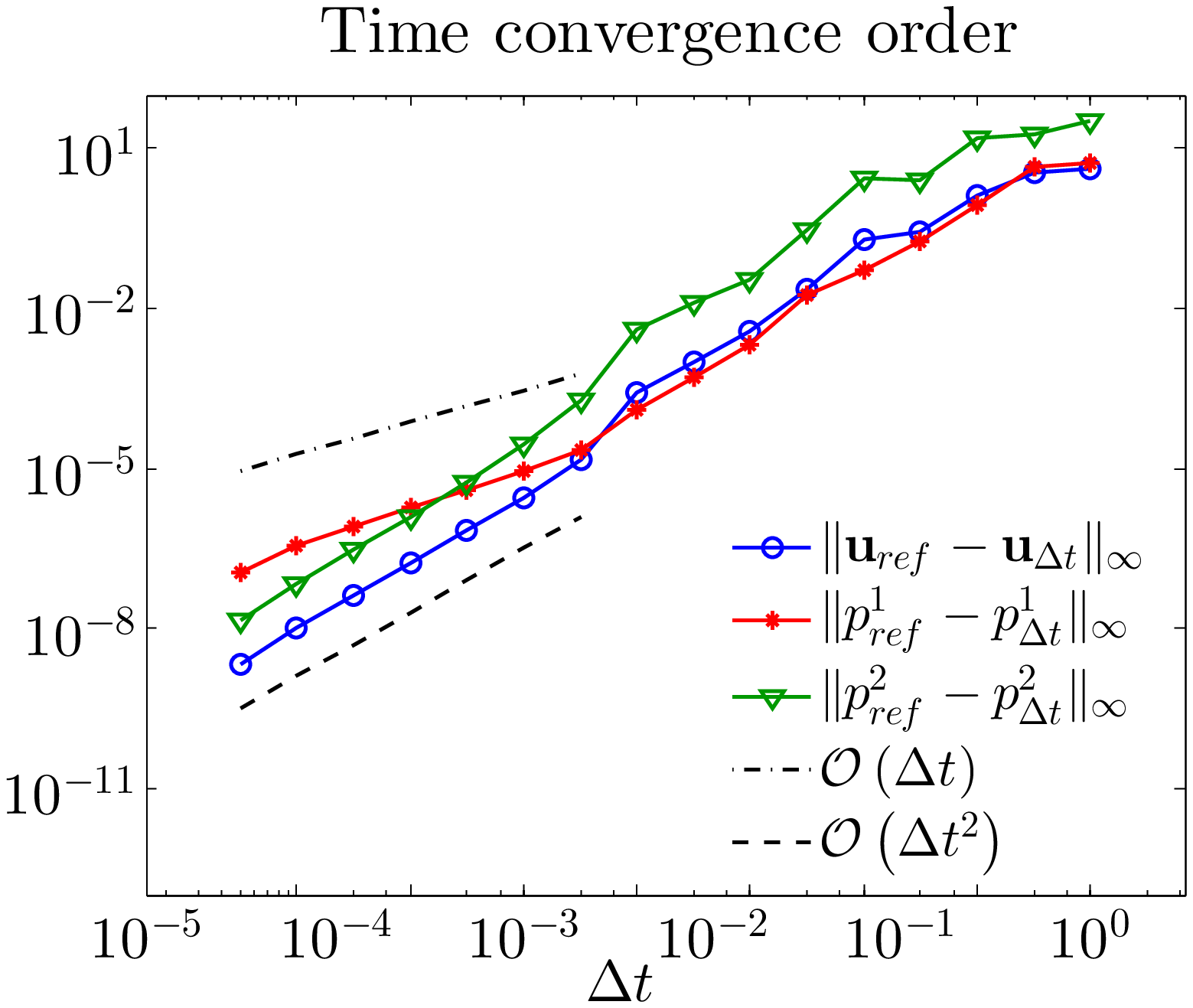} \label{fig:time_ff_rock2_PM1}} \\
\subfigure[RKC, PM1V, CP0.]{
\includegraphics[trim=0.0cm 0.0cm 0.0cm 0.0cm, clip=true, height=0.2\textheight]{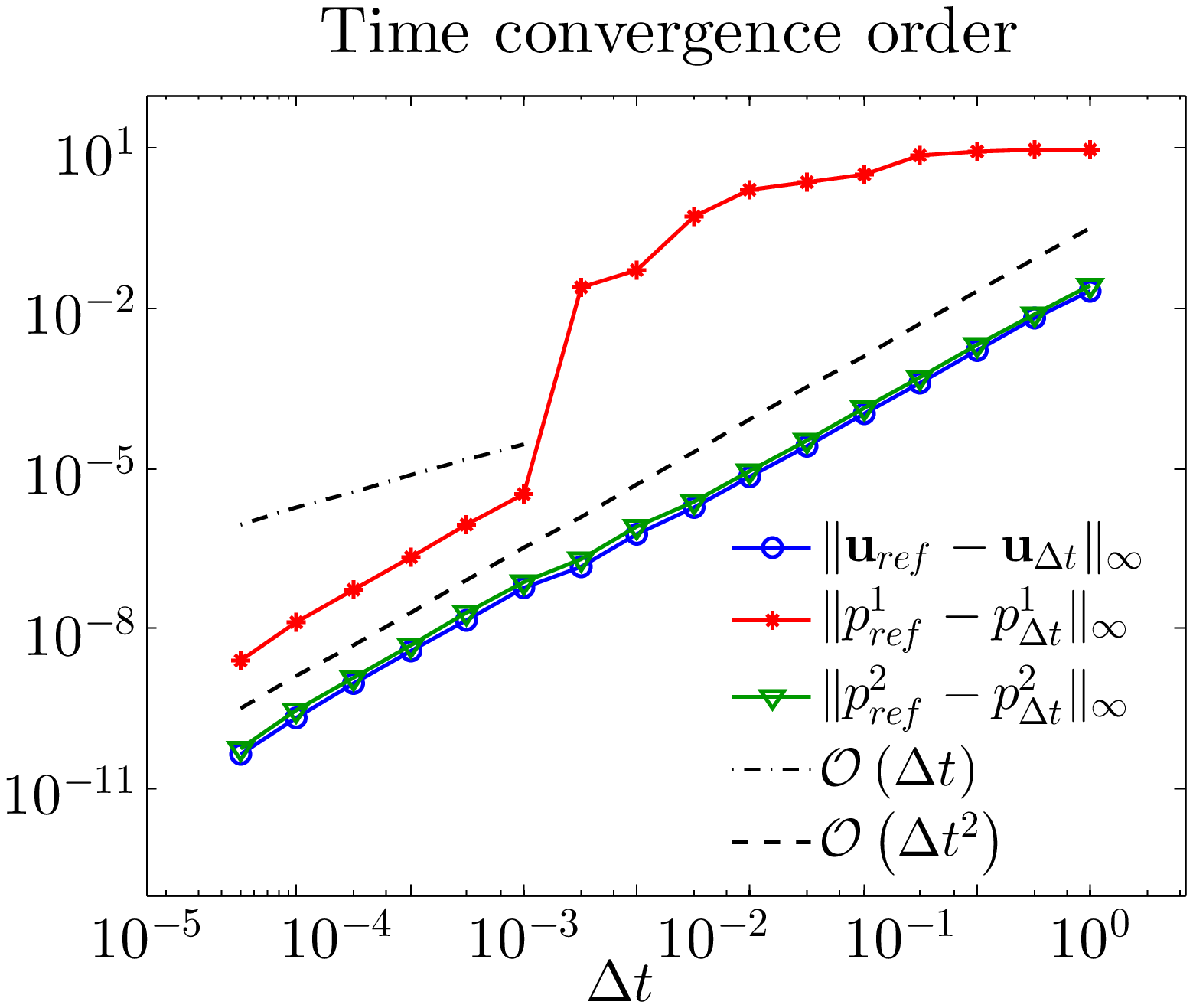} \label{fig:time_ff_rkc_PM1V}}
\subfigure[ROCK2, PM1V, CP0.]{
\includegraphics[trim=0.0cm 0.0cm 0.0cm 0.0cm, clip=true, height=0.2\textheight]{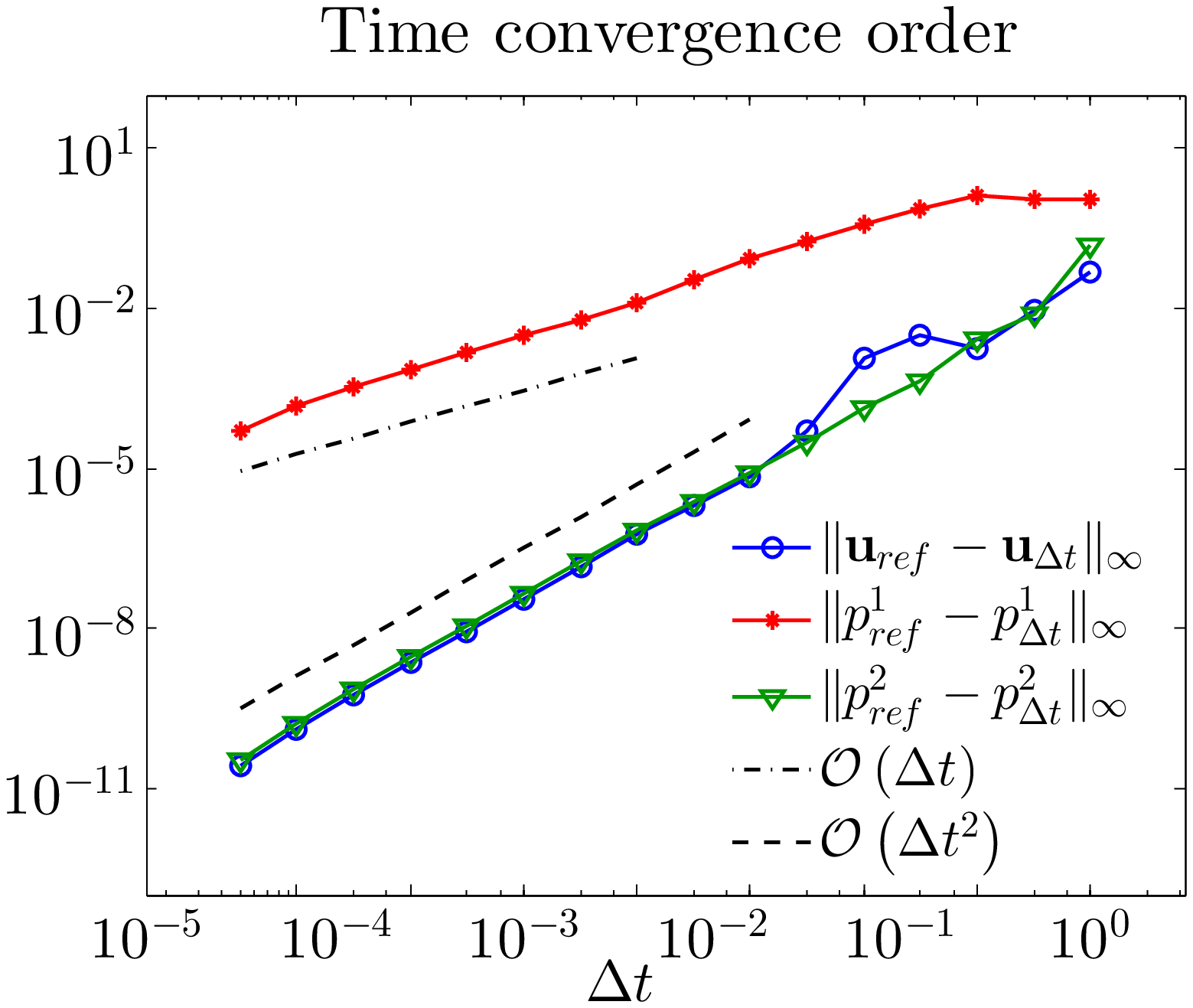} \label{fig:time_ff_rock2_PM1V}} \\
\subfigure[RKC, AP2 and AP2, CP0.]{
\includegraphics[trim=0.0cm 0.0cm 0.0cm 0.0cm, clip=true, height=0.2\textheight]{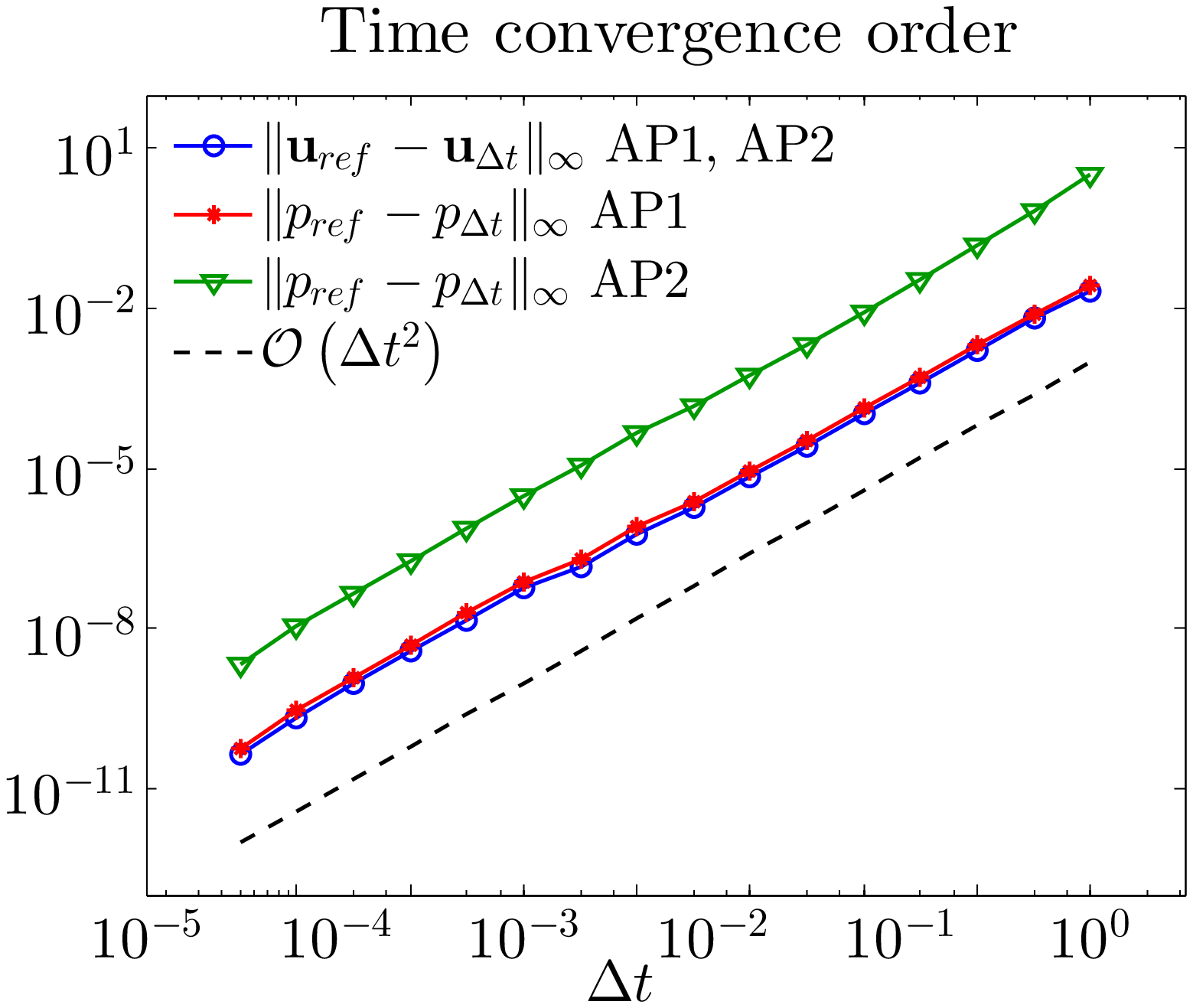} \label{fig:time_ff_rkc_AP12}}
\subfigure[ROCK2, AP1 and AP2W, CP0.]{
\includegraphics[trim=0.0cm 0.0cm 0.0cm 0.0cm, clip=true, height=0.2\textheight]{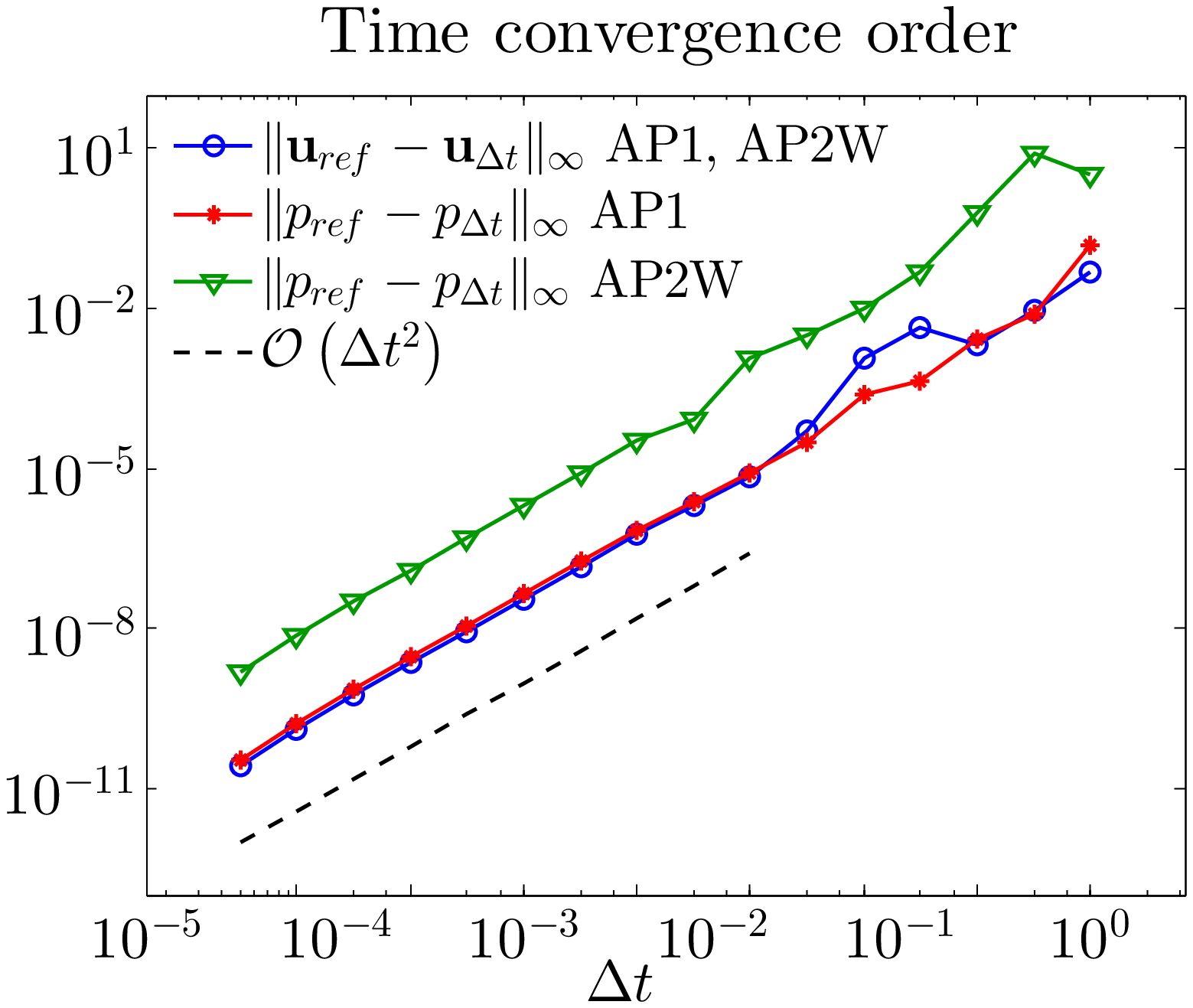} \label{fig:time_ff_rock2_AP12}}
\end{center}
\caption{Time convergence results of the forced flow problem.}
\label{fig:time_ff}
\end{figure}

In Figures \ref{fig:space_ff} and \ref{fig:time_ff} we see the convergence results of the methods for the forced flow problem. The results of the spatial convergence are visually the same for all the methods and are summarized in Figures \ref{fig:space_ff}(a,b), we see that second order space convergence is achieved for all the methods for both the velocity and the pressure. Figures \ref{fig:time_ff} show the results of the time convergence tests. In Figure \ref{fig:time_ff_pirock_PM1} we see that PIROCK, PM1, CP0 achieve the expected time order of convergence for all the quantities but we do not have results for $\Delta t \geq 10^{-2}$, this is because the method is not stable. In Figures \ref{fig:time_ff}(b-g) we show the results of RKC and ROCK2 for PM1, PM1V, AP1, AP2 and AP2W. We observe that in general the methods that project the velocity after each stage (PM1V, AP1, AP2, AP2W) have more regular results. Looking at Figures \ref{fig:time_ff}(f, g) we remark that the pressure computed with AP1 converges faster than the one computed with AP2 and AP2W. For RKC, PM1V, CP0 the pressure $p^1$ has second order of convergence, even if only first order was expected.

\clearpage

\subsection{Stability tests}\label{sec:stab_ff}
\begin{figure}[!hbtp]
\begin{center}
\subfigure[RKC with PM1, PM1V, AP1, AP2 and CP0.]{
\includegraphics[trim=0.0cm 0.0cm 0.0cm 0.0cm, clip=true, height=0.2\textheight]{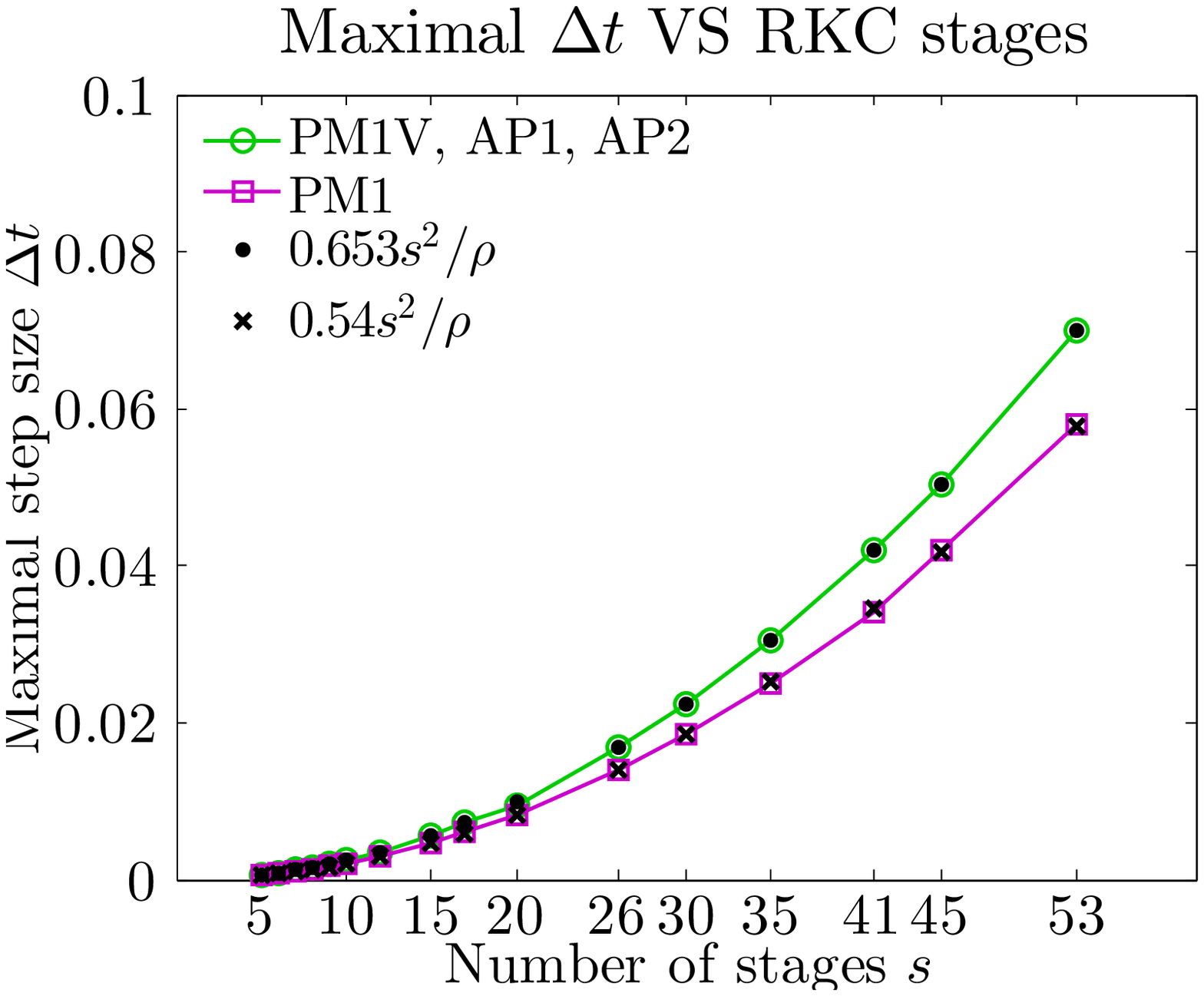} \label{fig:stab_ff_rkc_s}}
\subfigure[ROCK2 with PM1, PM1V, AP1, AP2W and CP0.]{
\includegraphics[trim=0.0cm 0.0cm 0.0cm 0.0cm, clip=true, height=0.2\textheight]{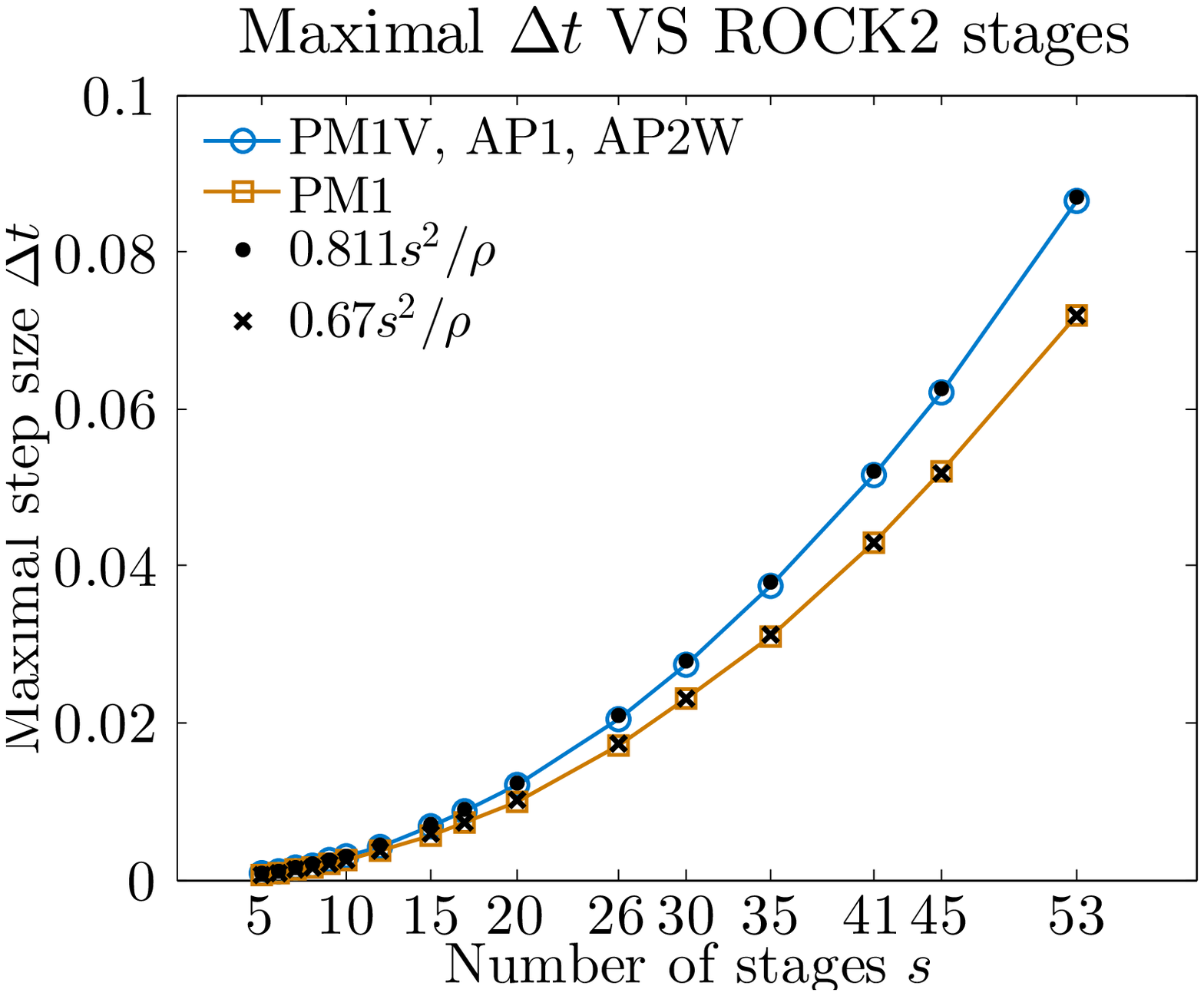} \label{fig:stab_ff_rock2_s}} \\
\subfigure[RKC with PM1, PM1V, AP1, AP2 and CP0.]{
\includegraphics[trim=0.0cm 0.0cm 0.0cm 0.0cm, clip=true, height=0.2\textheight]{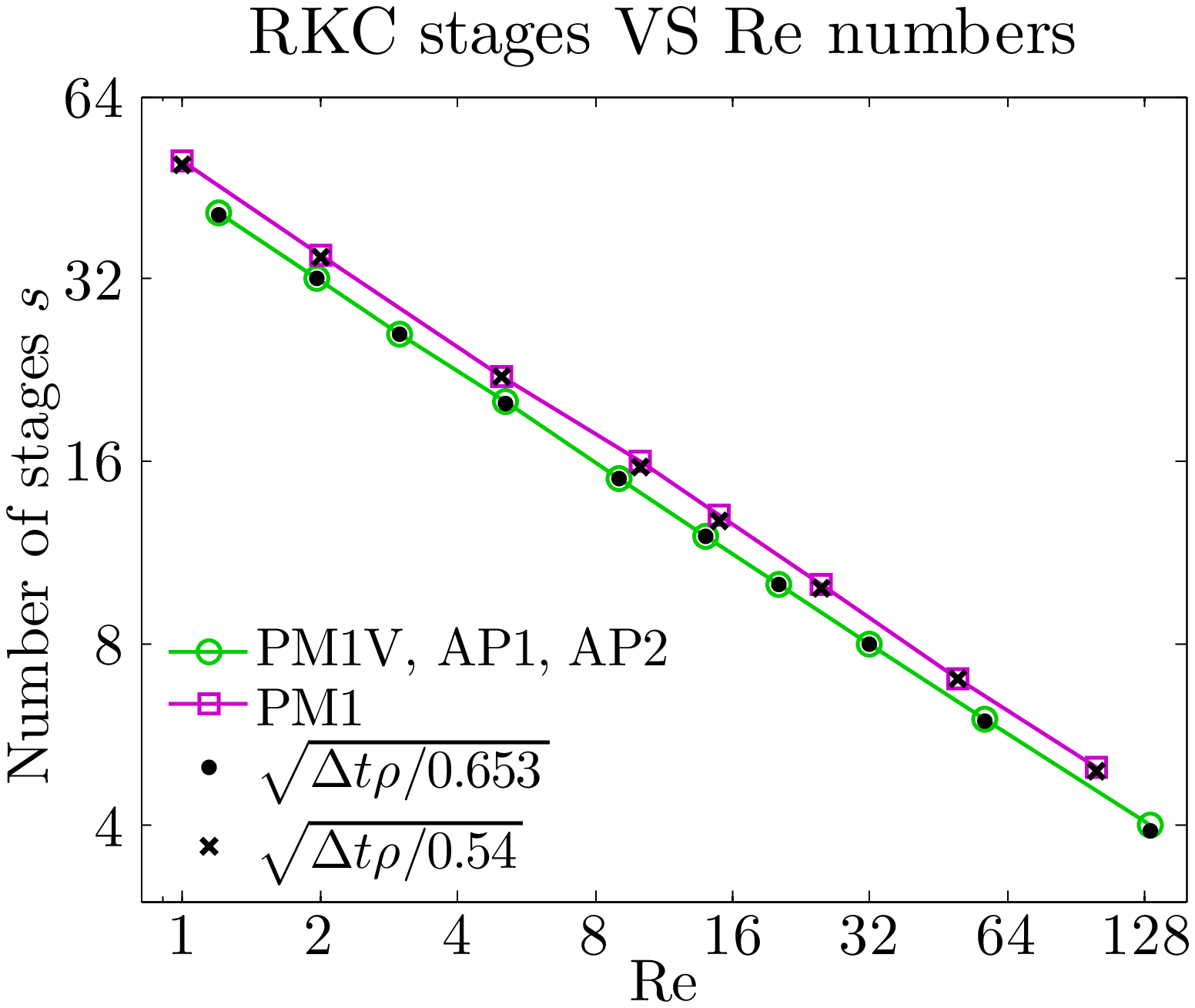} \label{fig:stab_ff_rkc_re}}
\subfigure[ROCK2 with PM1, PM1V, AP1, AP2W and CP0.]{
\includegraphics[trim=0.0cm 0.0cm 0.0cm 0.0cm, clip=true, height=0.2\textheight]{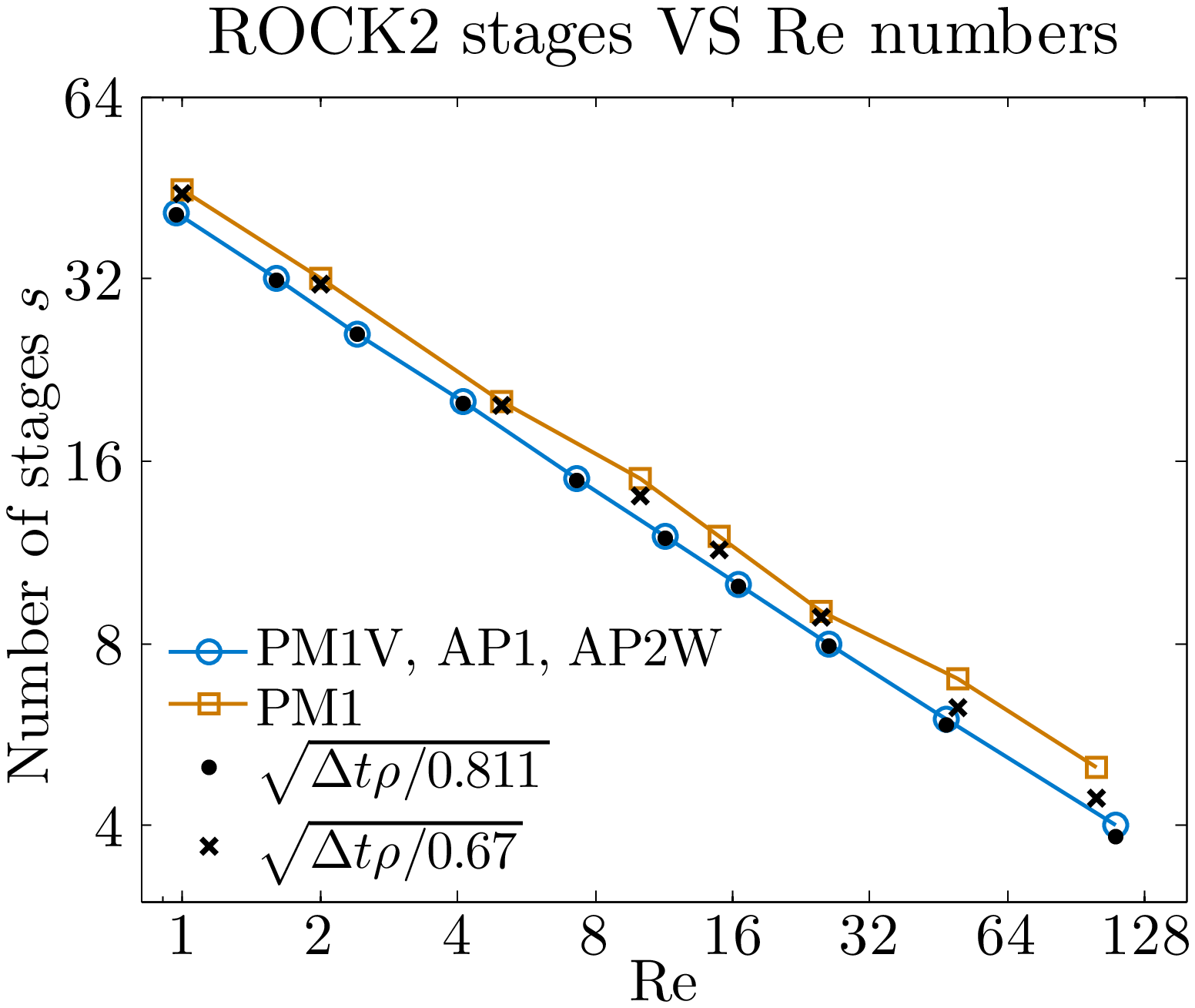} \label{fig:stab_ff_rock2_re}}
\end{center}
\caption{Stability test of the forced flow problem.}
\label{fig:stab_ff}
\end{figure}

In Figure \ref{fig:stab_ff} we present the results of the stability test done on the forced flow problem. From Figures \ref{fig:stab_ff}(a,b) we see that the methods that project the velocity after each stage (PM1V, AP1, AP2 and AP2W) maintain the same stability domain of the ODEs. On the other hand the stability domain grows only as $0.54 s^2$ and $0.67 s^2$ for RKC and ROCK2 respectively when we use PM1, so when we do just one projection per time step. Figures \ref{fig:stab_ff}(c,d) show that these properties does not depend on the Reynolds number.

\subsection{Errors at the boundaries, codes profiling, accuracy improvement by projections}\label{sec:errb_ff}
In this test we show the error of the tangential component of the velocity for different methods, show the codes profile and give a deeper analysis on the errors. 

\begin{figure}[!hbtp]
\begin{center}
\subfigure[$\dx{\phi_1}$ and error $u^*_1-u^e_1$ of PM1, PM1V and PM3.]{
\includegraphics[trim=0.0cm 0.0cm 0.0cm 0.0cm, clip=true, height=0.2\textheight]{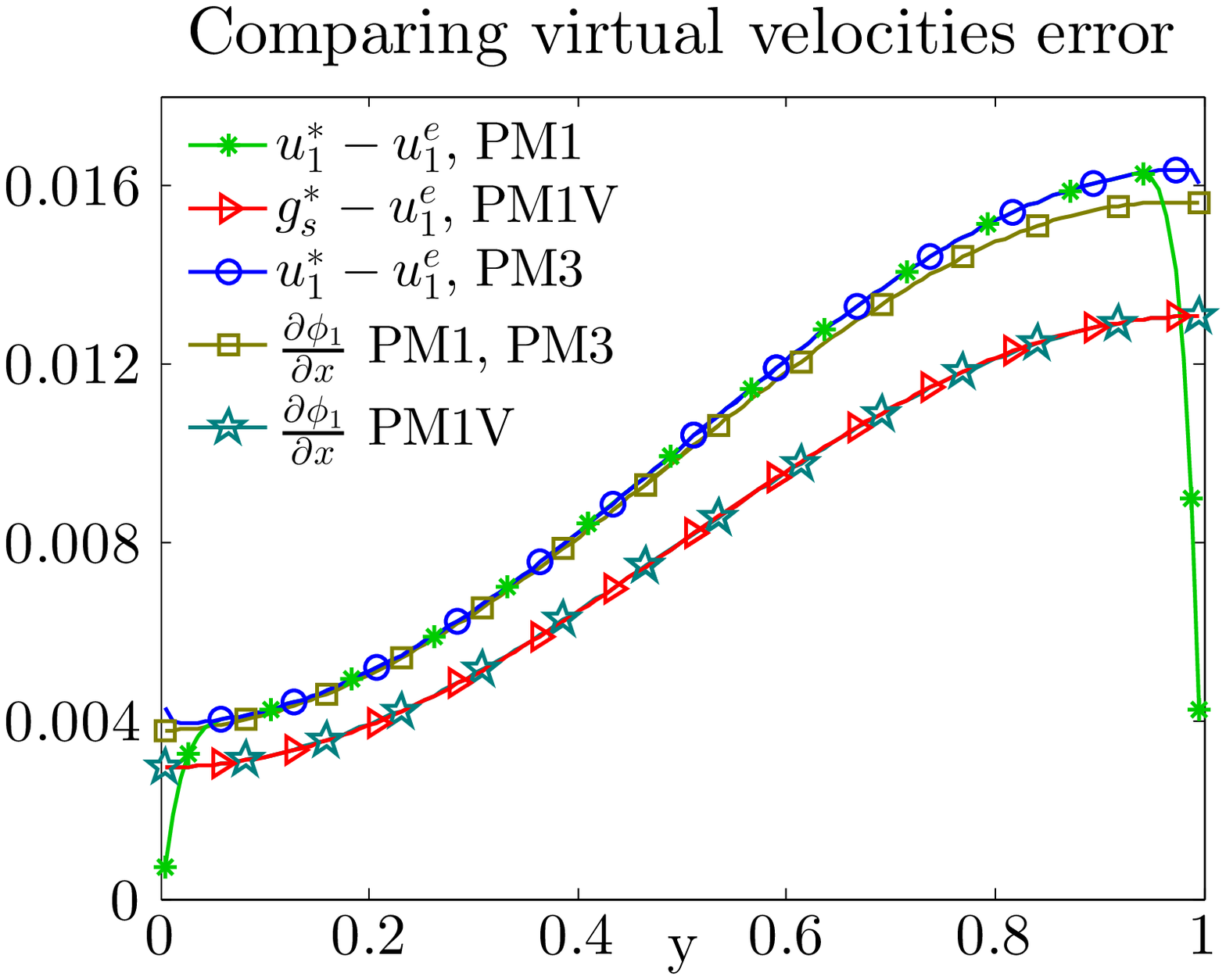} \label{fig:virt_1}}
\subfigure[Error $u_1-u^e_1$ of PM1, PM1V and PM3.]{
\includegraphics[trim=0.0cm 0.0cm 0.0cm 0.0cm, clip=true, height=0.2\textheight]{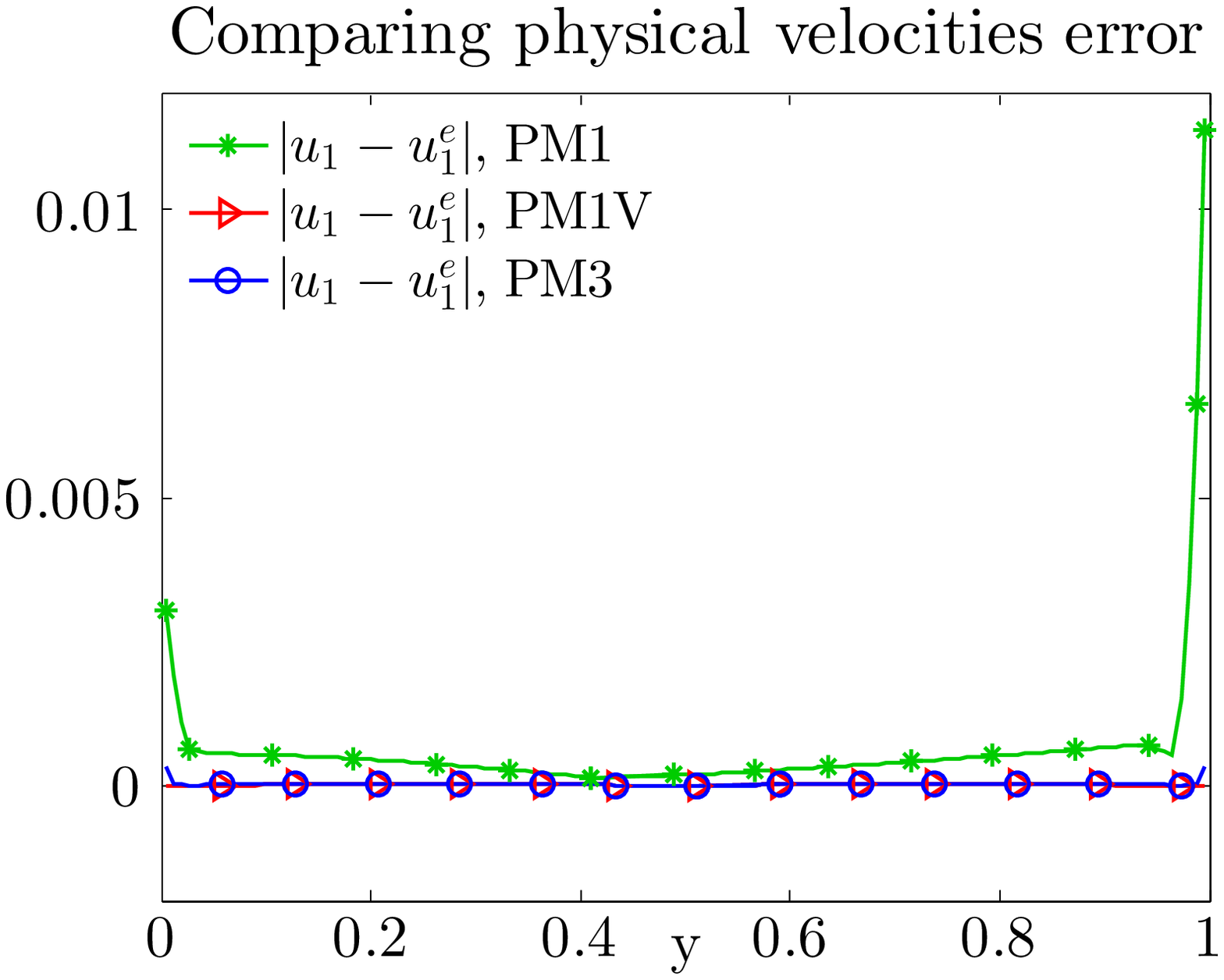} \label{fig:phys_1}} \\
\subfigure[Error $u_1-u^e_1$ of PM1V and PM3.]{
\includegraphics[trim=0.0cm 0.0cm 0.0cm 0.0cm, clip=true, height=0.2\textheight]{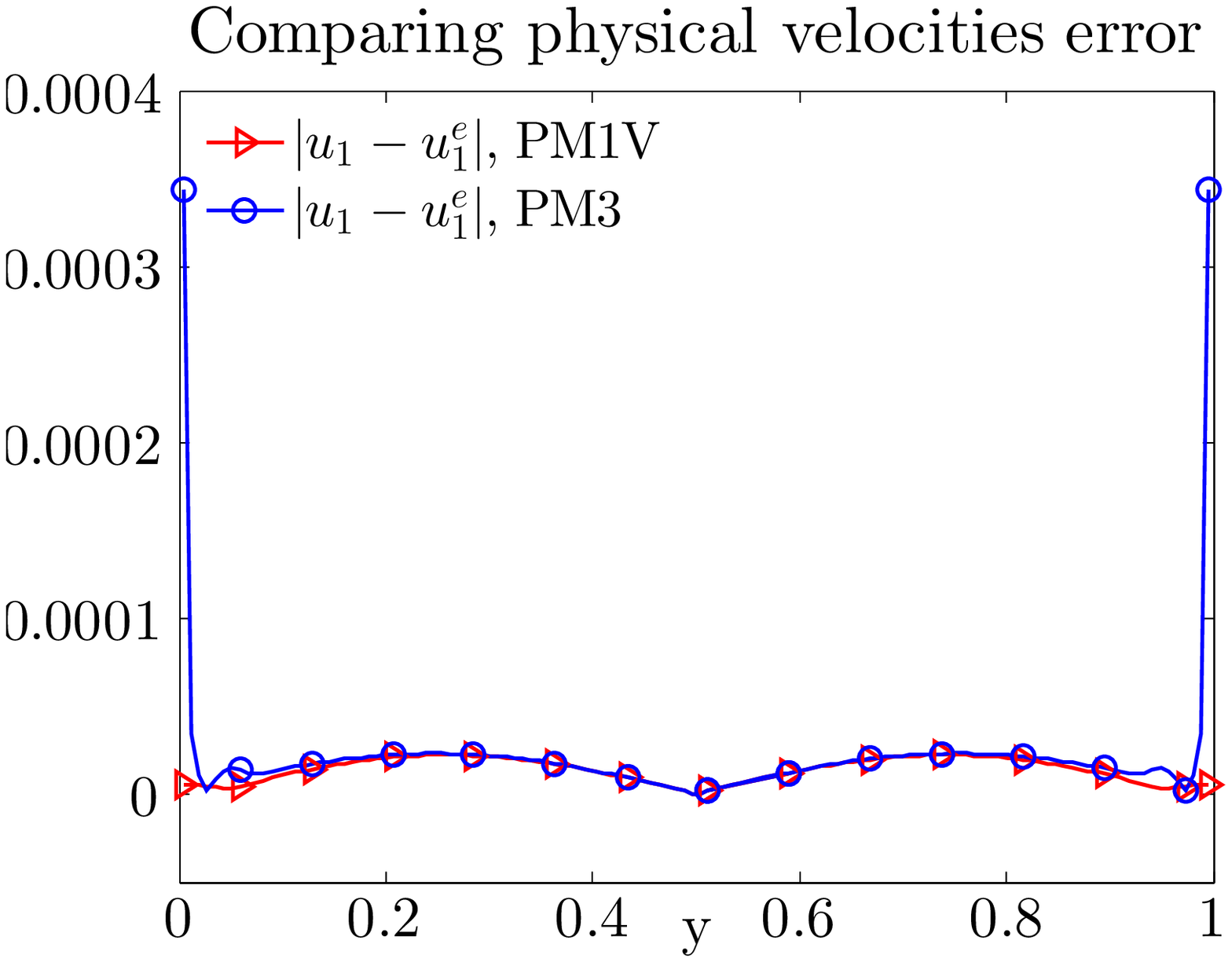} \label{fig:phys_1_s}}
\subfigure[Error $u^*_2-u^e_2$ of PM1, PM1V and PM3.]{
\includegraphics[trim=0.0cm 0.0cm 0.0cm 0.0cm, clip=true, height=0.2\textheight]{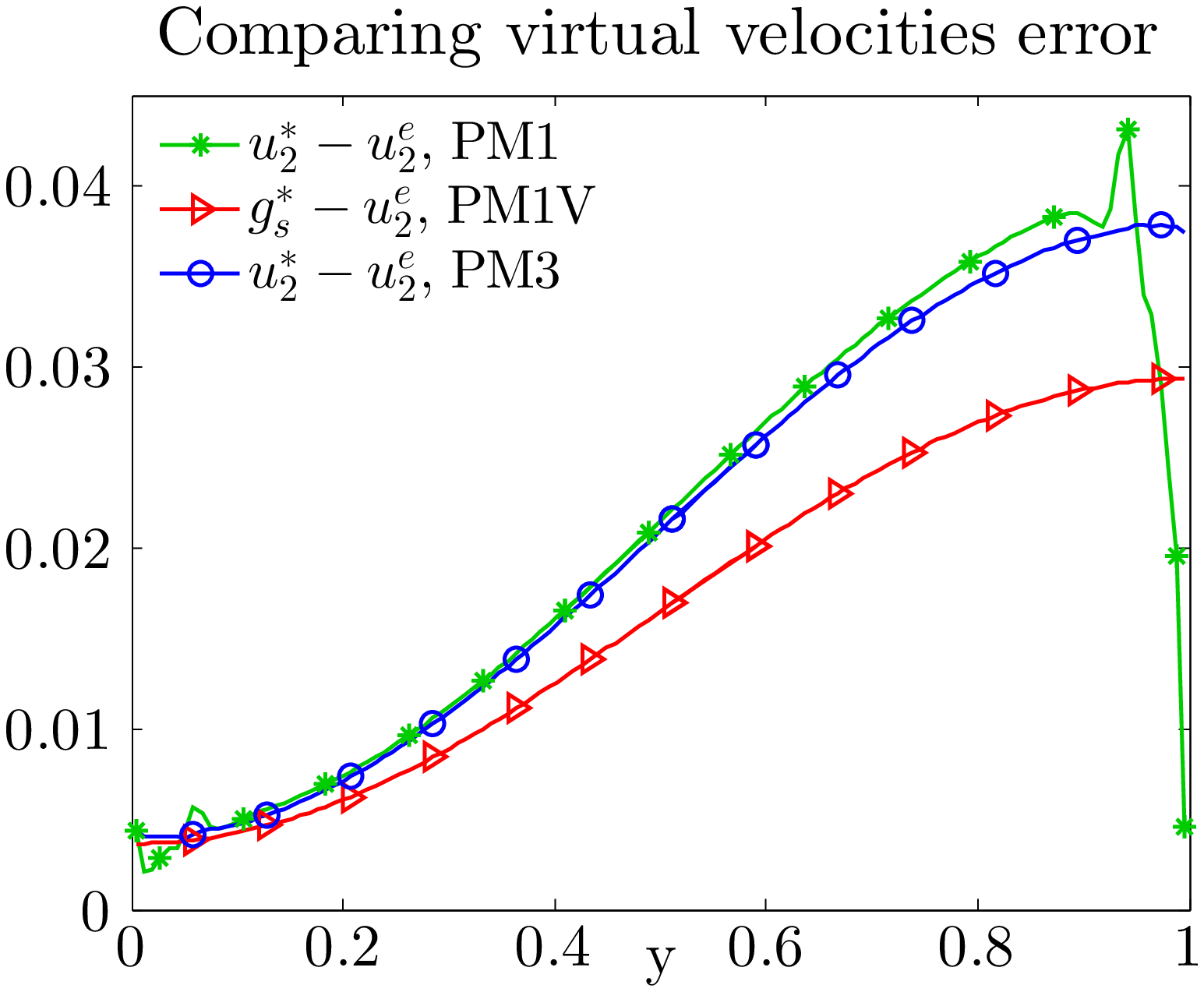} \label{fig:virt_2}}
\end{center}
\caption{Tangential component's error of different methods.}
\label{fig:errb_ff}
\end{figure}

In Figure \ref{fig:errb_ff} $u^*_n$ is the virtual velocity at time step $t_n$ of ROCK2, PM1 and ROCK2, PM3. $g_s^*$ is the last stage of ROCK2, PM1V before being projected (see \eqref{eq:pm1v_algo}). $u_n$ is the physical velocity and $u^e_n$ is the exact velocity given by equations \eqref{eq:ff}. The quantities in Figure \ref{fig:errb_ff} are defined in $[0,1]\times[0,1]$ but here we plot their maximum over the horizontal line, i.e. for $v: [0,1]\times [0,1]\rightarrow \Rb$ we plot $\max_{x\in [0,1]}v(x,y)$ for different $y$s.

In Figure \ref{fig:virt_1} we display the difference between the virtual velocity and the exact solution after one time step for PM1 and PM3. For PM1V we show the difference between the last stage $g_s^*$ before projection and the exact solution. We plot the correction term $\dx{\phi_1}$ as well for each method PM1, PM1V and PM3. For PM1 the imposition of the boundary conditions for $u^*_1$ nullify the error at the boundary $y=0$ and $y=1$. Moreover we see that the boundary conditions affect also the interior of the domain since several stages are computed without corrections. For PM1V the error also vanish in $y=0$ and $y=1$ but since the solution is projected after each stage there is no propagation inside the domain. For PM3 the error at the boundary satisfies $\dy{u^*_1-u^e_1}=0$ which corresponds to equations \eqref{eq:pm3b}. We see that PM1V naturally satisfies equations \eqref{eq:pm3b}, without being imposed (it is the case for PM3). Observe that the correction term $\dx{\phi_1}$ for PM1 and PM3 is much closer to the error $u^*_1-u^e$ of PM3  while for PM1 there is a larger gap between $\dx{\phi_1}$ and $u^*_1-u^e$ near the boundary. For PM1V the error and the correction term have an excellent match. In Figure \ref{fig:phys_1} we see the error $|u_1-u^e_1|$ between the physical velocity $u_1=u^*_1-\dx{\phi_1}$ and the exact solution. We see that the error of the physical velocity for PM1 is much larger at the boundary than the one of PM1V and PM3. 
In order to compare the errors of PM3 and PM1V we plot them in Figure \ref{fig:phys_1_s} using a different scale. We see that the error of PM1V is much smaller than the one of PM3 and that PM3 still have a boundary layer, something that is not present at all for PM1V since the maximal error is inside the domain. Figure \ref{fig:virt_2} shows the errors in the second time step, we see that the large error at the boundary of PM1 has propagated inside the domain in the second time step. The same is not true for PM1V and PM3. 

Analyzing the first two time steps we saw that changing the boundary conditions of PM1 (thus using PM3) can increase the accuracy of the solution. But more important we saw that projecting the velocity after each stage (so using PM1V) creates an excellent match between the last stage $g_s^*$ and the correction term, as a consequence the solution is much more accurate. We did not plot the internal stages $g_i^*$ for $i<s$ but we suppose that they also match the correction term as well as $g_s^*$. We have not done this test for AP1, AP2 and AP2W but we are confident that the results are the same of PM1V, this is confirmed by the rest of this section where we will analyze the overhead of the extra projections in PM1V, AP1, AP2W and see if it is worthwhile compared to the increased accuracy.

In Table \ref{tab:prof} we show part of the PM1, PM1V, AP1, AP2W codes profile when integrating from $t=0$ to $t=1$ again with $\Delta t=10^{-1}$, so doing then time steps. We would like to stress on the fact that the timings given in Table \ref{tab:prof} are adversely affected by the profiler tool (Score-P), in reality the code is much faster (about 4-5 times). What is interesting here is to see the relative difference between the methods, this quantity is not affected by the profiler. The first line of Table \ref{tab:prof} shows that PM1 is about $1.3$ times faster than the other methods. PM1 is the method which spends more time in computing the right hand side $f$, this is because it requires $15$ stages per step while the other methods just $13$ since they have a longer stability domain. On the other hand its time spent in projecting the velocity in negligible. For the other methods it is considerably high, looking at the third line we see that it takes about one third of the total time. The time spent in computing the pressure is negligible for all the methods since we do it in last time step only for PM1, PM1V and AP1. For AP2W it comes for free since no Poisson problems must be solved for the pressure.
\begin{table}[!htbp]
\begin{center}
\begin{tabular}{|l||c|c||c|c||c|c||c|c|}
\hline
 & \multicolumn{2}{|c||}{PM1} & \multicolumn{2}{|c||}{PM1V} & \multicolumn{2}{|c||}{AP1} & \multicolumn{2}{|c|}{AP2W} \\ \hline
 & sec & \% & sec & \% & sec & \% & sec & \% \\ \hline
Total time & 18.82  & 100  & 24.19  & 100  & 24.24  & 100  & 23.98  & 100  \\ \hline
$f$ time & 17.37  & 92.24  & 15.03  & 62.1  & 14.99  & 61.82  & 14.97  & 62.42  \\ \hline
Vel. Proj. time & 0.64  & 3.42  & 8.35  & 34.53  & 8.37  & 34.55  & 8.39  & 34.98 \\ \hline
Pre. Proj. time & 0.22  & 1.16  & 0.22  & 0.89  & 0.22  & 0.91  & 0  & 0  \\ \hline
\end{tabular}
\end{center}
\caption{Profile of ROCK2 with PM1, PM1V, AP1, AP2W codes.}
\label{tab:prof}
\end{table}

In Table \ref{tab:err} we write the errors of the methods when compared against the reference solution computed with the RK4 method (see \ref{sec:desc_eff}) at $t=1$. We see that PM1V, AP1, AP2W give exactly the same velocity errors and PM1V and AP1 give the same pressure error, in fact they compute it in the same way.
\begin{table}[!htbp]
\begin{center}
\begin{tabular}{|l|c|c|c|c|}
\hline
 & PM1 & PM1V & AP1 & AP2W \\ \hline
Velocity error & $1.88 \cdot 10^{-1}$ & $2.92 \cdot 10^{-4}$ & $2.92 \cdot 10^{-4}$ & $2.92 \cdot 10^{-4}$ \\ \hline
Pressure error & $2.22$ & $3.5 \cdot 10^{-1}$ & $3.5 \cdot 10^{-1}$ & $3.19 \cdot 10^{-2}$ \\ \hline
\end{tabular}
\end{center}
\caption{Errors of the PM1, PM1V, AP1, AP2W methods with $\Delta t = 10^{-1}$ at $t=1$.}
\label{tab:err}
\end{table}

Table \ref{tab:rat} shows the ratio between the errors of PM1 and the errors of the other methods. In Figure \ref{fig:errb_ff} we saw that after the first time step PM1V gives an error that is much smaller than the error of PM1. With this table we want to quantify how much it is smaller, taking into account also AP1 and AP2W. We see in the first line that for the velocity PM1V, AP1, AP2W are $644$ times more accurate than PM1. PM1V and AP1 are about $6$ times more accurate than PM1 for the pressure, AP2W is about $70$ times more accurate.
\begin{table}[!htbp]
\begin{center}
\begin{tabular}{|l|c|c|c|}
\hline
 & PM1/PM1V & PM1/AP1 & PM1/AP2W \\ \hline
Velocity error ratio & $644 \times$ & $644 \times$ & $644 \times$ \\ \hline
Pressure error ratio & $6.34 \times$ & $6.34 \times$ & $69.6 \times$ \\ \hline
\end{tabular}
\end{center}
\caption{Ratio between the PM1 and other methods errors.}
\label{tab:rat}
\end{table}

We saw that projecting the velocity after each stage makes the methods $1.3$ times slower but on the other hand they are $644$ times more accurate. Hence the additional projections are absolutely worthwhile.

\subsection{Numerical efficiency}\label{sec:eff_ff}
\begin{figure}[!hbtp]
\begin{center}
\subfigure[Vel. efficiency of PM1, PM1V and CP0.]{
\includegraphics[trim=0.0cm 0.0cm 0.0cm 0.0cm, clip=true, height=0.2\textheight]{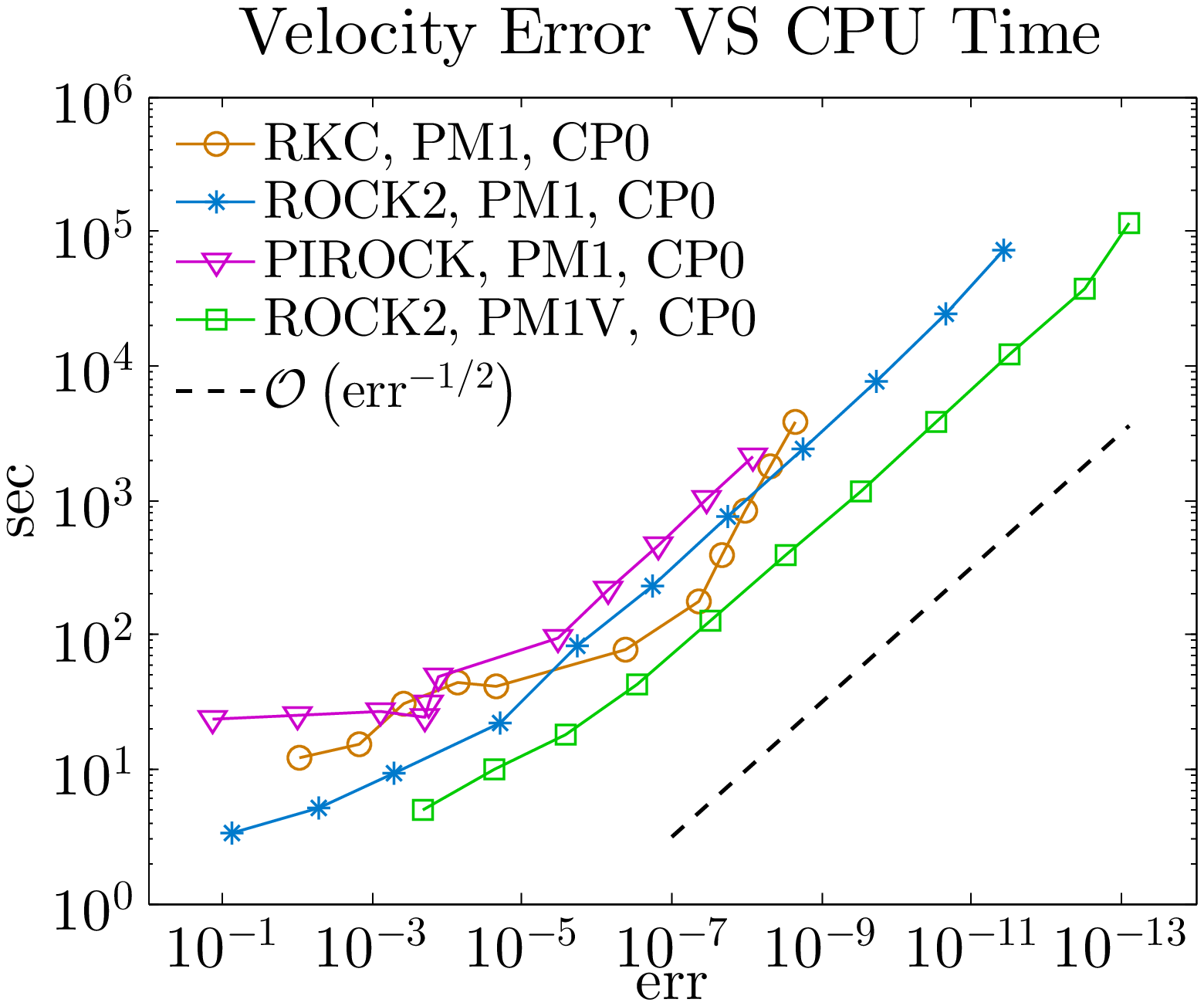} \label{fig:eff_ff_Proj_Vel_CP0}}
\subfigure[Vel. efficiency of PM1, PM1V and CP1.]{
\includegraphics[trim=0.0cm 0.0cm 0.0cm 0.0cm, clip=true, height=0.2\textheight]{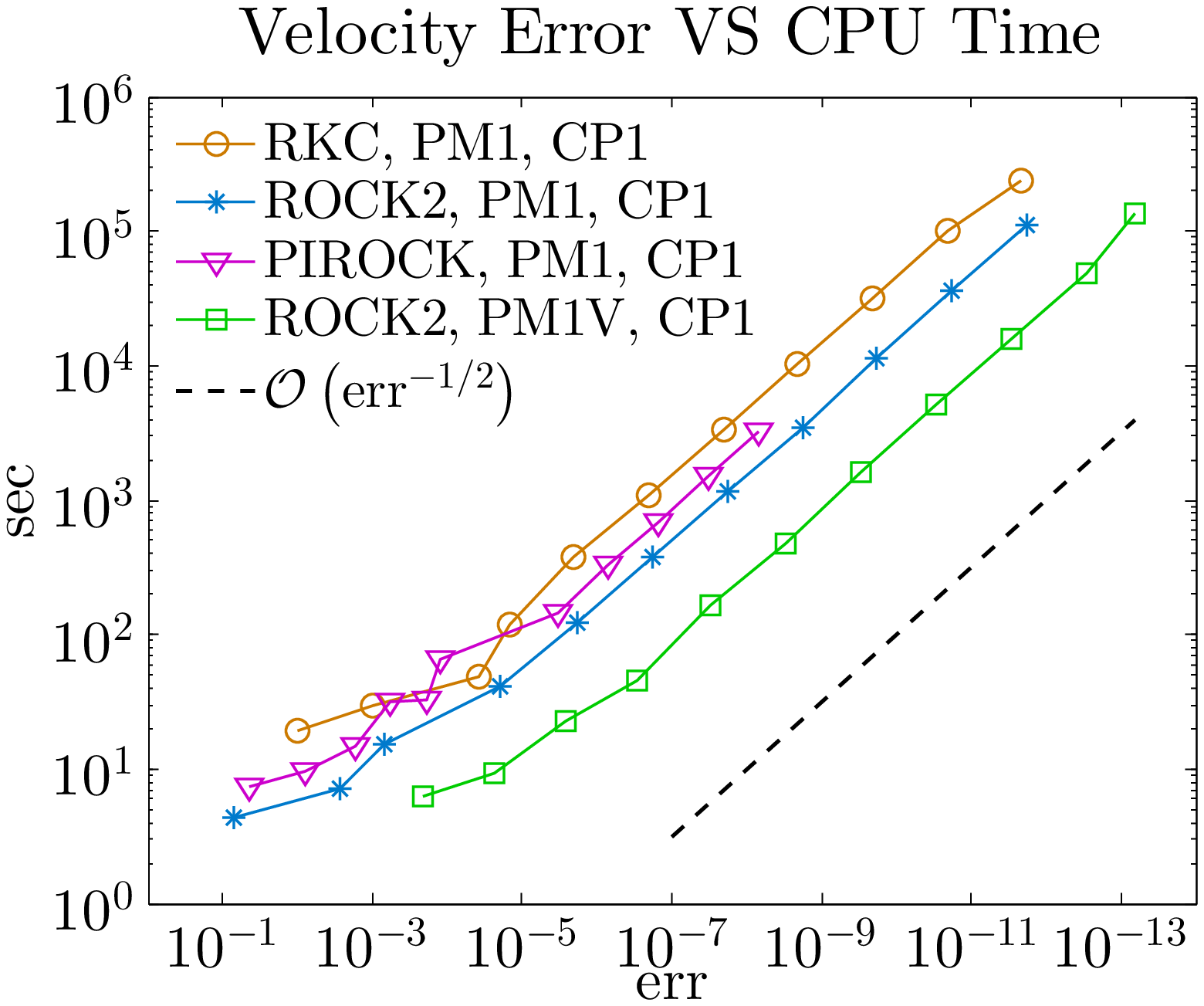} \label{fig:eff_ff_Proj_Vel_CP1}} \\
\subfigure[Pre. efficiency of PM1, PM1V and CP0.]{
\includegraphics[trim=0.0cm 0.0cm 0.0cm 0.0cm, clip=true, height=0.2\textheight]{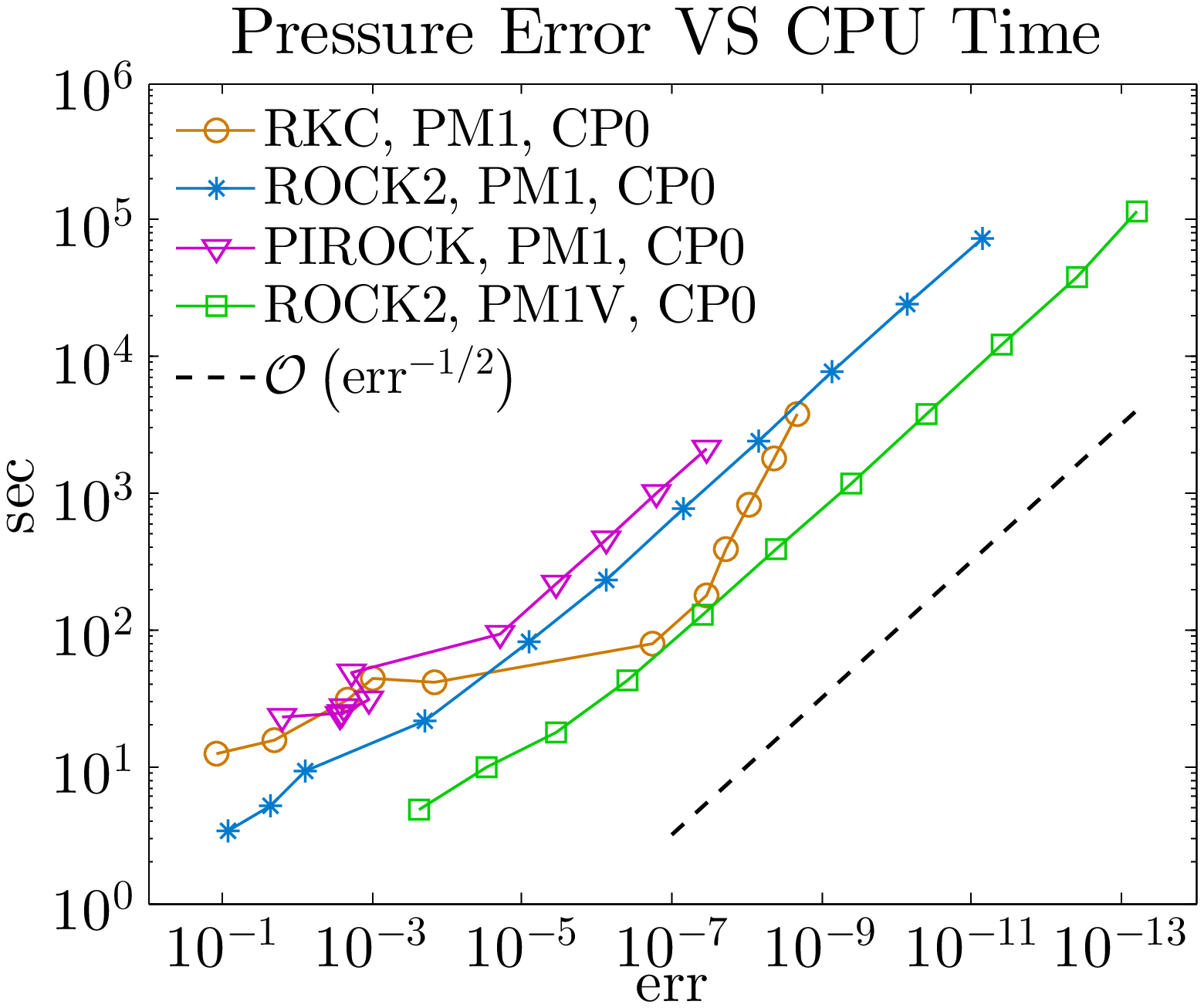} \label{fig:eff_ff_Proj_Pre_CP0}}
\subfigure[Pre. efficiency of PM1, PM1V and CP1.]{
\includegraphics[trim=0.0cm 0.0cm 0.0cm 0.0cm, clip=true, height=0.2\textheight]{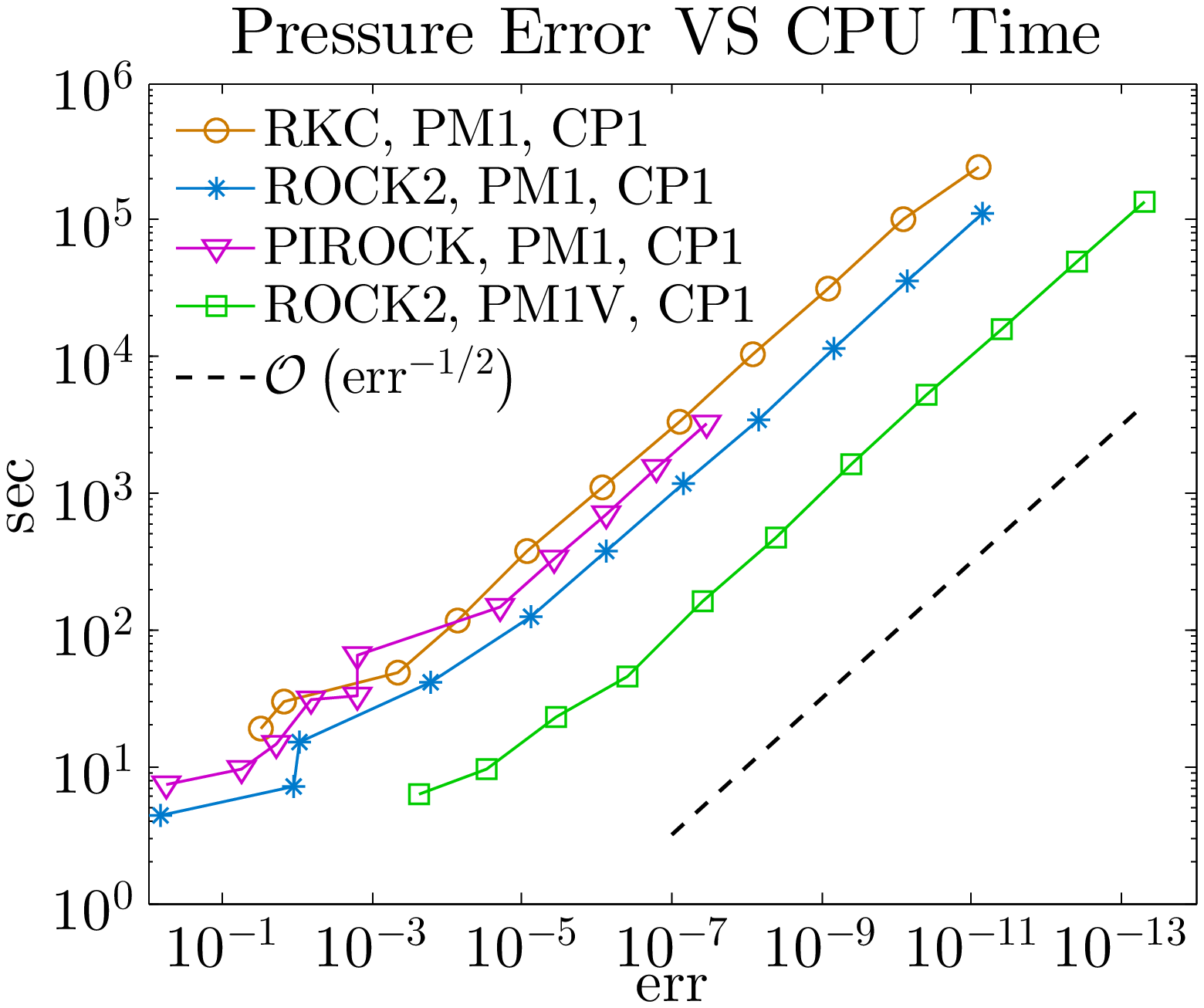} \label{fig:eff_ff_Proj_Pre_CP1}} \\
\subfigure[Vel. efficiency of AP1, AP2W and CP0.]{
\includegraphics[trim=0.0cm 0.0cm 0.0cm 0.0cm, clip=true, height=0.2\textheight]{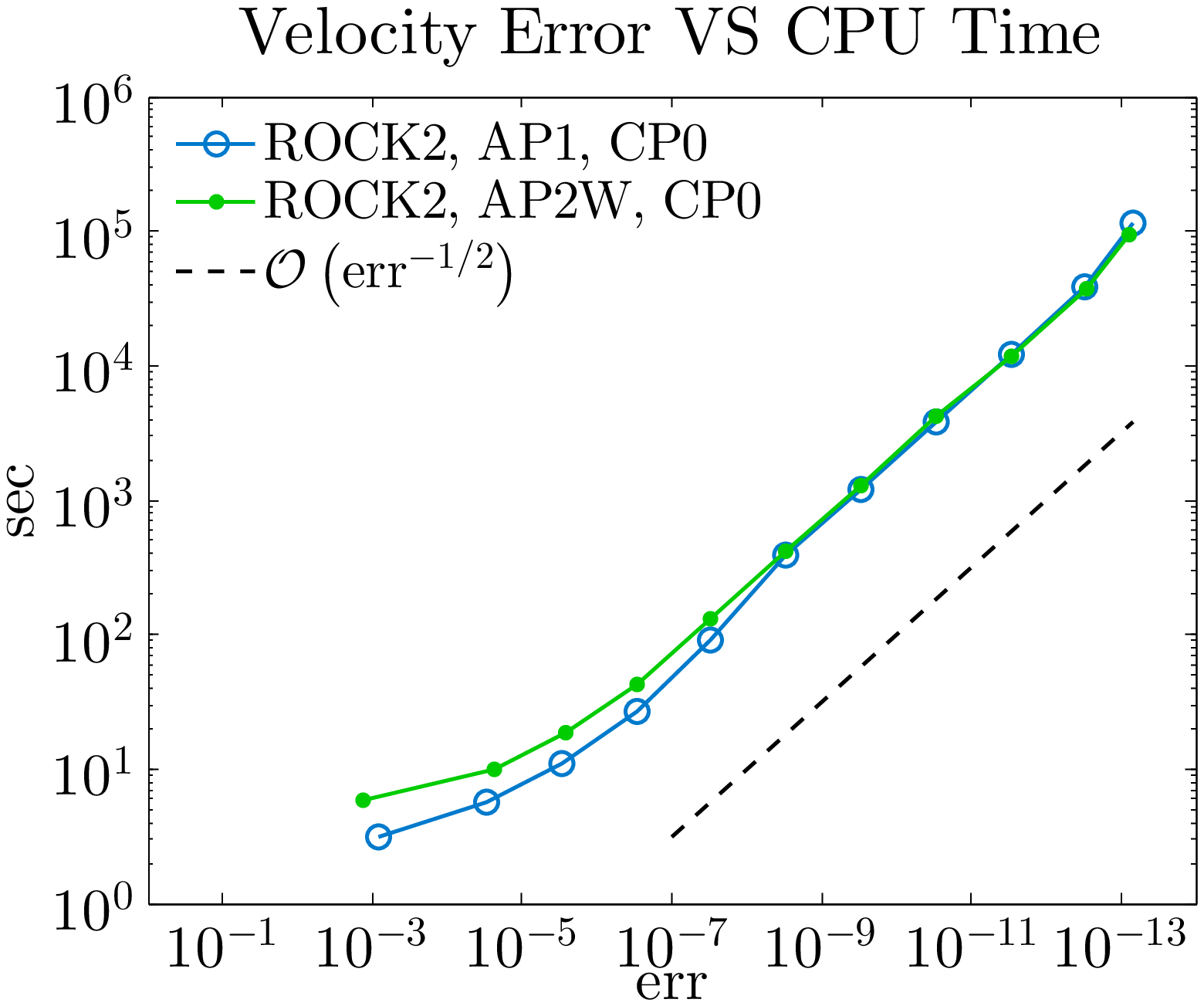} \label{fig:eff_ff_DAE_Vel_CP0}}
\subfigure[Vel. efficiency of AP1, AP2W and CP1.]{
\includegraphics[trim=0.0cm 0.0cm 0.0cm 0.0cm, clip=true, height=0.2\textheight]{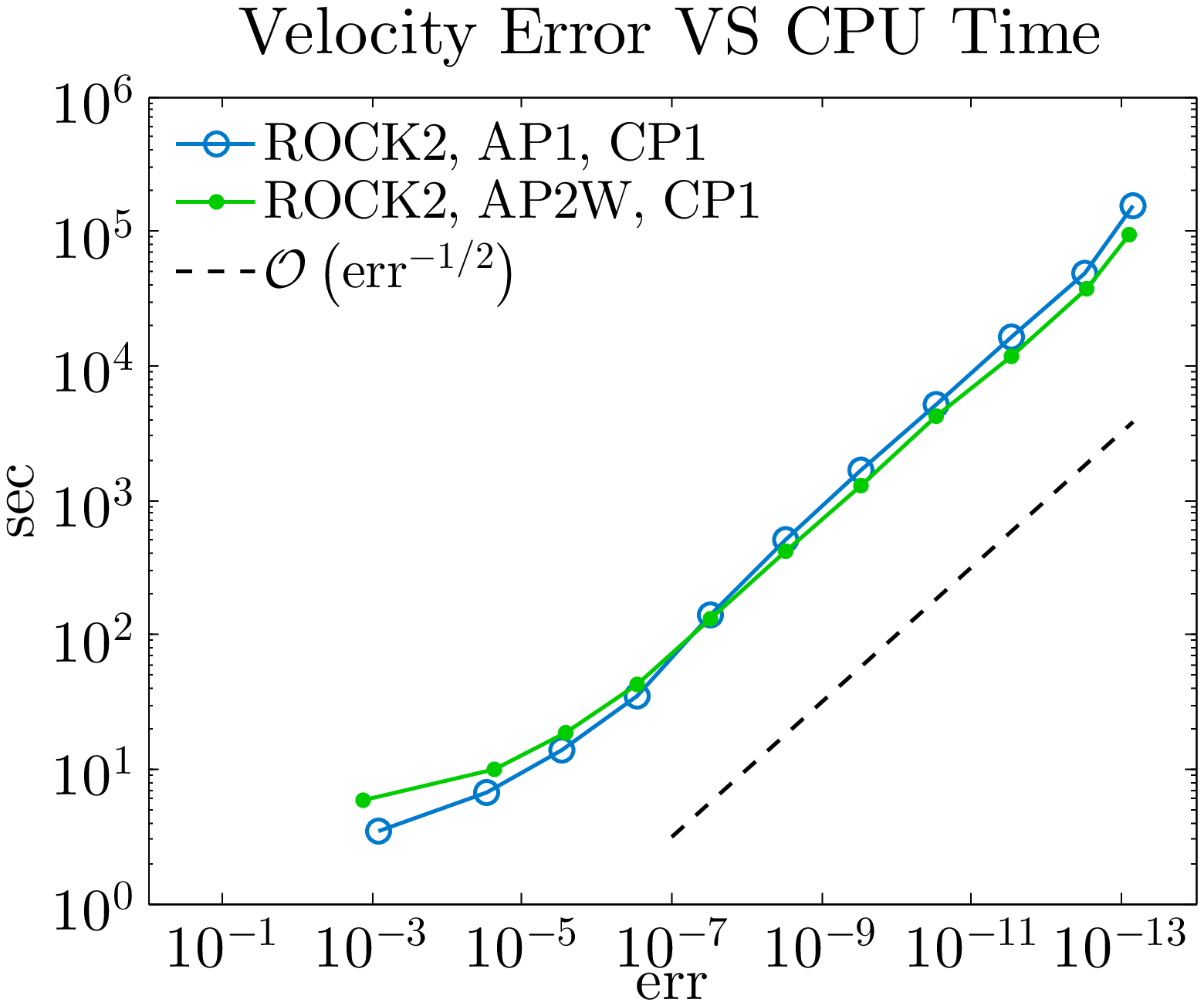} \label{fig:eff_ff_DAE_Vel_CP1}}\\
\subfigure[Pre. efficiency of AP1, AP2W and CP0.]{
\includegraphics[trim=0.0cm 0.0cm 0.0cm 0.0cm, clip=true, height=0.2\textheight]{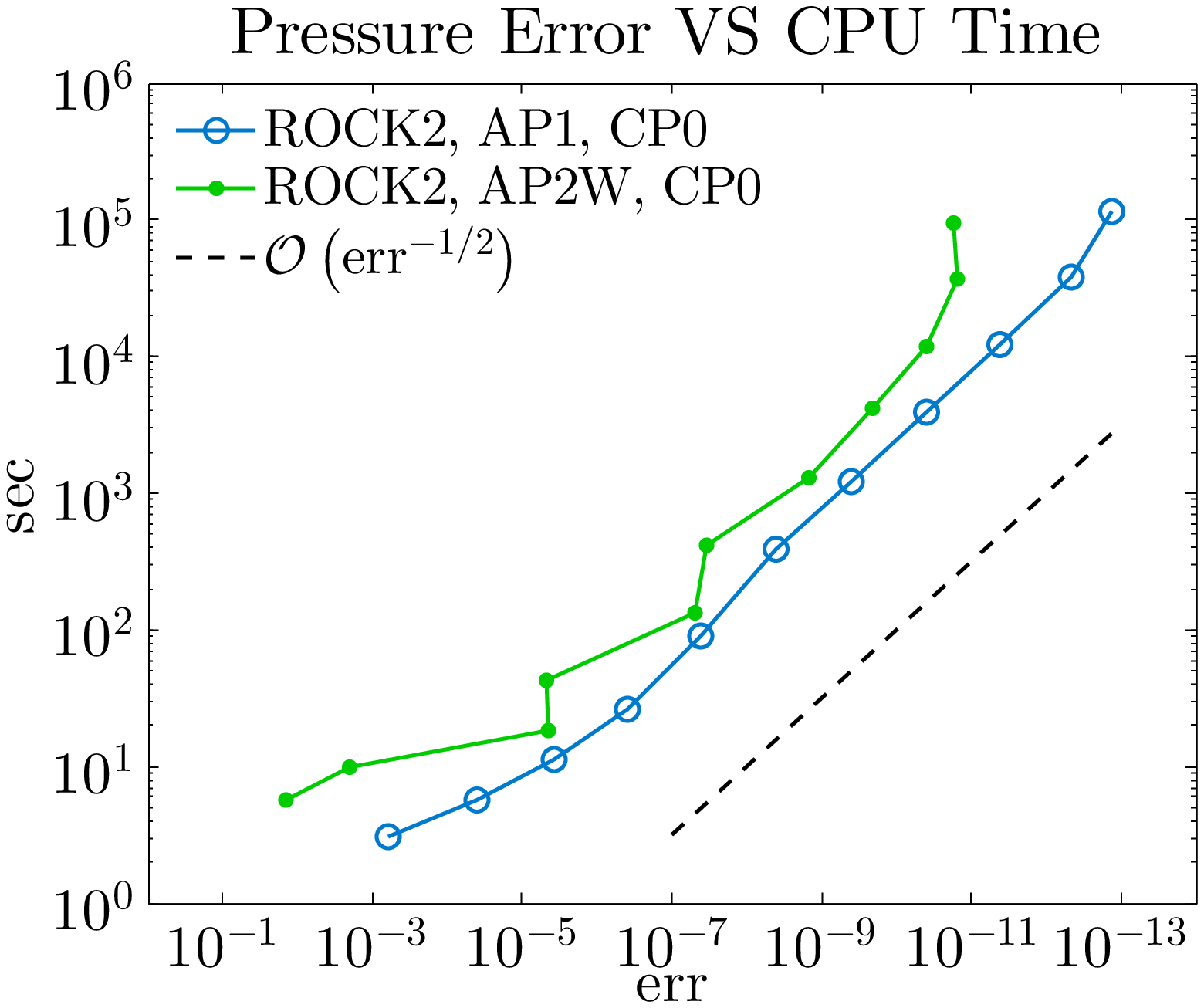} \label{fig:eff_ff_DAE_Pre_CP0}}
\subfigure[Pre. efficiency of AP1, AP2W and CP1.]{
\includegraphics[trim=0.0cm 0.0cm 0.0cm 0.0cm, clip=true, height=0.2\textheight]{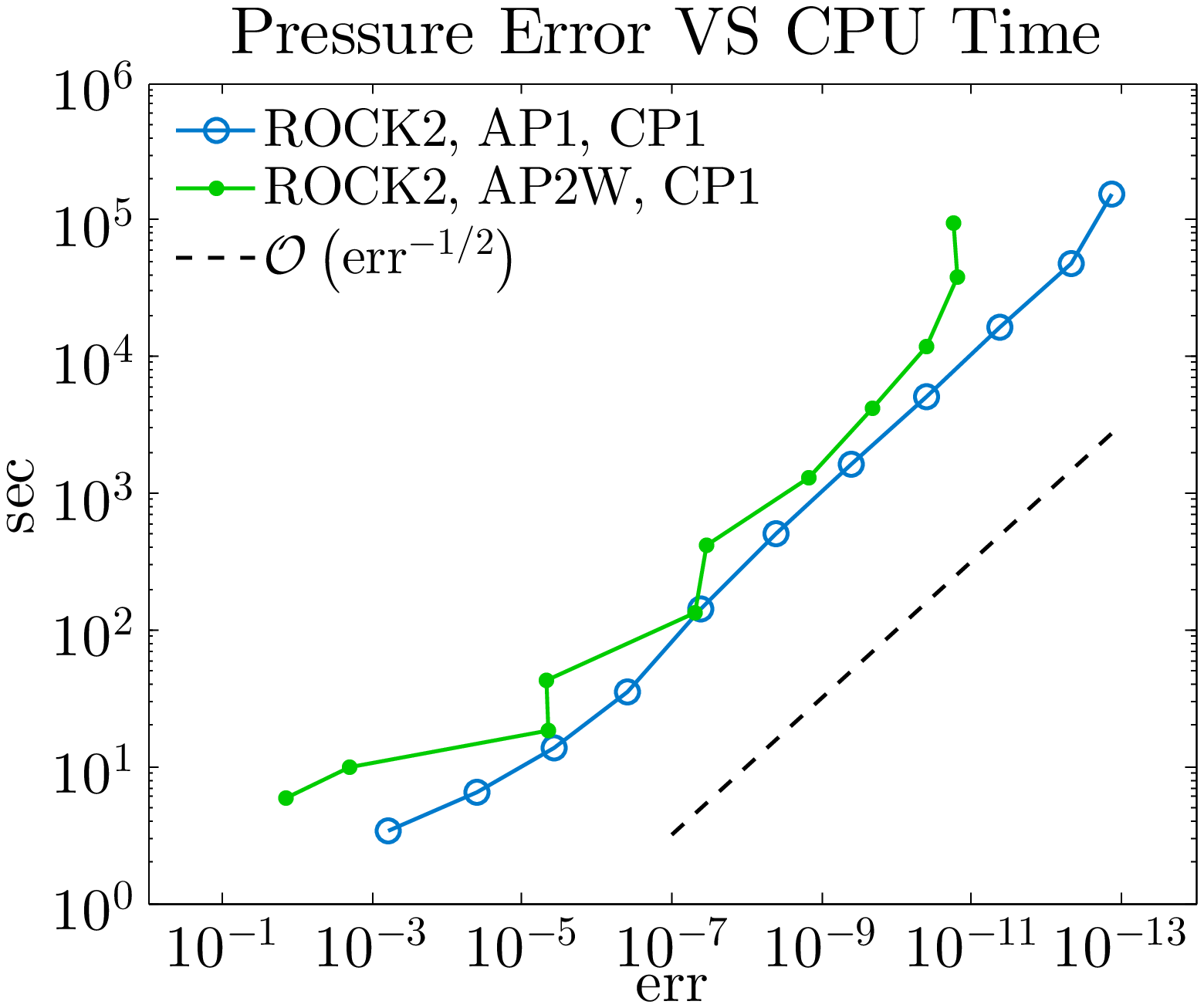} \label{fig:eff_ff_DAE_Pre_CP1}}
\end{center}
\caption{Numerical efficiency tests of the forced flow.}
\label{fig:eff_ff}
\end{figure}

In Figure \ref{fig:eff_ff} we compare the numerical efficiency of the methods for the forced flow problem. In Figures \ref{fig:eff_ff}(a,c) we see that among the PM1 and PM1V methods with CP0 the best one is ROCK2, PM1V while PIROCK is the less efficient. The efficiency lack in PIROCK can be associated to its stability issues. We also see that depending on the tolerance RKC, PM1, CP0 can beat ROCK2, PM1, CP0, especially for the pressure (see Figure \ref{fig:eff_ff_Proj_Pre_CP0}) but its accuracy increases slowly. When computing the second order pressure at each time step (see Figures \ref{fig:eff_ff}(b,d)) ROCK2, PM1V is again the best method and ROCK2, PM1, CP1 is better than RKC, PM1, CP1 for all the tolerances. In Figure \ref{fig:eff_ff}(e-h) we show the efficiency of ROCK2 with AP1 and AP2W. We see that for the largest tolerance there is an efficiency decrease. We think that it is not really an efficiency decreasing but instead it is the reference solution that is not accurate enough. This is even more visible for the pressure of AP2W since the pressure of the reference solution has been computed with AP1. However comparing Figures \ref{fig:eff_ff}(e,f) we see that switching from CP0 to CP1 does not change the computation time of AP2W while for AP1 there is a small increase because of the additional Poisson problem for the pressure. Looking at the pressure (Figures \ref{fig:eff_ff}(g,h)) we see that AP1 is more efficient than AP2W, specially for CP0. For CP1 the difference is smaller, it would be interesting to see what happens with a more accurate reference solution. When using CP0 AP1 is better than AP2W, when using CP1 it depends if the user is interested in the velocity or pressure accuracy.

As we already saw in section \ref{sec:errb_ff} projecting the velocity at each stage (PM1V, AP1, AP2, AP2W) is more efficient than projecting only once per time step. The overhead of the extra projections is compensated by the increased accuracy (see \ref{sec:errb_ff}), the better estimation of the local error (see \ref{sec:dae_tsa}) and the larger stability domain (see \ref{sec:stab_ff}). Yet, in this example a lot of computational time is spent computing the right hand side $f$ (see Table \ref{tab:prof}) because of the forcing terms, thus the time spent in projections is small in percentage. This means that in a problem without forcing term the overhead of the extra projections is more important and it could cause a decay in efficiency when projecting at each stage, we will investigate that in section \ref{sec:eff_gt}. The projection time depends strongly on the grid size $\Delta x$, it would be good to see what happens when a smaller grid size is used but we did not had time to do that. We are particularly interested in see how the results in Figures \ref{fig:eff_ff}(f,h) would change. It is also interesting to see that except for RKC, PM1, CP0 the computational time increases only as the square root of the required accuracy, while it increases linearly for RKC, PM1, CP0. Switching from CP0 to CP1 completely changes the behavior of RKC, PM1 (see Figures \ref{fig:eff_ff}(a-d)), it is much more sensible to the pressure accuracy than the other methods. In fact, comparing the first column of figures in \ref{fig:eff_ff} with the second column we remark that for ROCK2 and PIROCK computing a second order pressure after each time step slightly increases the computational time without improving accuracy at the final time step at $t=1$ (where the errors are measured). Hence, for ROCK2 and PIROCK one should use CP1 only when accurate pressures are needed at intermediate time steps. While for RKC, PM1 there's a large difference between CP0 and CP1.

With a side by side comparison ROCK2, AP1 is the best method when CP0 is used. For CP1 the choice is between ROCK2, AP1 and ROCK2, AP2W, it depends if the user is interested on velocity or pressure accuracy. If it is interest in velocity it choses AP2W because it is faster, if it is interested in pressure it choses AP1 because of its accuracy. ROCK2, PM1V performs almost as well as ROCK2, AP1. Except the pressure update the difference between these methods is the recurrence formulation. PM1V uses algorithm \eqref{eq:rkc_rec_proj} while AP1 uses \eqref{eq:rkc_rec_projf} (see section \ref{sec:rec_dae}). For this reason we prefer AP1 since its realization with recursive formulas is consistent with its realization with the coefficients $a_{ij}$, $b_i$ (see section \ref{sec:rec_dae}).

\subsection{Different Reynolds numbers behavior} \label{sec:rey_ff}
\begin{figure}[!hbtp]
\begin{center}
\subfigure[Errors in infinity norm.]{
\includegraphics[trim=0.0cm 0.0cm 0.0cm 0.0cm, clip=true, height=0.2\textheight]{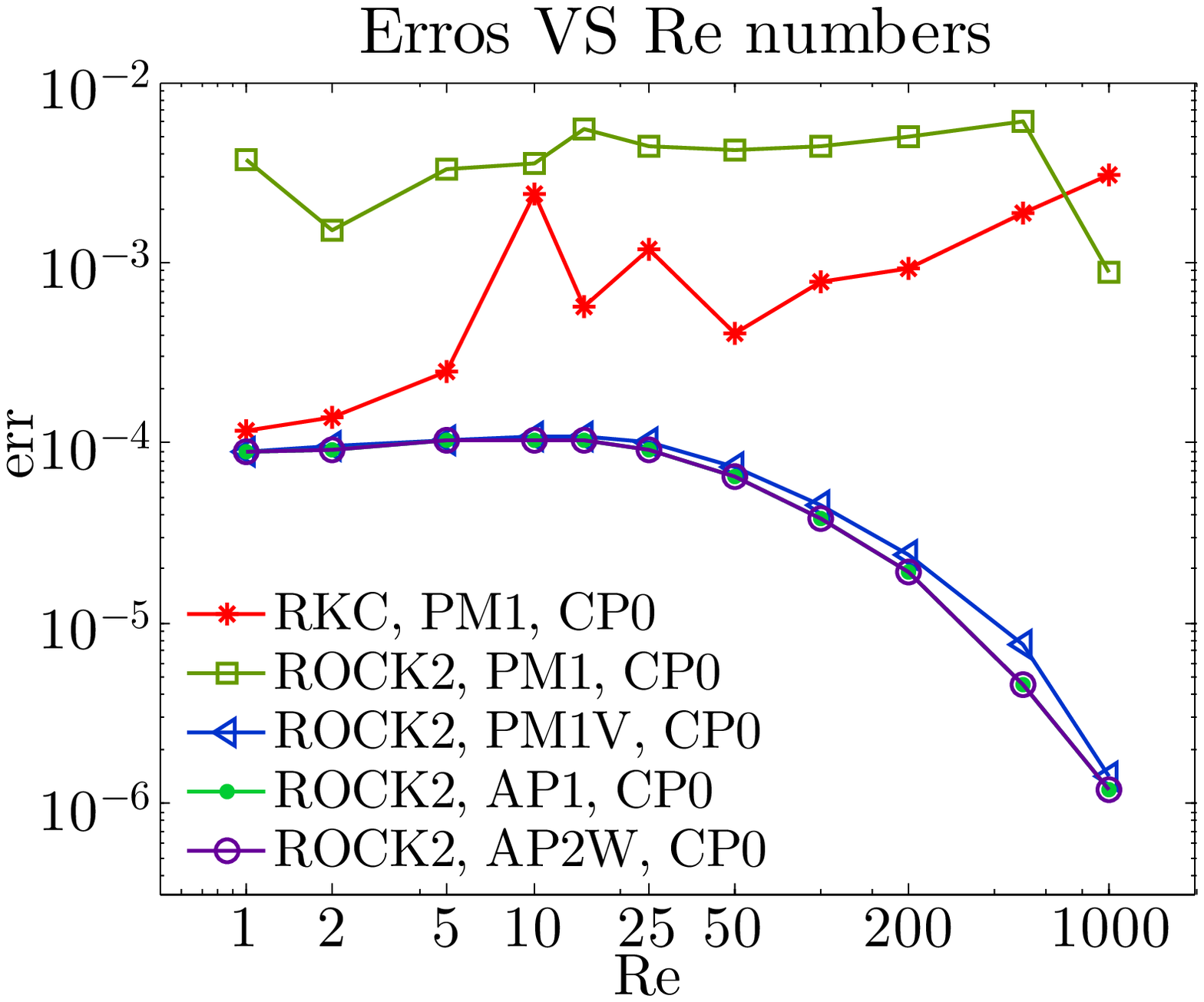} \label{fig:rey_err}}
\subfigure[CPU time.]{
\includegraphics[trim=0.0cm 0.0cm 0.0cm 0.0cm, clip=true, height=0.2\textheight]{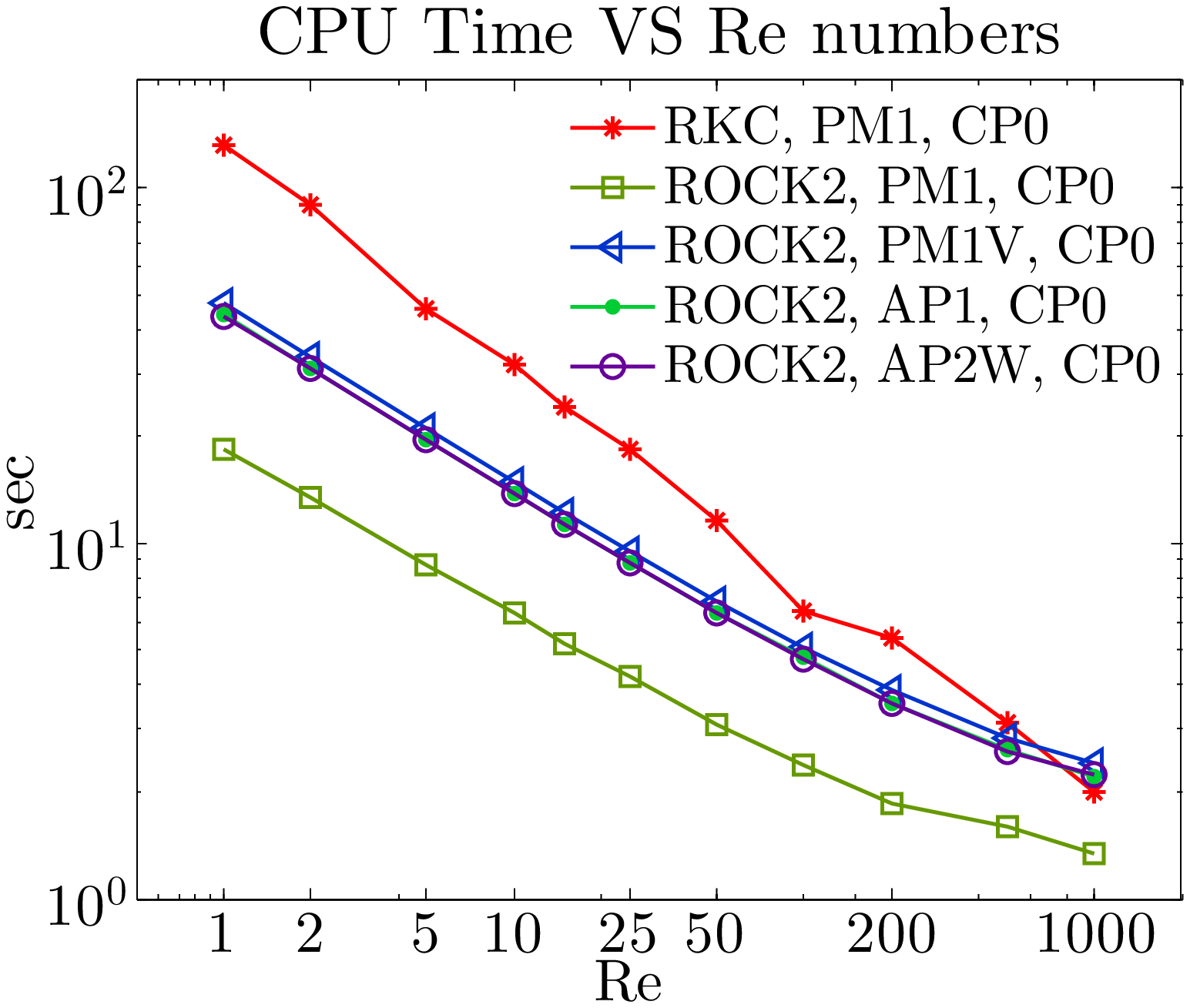} \label{fig:rey_time}} \\
\subfigure[Average number of stages per step.]{
\includegraphics[trim=0.0cm 0.0cm 0.0cm 0.0cm, clip=true, height=0.2\textheight]{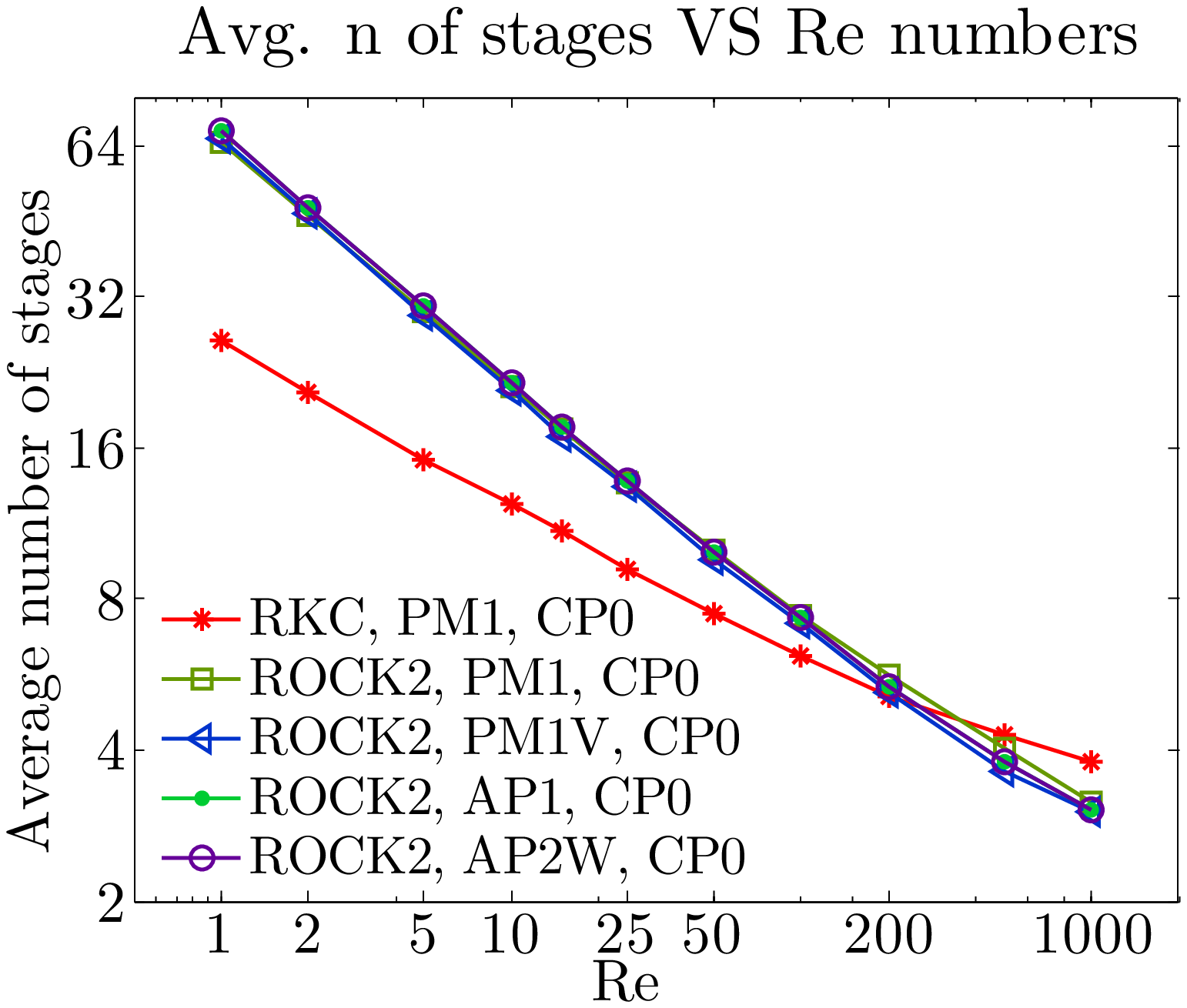} \label{fig:rey_avg}}
\subfigure[Total number of stages.]{
\includegraphics[trim=0.0cm 0.0cm 0.0cm 0.0cm, clip=true, height=0.2\textheight]{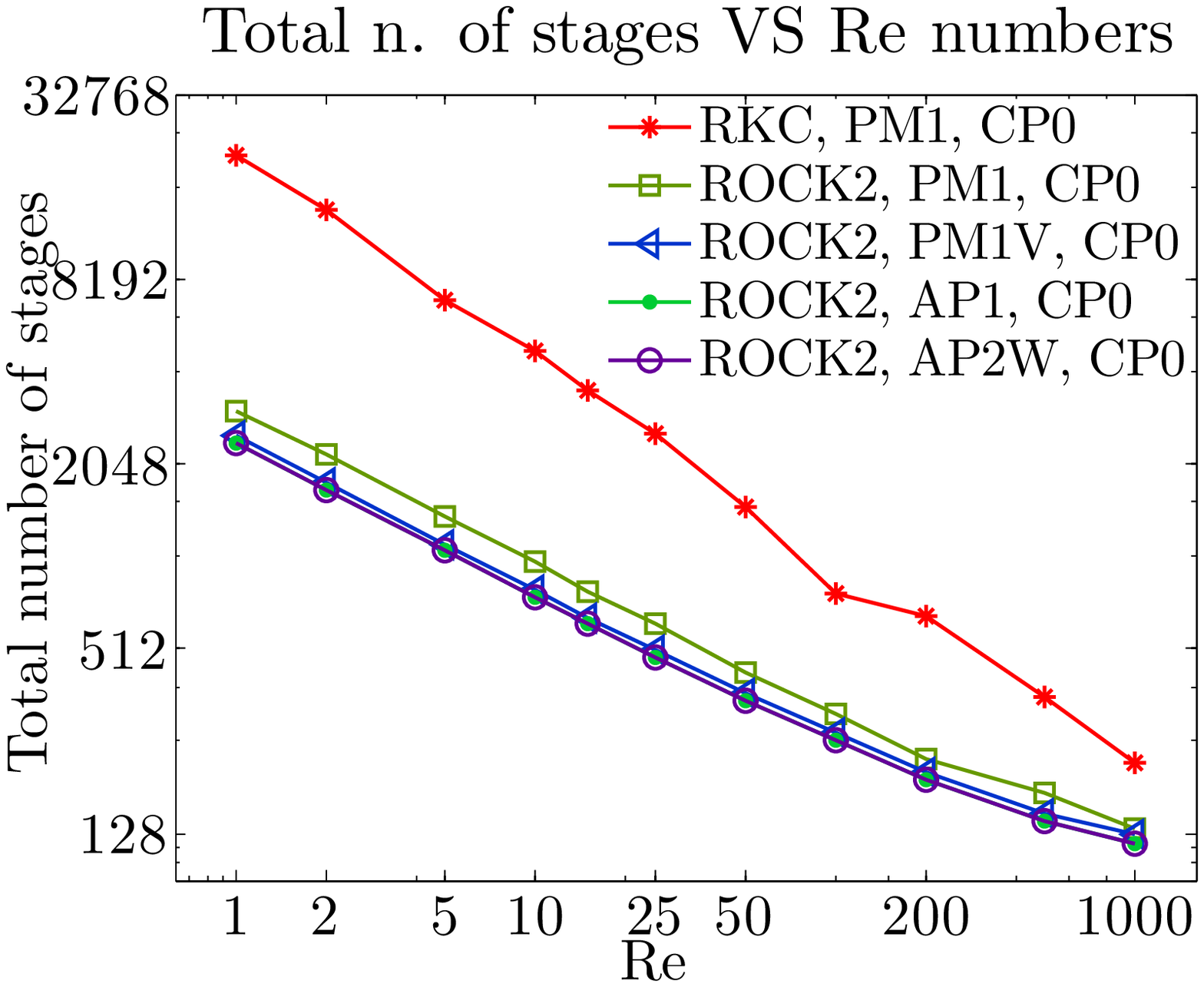} \label{fig:rey_sum}} \\
\subfigure[Number of time steps.]{
\includegraphics[trim=0.0cm 0.0cm 0.0cm 0.0cm, clip=true, height=0.2\textheight]{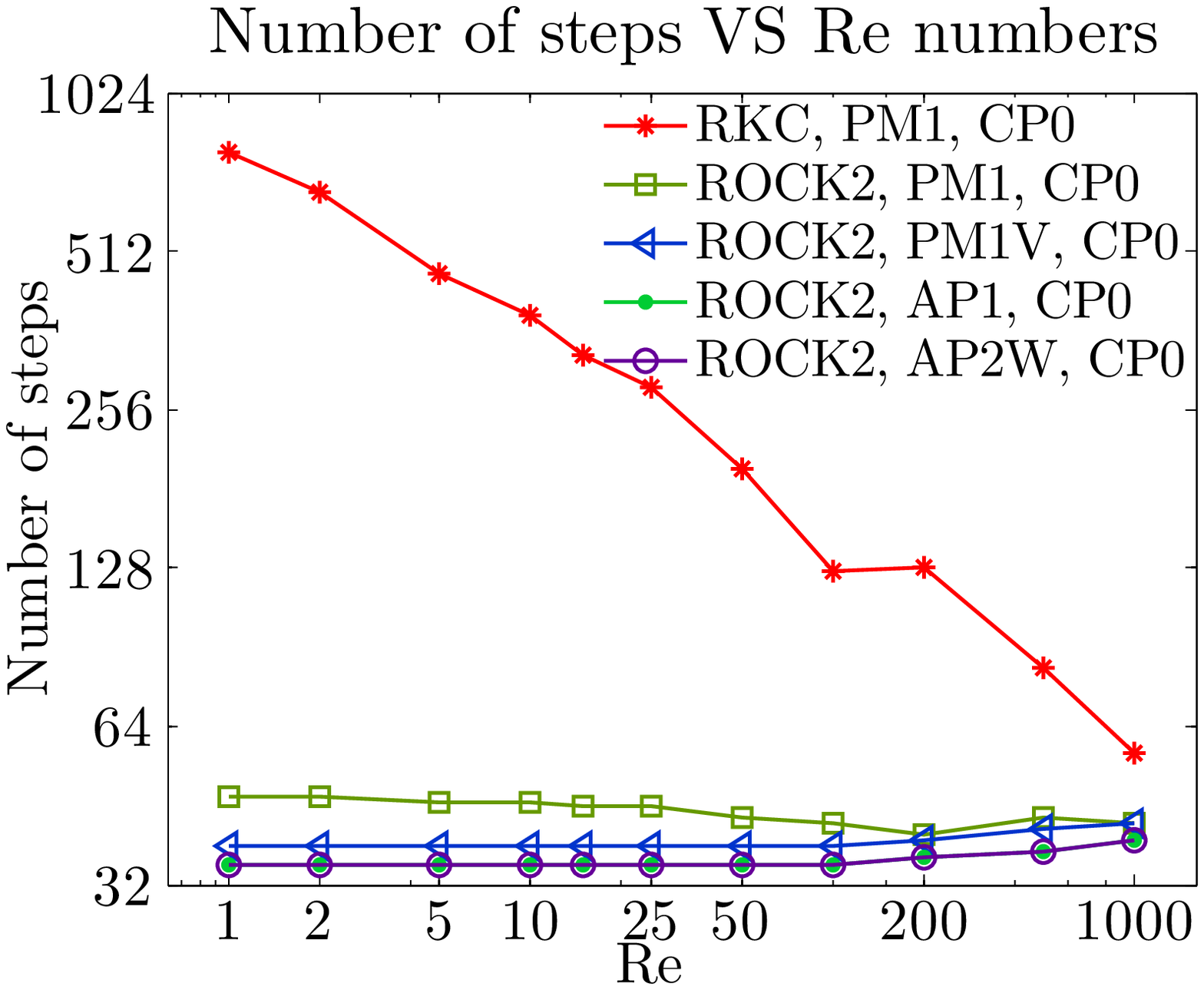} \label{fig:rey_tsteps}}
\subfigure[Number of rejected steps.]{
\includegraphics[trim=0.0cm 0.0cm 0.0cm 0.0cm, clip=true, height=0.2\textheight]{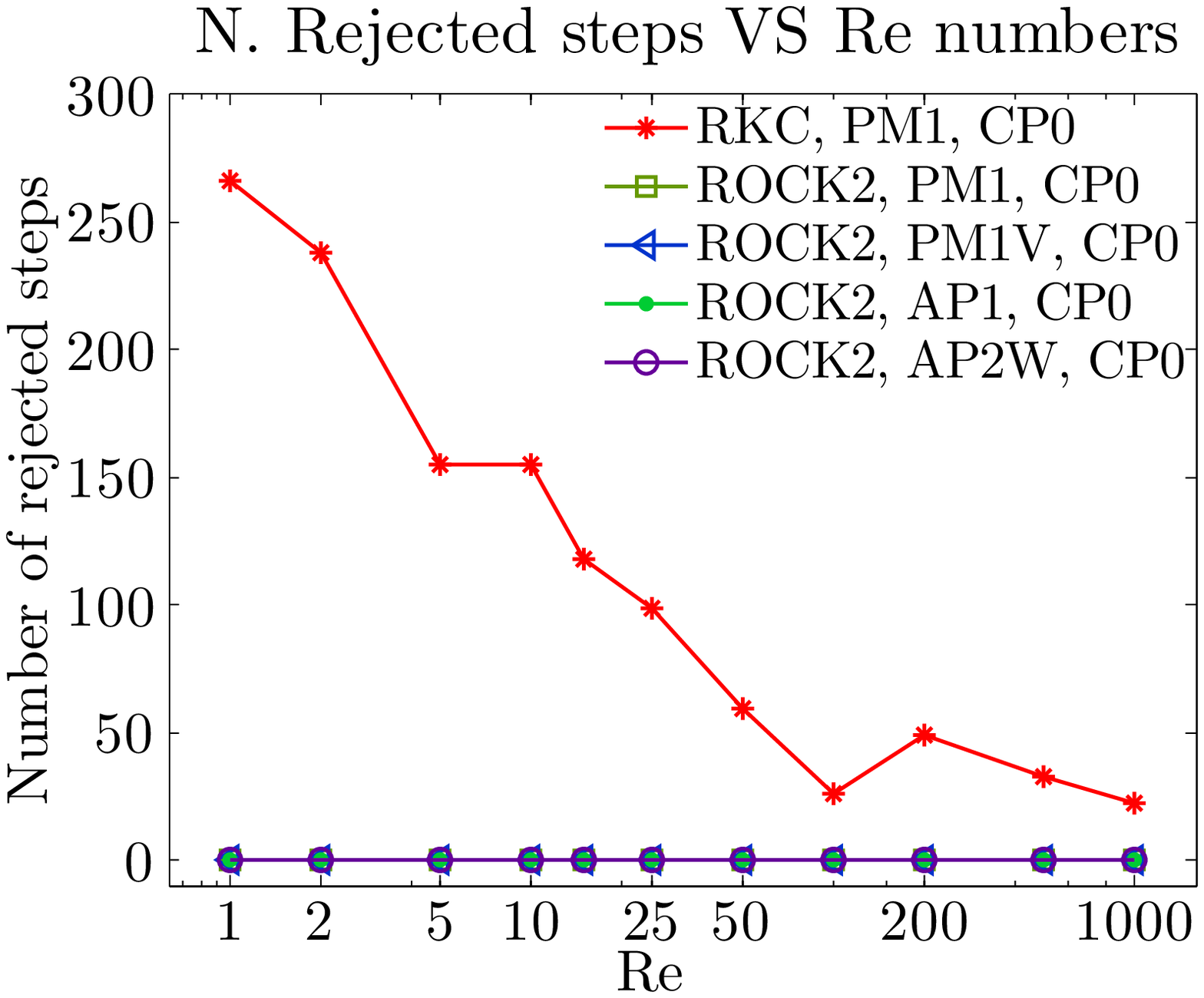} \label{fig:rey_rsteps}} 
\end{center}
\caption{Different Reynolds number behaviour.}
\label{fig:rey_ff}
\end{figure}

In this test we compare the different methods when changing the Reynolds number, the results are in Figure \ref{fig:rey_ff}. Figure \ref{fig:rey_err} shows the error of the different methods against the exact solution, we see that again projecting the velocity after each stage (PM1V, AP1, AP2W) is much more accurate than PM1. Figure \ref{fig:rey_time} shows that RKC, PM1 is not only badly accurate but also the most computationally expensive method. ROCK2, PM1 is less accurate than ROCK2 with PM1V, AP1, AP2W but at least it is faster. However we saw in \ref{sec:errb_ff} and \ref{sec:eff_ff} that this velocity is not worthwhile. We see that the computational time of ROCK2, PM1V is higher than ROCK2 with AP1 and AP2W. Figure \ref{fig:rey_avg} shows the average number of stages used in each time step, we see that RKC, PM1 uses much less stages per time step but on the other hand in Figure \ref{fig:rey_sum} we see that in total it uses more stages than the other methods. In fact it has a smaller average $\Delta t$, see Figure \ref{fig:rey_tsteps} which displays the total number of time steps. From this figure we see that PM1 and PM1V require more steps than AP1 and AP2W without gains in accuracy (see Figure \ref{fig:rey_err}), proving that algorithm \eqref{eq:rkc_rec_projf} is better than \eqref{eq:rkc_rec_proj}. Figure \ref{fig:rey_rsteps} shows that RKC has a considerable number of rejected time steps while ROCK2 do not have any rejected step.

\section{The Green-Taylor vortex} \label{sec:gt}
The Green-Taylor vortex is an unsteady flow of a decaying vortex. Its analytical solution in the domain $\Omega=[0,1]\times[0,1]$ is given by
\begin{align} \label{eq:gtv}
\begin{aligned}
u(t,x,y) &= -e^{-2\pi^2 t/\Rey}\sin(\pi x)\cos(\pi y) ,\\
v(t,x,y) &= e^{-2\pi^2 t/\Rey}\cos(\pi x)\sin(\pi y), \\
p(t,x,y) &= \frac{1}{4}e^{-4\pi^2 t/\Rey} \left(\cos(2\pi x)+\cos(2\pi y)\right). 
\end{aligned}
\end{align}
This test has been taken from \cite{acc_an}, here no forcing terms are needed. In the following we will present the convergence results and compare the numerical efficiency of the different methods.

\subsection{Convergence order}\label{sec:conv_gt}
\begin{figure}[!hbtp]
\begin{center}
\subfigure[PIROCK, PM1, CP0.]{
\includegraphics[trim=0.0cm 0.0cm 0.0cm 0.0cm, clip=true, height=0.2\textheight]{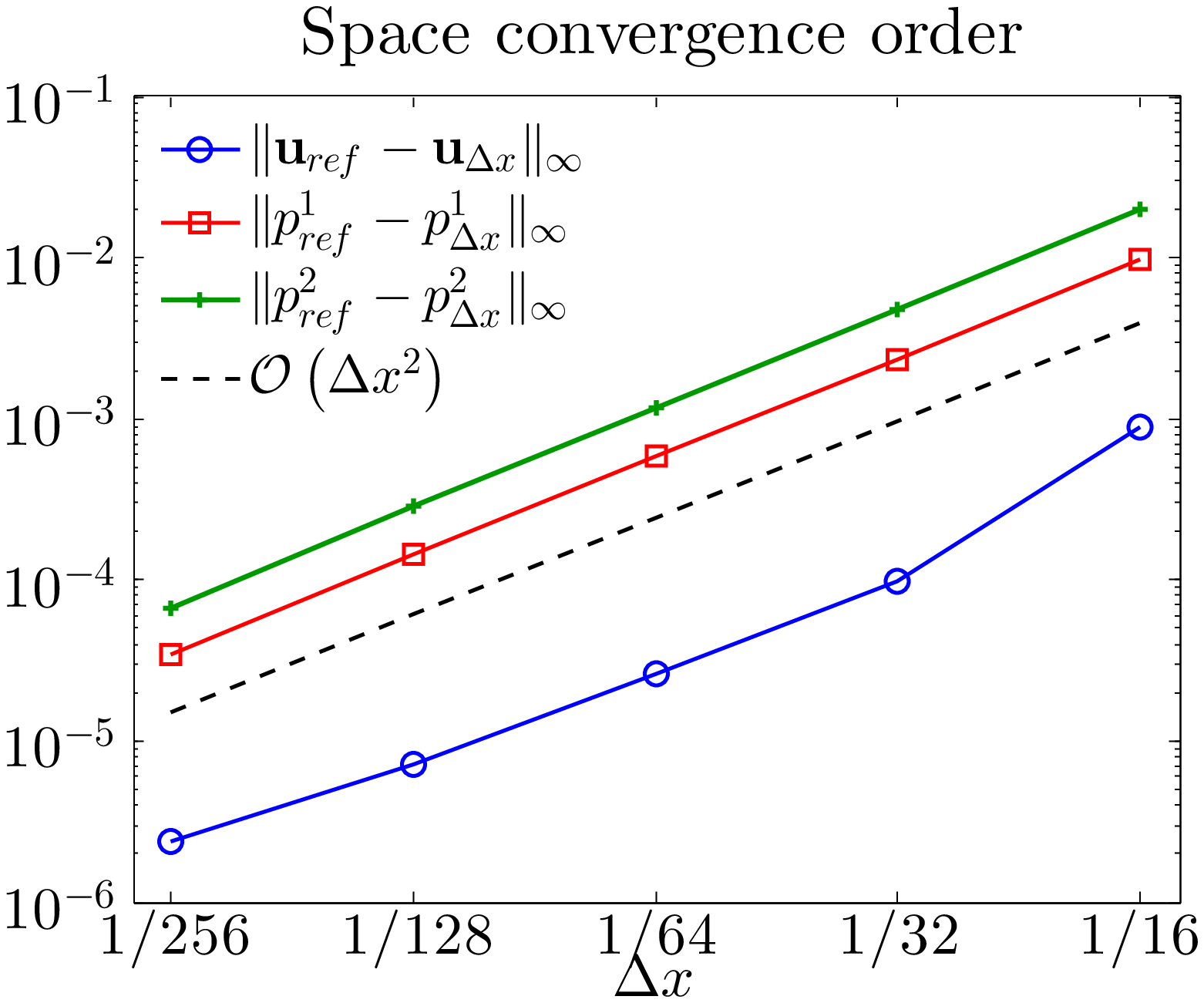} \label{fig:space_gt_pirock_PM1}} \\
\subfigure[RKC, PM1, CP0.]{
\includegraphics[trim=0.0cm 0.0cm 0.0cm 0.0cm, clip=true, height=0.2\textheight]{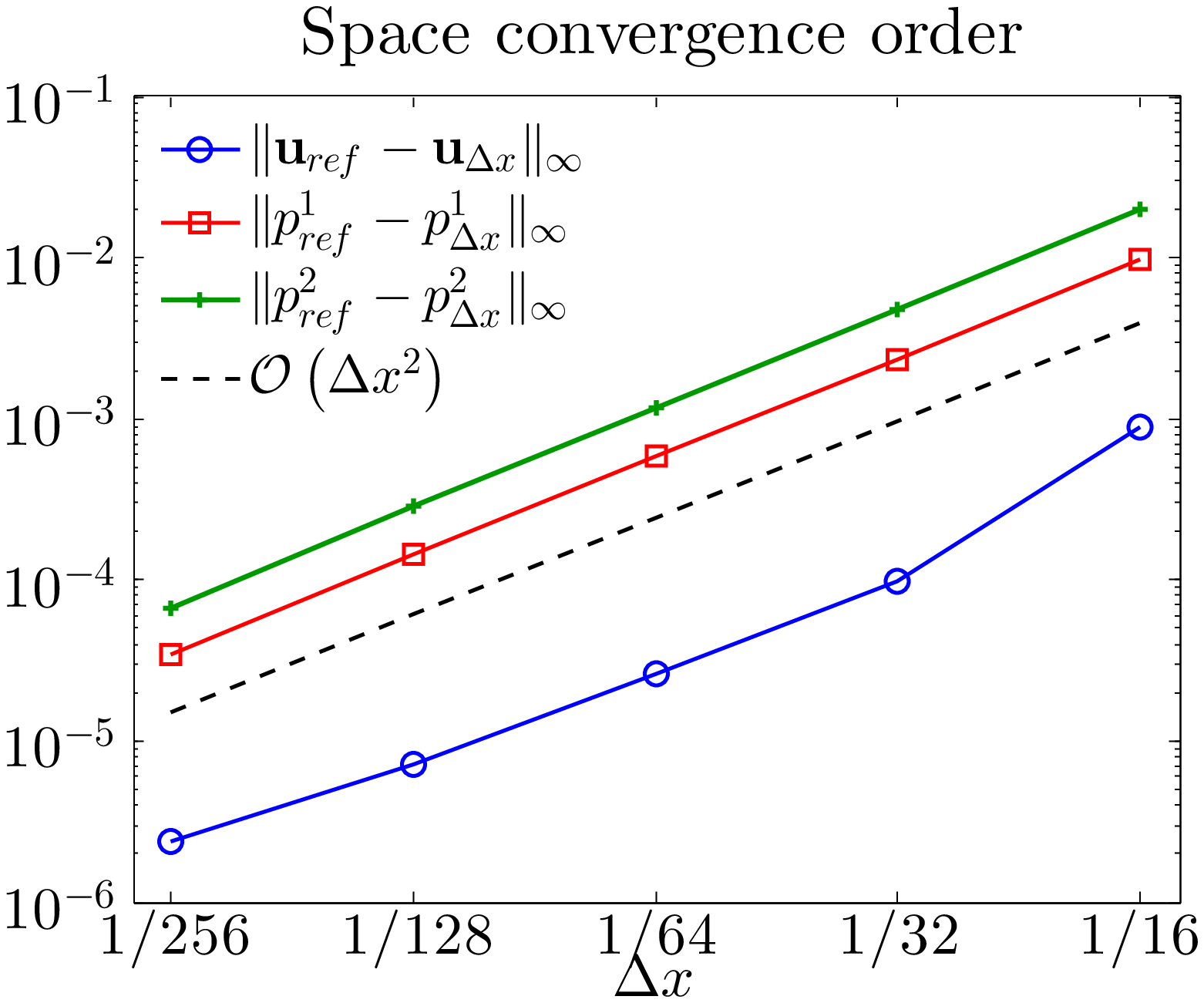} \label{fig:space_gt_rkc_PM1}}
\subfigure[ROCK2, PM1, CP0.]{
\includegraphics[trim=0.0cm 0.0cm 0.0cm 0.0cm, clip=true, height=0.2\textheight]{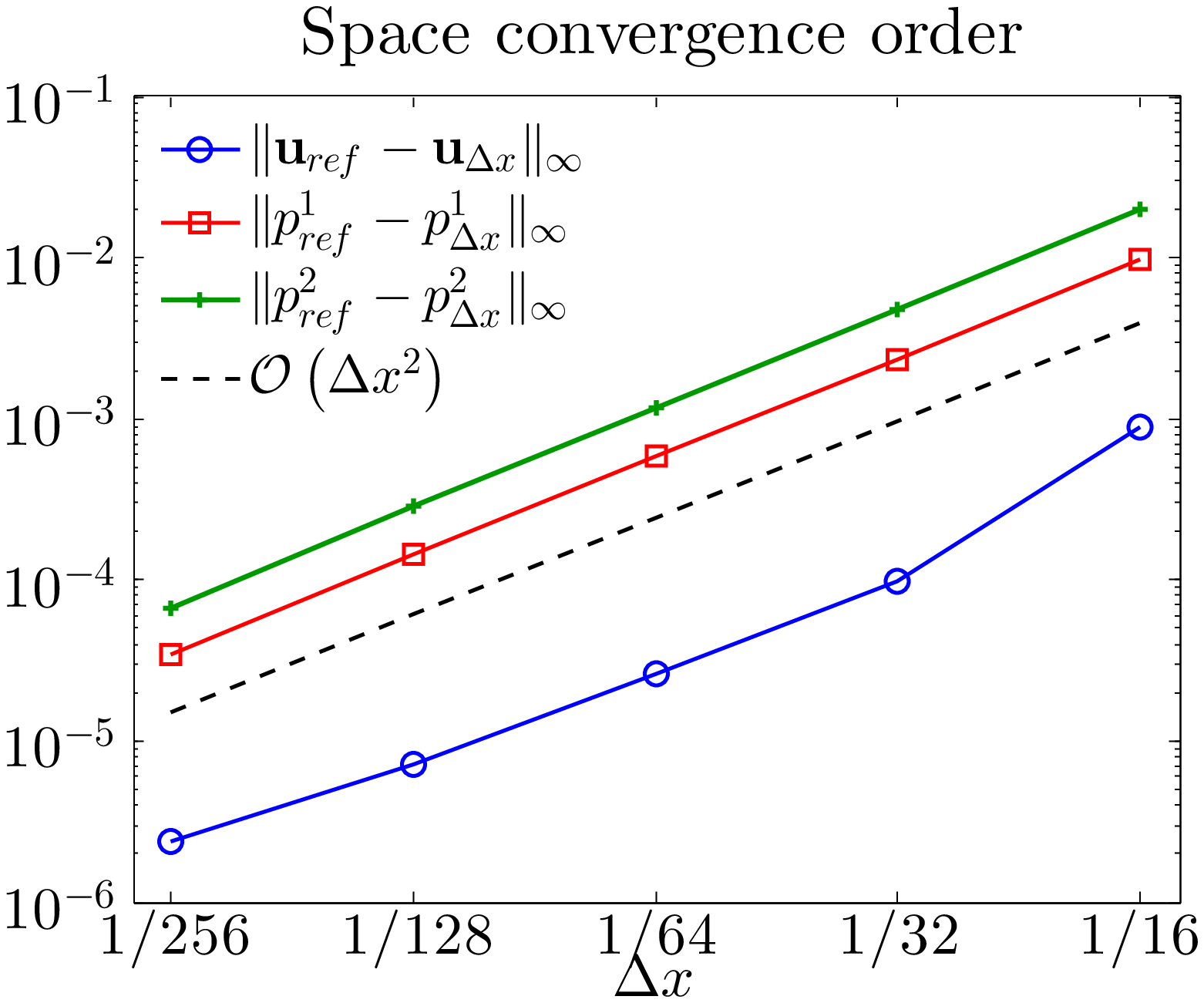} \label{fig:space_gt_rock2_PM1}} \\
\subfigure[RKC, PM1V, CP0.]{
\includegraphics[trim=0.0cm 0.0cm 0.0cm 0.0cm, clip=true, height=0.2\textheight]{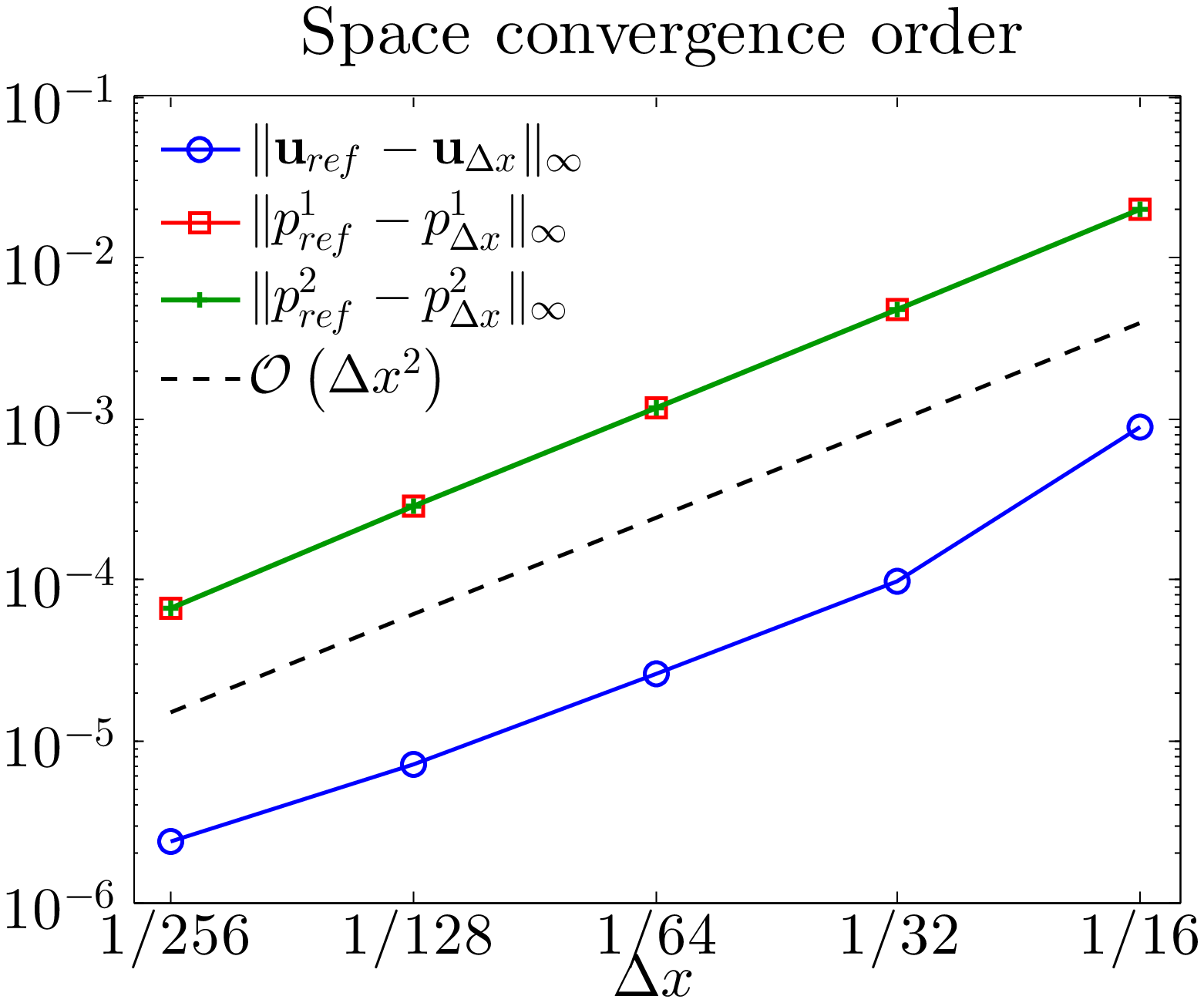} \label{fig:space_gt_rkc_PM1V}}
\subfigure[ROCK2, PM1V, CP0.]{
\includegraphics[trim=0.0cm 0.0cm 0.0cm 0.0cm, clip=true, height=0.2\textheight]{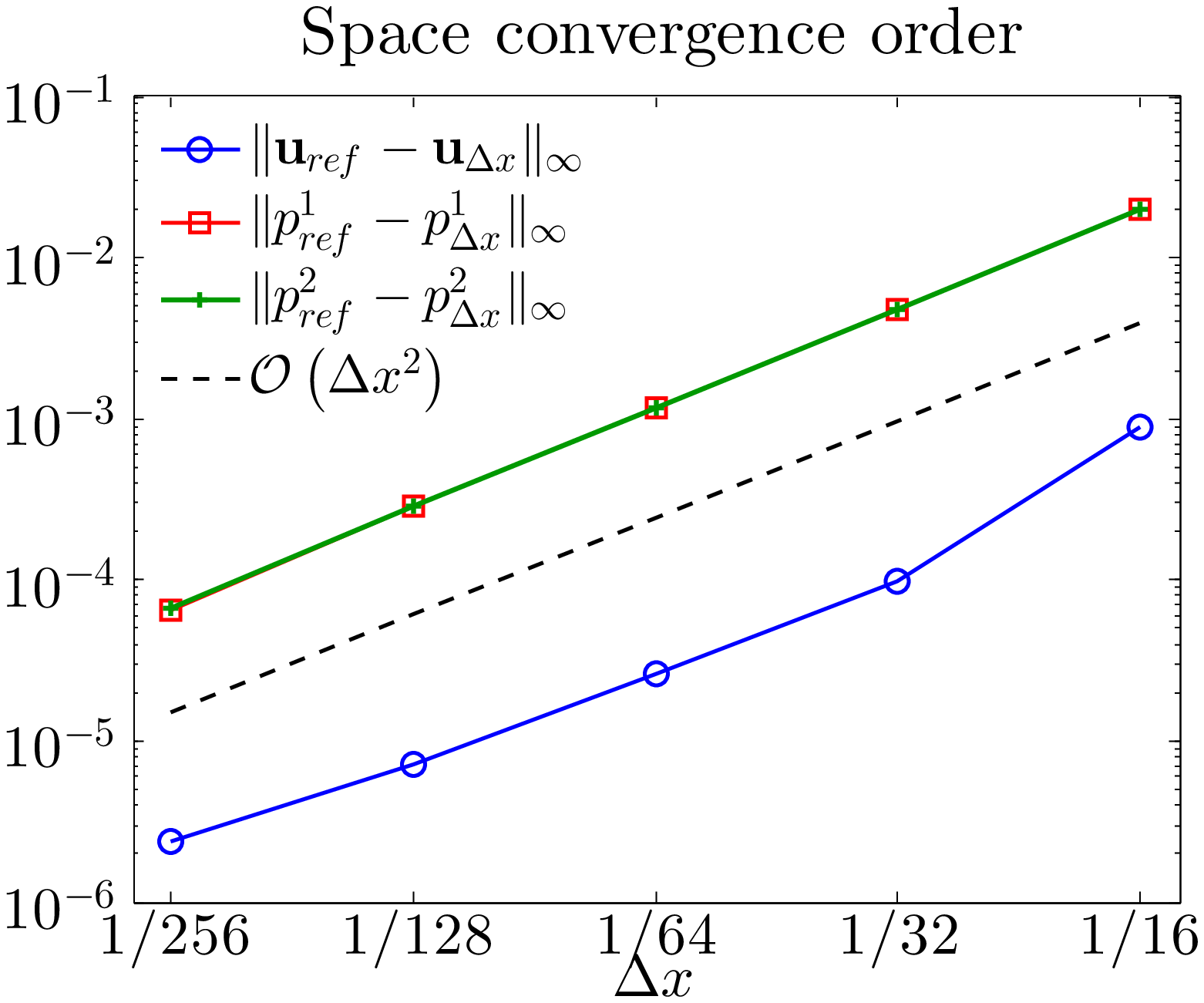} \label{fig:space_gt_rock2_PM1V}} \\
\subfigure[RKC, AP2 and AP2, CP0.]{
\includegraphics[trim=0.0cm 0.0cm 0.0cm 0.0cm, clip=true, height=0.2\textheight]{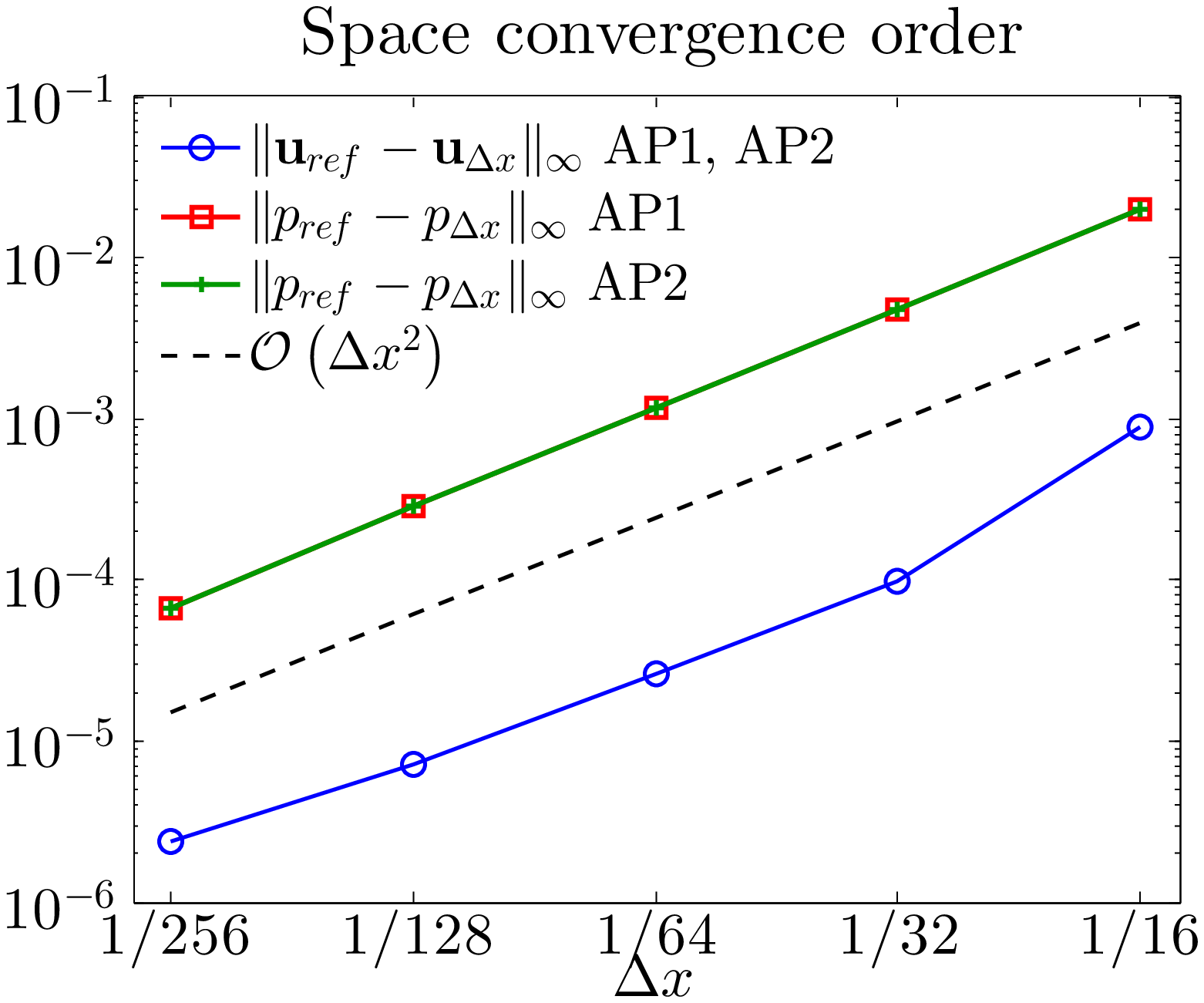} \label{fig:space_gt_rkc_AP12}}
\subfigure[ROCK2, AP1 and AP2W, CP0.]{
\includegraphics[trim=0.0cm 0.0cm 0.0cm 0.0cm, clip=true, height=0.2\textheight]{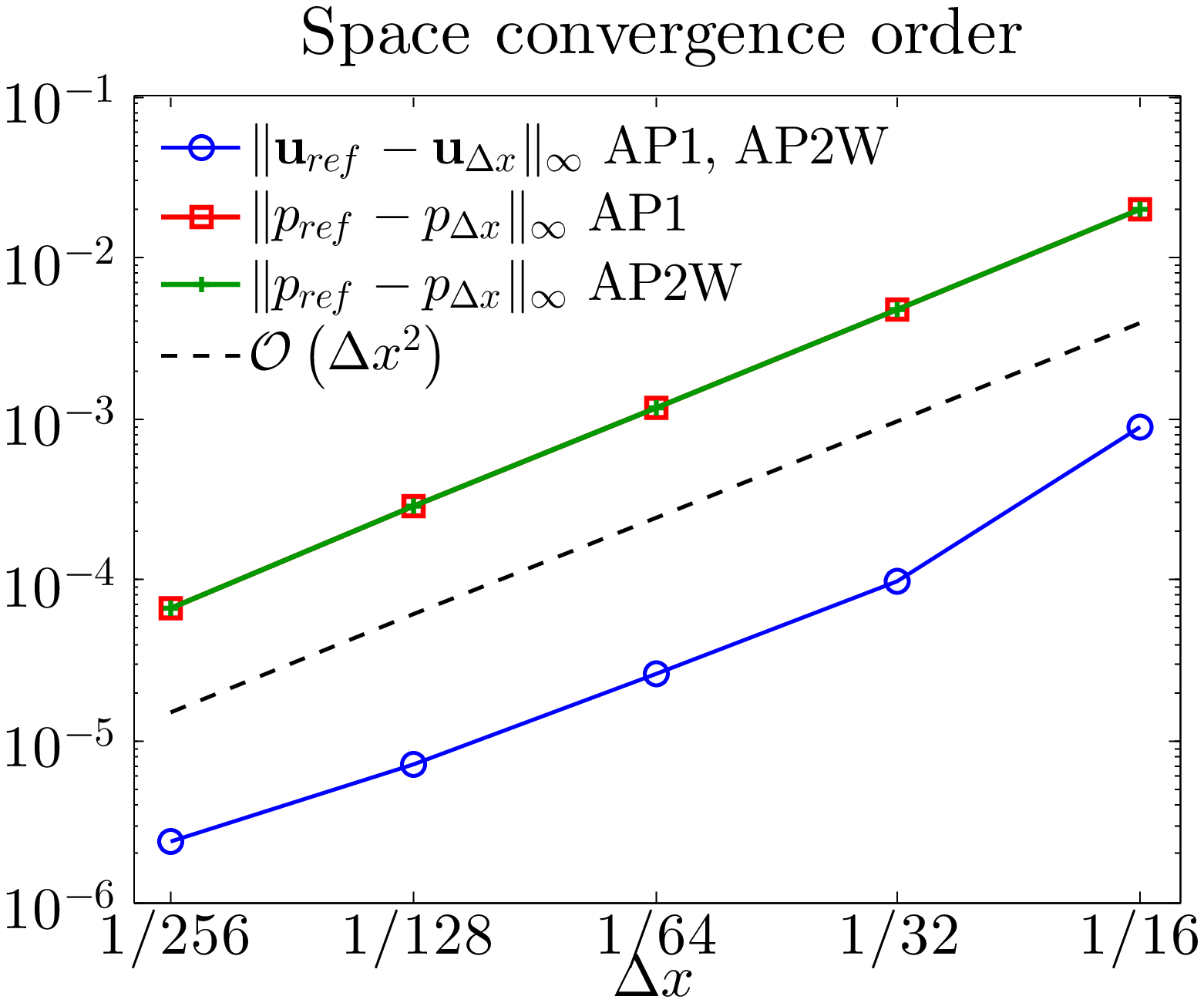} \label{fig:space_gt_rock2_AP12}}
\end{center}
\caption{Space convergence results of the Green-Taylor vortex.}
\label{fig:conv_space_gt}
\end{figure}

\begin{figure}[!hbtp]
\begin{center}
\subfigure[PIROCK, PM1, CP0.]{
\includegraphics[trim=0.0cm 0.0cm 0.0cm 0.0cm, clip=true, height=0.2\textheight]{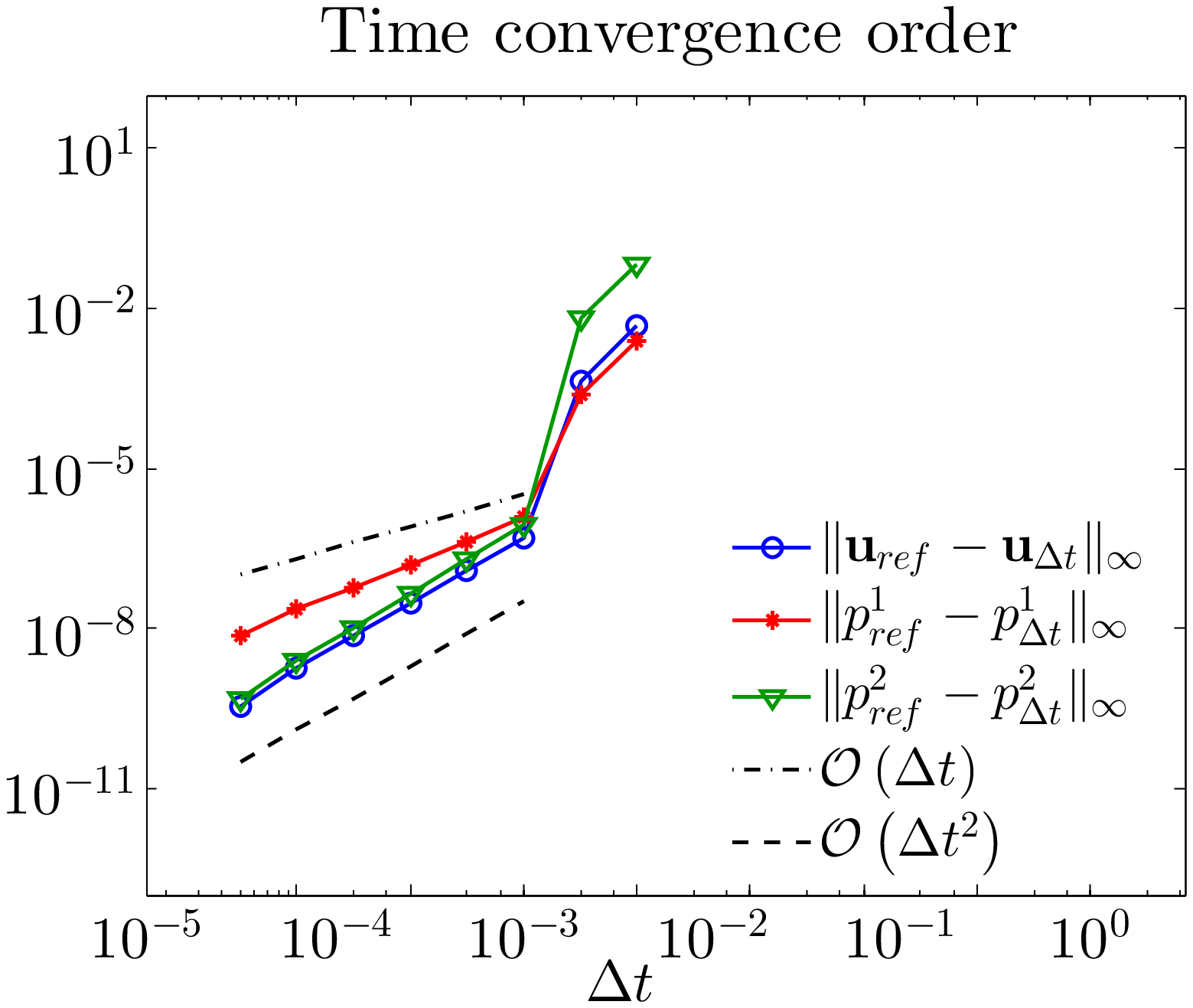} \label{fig:time_gt_pirock_PM1}} \\
\subfigure[RKC, PM1, CP0.]{
\includegraphics[trim=0.0cm 0.0cm 0.0cm 0.0cm, clip=true, height=0.2\textheight]{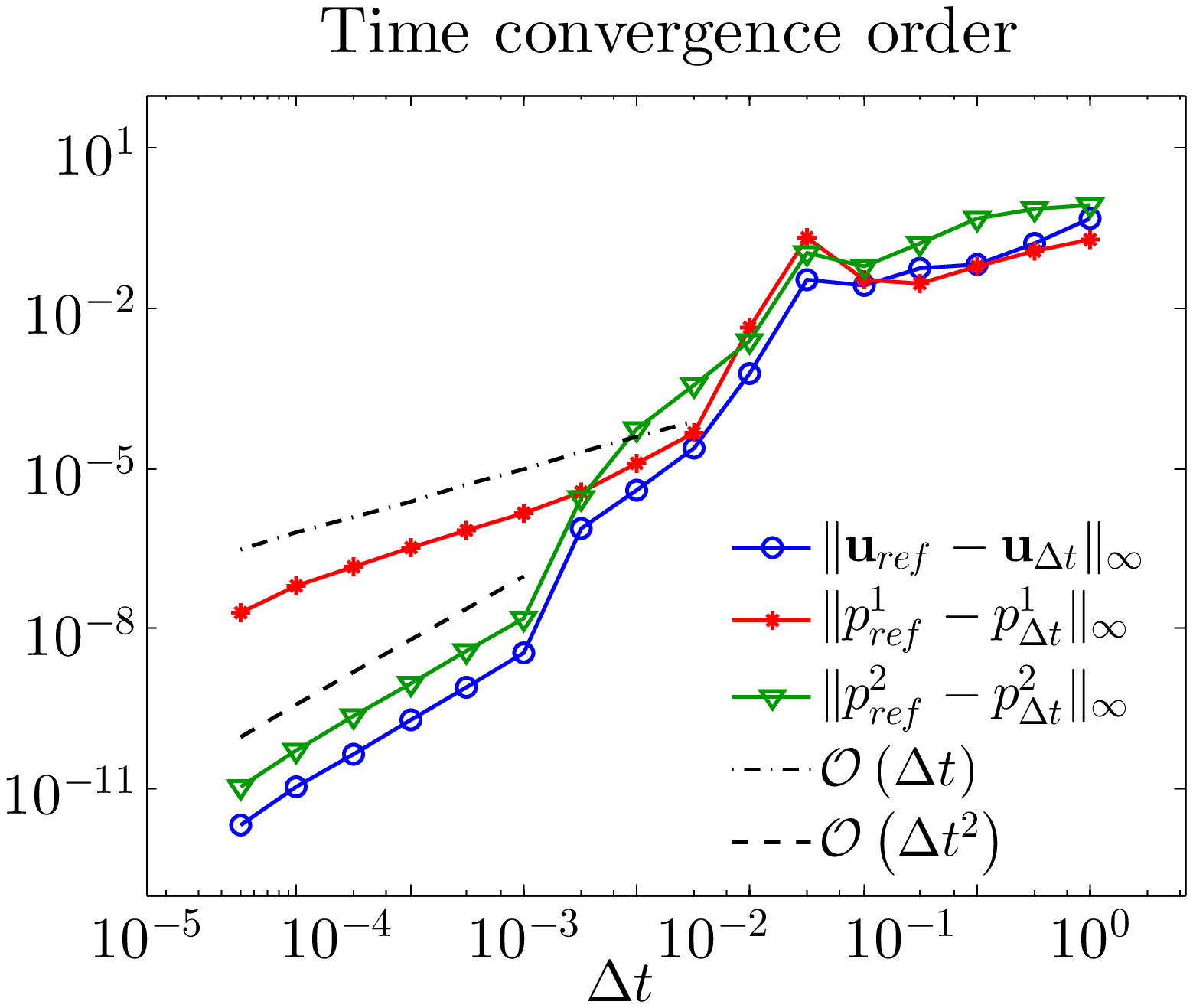} \label{fig:time_gt_rkc_PM1}}
\subfigure[ROCK2, PM1, CP0.]{
\includegraphics[trim=0.0cm 0.0cm 0.0cm 0.0cm, clip=true, height=0.2\textheight]{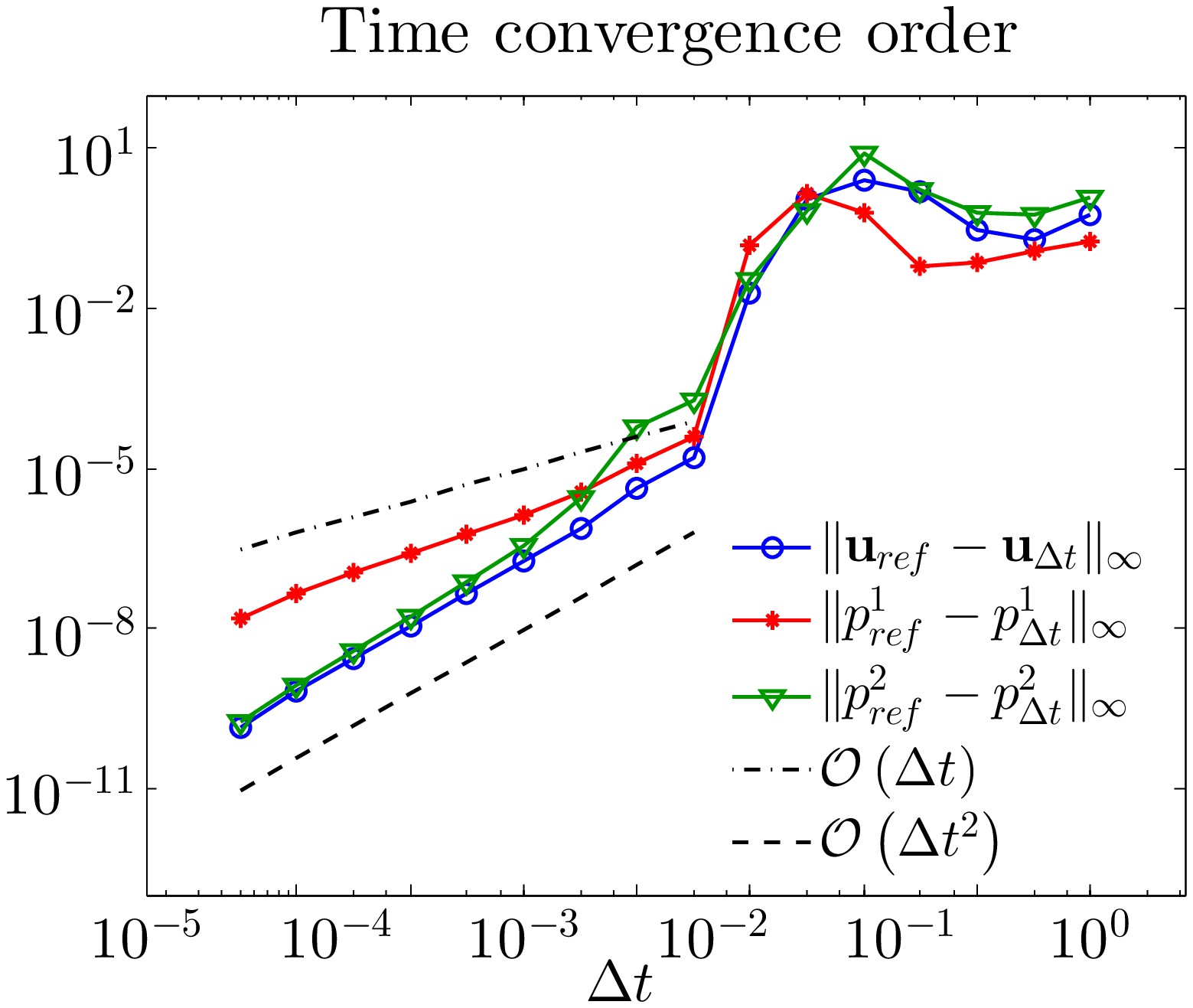} \label{fig:time_gt_rock2_PM1}} \\
\subfigure[RKC, PM1V, CP0.]{
\includegraphics[trim=0.0cm 0.0cm 0.0cm 0.0cm, clip=true, height=0.2\textheight]{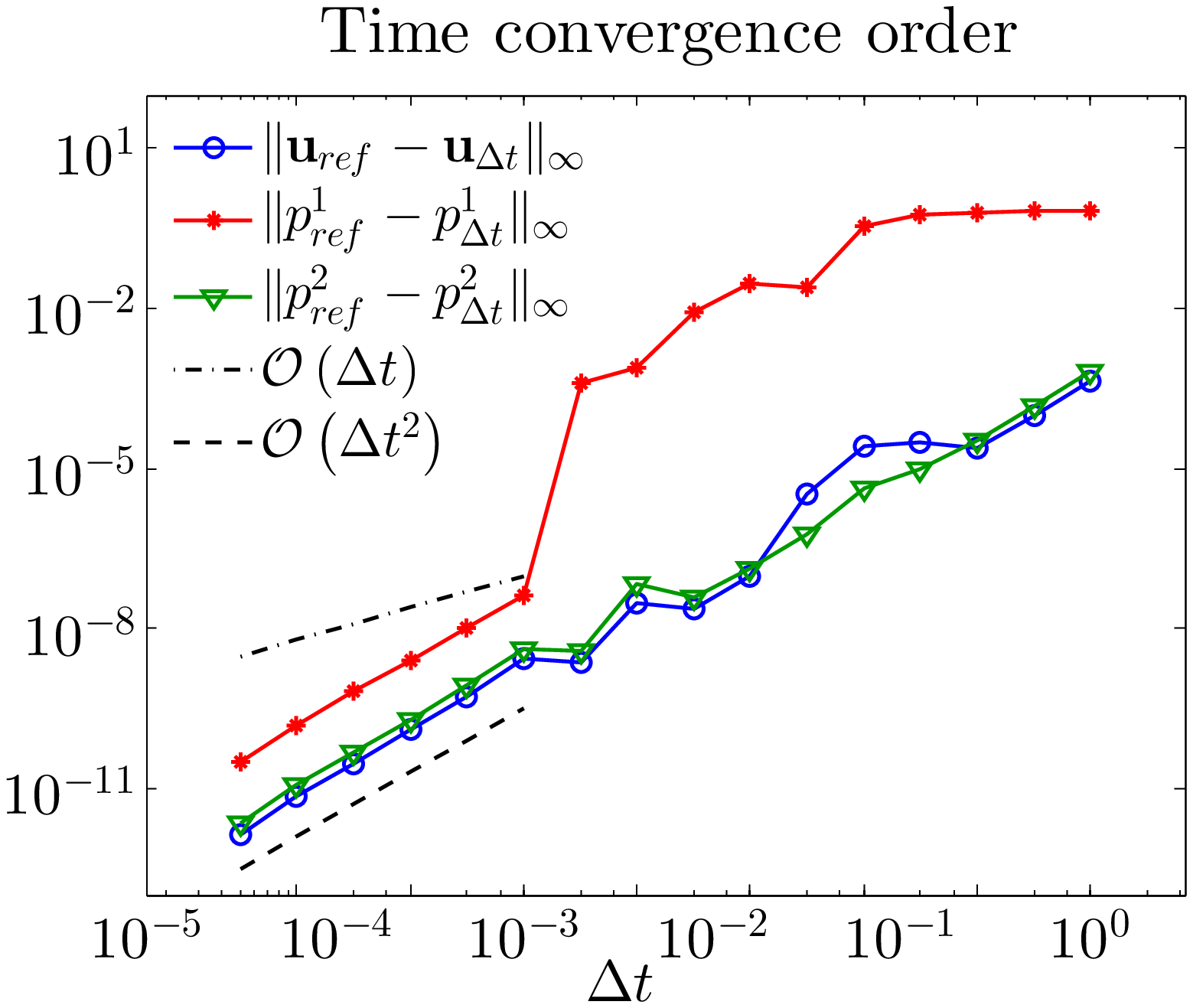} \label{fig:time_gt_rkc_PM1V}}
\subfigure[ROCK2, PM1V, CP0.]{
\includegraphics[trim=0.0cm 0.0cm 0.0cm 0.0cm, clip=true, height=0.2\textheight]{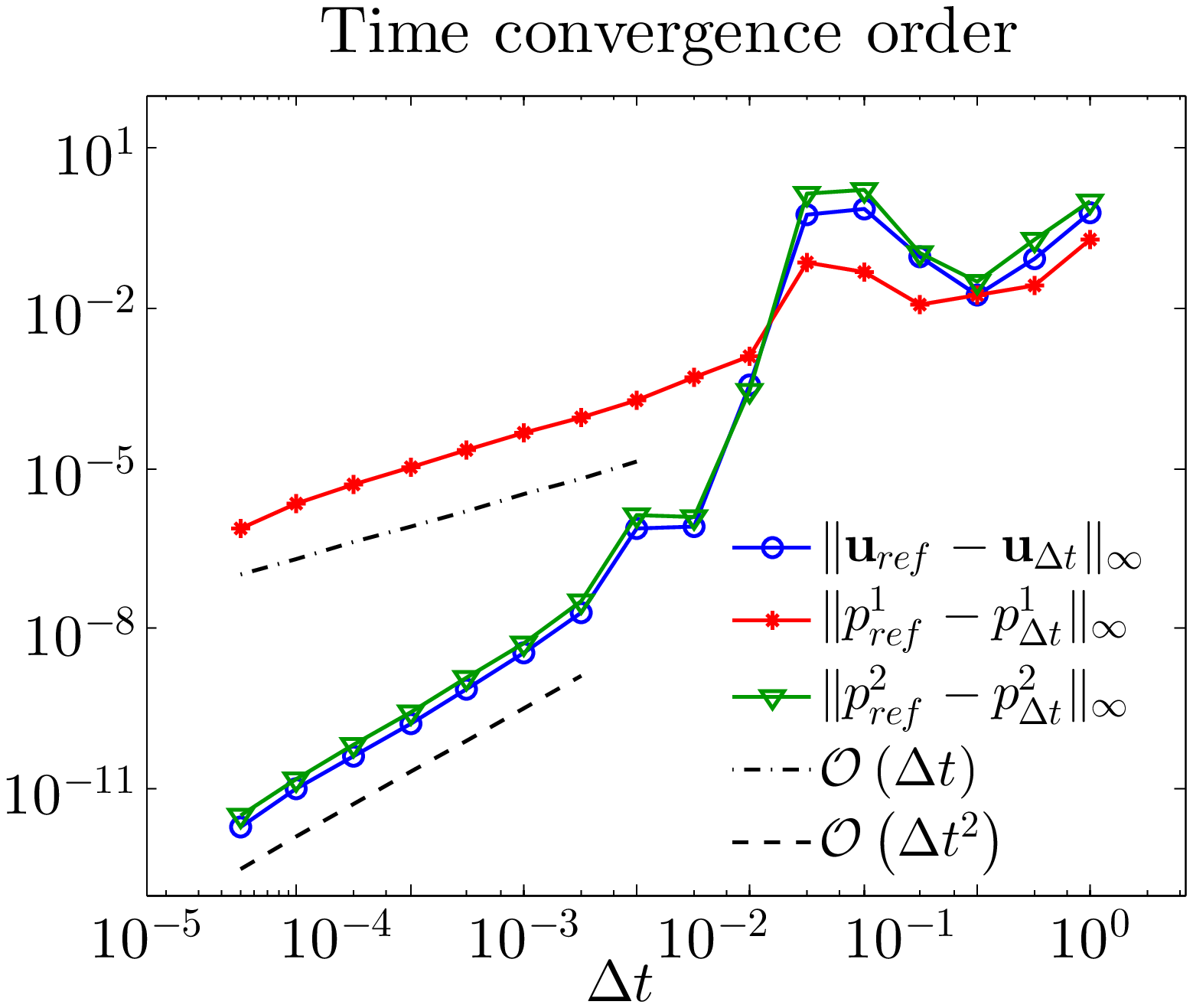} \label{fig:time_gt_rock2_PM1V}} \\
\subfigure[RKC, AP2 and AP2, CP0.]{
\includegraphics[trim=0.0cm 0.0cm 0.0cm 0.0cm, clip=true, height=0.2\textheight]{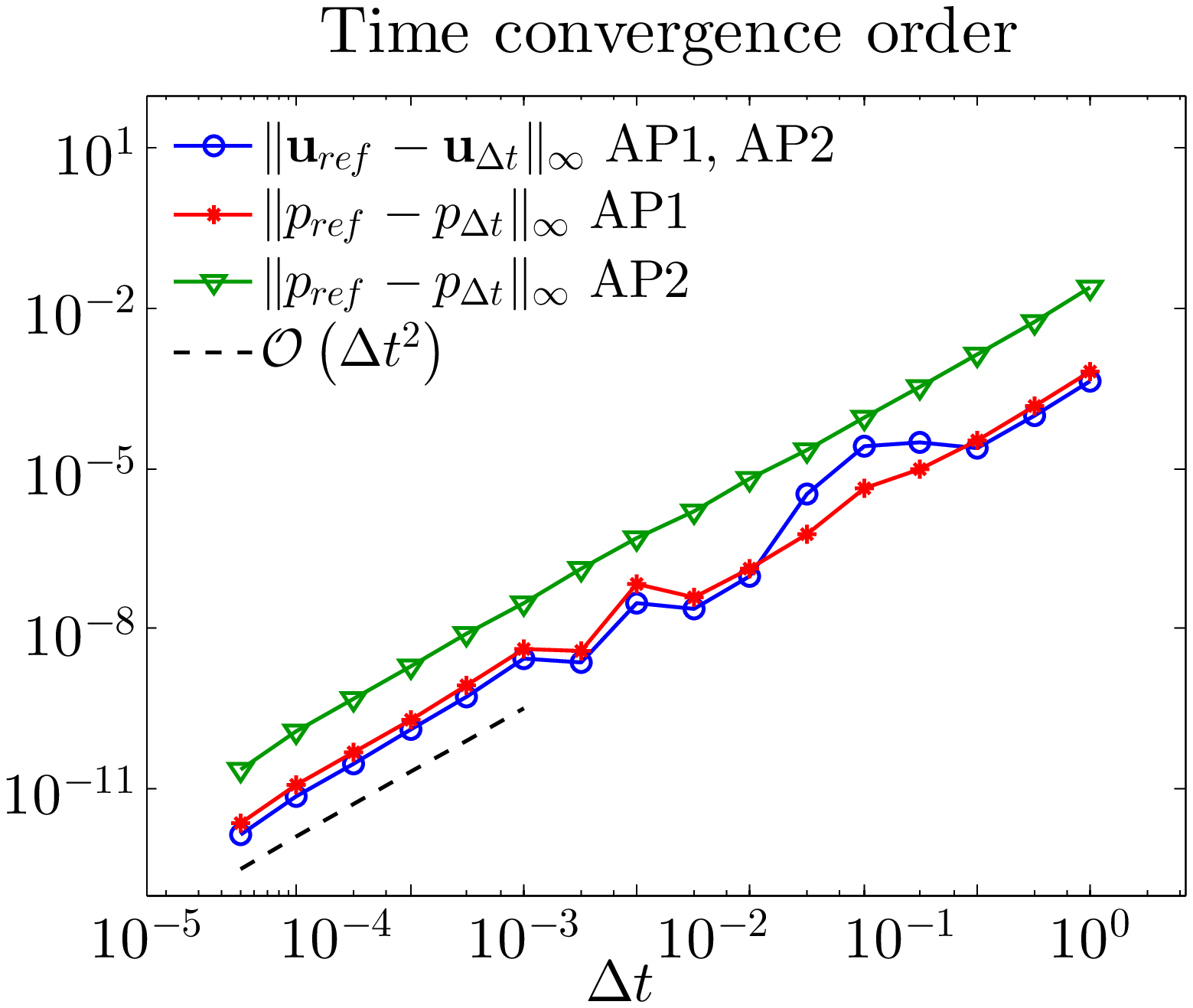} \label{fig:time_gt_rkc_AP12}}
\subfigure[ROCK2, AP1 and AP2W, CP0.]{
\includegraphics[trim=0.0cm 0.0cm 0.0cm 0.0cm, clip=true, height=0.2\textheight]{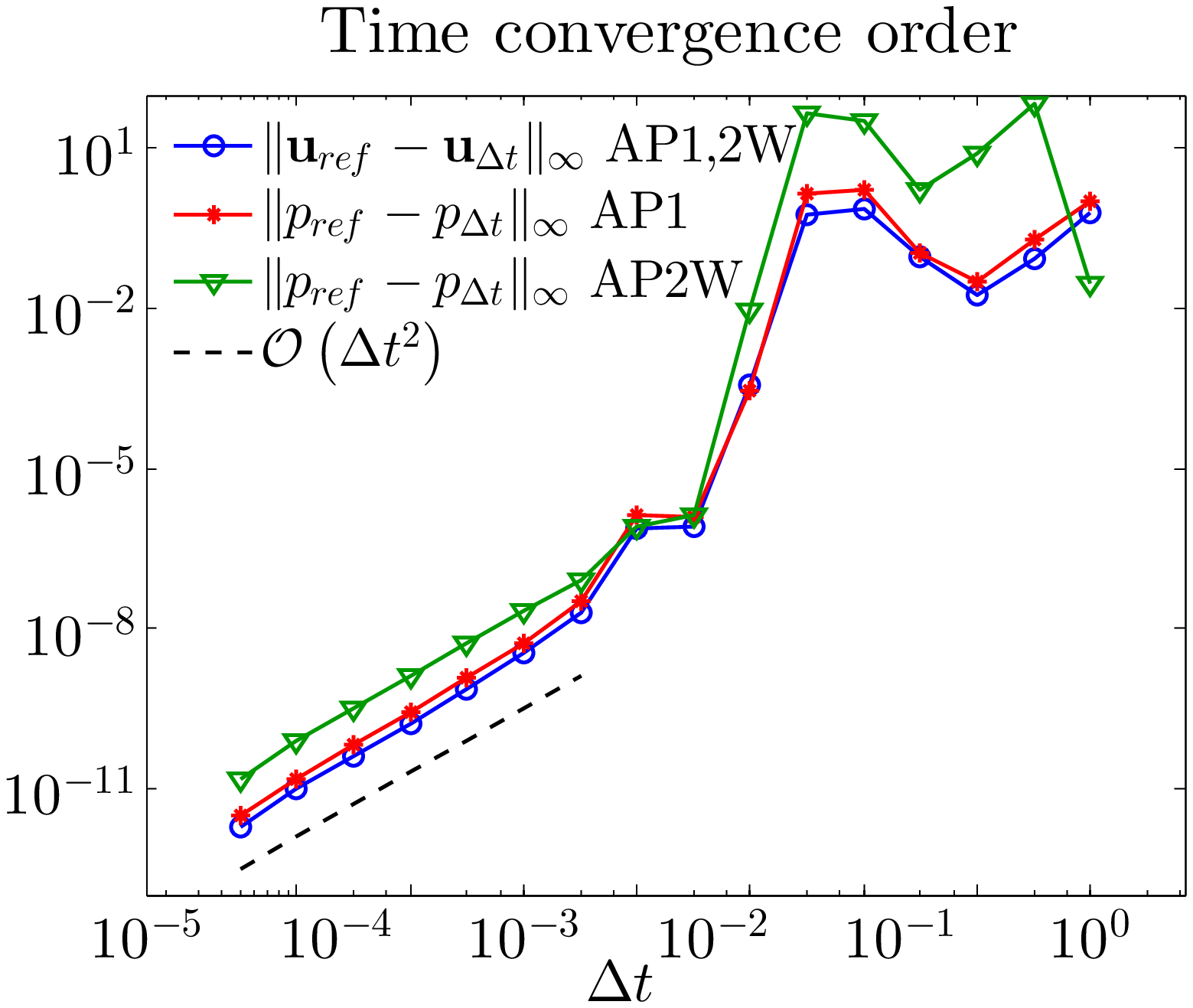} \label{fig:time_gt_rock2_AP12}}
\end{center}
\caption{Time convergence results of the Green-Taylor vortex.}
\label{fig:conv_time_gt}
\end{figure}

In Figures \ref{fig:conv_space_gt} and \ref{fig:conv_time_gt} we see the space and time convergence results for the Green-Taylor vortex. 
It is shown in Figure \ref{fig:conv_space_gt} that velocity and pressure achieves second order convergence in space for all the methods. For this problem we had to compute a reference solution with $\Delta x=1/1024$ in order to see a good convergence. In Figure \ref{fig:time_gt_pirock_PM1} we see that PIROCK, PM1, CP0 achieve the expected order of convergence for all the quantities but again it is not stable for $\Delta t \geq 10^{-2}$. In Figures \ref{fig:conv_time_gt}(b-g) we show the results of RKC and ROCK2 for PM1, PM1V, AP1, AP2 and AP2W. All of them have the expected order of convergence and as for the forced flow problem convergence results (see section \ref{sec:conv_ff}) RKC, PM1V, CP0 exhibits second order of convergence for the pressure $p^1$ even if first order was expected. Looking at Figures \ref{fig:conv_time_gt}(f, g) we notice that the pressure computed with AP1 is much more related to the velocity accuracy than the one computed with AP2 or AP2W. We believe that the reason is that using AP2 or AP2W the pressure is computed using the average values $\phi_i$s, hence large but very local errors are damped. Using AP1 a point value (the velocity) is used, so no damping occurs.

\subsection{Numerical efficiency}\label{sec:eff_gt}
\begin{figure}[!hbtp]
\begin{center}
\subfigure[Vel. efficiency of PM1, PM1V and CP0.]{
\includegraphics[trim=0.0cm 0.0cm 0.0cm 0.0cm, clip=true, height=0.2\textheight]{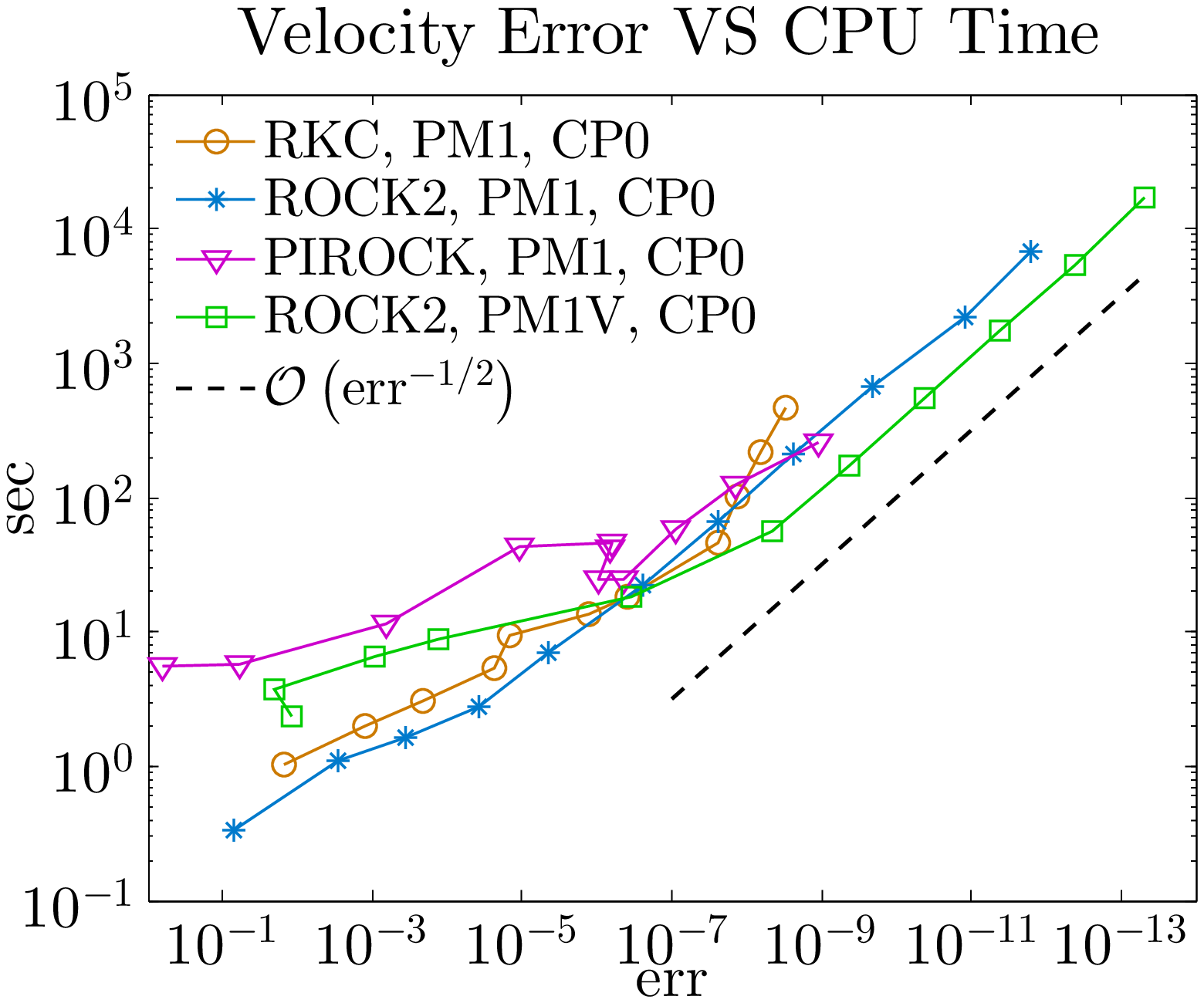} \label{fig:eff_gt_Proj_Vel_CP0}}
\subfigure[Vel. efficiency of PM1, PM1V and CP1.]{
\includegraphics[trim=0.0cm 0.0cm 0.0cm 0.0cm, clip=true, height=0.2\textheight]{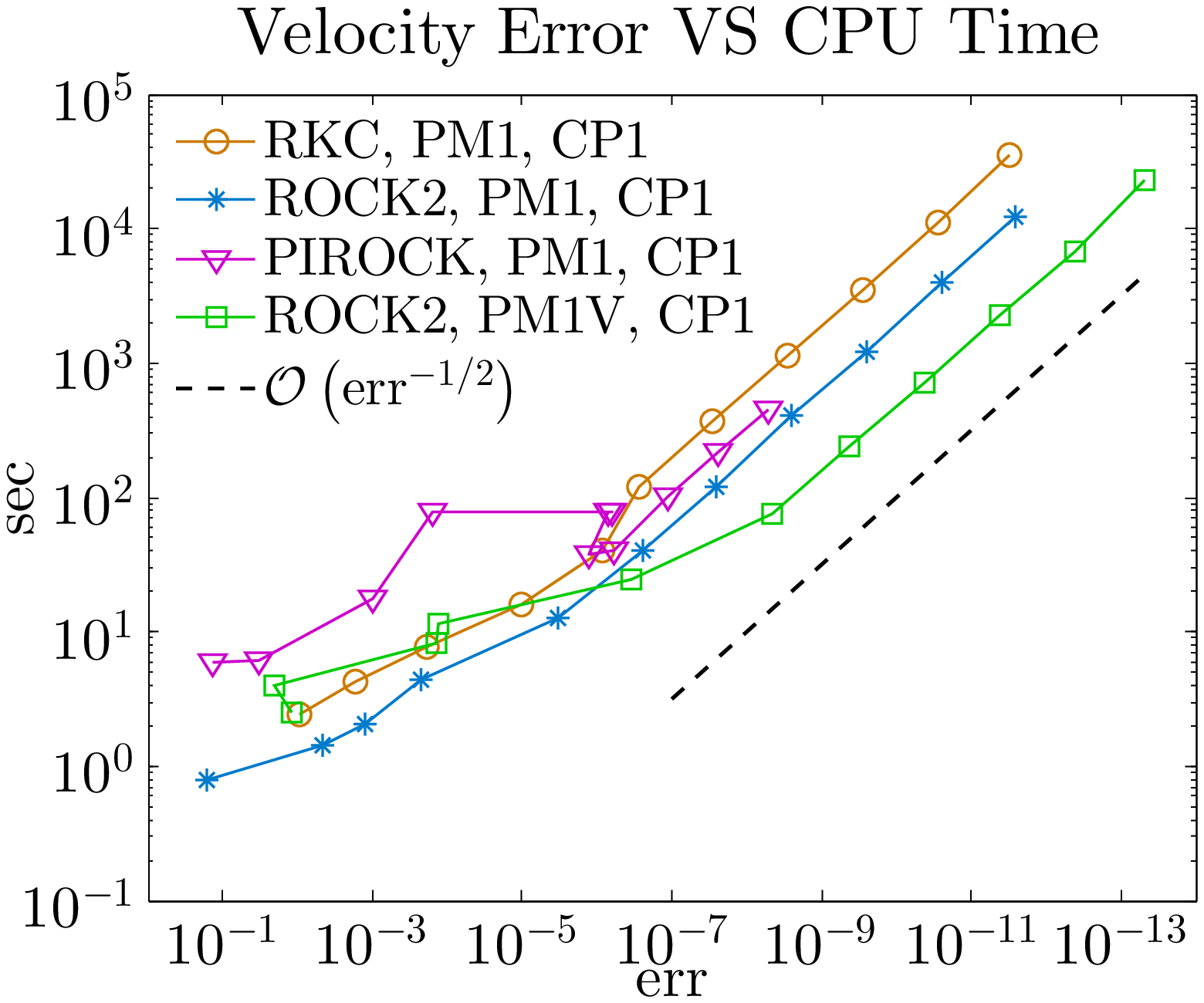} \label{fig:eff_gt_Proj_Vel_CP1}} \\
\subfigure[Pre. efficiency of PM1, PM1V and CP0.]{
\includegraphics[trim=0.0cm 0.0cm 0.0cm 0.0cm, clip=true, height=0.2\textheight]{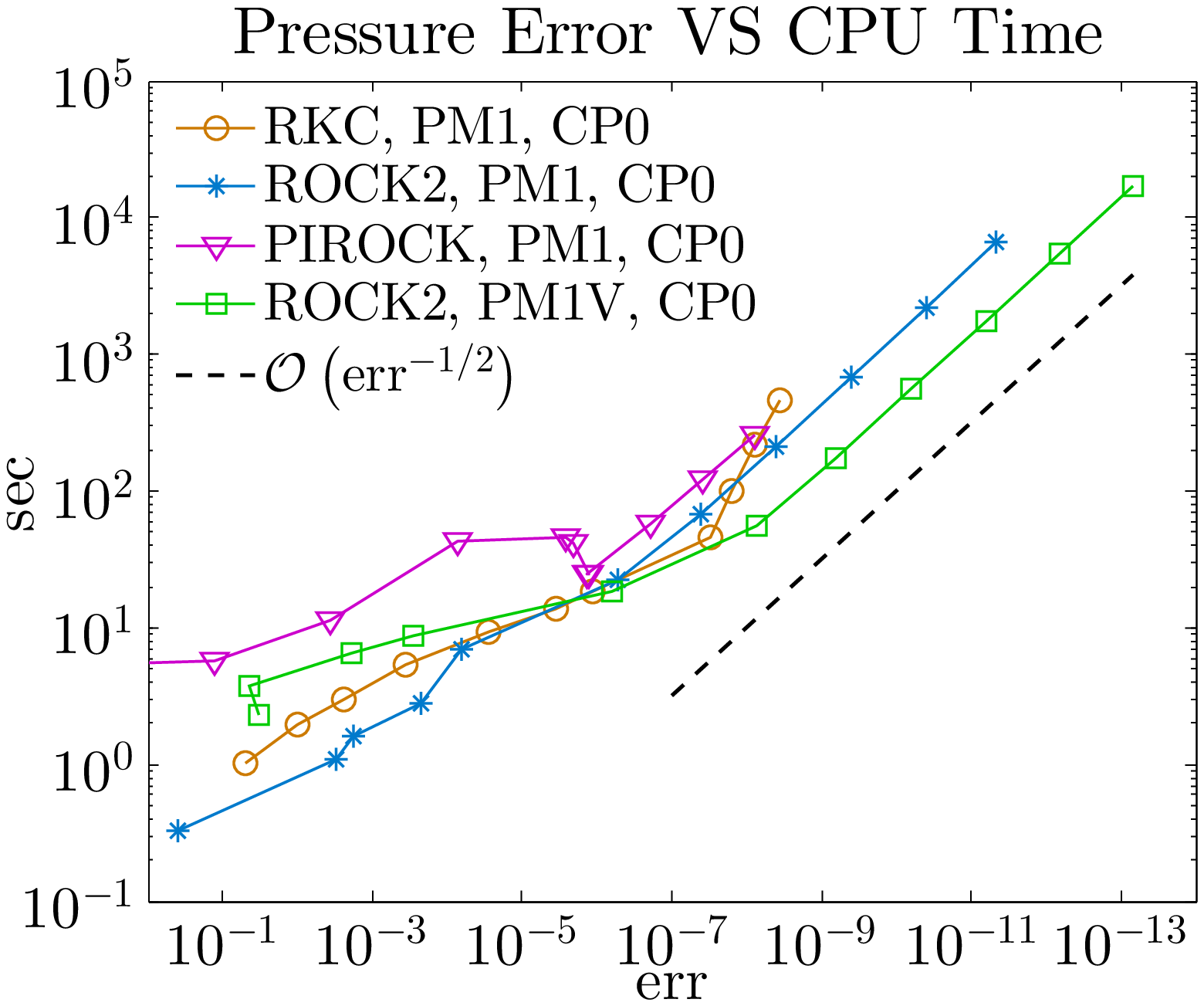} \label{fig:eff_gt_Proj_Pre_CP0}}
\subfigure[Pre. efficiency of PM1, PM1V and CP1.]{
\includegraphics[trim=0.0cm 0.0cm 0.0cm 0.0cm, clip=true, height=0.2\textheight]{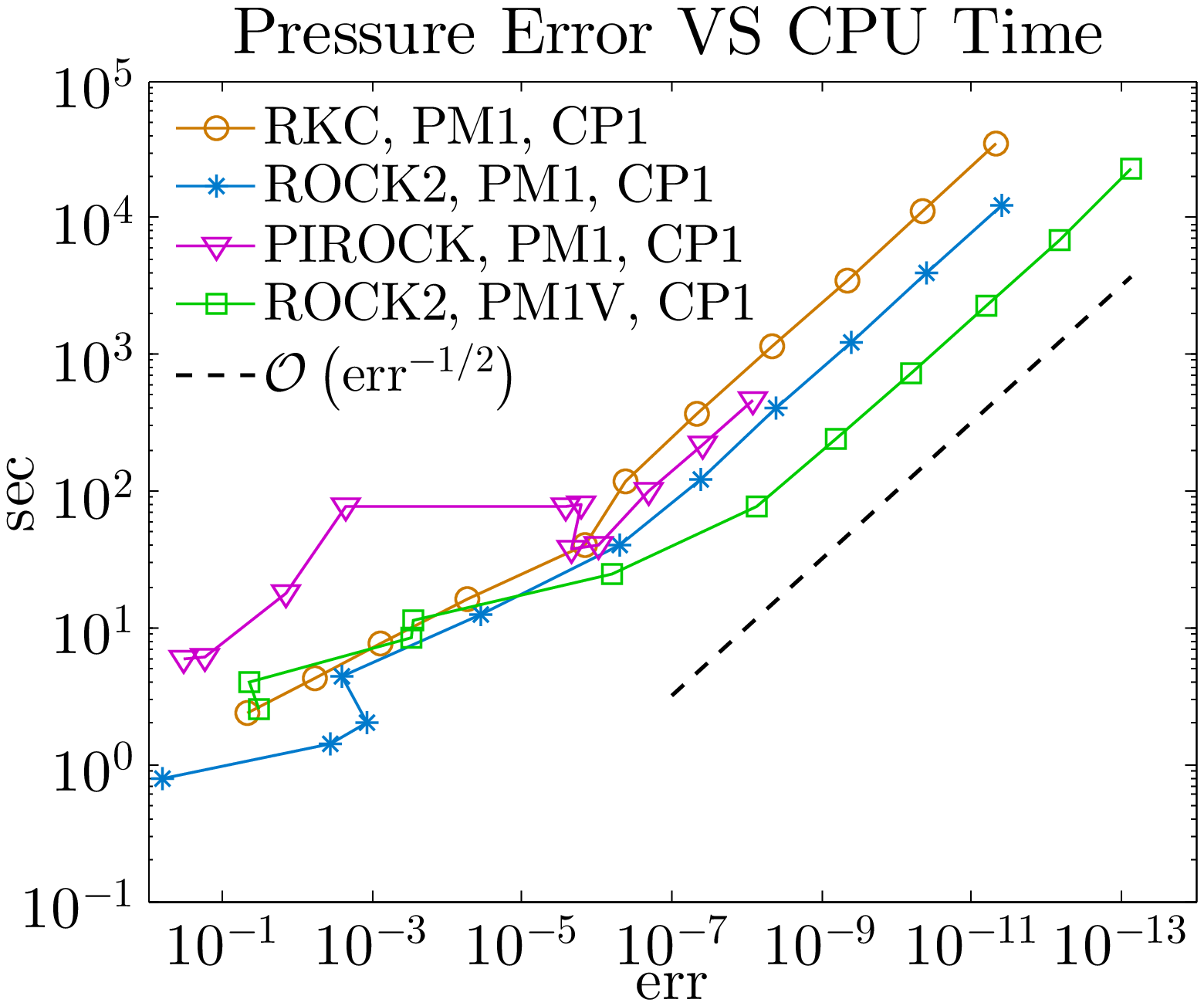} \label{fig:eff_gt_Proj_Pre_CP1}} \\
\subfigure[Vel. efficiency of AP1, AP2W and CP0.]{
\includegraphics[trim=0.0cm 0.0cm 0.0cm 0.0cm, clip=true, height=0.2\textheight]{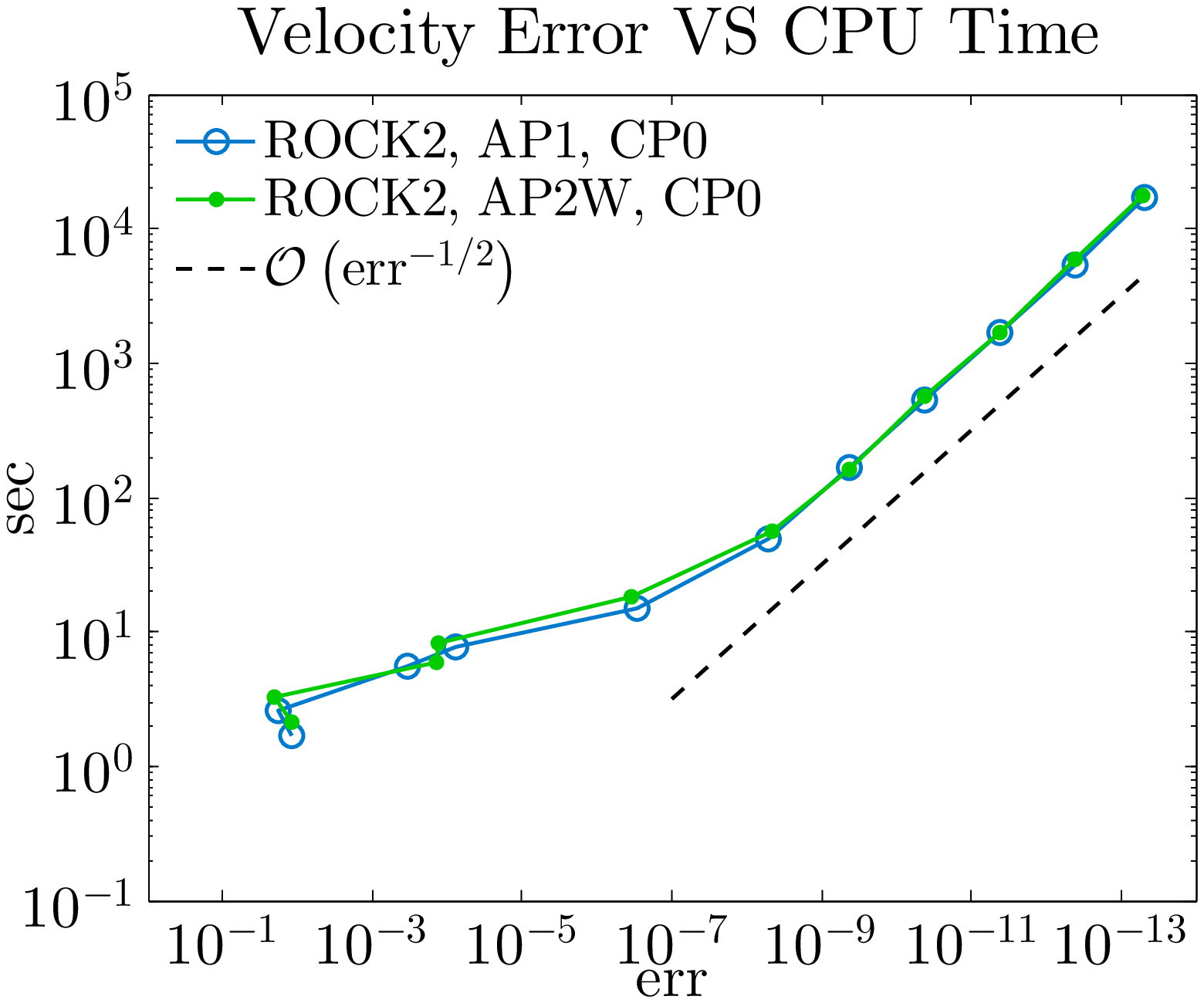} \label{fig:eff_gt_DAE_Vel_CP0}}
\subfigure[Vel. efficiency of AP1, AP2W and CP1.]{
\includegraphics[trim=0.0cm 0.0cm 0.0cm 0.0cm, clip=true, height=0.2\textheight]{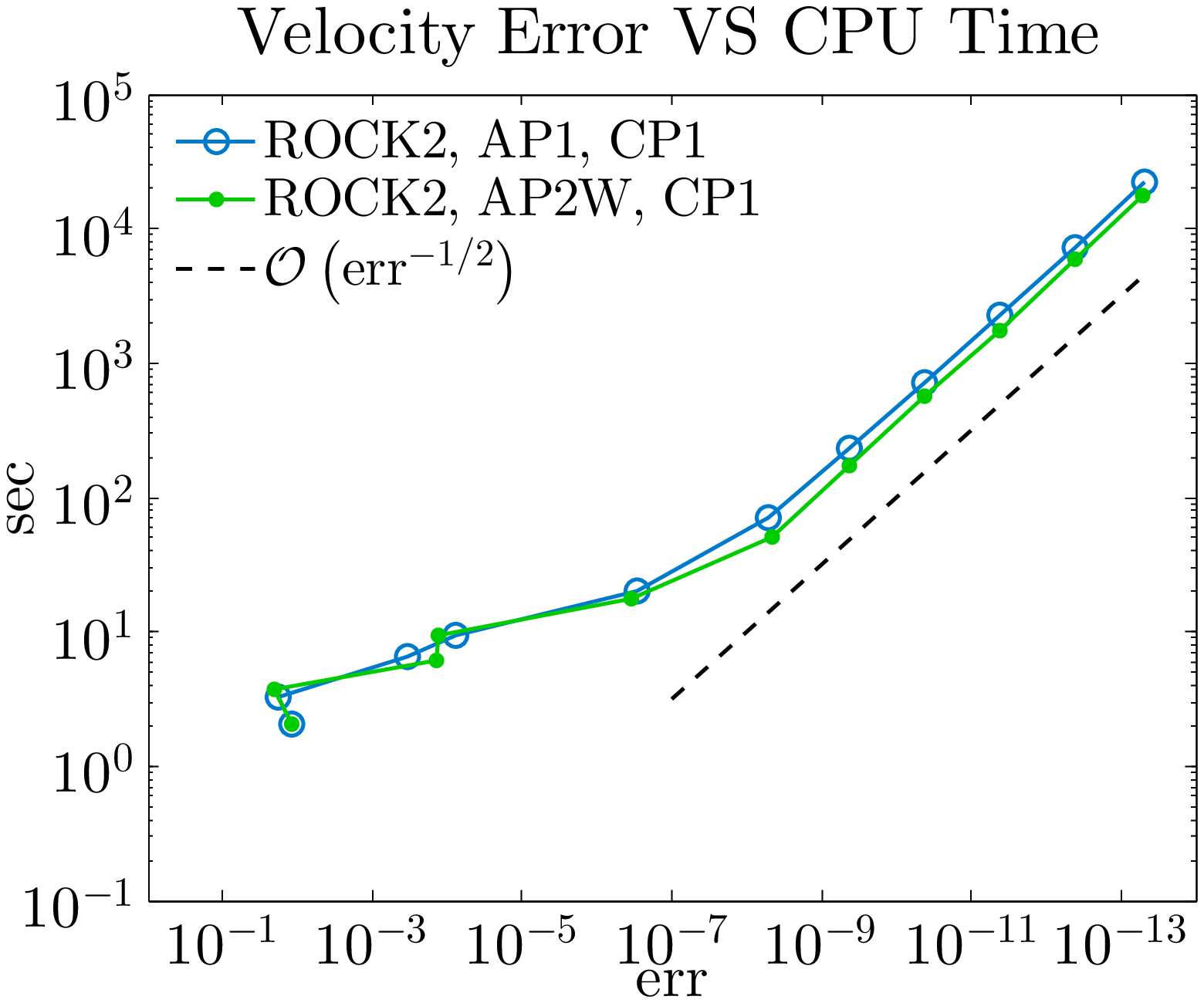} \label{fig:eff_gt_DAE_Vel_CP1}}\\
\subfigure[Pre. efficiency of AP1, AP2W and CP0.]{
\includegraphics[trim=0.0cm 0.0cm 0.0cm 0.0cm, clip=true, height=0.2\textheight]{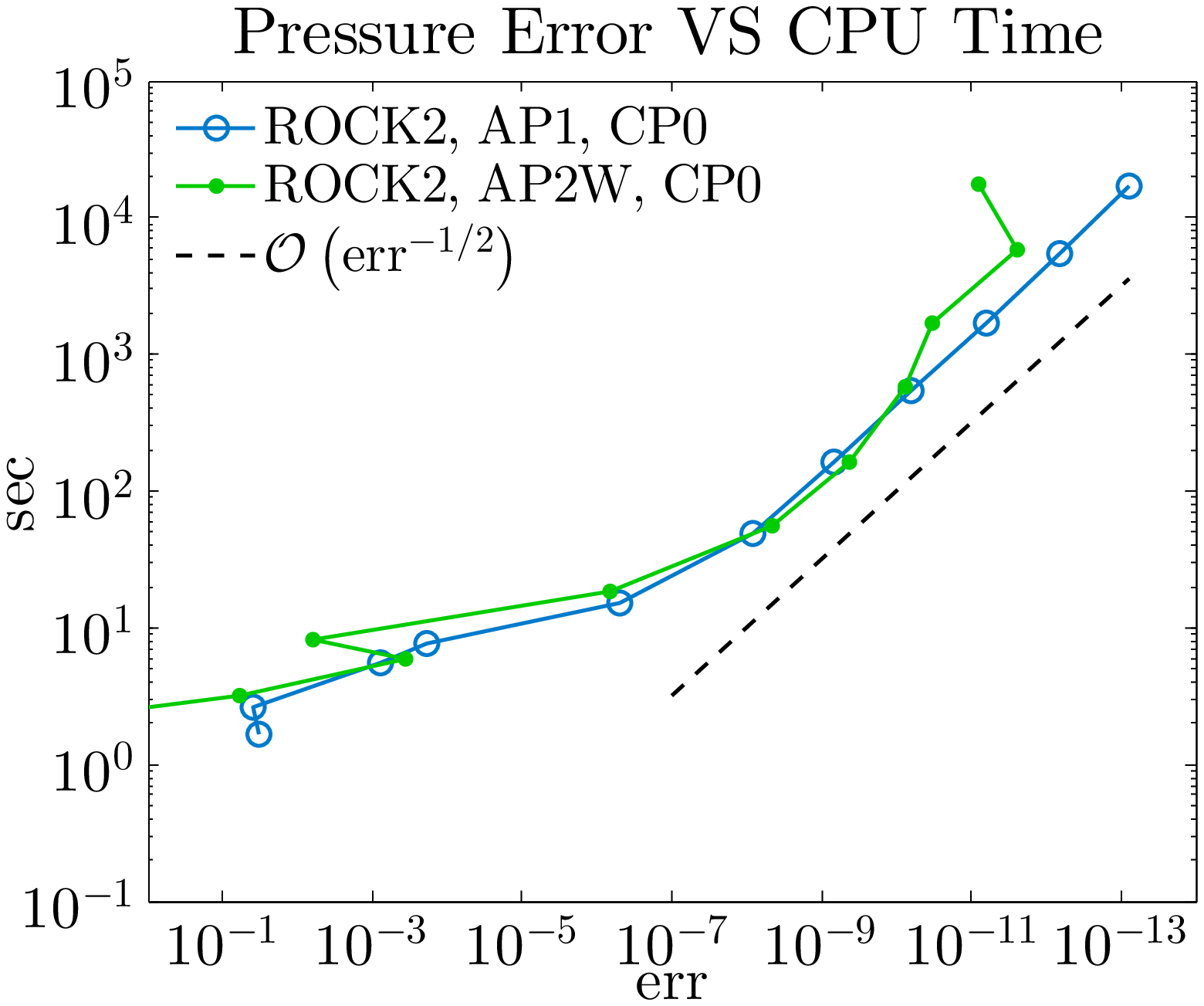} \label{fig:eff_gt_DAE_Pre_CP0}}
\subfigure[Pre. efficiency of AP1, AP2W and CP1.]{
\includegraphics[trim=0.0cm 0.0cm 0.0cm 0.0cm, clip=true, height=0.2\textheight]{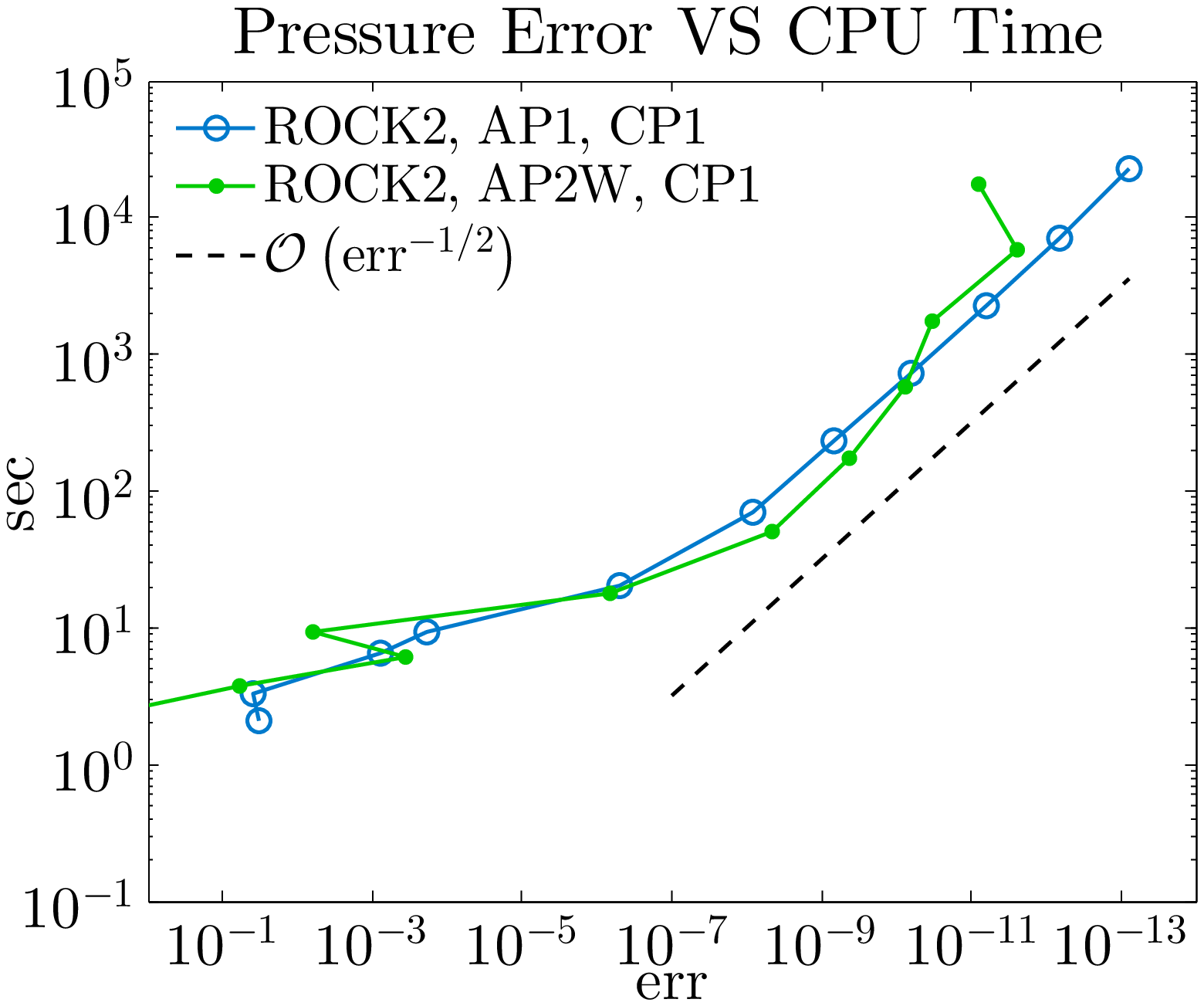} \label{fig:eff_gt_DAE_Pre_CP1}}
\end{center}
\caption{Numerical efficiency tests of the Green-Taylor vortex.}
\label{fig:eff_gt}
\end{figure}

In Figure \ref{fig:eff_gt} we compare the numerical efficiency of the methods for the Green-Taylor vortex. From Figures \ref{fig:eff_gt}(a-d) we see that again among the PM1, PM1V methods ROCK2, PM1V is the best one. RKC, PM1, CP0 is better than ROCK2, PM1, CP0 only for a few tolerances, less than in the case of the forced flow. The gain in efficiency when using PM1V is less than in the forced flow problem (see Figures \ref{fig:eff_ff}(a-d)) because in this example there is no forcing term, thus the overhead of the extra projections is more significant. Nonetheless PM1V remains more efficient than PM1. Comparing Figures \ref{fig:eff_gt}(a,c) with Figures \ref{fig:eff_gt}(b,d) we remark that again RKC has a much better performance when a second order pressure is computed after each step while ROCK2 and PIROCK are less sensible to this. Figures \ref{fig:eff_gt}(g,h) show that again AP1 is more reliable than AP2W for the pressure but on the other hand Figure \ref{fig:eff_gt}(f) shows that AP2W is faster than AP1 when the pressure is computed after each time step. Comparing the first column of figures in \ref{fig:eff_gt} with the second column we can conclude that for ROCK2 and PIROCK one should use CP1 only when accurate pressures are needed at intermediate time steps because at $t=1$ there is no gain in accuracy using CP1 instead of CP0, as for the forced flow.

Again the best methods are ROCK2, AP1 and ROCK2, AP2W. The choice depends on the preference between the velocity and the pressure accuracy.

\section{The Lid-driven cavity problem}
The lid-driven cavity problem is considered as the classical test case for the validation of new Navier-Stokes simulation codes and algorithms. Its geometry and boundary conditions are very simple. It consists of a square domain with Dirichlet boundary conditions on all sides. Three sides are stationary while the last one is moving. The mathematical formulation is equations \eqref{eq:ns} with boundary conditions
\begin{align}
\begin{aligned}
\bu(t,x,1) &= (1,0)^\top  \quad x\in [0,1], \\
\bu(t,x,0) &= (0,0)^\top \quad x\in [0,1], \\
\bu(t,x,y) &= (0,0)^\top \quad x\in\{0,1\}, y\in [0,1] .
\end{aligned}
\end{align}

It is known \cite{bench-cavity} that this problem presents singularities at the corners $(0,1),(1,1)$ where the velocity is discontinuous. This property makes difficult to properly evaluate the accuracy of the numerical results, mainly in the neighborhood of these points where pressure and vorticity are not finite. Since the space accuracy is associated to the smoothness and boundedness of the of the solution derivatives it is then completely lost at these corners. This is the reason why one would use the regularized driven cavity problem \cite{reg-cavity} in order to evaluate the solution accuracy. In this regularized problem the velocity is smoothed, so that the above problems are avoided. On the other hand these difficulties makes the driven cavity problem a widely used test case for the evaluation of incompressible flow solvers. 

In the following we will show the time convergence results and compare our results with tabular data given in \cite{ghia}.

\subsection{Convergence order}
\begin{figure}[!hbtp]
\begin{center}
\subfigure[PIROCK, PM1, CP0.]{
\includegraphics[trim=0.0cm 0.0cm 0.0cm 0.0cm, clip=true, height=0.2\textheight]{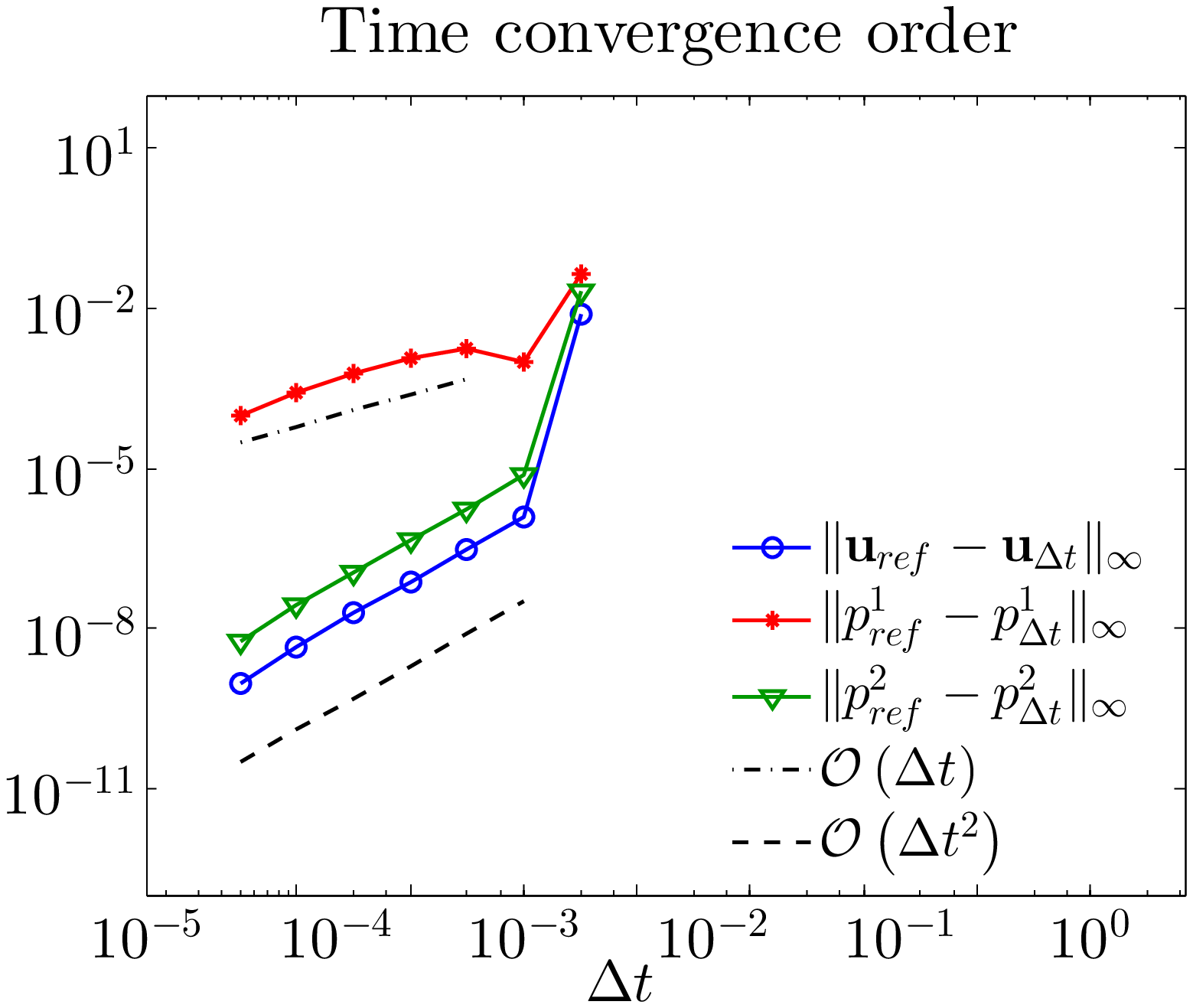} \label{fig:time_dc_pirock_PM1}} \\
\subfigure[RKC, PM1, CP0.]{
\includegraphics[trim=0.0cm 0.0cm 0.0cm 0.0cm, clip=true, height=0.2\textheight]{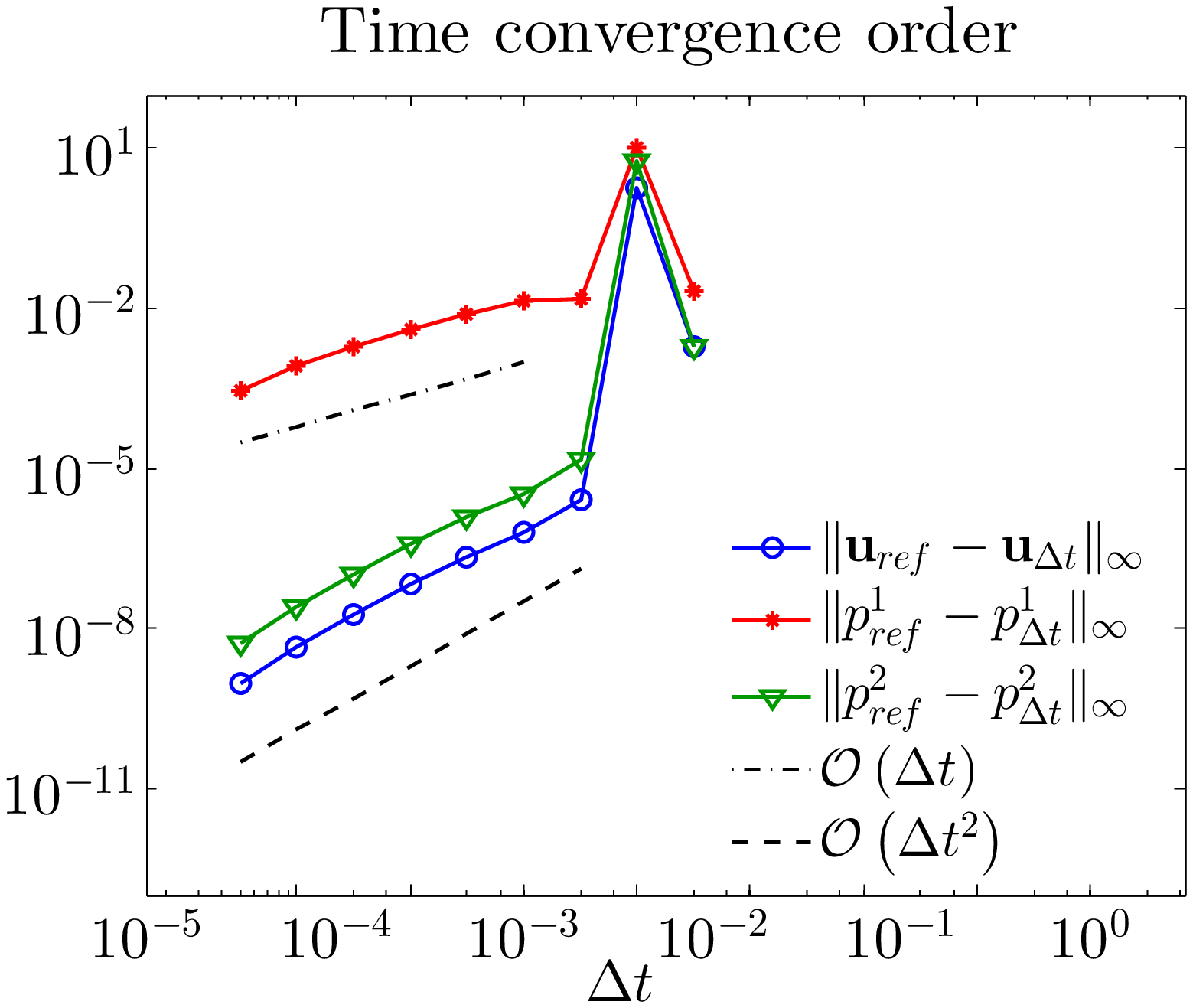} \label{fig:time_dc_rkc_PM1}}
\subfigure[ROCK2, PM1, CP0.]{
\includegraphics[trim=0.0cm 0.0cm 0.0cm 0.0cm, clip=true, height=0.2\textheight]{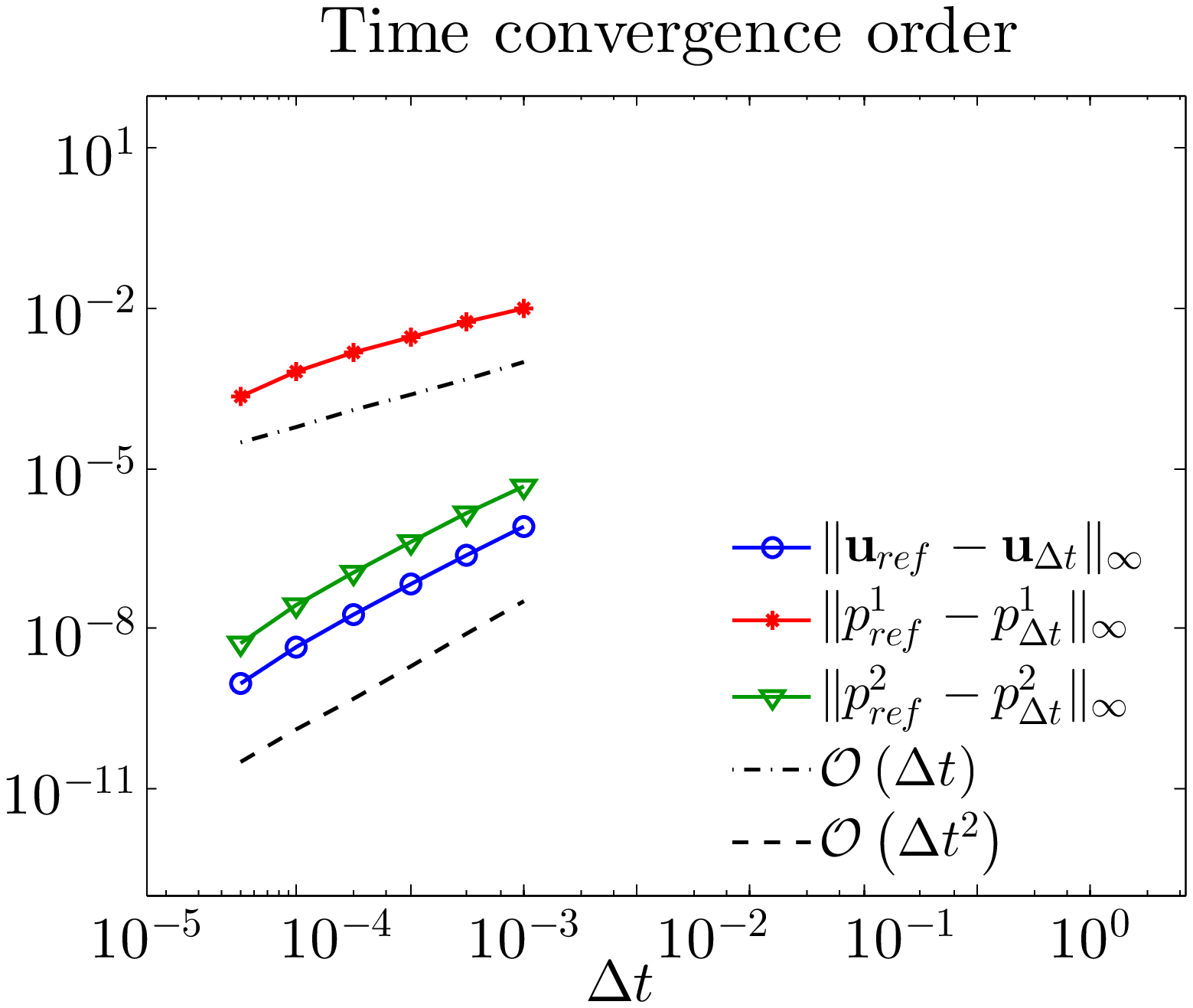} \label{fig:time_dc_rock2_PM1}} \\
\subfigure[RKC, PM1V, CP0.]{
\includegraphics[trim=0.0cm 0.0cm 0.0cm 0.0cm, clip=true, height=0.2\textheight]{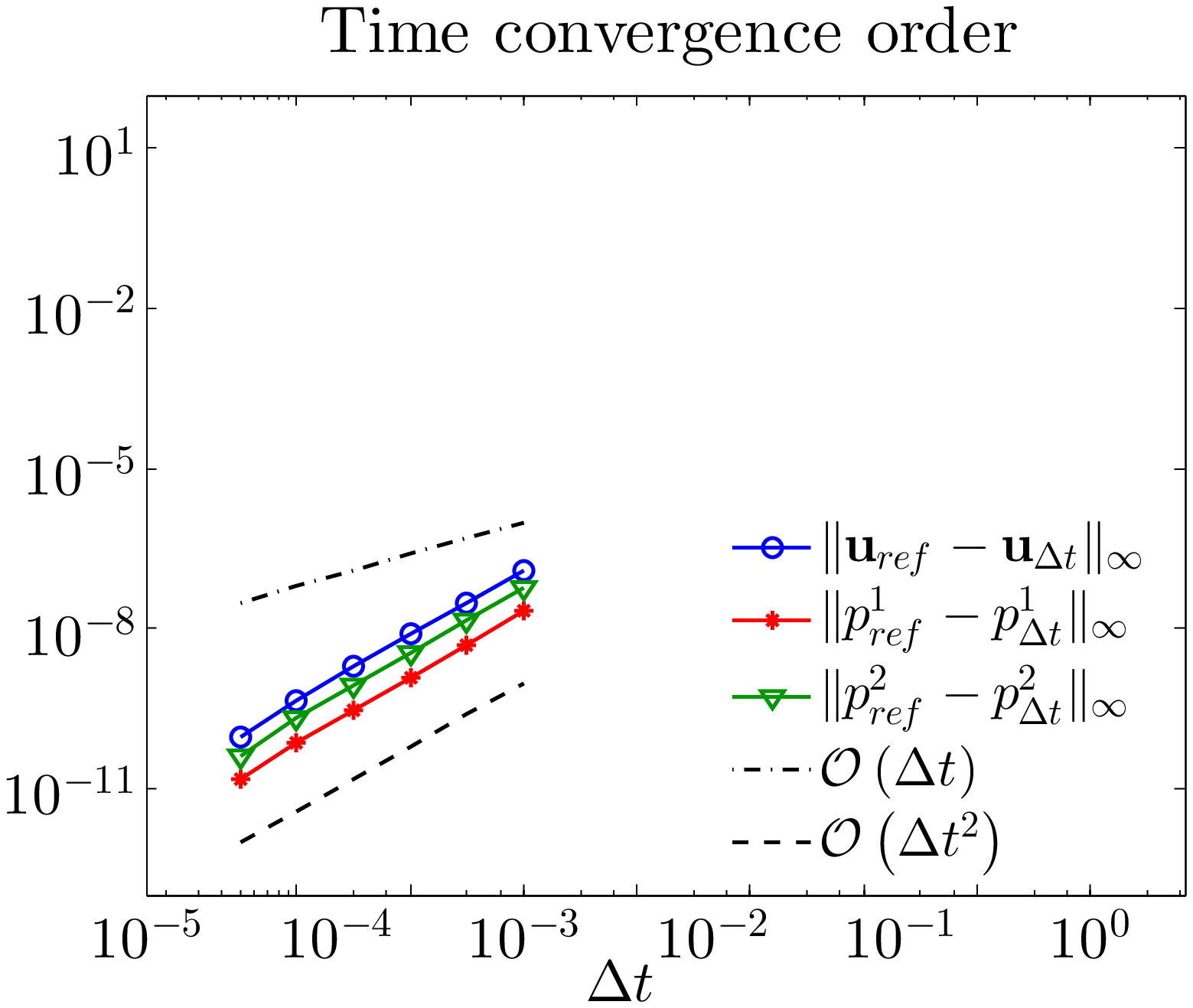} \label{fig:time_dc_rkc_PM1V}}
\subfigure[ROCK2, PM1V, CP0.]{
\includegraphics[trim=0.0cm 0.0cm 0.0cm 0.0cm, clip=true, height=0.2\textheight]{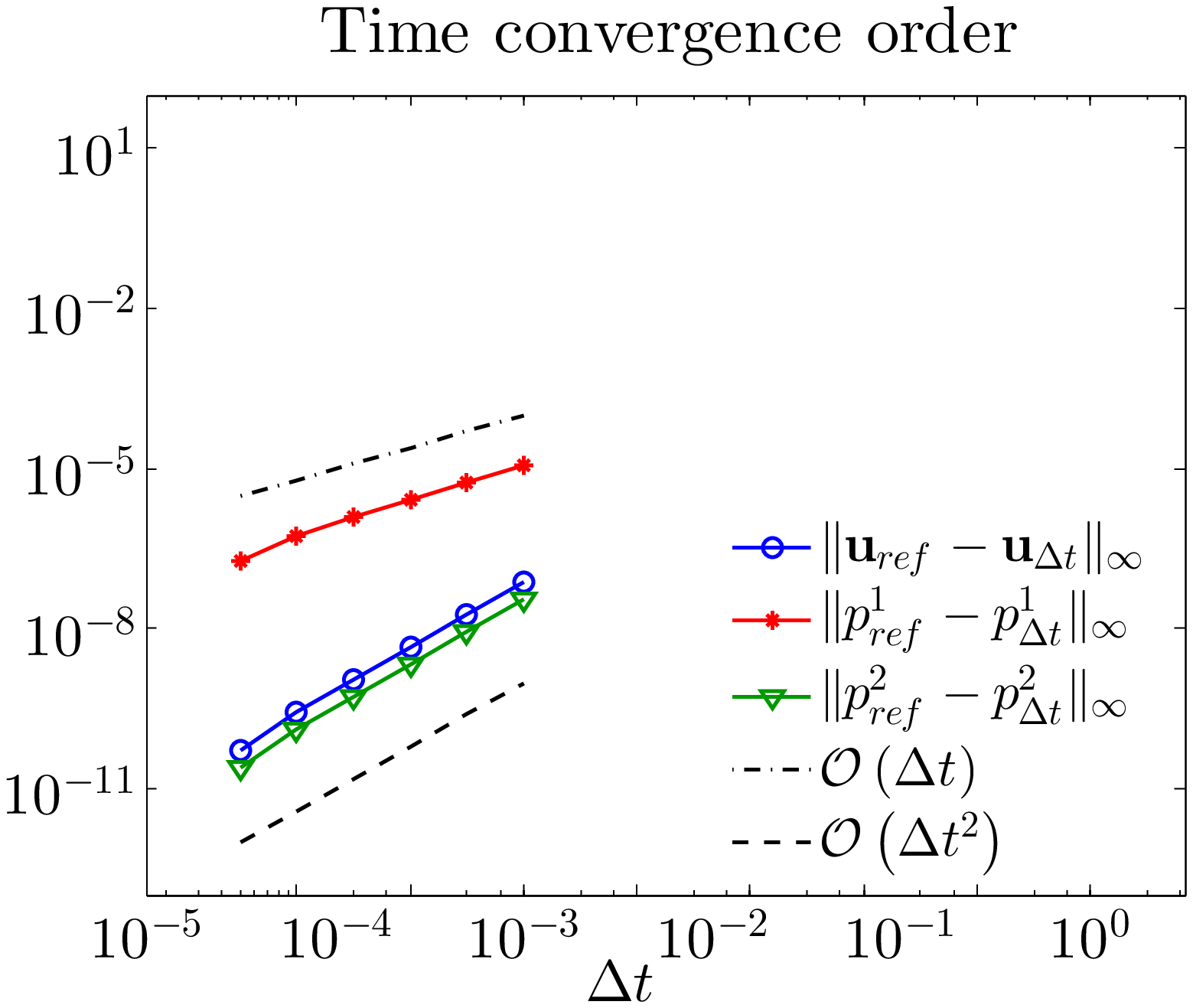} \label{fig:time_dc_rock2_PM1V}} \\
\subfigure[RKC, AP2 and AP2, CP0.]{
\includegraphics[trim=0.0cm 0.0cm 0.0cm 0.0cm, clip=true, height=0.2\textheight]{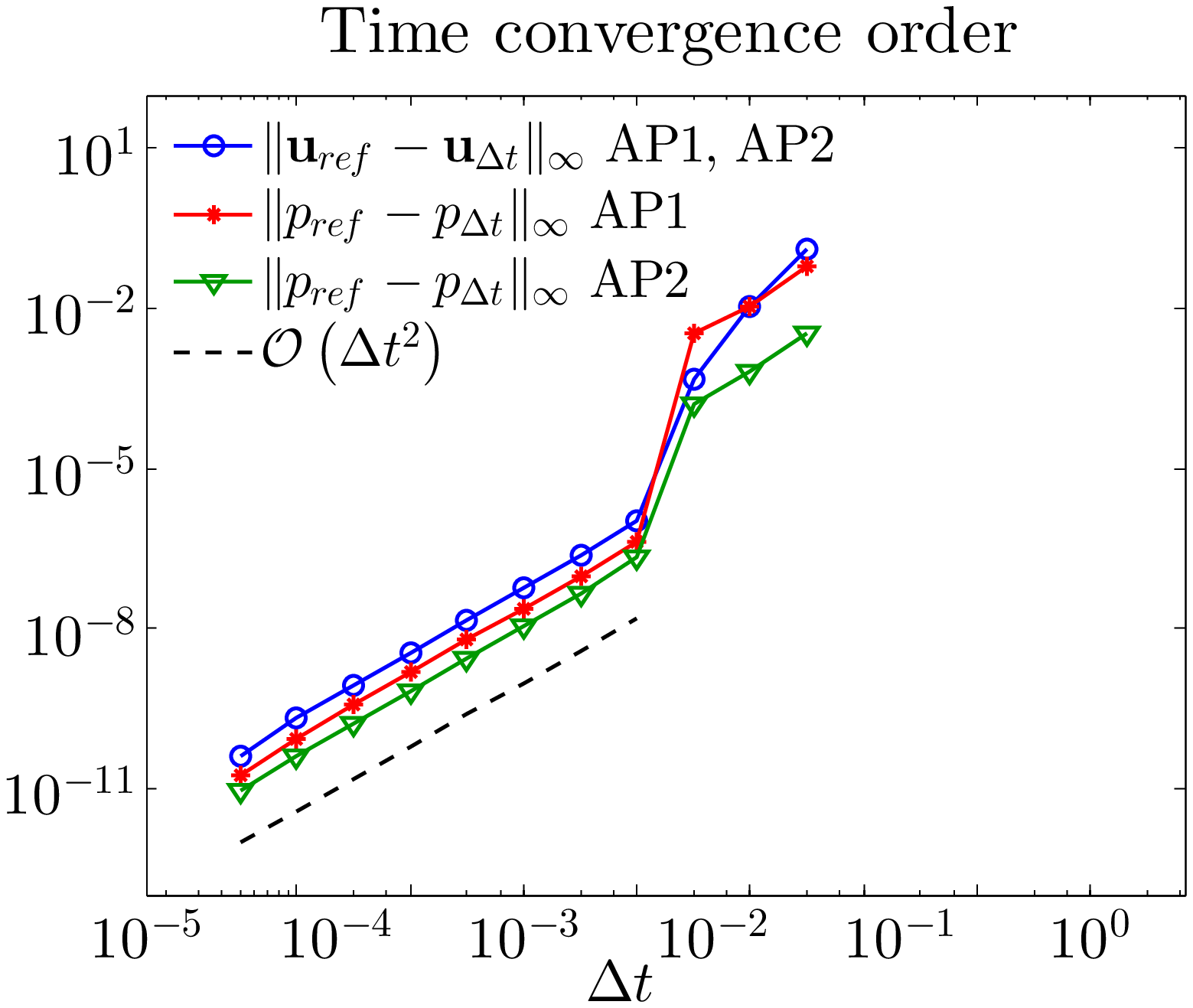} \label{fig:time_dc_rkc_AP12}}
\subfigure[ROCK2, AP1 and AP2W, CP0.]{
\includegraphics[trim=0.0cm 0.0cm 0.0cm 0.0cm, clip=true, height=0.2\textheight]{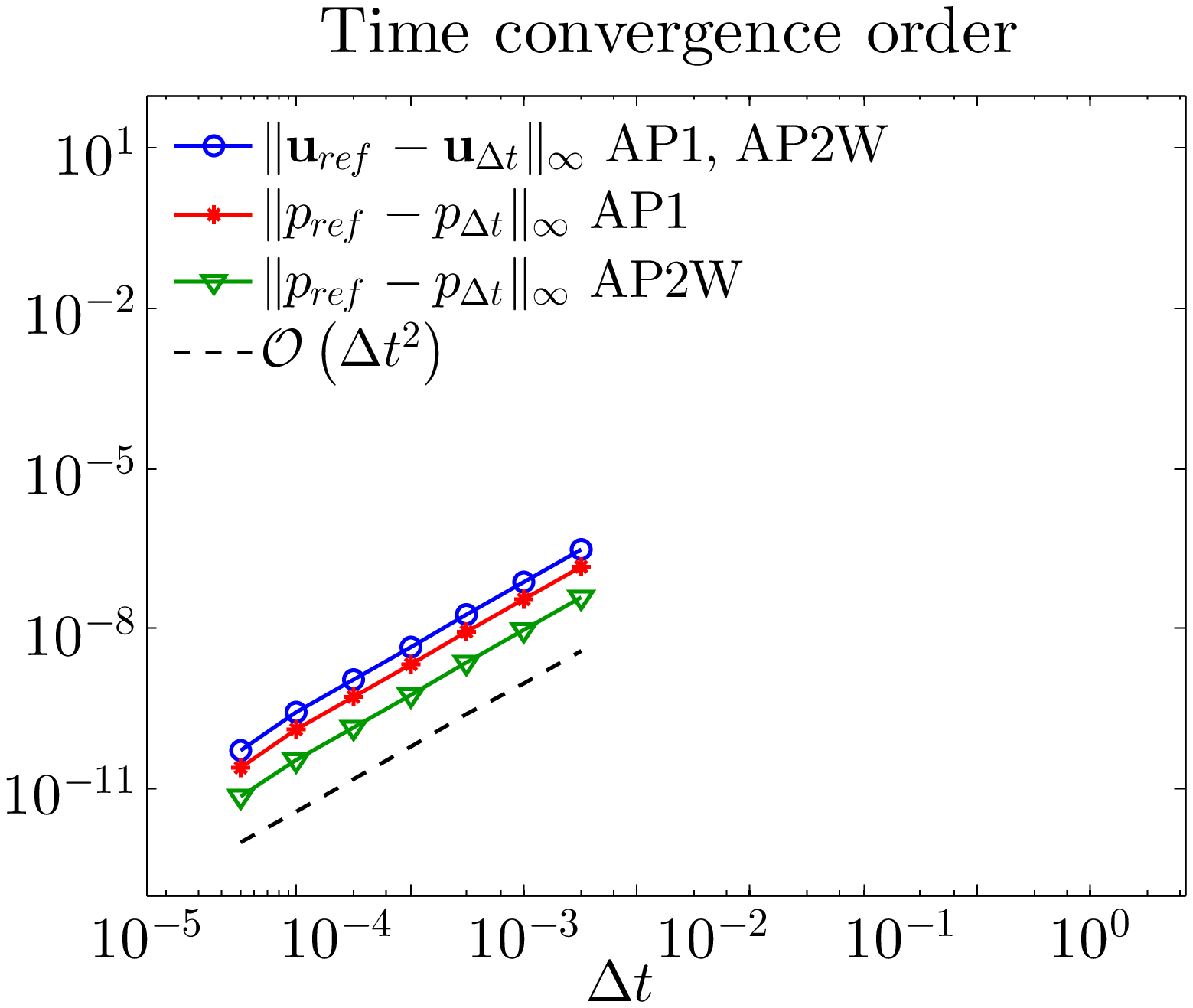} \label{fig:time_dc_rock2_AP12}}
\end{center}
\caption{Time convergence results of the lid-driven cavity problem.}
\label{fig:conv_time_dc}
\end{figure}

As we said the solution is non smooth at the corners and it follows that the expected space convergence order is not achieved. However time convergence is not affected by these singularities. In Figure \ref{fig:conv_time_dc} we show the time convergence of the lid-driven cavity problem for the methods listed in section \ref{sec:desc_conv}. We see that all the methods have the expected convergence order. In this problem not only PIROCK but also RKC and ROCK2 are not stable for $\Delta t \geq 10^{-2}$. From Figures \ref{fig:conv_time_dc}(f,g) we see that RKC is more stable than ROCK2, we suspect that it is because its stability region covers a larger area in the complex plane. Looking at Figures \ref{fig:conv_time_dc}(d,f) we see that one more time algorithm \eqref{eq:rkc_rec_projf} performs better than \eqref{eq:rkc_rec_proj}, in this case from the stability point of view.

\subsection{Comparing our results with an established reference}
Since this problem has been solved many times there is a great deal of data to compare with in the literature. An established reference is given by Ghia et al. in \cite{ghia} since it contains tabular results for various Reynolds numbers. We have compared the results given by the methods that allow time step adaptivity at the stationary point $t=500$ with the data given in \cite{ghia} for $\Rey=1000$. Like in \cite{ghia} we did the simulation with a $128\times 128$ grid. In our simulation we used an adaptive time step with $rtol=atol=10^{-3}$. 

\begin{figure}[!hbtp]
\begin{center}
\subfigure[Solution $u$.]{
\includegraphics[trim=0.0cm 0.0cm 0.0cm 0.0cm, clip=true, width=0.45\textwidth]{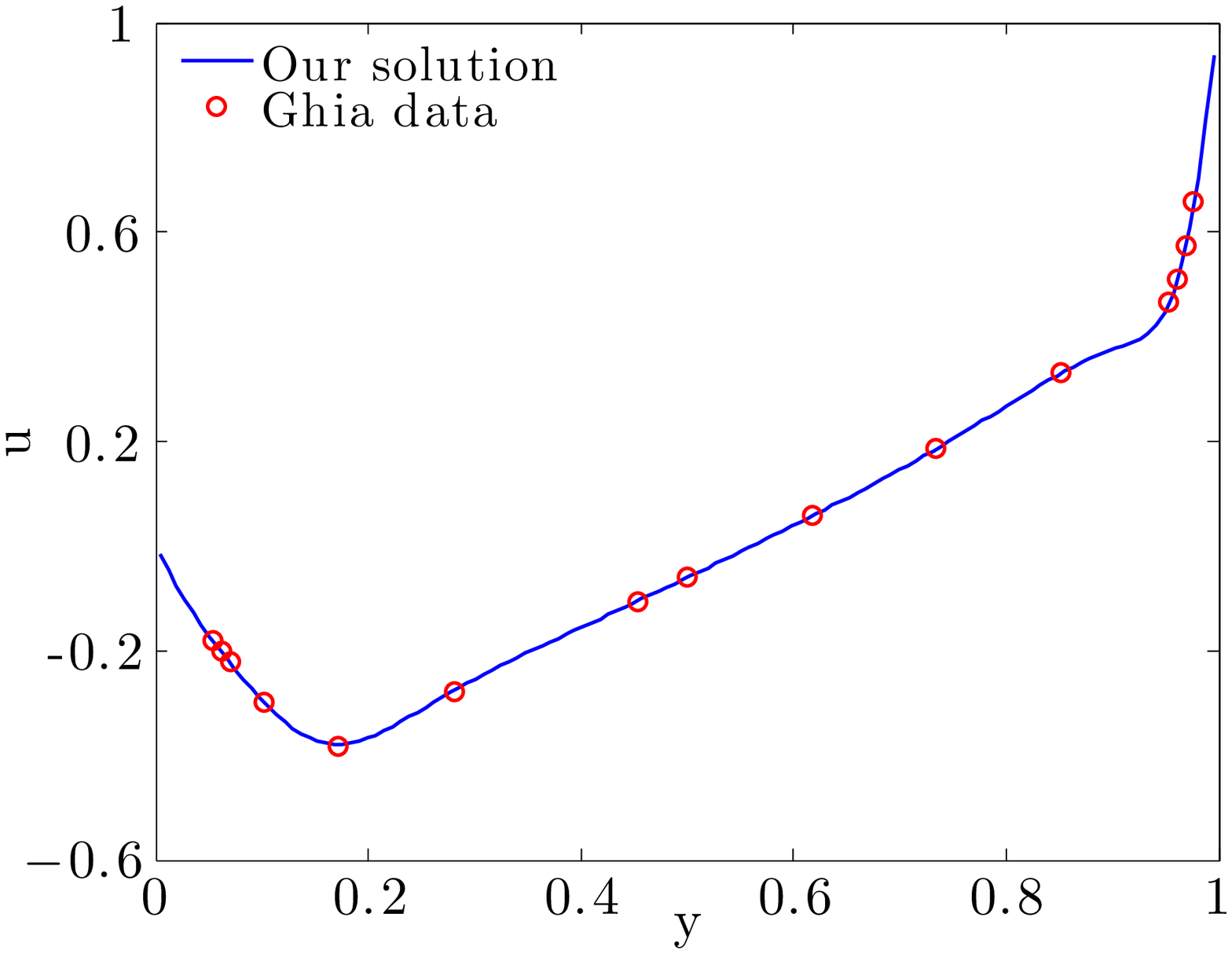} \label{fig:all_u}}
\subfigure[Solution $v$.]{
\includegraphics[trim=0.0cm 0.0cm 0.0cm 0.0cm, clip=true, width=0.45\textwidth]{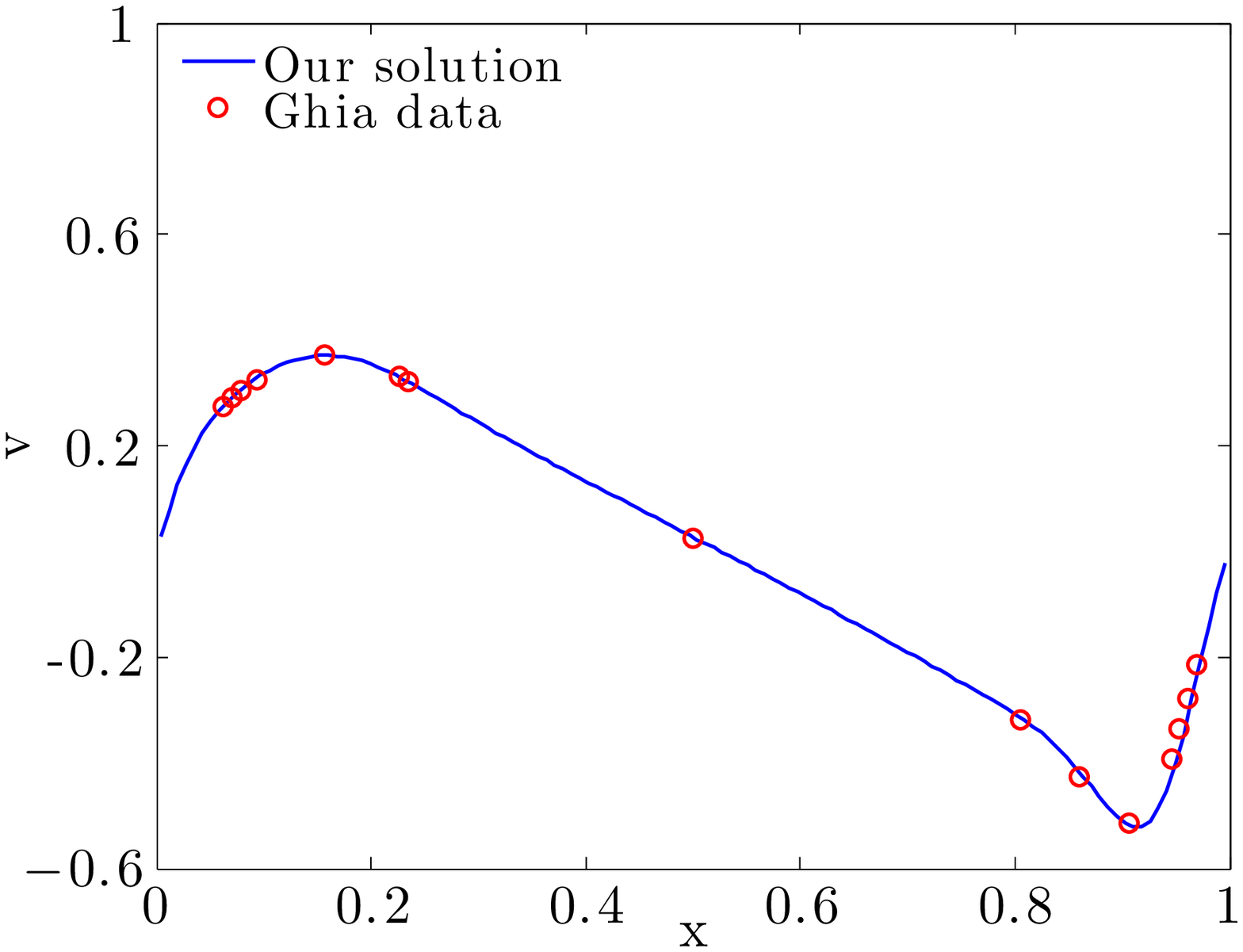} \label{fig:all_v}}
\end{center}
\caption{Comparing our results with Ghia et al.}
\label{fig:ghia}
\end{figure}

In the Figures the velocity $u$ is taken at the vertical centerline of the cavity, the velocity $v$ at the horizontal centerline. We did the test for the methods listed in section \ref{sec:desc_comp} and the results look the same, so they are summarized in Figure \ref{fig:ghia} under the label 'Our solution'. We see a very good matching between our solutions on the ones obtained by Ghia et al. in \cite{ghia}

\newpage
\thispagestyle{empty}
\mbox{}
\clearpage

\chapter{Conclusions}

In the previous chapters we have presented the stabilized explicit Runge-Kutta methods RKC, ROCK2 and PIROCK and two approaches for handling the Navier-Stokes equations incompressibility constraint: the projection method and the differential algebraic approach. Then we described the method used for the spatial discretization and for the solution of the Poisson equation on this particular grid. At the end we showed the numerical results with discussions.

Already in the first chapters we noticed that RKC has a limitation concerning time step adaptivity which is not present in ROCK2 and PIROCK. This issue can be solved easily but it was not our purpose. Due to this limit when using RKC we could enable time step adaptivity only with the projection method.

The main results obtained with numerical experiments are the following. We showed that all the methods work and attain the expected second order of convergence, in space and time. The stability tests showed that the projection method has a stability domain shorter than the one obtained for ordinary differential equations. On the contrary the differential algebraic approach has exactly the same stability domain as the ordinary differential equations case. 
We did also a very close comparison of the projection method and the differential algebraic approach using ROCK2. It turned out that with the used time step the projection method has a non negligible boundary layer produced by the boundary conditions for the virtual velocity. And, more important, the boundary layer is not present at all when the differential algebraic approach is used. Using a fixed time step we saw that the differential algebraic approach was $1.3$ times slower but on the other hand the solution was $644$ times more accurate! In fact from the numerical efficiency tests it turned out that the differential algebraic approach is much more efficient than the projection method.

The comparison between RKC, ROCK2 and PIROCK showed that ROCK2 is the most efficient and PIROCK is braked by his stability issues. We showed that the accuracy of the pressure does not affects the velocity, so accurate pressures should be computed only if needed. When comparing the methods with different Reynolds numbers again ROCK2 is way more efficient than RKC and the differential algebraic approach is much better than the projection method. 

For the differential algebraic approach two different methods have been used for the pressure computation: one uses a Poisson problem the other an interpolation with Lagrange polynomials. The latter is limited by the stages' order of the method but we showed how this limit can be circumvented. It turned out that the method which uses the Poisson problem is slightly more expensive but also more accurate and reliable. The one which uses the interpolation computes the pressure with a negligible cost. So the choice between them depends on how much the user is interested in the pressure accuracy. Anyway both of them give a second order accurate pressure.

As a last experiment we compared our results with an established reference and we found a very good match.

The differential algebraic approach resulted to be highly superior than the projection method: it his more accurate, it has a longer stability region and it permits to have a better estimation of the local error. Moreover the projection method cannot provide solutions of order higher than two. Differently, the differential algebraic approach accuracy depends only on the order of the Runge-Kutta method. It follows that it can be applied to higher order stabilized explicit methods, as ROCK4 \cite{rock4}, in a straightforward way.

In future we would like to compare the different methods with a finer grid, where the Poisson problems are more expensive. This might change the difference in efficiency between the projection method and the differential algebraic approach. Also the difference between the two methods for computing the pressure might change.

The main limit of RKC and ROCK2 is that their stability region is in the neighborhood of the negative real axis, hence they are unstable when the Navier-Stokes equations are advection dominated. In the future our main goal is to apply the differential algebraic approach to PIROCK (it has not been done here). We expect that it behaves better than the projection method. We saw that the differential algebraic approach has better stability properties, the same of ordinary differential equations, so it could be that the stability problems are also fixed. We expect that PIROCK behaves similarly to ROCK2 for diffusion dominated flows. For advection dominated flows we expect from it a better behavior due to its large stability domain on the complex plane.

\addcontentsline{toc}{chapter}{Bibliography}

\end{document}